\UseRawInputEncoding
\documentclass{memo-l}

\usepackage{amsfonts, latexsym, amscd, float,amsmath,amssymb,mathrsfs,bm,multirow,graphics}
\usepackage[dvips]{graphicx}
\usepackage[percent]{overpic}

\newtheorem{theorem}{Theorem}[chapter]
\newtheorem{lemma}[theorem]{Lemma}
\newtheorem{corollary}[theorem]{Corollary}

\theoremstyle{definition}

\theoremstyle{remark}
\newtheorem{remark}[theorem]{Remark}

\numberwithin{section}{chapter}
\numberwithin{equation}{chapter}
\numberwithin{figure}{chapter}

\newcommand{\R}{{\Bbb R}}

\newcommand{\C}{{\Bbb C}}

\newcommand{\Z}{{\Bbb Z}}
\newcommand{\re}{\text{\upshape Re\,}}
\newcommand{\im}{\text{\upshape Im\,}}
\newcommand{\Li}{\text{\upshape Li}}
\newcommand{\dist}{\text{\upshape dist}}
\newcommand{\sgn}{\text{\upshape sgn\,}}

\def\Xint#1{\mathchoice
{\XXint\displaystyle\textstyle{#1}}%
{\XXint\textstyle\scriptstyle{#1}}%
{\XXint\scriptstyle\scriptscriptstyle{#1}}%
{\XXint\scriptscriptstyle\scriptscriptstyle{#1}}%
\!\int}
\def\XXint#1#2#3{{\setbox0=\hbox{$#1{#2#3}{\int}$}
\vcenter{\hbox{$#2#3$}}\kern-.5\wd0}}

\def\dashint{\Xint-}

\makeindex

\begin{document}

\frontmatter

\title[On the Asymptotics to all Orders of the Riemann Zeta Function]{On the Asymptotics to all Orders of the Riemann Zeta Function and of a Two-Parameter Generalization of the Riemann Zeta Function}


\author{Athanassios S. Fokas}
\address{Department of Applied Mathematics and Theoretical Physics, University of Cambridge, Cambridge, CB3 0WA, UK}
\curraddr{}
\email{T.Fokas@damtp.cam.ac.uk}
\thanks{The authors are grateful to Professor D. Bump for several useful suggestions that led to a significant improvement of the manuscript. They are also grateful to Kostis Kalimeris for many helpful remarks.}
\thanks{The first author was supported by the EPSRC, UK via a senior fellowship and by the Guggenheim Foundation, USA. He is deeply grateful to his students Michail Dimakos (partially supported by the Onassis Foundation) and Dionyssis Mantzavinos (partially supported by the EPSRC, UK) for performing a large number of numerical experiments during the early stages of this project. He also expresses his sincere gratitude to Anthony Ashton, Bryce McLeod, Joe Keller, David Stuart, and Eugene Shargorodsky, for extensive discussions and important suggestions.
}

\author{Jonatan Lenells}
\address{Department of Mathematics, KTH Royal Institute of Technology, 100 44 Stockholm, Sweden}
\curraddr{}
\email{jlenells@kth.se}
\thanks{The second author acknowledges support from the EPSRC, UK, the G\"oran Gustafsson Foundation, Sweden, the European Research Council, Grant Agreement No. 682537, and the Swedish Research Council, Grant No. 2015-05430.}

\date{}



\dedicatory{Dedicated to Trevor Stuart with deep gratitude}

\maketitle

\tableofcontents

\begin{abstract}
We present several formulae for the large $t$ asymptotics of the Riemann zeta function $\zeta(s)$, $s=\sigma+i t$, $0\leq \sigma \leq 1$, $t>0$, which are valid to all orders. A particular case of these results coincides with the classical results of Siegel. Using these formulae, we derive explicit representations for the sum $\sum_a^b n^{-s}$ for certain ranges of $a$ and $b$. In addition, we present precise estimates relating this sum with the sum $\sum_c^d n^{s-1}$ for certain ranges of $a, b, c, d$. 
We also study a two-parameter generalization of the Riemann zeta function which we denote by $\Phi(u,v,\beta)$, $u\in \mathbb{C}$, $v\in \mathbb{C}$, $\beta \in \mathbb{R}$. Generalizing the methodology used in the study of $\zeta(s)$, we derive asymptotic formulae for $\Phi(u,v,\beta)$.\end{abstract}


\mainmatter


\part{Asymptotics to all Orders of the Riemann Zeta Function}

\chapter{Introduction}
It is well known, see for example theorem 4.11 in Titchmarsh \cite{T1986}, that the Riemann zeta function $\zeta(s)$, $s=\sigma+ i t$, $\sigma, t \in \R$, satisfies the equation
\begin{align}\label{1.1}
\zeta(s)=\sum_{n \leq x}\frac{1}{n^s} -\frac{x^{1-s}}{1-s}+O \left( \frac{1}{x^\sigma} \right), \qquad |t| < \frac{2 \pi x}{c}, \quad c>1, \quad \sigma>0, \quad t \to \infty.
\end{align}
Here we present a systematic investigation of asymptotic formulae for $\zeta(s)$. In this connection, in chapter \ref{sec2} we derive the following exact formula (see theorem \ref{theorem2.1}):
\begin{align}\nonumber
& \zeta(1-s) =  \sum^{\left[ \frac{\eta}{2 \pi} \right]}_{n=1} n^{s-1}+\frac{1}{(2\pi)^s} \left[- \frac{\eta^s}{s} +
\left( e^{\frac{i\pi s}{2}} \int^{\infty e^{i\phi_1}}_{\eta e^{- \frac{i\pi}{2}}} + e^{-\frac{i\pi s}{2}} \int^{\infty e^{i\phi_2}}_{\eta e^{\frac{i\pi}2}} \right) \frac{z^{s-1}dz}{e^z - 1}\right],
	\\ \label{1.2}
&\hspace{5cm} 0 <\eta <\infty, \quad -\frac{\pi}{2}<\phi_j< \frac{\pi}{2}, \quad j=1,2,
\end{align}
where $[\frac{\eta}{2 \pi}]$ denotes the integer part of $\frac{\eta}{2 \pi}$ and the contours of integration in the first and second integrals in the rhs of equation (\ref{1.2}) are the rays from $\eta  \exp(-i\pi/2)$ to $\infty \exp (i\phi_1)$ and from $\eta  \exp(i\pi/2)$ to $\infty \exp(i \phi_2)$, respectively.

\section{The large $t$ asymptotics of $\zeta(s)$ valid to all orders}
Equation (\ref{1.2}) suggests a separate analysis for the cases $t<\eta $, $t=\eta $ and $t>\eta $. The first two cases are analysed in theorems \ref{th3.1} and \ref{th3.2} of chapter \ref{sec3}. Theorem \ref{th3.1} presents the large $t$ asymptotics of $\zeta(s)$ valid to all orders in the case of $(1+\epsilon)t < \eta $. For example, the formulae of theorem \ref{th3.1} yield the following equation for the leading asymptotic terms in the case when $(1+\epsilon)t < \eta$ for some $\epsilon > 0$:
\begin{align}\nonumber
\zeta(s) = &\; \sum_{n =1}^{[\frac{\eta}{2\pi}]} n^{-s} - \frac{1}{1-s} \left(\frac{\eta}{2\pi}\right)^{1-s} 
+ \frac{2i\eta^{-s}}{(2\pi)^{1-s}} \biggl\{
-i\arg(1 - e^{i \eta})
+ \frac{t - i \sigma}{\eta} \re \Li_{2}(e^{i\eta}) \biggr\}
	\\ \label{zetaformula2.1intro}
&    + O\biggl(\frac{t^2}{\eta^{2 + \sigma}}\biggr),
\qquad (1+\epsilon)t < \eta < \infty, \quad 0 \leq \sigma \leq 1, \quad t \to \infty,	
\end{align}
where the error term is uniform for $\eta, \sigma$ in the above ranges and the polylogarithm $\Li_m(z)$ is defined by
\begin{align*}
  \Li_m(z) = \sum_{k=1}^\infty \frac{z^k}{k^m}, \qquad m \geq 1.
\end{align*}  
Theorem \ref{th3.2} presents the large $t$ asymptotics of $\zeta(s)$ valid to all orders in the case of $\eta=t$. The formula of theorem \ref{th3.2} yields the following equation for the leading asymptotic terms:
\begin{align}\nonumber
\zeta(s) = & \sum_{n =1}^{[\frac{t}{2\pi}]} n^{-s} - \frac{1}{1-s} \left(\frac{t}{2\pi}\right)^{1-s} 
	\\ \nonumber
&- \frac{2 it^{-s}e^{it}}{(2\pi)^{1-s}} \biggl\{-\frac{e^{-it}}{2t}[t + i(\sigma -1)]
- i \im \Li_1(e^{it})
+ \frac{i\sigma}{t} \re \Li_2(e^{it})
	\\\nonumber
& \hspace{1.9cm} - \frac{1}{t} \im \Li_3(e^{it})\biggr\}
	\\\nonumber
&+ \frac{ t^{-s}e^{it}}{(2\pi)^{1-s}}\biggl\{ \frac{1+i}{2} \sqrt{\pi t} -\frac{1}{3} i (3 \sigma -2)
	\\\nonumber
&\hspace{1.9cm}
+ \frac{i-1}{24 \sqrt{t}} \sqrt{\pi } \left(6\sigma ^2-6 \sigma +1\right)
+ \frac{1}{135 t } \left(45 \sigma^3-45 \sigma ^2+4\right)
	\\\nonumber
&\hspace{1.9cm}
-\frac{1+i}{576 t^{\frac{3}{2}}}\sqrt{\pi } \left(36 \sigma^4-24 \sigma ^3-24 \sigma ^2+12 \sigma +1\right)
\biggr\}
 + O\bigl(t^{-\sigma - 2}\bigr),
	\\ \label{zetaformula2.2intro}
& \hspace{5cm}  0 \leq \sigma \leq 1, \quad t \to \infty,	
\end{align}
where the error term is uniform for $\sigma$ in the above range.

The best estimate for the growth of $\zeta(s)$ as $t \to \infty$, is not based on equation (\ref{1.1}) but on the well known approximate functional equation, see for example equation (4.12.4) on page 79 of Titchmarsh \cite{T1986}:
\begin{align}\label{1.6}
\zeta(s)=\sum_{n \leq x}\frac{1}{n^s}+\chi(s) \sum_{n \leq y} \frac{1}{n^{1-s}} +O \left(x^{-\sigma} + t^{\frac{1}{2}-\sigma} y^{\sigma-1} \right), \qquad xy=\frac{t}{2\pi},
\end{align}
where $0<\sigma<1$ and $\chi(s)$ is defined by
$$\chi(s) = \frac{(2\pi)^s}{\pi} \sin\left( \frac{\pi s}{2} \right) \Gamma(1-s),$$
with $\Gamma(s)$ denoting the Gamma function.
Building on Riemann's unpublished notes, Siegel in his classical paper \cite{S1932} presented the error terms of the rhs of equation (\ref{1.6}) to \textit{all} orders, only in the important case of $x=y=\sqrt{t/2\pi}$. In theorems \ref{zetath2} and \ref{ZETATH} we present analogous results for \textit{any} $x$ and $y$ valid to \textit{all} orders. Theorem \ref{zetath2} yields the following explicit result for the subleading term in the rhs of equation (\ref{1.6}) when $1 < \eta = 2\pi y < \sqrt{t}$. For every $\epsilon > 0$, there exists an $A > 0$ such that
\begin{align*}
\zeta(s) 
= &\; \sum_{n=1}^{[\frac{t}{\eta}]} \frac{1}{n^s} 
+ \chi(s) \sum_{n=1}^{[\frac{\eta}{2\pi}]} \frac{1}{n^{1-s}}
 - \frac{e^{-i\pi s}\Gamma(1-s)}{i\sqrt{2\pi t} }e^{-([\frac{t}{\eta}]+1)i\eta} e^{\frac{i \pi}{2}(s-1)} \eta^{s}
\frac{e^{\frac{i\pi}{4}} }{1 - e^{-i\eta}}  
	\\ \nonumber
& + 
e^{-i\pi s}\Gamma(1-s) e^{-\frac{\pi t}{2}}\eta^{\sigma-1}
\times \begin{cases}	
O\bigl(\frac{\eta}{t^{5/6}}\bigr), & 1 < \eta < A t^{\frac{1}{3}} < \infty, \\
  O\bigl(e^{- \frac{At}{\eta^2}} + \frac{\eta^3}{t^{3/2}}\bigr), & t^{\frac{1}{3}} < \eta < A\sqrt{t} < \infty, \end{cases}  
 	\\ \nonumber
& \hspace{4cm} \dist(\eta, 2\pi \Z) > \epsilon, \quad 0 \leq \sigma \leq 1, \quad t \to \infty, 
\end{align*}
where the error term is uniform for all $\eta, \sigma$  in the above ranges. In corollary \ref{zetacor2}, similar formulas are presented for the case $2\pi\sqrt{t} < \eta < \frac{2\pi}{\epsilon} t$.

The case of $\epsilon \sqrt{t} < \eta < t$ is analyzed in theorem \ref{ZETATH}, which provides a formula valid to all orders in terms of the function $\Phi(\tau, u)$ defined by
\begin{align}\label{Phidefintro}
\Phi(\tau, u) = \int_{0 \nwarrow 1} \frac{e^{\pi i \tau x^2 + 2\pi i u x}}{e^{\pi i x} - e^{-\pi i x}} dx, \qquad \tau < 0, \quad u \in \C,
\end{align}
where $0 \nwarrow 1$ denotes a straight line parallel to $e^{3\pi i/4}$ which crosses the real axis between $0$ and $1$. The formula of theorem \ref{ZETATH} yields the following result for the leading asymptotic terms: For every $\epsilon > 0$, 
\begin{align}\label{1.7}
\zeta(s) 
=& \sum_{n=1}^{[\frac{t}{\eta}]} \frac{1}{n^s} 
+ \chi(s) \sum_{n=1}^{[\frac{\eta}{2\pi}]} \frac{1}{n^{1-s}}
 	\\ \nonumber
& + e^{-i\pi s}\Gamma(1-s)\biggl\{ e^{\frac{i\pi(s-1)}{2}}\eta^{s-1}  e^{\frac{2t}{\eta}[\frac{\eta}{2\pi}]\pi i - it 
- \frac{it}{2\eta^2}(2[\frac{\eta}{2\pi}]\pi - \eta)^2} 
	\\ \nonumber
& \qquad \times \biggl[\Phi + \frac{\sigma -1}{i\eta} \biggl(
  \partial_2 \Phi  
  + 
\biggl(2\Bigl[\frac{\eta}{2\pi}\Bigr]\pi i - i\eta\biggr)   \Phi  \biggr) \biggr]
 + O\biggl( e^{-\frac{\pi t}{2}} \frac{\eta^\sigma}{t^{5/6}} \biggr)\biggr\},
	\\ \nonumber
& \hspace{3cm} 
 \epsilon \sqrt{t} < \eta < t , \quad 0 \leq \sigma \leq 1, \quad t \to \infty,
\end{align}
where the error term is uniform for $\eta, \sigma$ in the above ranges and $\Phi$ and $\partial_2\Phi$ are evaluated at the point 
\begin{align}\label{Phievaluationpoint}
\left(-\frac{2\pi t}{\eta^2}, \frac{2t}{\eta} - \frac{2\pi t}{\eta^2}\Bigl[\frac{\eta}{2\pi}\Bigr] - \Bigl[\frac{t}{\eta}\Bigr] - \frac{1}{2}\right).
\end{align}

\begin{remark}\upshape
1. In the particular case of $\eta = \sqrt{2\pi t}$, equation (\ref{1.7}) agrees with the analogous formula of Siegel \cite{S1932}.

2. Although Siegel suggested already in \cite{S1932} the possibility of deriving a formula such as (\ref{1.7}) based on the function (\ref{Phidefintro}), the authors have not been able to locate such a formula in the literature. In the case of $\sigma=1/2$, an alternative asymptotic representation for $\zeta(s)$ involving a sum which is smoothly rather than sharply truncated is presented in \cite{BK1992}. An introduction to the Riemann-Siegel formula can be found in chapter 7 of \cite{E2001}.

3. Equation (\ref{1.7}) is particularly useful in the case when $\eta = \sqrt{\frac{2\pi t}{b}}$ where $b> 0$ is a rational number. Indeed, in this case the function $\Phi$ and its derivatives evaluated at the point (\ref{Phievaluationpoint}) can be computed explicitly (see equation (\ref{Phiformula}) below). Therefore, in this case theorem \ref{ZETATH} yields an {\it explicit} asymptotic expansion of $\zeta(s)$ to all orders. 

4. Equation (\ref{zetaformula2.1intro}) is only useful as an asymptotic formula in the case when $t/\eta \to 0$ as $t \to \infty$ (otherwise the error term is as large as the retained terms). Therefore the asymptotic range where $t < \eta$ and $t/\eta = O(1)$ is not covered by this equation. However, this asymptotic sector is covered by corollary \ref{zetacor2}. 
\end{remark}

\section{The explicit form of certain sums}
Theorems \ref{th3.1} and \ref{th3.2} allow us to evaluate the sum $\sum_{a}^{b} n^{-s}$, for certain $a$ and $b$, to \textit{all} orders, see theorems \ref{th5.1} and \ref{th5.2}. For example, the leading order of such sums follows immediately from equations (\ref{zetaformula2.1intro}) and (\ref{zetaformula2.2intro}): For every $\epsilon > 0$,
\begin{align}\label{1.8}
& \sum_{n=\left[ \frac{\eta_1}{2\pi} \right]+1}^{\left[ \frac{\eta_2}{2\pi} \right]}n^{-s}=\frac{1}{1-s} \left[ \left(\frac{\eta_2}{2\pi} \right)^{1-s}- \left(\frac{\eta_1}{2\pi} \right)^{1-s} \right]+O(\eta_1^{-\sigma}), 
	\\ \nonumber
& \hspace{3cm} (1+\epsilon)t<\eta_1<\eta_2<\infty, \quad 0 \leq \sigma \leq 1, \quad t \to \infty,
	\\ \label{1.9}
& \sum_{n=\left[ \frac{t}{2\pi} \right]+1}^{\left[ \frac{\eta}{2\pi} \right]}n^{-s}
= \frac{1}{1-s} \left(\frac{\eta}{2\pi} \right)^{1-s} 
+ \frac{1}{2}e^{it}e^{\frac{\pi i}{4}} \bigg(\frac{t}{2\pi}\bigg)^{\frac{1}{2} -s}
+O(t^{-\sigma}), 
	\\ \nonumber
&\hspace{3cm} (1+\epsilon)t<\eta<\infty, \quad 0 \leq \sigma \leq 1, \quad t \to  \infty,
\end{align}
uniformly with respect to $\eta_1, \eta_2, \eta$, and $\sigma$.

\section{The explicit form of the difference of certain sums}
It is well known, see for example page 78 of Titchmarsh \cite{T1986}, that for large $t$, the following sums coincide to the leading order:
$$\sum_{x < n \leq N} n^{-s} \qquad \text{and} \qquad \chi(s)\!\!\!\!\!\! \sum_{\frac{t}{2\pi N} < n \leq \frac{t}{2 \pi x}} n^{s-1}.$$
Theorems \ref{zetath2} and \ref{ZETATH} allow us to evaluate the difference of  similar sums to all orders, see theorem \ref{th5.3}. In particular, equation (\ref{1.7}) immediately  implies the following equation:
\begin{align}\label{1.10}
\sum_{n=[\frac{t}{\eta_2}]+1}^{[\frac{t}{\eta_1}]} \frac{1}{n^s}  = & \; \chi(s) \sum_{n=[\frac{\eta_1}{2\pi}]+1}^{[\frac{\eta_2}{2\pi}]} \frac{1}{n^{1-s}}
 	\\ \nonumber
& + e^{-i\pi s}\Gamma(1-s)\biggl\{ e^{\frac{i\pi(s-1)}{2}}\eta^{s-1}  e^{\frac{2t}{\eta}[\frac{\eta}{2\pi}]\pi i - it 
- \frac{it}{2\eta^2}(2[\frac{\eta}{2\pi}]\pi - \eta)^2} 
	\\ \nonumber
&  \times \biggl[\Phi + \frac{\sigma -1}{i\eta} \biggl(
  \partial_2 \Phi  
  + 
\biggl(2\Bigl[\frac{\eta}{2\pi}\Bigr]\pi i - i\eta\biggr)   \Phi  \biggr) \biggr]
 + O\biggl( e^{-\frac{\pi t}{2}} \frac{\eta^\sigma}{t^{5/6}} \biggr)\biggr\}\Biggr|_{\eta = \eta_1}^{\eta_2},
	\\ \nonumber
& \hspace{3cm} 
 \epsilon \sqrt{t} < \eta_1 < \eta_2 < t , \quad 0 \leq \sigma \leq 1, \quad t \to \infty.
\end{align}
These results are used in \cite{KFpreprint} in order to obtain two asymptotic identities, which are then used for the derivation of estimates of certain exponential sums. In particular, the first of these identities is derived from (\ref{1.8}) and has applications in the asymptotic estimation of several single and double exponential sums analysed in \cite{KFpreprint}; furthermore, a variant of (\ref{1.10}) is used for the derivation
of the second important identity employed in \cite{KFpreprint}. These new identities are used
either by themselves or in combination with classical techniques \cite{T1986} in order to obtain estimates which are either derived in a simpler way, or are sharper than the estimates obtained using only the classical techniques.

\section[Asymptotics of a two-parameter generalization]{Asymptotics of a two-parameter generalization of Riemann's zeta function}
\begin{figure}
\begin{center}
\bigskip
 \begin{overpic}[width=.6\textwidth]{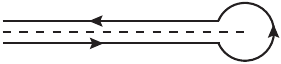}
 \put(86,6.5){$0$}
   \end{overpic}
   \caption{The Hankel contour $H_\alpha$.} \label{XHankel.pdf}
   \end{center}
 \end{figure}

In chapters \ref{sec6}-\ref{sec8}, we study the following two-parameter generalization of Riemann's zeta function:
\begin{align}\label{X1.1}
  \Phi(u,v,\beta) = \int_{H_\alpha} z^{u-1} (z-2i\pi \beta)^{v-1} \frac{dz}{e^{-z} - 1}, \qquad u \in \C, \quad v \in \C, \quad \beta \in \R \setminus \{0\},
\end{align}
where $H_\alpha$  denotes the usual Hankel contour appearing in the definition of Riemann's zeta function, namely the union of the following curves in the complex $z$-plane cut along the negative axis with the orientation shown in figure \ref{XHankel.pdf}:
\begin{align}\label{X1.2}
& H_\alpha = \left\{re^{-i\pi} \, | \, \alpha < r <\infty\right\} \cup \left\{\alpha e^{i\theta}\, | \, -\pi <\theta <\pi\right\}  \cup \left\{re^{i\pi}\, | \, \alpha< r <\infty\right\},
\end{align}
and we require that $0 < \alpha < 2\pi\min(|\beta|, 1)$ so that the contour surrounds neither of the points $\pm2i\pi$ or $\pm2i\pi \beta$. It follows immediately from the definition (\ref{X1.1}) that $\Phi(u,v,\beta)$ is an entire function of both $u$ and $v$ for each $\beta \in \R \setminus \{0\}$.

Recall that the Riemann function $\zeta(s)$ can be defined by 
\begin{align}\label{2.4}
\zeta(s) = \frac{\Gamma(1-s)}{2i\pi} \int_{H_\alpha} \frac{z^{s-1}}{e^{-z}-1} dz, \qquad s \in \C \setminus \{1\},
\end{align}
where $0 < \alpha < 2\pi$.
Thus, in the particular case that $v = 1$, the function  $\Phi$ is proportional to the Riemann zeta function $\zeta(u)$:
\begin{align}\label{X1.3}
  \Phi(u,1, \beta) = \frac{2i\pi}{\Gamma(1-u)} \zeta(u), \qquad u \in \C.
\end{align}

The analog of the identity (\ref{1.2}) for $\Phi$ is derived in chapter \ref{sec6}. The large $t$ asymptotics of $\Phi$, where $u$ and $v$ are given by
$$u = \sigma_1 + it, \qquad v = \sigma_2 - it,$$
is computed in chapters \ref{sec7} and \ref{sec8}.

\section{Fourier coefficients of the product of two Hurwitz zeta functions}
One reason for studying the function $\Phi(u,v,\beta)$ is that it embeds the Riemann zeta function into a larger family of functions, in the same way that the Hurwitz and Lerch zeta functions embed $\zeta(s)$ into a larger family depending on one and two additional parameters, respectively.
Another reason (which was also our original motivation for considering $\Phi$) is that the Fourier coefficients of the absolute value squared of the Hurwitz zeta function (or more generally of the product of two Hurwitz zeta functions) can be expressed in terms of $\Phi$. In what follows, we explain this connection; see chapter \ref{hurwitzsec} for details.

Recall that the Hurwitz zeta function $\zeta(s,\alpha)$ is the analytic continuation of the sum
$$\zeta(s,\alpha) = \sum_{n=0}^\infty \frac{1}{(n + \alpha)^s}, \qquad \re s > 1, \quad \re \alpha > 0,$$
and that the modified Hurwitz zeta function $\zeta_1(s,\alpha)$ is defined by
$$\zeta_1(s, \alpha) = \zeta(s, \alpha) - \alpha^{-s} = \zeta(s, \alpha +1).$$
For each choice of $u,v \in \C \setminus\{1\}$, the product $\zeta_1(u,\alpha)\zeta_1(v,\alpha)$ is a smooth function of $\alpha \in [0,1]$; hence it can be represented by its Fourier series
\begin{align}\label{introzeta1zeta1Fourierseries}
\zeta_1(u,\alpha)\zeta_1(v,\alpha) = \sum_{n \in \Z} q_n(u,v) e^{2\pi i n \alpha},
\end{align}
where the Fourier coefficients $q_n(u,v)$ are defined by
$$q_n(u,v) = \int_0^1 \zeta_1(u,\alpha) \zeta_1(v,\alpha) e^{-2\pi i n \alpha} d\alpha, \qquad n \in \Z.$$

For $u = \bar{v} = \sigma + it$, the zeroth Fourier coefficient $q_0(u,v)$ is given by the mean square average of the modified Hurwitz function:
$$q_0(\sigma + it,\sigma - it) = \int_0^1 |\zeta_1(\sigma + it,\alpha)|^2 d\alpha.$$
Attempts to compute the large $t$ asymptotics of $q_0(\sigma + it,\sigma - it)$ have a long history, especially in the case when $\sigma = 1/2$. Increasingly refined asymptotic estimates were derived in an extensive series of papers (see \cite{KL1952, B1979, R1983, Z1990, Z1991, KM1993}) before the following formula was finally obtained independently in \cite{A1992} and \cite{Z1994}:
\begin{align}\label{introint01hurwitz}
\int_0^1 \bigg|\zeta_1\bigg(\frac{1}{2} + it, \alpha\bigg)\bigg|^2 d\alpha = \ln \frac{t}{2\pi}  + \gamma - 2 \re \frac{\zeta(\frac{1}{2} + it)}{\frac{1}{2} + it} + O(t^{-1}),
\end{align}
where $\gamma$ denotes the Euler constant. 
An alternative derivation of (\ref{introint01hurwitz}) was presented in \cite{KM1996}. 
In the derivation of \cite{KM1996}, the asymptotic formula (\ref{introint01hurwitz}) is obtained as an easy consequence of the following identity \cite[Eq. (2.1) with $N = 1$]{KM1996}:
\begin{align}\label{q0KMidentity}
q_0(u,v) = &\; \frac{1}{u+v-1} + R_0(u,v) + R_0(v,u) - T_0(u,v) - T_0(v,u), 
\end{align}
where $0 < \re u, \re v < 2$ and $R_0(u,v)$ and $T_0(u,v)$ are defined by
\begin{align}\label{introR0def}
& R_0(u,v) = \Gamma(u+v-1)\zeta(u+v-1) \frac{\Gamma(1-v)}{\Gamma(u)},
	\\\label{introT0def}
& T_0(u,v) = \frac{\zeta(u) - 1}{1 - v} 
+ \frac{u}{1-v} \sum_{l=1}^\infty l^{1-u-v} \int_l^\infty \beta^{u+v-2}(1+\beta)^{-u-1} d\beta.
\end{align}
Actually, a version of the identity (\ref{q0KMidentity}) is presented as the main result of \cite{KM1996}, because from this identity many asymptotic results (such as (\ref{introint01hurwitz})) can easily be obtained. A key point here is that, for $u = \sigma + it$ and $v = \sigma - it$, the sum appearing in the definition (\ref{introT0def}) of $T_0(u,v)$ is small as $t \to \infty$.  

In chapter \ref{hurwitzsec}, we will analyze the large $t$ asymptotics of the $n$th Fourier coefficient $q_n$ for any $n \in \Z$ by establishing a generalization of (\ref{q0KMidentity}).
We will first show that $R_0$ and $T_0$ are the analytic continuations of the integrals
\begin{subequations}\label{introR0T0}
\begin{align}
& R_0(u,v) = \int_0^\infty \alpha^{-v} \zeta_1(u, \alpha) d\alpha, \qquad \re (u+v) > 2, \quad \re v < 1,
	\\ 
& T_0(u,v) = \int_0^1 \alpha^{-v}  \zeta_1(u, \alpha)  d\alpha, \qquad  \re v < 1.
\end{align}
\end{subequations}
Then, we will prove (see Theorem \ref{QNTH}) the following generalization of (\ref{q0KMidentity}) which is valid for $n \in \Z$ and $\re u, \re v < 2$: 
\begin{align}\label{introqnidentity}
q_n(u,v) = &\;  b_n(u+v) + R_n(u,v) + R_n(v,u) - T_n(u,v) - T_n(v,u), 
\end{align}
where, for any $n \in \Z$, $b_n(s)$, $R_n(u,v)$, and $T_n(u,v)$ are the analytic continuations of the following functions:
\begin{align*}
& b_n(s) =  \int_1^\infty \alpha^{-s}  e^{-2\pi i n \alpha}  d\alpha, \qquad \re s > 1,
	\\
& R_n(u,v) = \int_0^\infty \alpha^{-v} \zeta_1(u,\alpha) e^{-2i\pi n \alpha} d\alpha, \qquad \re (u+v) > 2, \quad \re v < 1,
	\\
& T_n(u,v) = \int_0^1 \alpha^{-v}  \zeta_1(u, \alpha) e^{-2\pi i n\alpha} d\alpha, \qquad \re v < 1.
\end{align*}
Since $b_0(s) = \frac{1}{s-1}$, it follows immediately from (\ref{introR0T0}) that equation (\ref{introqnidentity}) reduces to (\ref{q0KMidentity}) when $n = 0$.

Just like the identity (\ref{q0KMidentity}) can be used to find the large $t$ asymptotics of the zeroth Fourier coefficient $q_0$, the identity (\ref{introqnidentity}) can be used to find the large $t$ behavior of the $n$th Fourier coefficient $q_n$.
The function  $\Phi(u,v,\beta)$ enters in this formulation, because it turns out (see equation (\ref{Rncontinuation})) that for $n \neq 0$ the function  $R_n(u,v)$ is given by the following analog of (\ref{introR0def}):
\begin{align*}
R_n(u,v) = \frac{\Gamma(1-v)}{\Gamma(u)} 
   \Phi(u,v, n ) \times \begin{cases} 
  \frac{ie^{\pi i v}}{2\sin(\pi u)}, \quad &  n \geq 1, \\
 \frac{ie^{-\pi i v}}{2\sin(\pi u)}, & n \leq -1.
  \end{cases}
\end{align*}
We will show in chapter \ref{hurwitzsec} that the asymptotics of $T_n$ can be obtained (as in the case $n = 0$) for any $n \neq 0$. 
Thus, the computation of the leading asymptotics of the Fourier coefficients $q_n$, $n \neq 0$, reduces to the calculation of the large $t$ asymptotics of $\Phi(u,v,n)$ for integer values of $n$. 
By substituting the asymptotic results for $\Phi$ obtained in chapter \ref{sec7} into (\ref{introqnidentity}), we obtain in chapter \ref{hurwitzsec} an asymptotic formula for the $n$th Fourier coefficient of $|\zeta_1(\sigma + it,\alpha)|^2$.
These considerations also lead to Lindel\"of type asymptotic bounds on the function $\Phi(u,v,n)$ for large $t$ (see Corollary \ref{Phiboundcor}):
\begin{align*}
\Phi(\sigma + it, \sigma - it,n) = \begin{cases} O(1), & \sigma \in [0, 1/2), \\
O(\ln t), & \sigma = 1/2, \\
O(t^{2\sigma - 1}), & \sigma \in (1/2, 1],
\end{cases} \quad n \geq 1, \quad t \to \infty.
\end{align*}

\section{Several representations for the basic sum}
The remarkable fact about the large $t$-asymptotics of the Riemann zeta function is that, whereas the higher order terms in the asymptotic expansion can be computed {\it explicitly}, the sum appearing in the leading order is {\it transcendental}. In chapter \ref{sec9} we present several integral representations of this basic sum with the hope that these representations may be useful for the estimation of this fundamental sum (some results in this direction can be found in \cite{AFpreprint}).

Let $C_\eta^t$ denote the semicircle from $i\eta$ to $it$ with $\re z \geq 0$:
\begin{align}\label{Cetatdef}
  C_\eta^t = \left\{ \frac{i(\eta + t)}{2} + \frac{t - \eta}{2}e^{i\theta} \; \middle| \; -\frac{\pi}{2} < \theta < \frac{\pi}{2} \right\}.
\end{align}
The following representations for the basic sum $ \sum n^{s-1} $ are derived in chapter \ref{sec9}:
\begin{align*}
 \sum_{n=[\frac{\eta}{2\pi}]+1}^{[\frac{t}{2\pi}]} n^{s-1} 
&  = \frac{e^{-\frac{i\pi s}{2}}}{(2\pi)^s} \int_{C_\eta^t} \frac{z^{s-1}}{e^{z} - 1} dz + O(t^{\sigma -1})
  	\\
&  = \frac{2}{(2\pi)^s} \dashint_{\eta}^t \frac{\rho^{s-1}}{e^{i\rho} - 1} d\rho + O(t^{\sigma -1})
	\\
& = \frac{1}{(2\pi)^s} \lim_{\substack{\epsilon \to 0 \\ \epsilon t \to 0}} \int_{\eta}^t \frac{u^{s-1}}{e^{iue^{-i\epsilon}} - 1} du + O(t^{\sigma -1}), \qquad 0 < \eta < t, \quad t \to \infty,
\end{align*}
where it is assumed that $\dist(\eta, 2\pi \Z) > \delta$, $\dist(t, 2\pi \Z) > \delta$, and the contour in the second equality above denotes the principal value integral with respect to the points 
$$\left\{2\pi n \; \middle| \; n \in \Z, \; \Bigl[\frac{\eta}{2\pi}\Bigr] + 1\leq n \leq \Bigl[\frac{t}{2\pi}\Bigr]\right\}.$$

\smallskip\bigskip\noindent
{\bf Basic Assumption:} In the statements of all results throughout this work, it will be assumed (also when not explicitly mentioned) that 
\begin{align}\label{basicassumption}
\boxed{\eta > 0, \quad t>0, \quad \eta \notin 2 \pi \Z, \quad t \notin 2\pi \Z.}
\end{align}

\chapter{An Exact Representation for $\zeta(s)$}\label{sec2}

In chapters \ref{sec2} and \ref{sec3}, we assume that the branch cut for the logarithm runs along the negative real axis. By the standing assumption (\ref{basicassumption}), we always have $\eta, t \notin 2 \pi \Z$.

\begin{theorem}[An exact representation for Riemann's zeta function]\label{theorem2.1}
Let $\zeta(s)$, $s=\sigma+ i t$, $\sigma, t \in \R$, denote the Riemann zeta function. Then,
\begin{align}\nonumber
\zeta(s) = &\; \chi(s) \Bigg\{ \sum^{\left[ \frac{\eta}{2 \pi} \right]}_{n=1} n^{s-1} 
	\\\nonumber
&\hspace{1cm}  + \frac{1}{(2\pi)^s} \left[ - \frac{\eta^s}{s} + \left( e^{ \frac{i\pi s}{2}} \int^{\infty e^{i\phi_1}}_{\eta e^{-\frac{i\pi}{2}}} + e^{- \frac{i\pi s}{2}} \int^{\infty e^{i\phi_2}}_{\eta e^{\frac{i\pi}{2}}}\right) \frac{e^{-z}z^{s-1}dz}{1-e^{-z}}  \right] \Bigg\},
	\\\label{2.1}
&\hspace{5cm} 0<\eta<\infty, \ - \frac{\pi}{2} < \phi_j < \frac{\pi}{2}, \ \ j = 1,2, 
\end{align}
where
\begin{align}\label{2.2}
  \chi(s) = \frac{(2\pi)^s}{\pi} \sin\left( \frac{\pi s}{2}\right) \Gamma(1-s),
\end{align}
with $\Gamma(s)$, $s\in \C$, denoting the Gamma function, and where the contours of the first and second integrals in the rhs of (\ref{2.1}) are the rays from $\eta \exp(-i\pi/2) $ to $\infty \exp(i\phi_1)$ and from $\eta\exp(i\pi/2)$ to $\infty \exp(i\phi_2)$ respectively.
An equivalent formula is
\begin{align}\nonumber
&\zeta(1-s) =  \sum^{\left[ \frac{\eta}{2 \pi} \right]}_{n=1} n^{s-1}+\frac{1}{(2\pi)^s} \left[- \frac{\eta^s}{s} +
\left( e^{\frac{i\pi s}{2}} \int^{\infty e^{i\phi_1}}_{\eta e^{- \frac{i\pi}{2}}} + e^{-\frac{i\pi s}{2}} \int^{\infty e^{i\phi_2}}_{\eta e^{\frac{i\pi}2}} \right) \frac{e^{-z}z^{s-1}dz}{1-e^{-z}}\right],
	\\ \label{2.3}
&\hspace{5cm} 0<\eta<\infty, \ - \frac{\pi}{2} < \phi_j < \frac{\pi}{2}, \ \ j = 1,2.	
\end{align}
\end{theorem}
\begin{proof}
We will first derive equation (\ref{2.1}) in the particular case $0<\eta<2\pi$. In this case, replacing $\eta$ with $\alpha$, (\ref{2.1}) can be written as
\begin{align}\nonumber
& \zeta(s) = \frac{\chi(s)}{(2\pi)^s} \left \{- \frac{\alpha^s}{s} +\left( e^{\frac{i\pi s}{2}}
\int^{\infty e^{i\phi_1}}_{\alpha e^{- \frac{i\pi}{2}}} + e^{- \frac{i\pi s}{2}} \int^{\infty e^{i\phi_2}}_{\alpha e^{\frac{i\pi}{2}}}\right)\frac{e^{-z}z^{s-1}}{1-e^{-z}} dz\right \}, 
	\\ \label{2.6}
& \hspace{4cm} 0<\alpha <2\pi; \quad \frac{-\pi}{2} < \phi_j < \frac{\pi}{2}, \quad j = 1,2.
\end{align}

We decompose the Hankel contour $H_\alpha$ in the expression (\ref{2.4}) for $\zeta(s)$ into the union of  three different contours (see figure \ref{L1L2L3.pdf}), namely
\begin{figure}
\begin{center}
 \begin{overpic}[width=.7\textwidth]{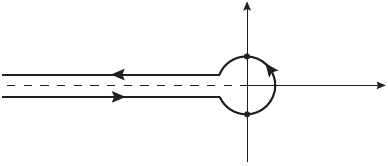}
 \put(28.5,26.5){$L_1$}
 \put(28.5,12.5){$L_2$}
 \put(71,26){$L_3$}
 \put(101,19.5){$\re z$}
  \put(59.3,29.7){$i\alpha$}
 \put(56.2,9.5){$-i\alpha$}
   \end{overpic}
   \caption{The decomposition of $H_\alpha$ into $L_1+L_2+L_3$.}\label{L1L2L3.pdf}
   \end{center}
\end{figure}
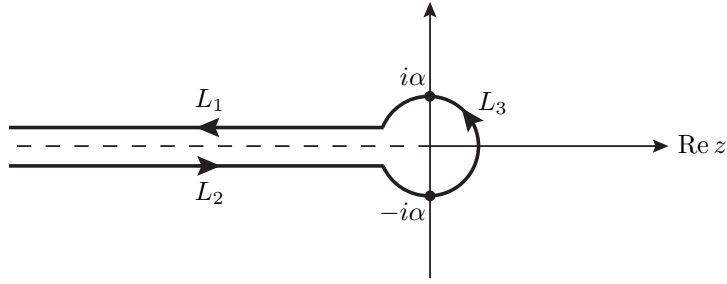
$$H_\alpha = L_1 \cup L_2 \cup L_3,$$
where
\begin{align}\nonumber
& L_1 = \left\{ \alpha e^{i\theta} \, \Big| \, \frac{\pi}{2} < \theta < \pi \right\} \cup \{ |z|e^{i\pi} \, | \, \alpha  < |z| < \infty)\},
	\\ \nonumber
& L_2 = \{ |z|e^{-i\pi} \, | \, \alpha  <|z| <\infty\} \cup \left\{\alpha e^{i\theta} \, \Big| \, -\pi <\theta < - \frac{\pi}{2} \right\},
	\\ \label{2.7}
& L_3 = \left\{ \alpha e^{i\theta} \, \Big| \, -\frac{\pi}{2} < \theta < \frac{\pi}{2} \right\}.
\end{align}
The integral along $L_1$ can be written as follows:
\begin{align}\label{2.8}
\int_{L_1} \frac{z^{s-1}}{e^{-z}-1} dz = \int^{\infty e^{i\pi}}_{\alpha e^{\frac{i\pi}{2}}} \frac{e^z}{1-e^z} z^{s-1} dz = e^{i\pi s} \int^\infty_{\alpha e^{- \frac{i\pi}{2}}} \frac{e^{-u}u^{s-1}du}{1-e^{-u}},
\end{align}
where the contours of the second and third integrals in (\ref{2.8}) are the rays in the cut complex $z$-plane from $\alpha  \exp(i\pi/2)$ to $\infty\exp(i\pi)$ and from $\alpha \exp(-i\pi/2)$ to $\infty$ respectively.

Indeed, in the domain enclosed by $L_1$ and by the ray from $\alpha \exp(i\pi/2)$ to $\infty \exp(i\pi)$, i.e., in the shaded domain 1 of figure \ref{domains123.pdf}, we have
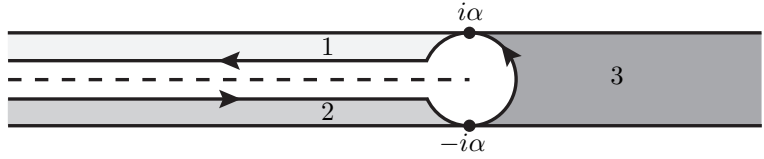
\begin{figure}
\begin{center}
 \begin{overpic}[width=.8\textwidth]{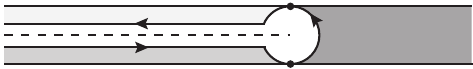}
 \put(42,10.4){$1$}
 \put(42,1.8){$2$}
 \put(80,6.5){$3$}
 \put(59.9,14.9){$i\alpha$}
 \put(57.5,-2.3){$-i\alpha$}
   \end{overpic}
   \caption{The domains $1$, $2$, and $3$.}\label{domains123.pdf}
   \end{center}
\end{figure}
$$ z = re^{i\theta}, \quad \alpha \leq r \leq \infty, \quad \frac{\pi}{2} \leq \theta \leq \pi;$$
thus,
$$ |e^z| = e^{r\cos \theta}, \quad \alpha \leq r \leq \infty, \quad 0 \leq \cos \theta \leq -1.$$
Hence, Cauchy's theorem implies the first equality of (\ref{2.8}).  Then, the substitution $u=z\exp(-i\pi)$, yields the second equality in (\ref{2.6}).
Similarly,
\begin{align}\label{2.9}
\int_{L_2} \frac{z^{s-1}}{e^{-z}-1} dz = \int^{\alpha e^{- \frac{i\pi}{2}}}_{\infty e^{-i\pi}} \frac{e^z}{1-e^z} z^{s-1} dz = e^{-i\pi s} \int^{\alpha e^{\frac{i\pi}{2}}}_\infty \frac{e^{-u}u^{s-1}}{1-e^{-u}} du,
\end{align}
where the curve $L_2$ is defined in (\ref{2.7}) and the contours of the second and third integrals in (\ref{2.9}) are the rays in the cut complex $z$-plane from $\infty\exp(-i\pi)$ to $\alpha \exp(-i\pi/2)$ and from $\infty$ to $\alpha \exp(i\pi/2)$ respectively.

Indeed, the first equality in (\ref{2.9}) is the consequence of Cauchy's theorem applied in the domain enclosed by $L_2$ and by the ray from $\infty \exp(-i\pi)$ to $ \alpha \exp(-i\pi/2)$ (i.e., in the shaded domain 2 of figure \ref{domains123.pdf}), whereas the second equality in (\ref{2.9}) follows from the substitution $u= z\exp(i\pi)$.

The integral along $L_3$ can be written as follows:
\begin{align}\label{2.10}
 \int_{L_3} \frac{z^{s-1}}{e^{-z}-1}dz = - \left( \int^\infty_{\alpha e^{- \frac{i\pi}{2}}} + \int^{\alpha e^{\frac{i\pi}{2}}}_\infty\right) \frac{e^{-z}}{1-e^{-z}} z^{s-1}dz - 2i\sin \left( \frac{\pi s}{2}\right) \frac{\alpha^s}{s},
\end{align}
where the curve $L_3$ is defined in (\ref{2.7}) and the contours in the first and second integrals in (\ref{2.10}) are the rays in the cut complex $z$-plane from $\alpha \exp(-i\pi/2)$ to $\infty$ and from $\infty$ to $\alpha \exp (i\pi/2)$.

Indeed,
\begin{align}\label{2.11}
\int_{L_3} \frac{z^{s-1}}{e^{-z}-1} dz = - \int_{L_3} \frac{e^{-z}}{1-e^{-z}} z^{s-1}dz - \int_{L_3} z^{s-1}dz.
\end{align}
The second term in the rhs of (\ref{2.11}) yields the second term in the rhs of (\ref{2.10}).  Furthermore, in the domain enclosed by $L_3$ and by the two rays $(\alpha \exp(-i\pi/2),\infty)$ and $(\infty, \alpha\exp (i\pi/2))$, i.e., in the shaded domain 3 of figure \ref{domains123.pdf}, we have
$$ z = re^{i\theta}, \quad \alpha \leq r \leq \infty, \quad - \frac{\pi}{2} <\theta < \frac{\pi}{2};$$
thus
$$ |e^{-z}| = e^{-r\cos \theta}, \quad \alpha \leq r \leq \infty, \quad \cos \theta \geq 0.$$
Hence, Cauchy's theorem implies that the first term of the rhs of (\ref{2.11}) equals the first term in the rhs of (\ref{2.10}).

Adding equations (\ref{2.8})-(\ref{2.10}) we find equation (\ref{2.6}) but with $\phi_1=\phi_2=0$. However, Cauchy's theorem implies that the rays $(\alpha \exp(-i\pi/2),\infty)$ and $(\infty, \alpha \exp(i\pi/2))$ can be replaced with the rays 
$$(\alpha\exp(-i\pi/2), \infty\exp(i\phi_1)) \quad \text{and}\quad (\infty\exp(i\phi_2), \alpha\exp(i\pi/2)),$$ 
respectively, and then we find (\ref{2.6}).

In order to derive equation (\ref{2.1}) with $2\pi<\eta<\infty$, we introduce the curves $C^\eta_\alpha$, $C^{-\alpha}_{-\eta}$ and $\hat C^\alpha_\eta$, which are the following semi-circles in the complex $z$-plane, see figure \ref{Calphaa.pdf}:
\begin{figure}
\begin{center}
 \begin{overpic}[width=.5\textwidth]{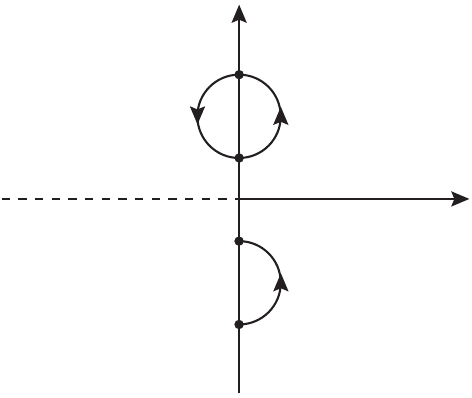}
 \put(44,70){$i\eta$}
 \put(43,47){$i\alpha$}
 \put(38.5,32){$-i\alpha$}
 \put(39,14){$-i\eta$}
 \put(62.5,58){$C_\alpha^\eta$}
 \put(31.5,58.5){$\hat{C}_\eta^\alpha$}
 \put(62.5,22.5){$C_{-\eta}^{-\alpha}$}
 \put(101.5,40.5){$\re z$} 
   \end{overpic}
   \caption{The contours $C^\eta_\alpha$, $C^{-\alpha}_{-\eta}$ and $\hat C^\alpha_\eta$.} \label{Calphaa.pdf}
   \end{center}
\end{figure}
\begin{subequations}\label{Cetaalphadef}
\begin{align}\label{2.12a}
 C^\eta_\alpha & = \left\{ \frac{i(\eta+\alpha)}{2} + \frac{(\eta-\alpha)}{2}e^{i\theta} \, \middle| \, - \frac{\pi}{2} <\theta < \frac{\pi}{2}\right\}, 
	\\ \label{2.12b}
 C^{-\alpha}_{-\eta} &  = \left\{ \frac{-i(\eta+\alpha)}{2} + \frac{(\eta-\alpha)}{2} e^{i\theta}\, \middle| \,  - \frac{\pi}{2} < \theta < \frac{\pi}{2} \right\},
	\\ \label{2.12c}
 \hat C^\alpha_\eta &  = \left\{ \frac{i(\eta+\alpha)}{2} + \frac{(\eta-\alpha)}{2} e^{i\theta}\, \middle| \,  \left(-\pi <\theta < \frac{\pi}{2}\right) \cup \left(\frac{\pi}{2} < \theta < \pi\right) \right\}.
\end{align}
\end{subequations}
The first integral in the rhs of equation (\ref{2.6}) can be rewritten as an integral along the curve $-C^{-\alpha}_{-\eta}$ plus an integral along the ray $(\eta \exp(i\pi/2),\infty \exp(i \phi_2) )$. Similarly, the second integral in the rhs of equation (\ref{2.6}) can be rewritten as an integral along the curve $C^\eta_\alpha$ plus an integral along the ray $(\eta \exp(i\pi/2),\infty \exp(i \phi_1) )$. Hence, the sum of the two integrals in the rhs of equation (\ref{2.6}) yields the sum of the two integrals in equation (\ref{2.1}) plus the sum of the following two integrals:
\begin{align}\label{2.13}
I_1=-e^{\frac{i \pi s}{2}} \int_{C^{-\alpha}_{-\eta}}\frac{e^{-z}}{1-e^{-z}}z^{s-1}dz, \quad  I_2=e^{-\frac{i \pi s}{2}}\int_{C_{\alpha}^{\eta}}\frac{e^{-z}}{1-e^{-z}}z^{s-1}dz. 
\end{align} 
Making the change of variables $u = z \exp(i\pi)$ in the integral $I_1$, we find
\begin{align} \nonumber
I_1 & =e^{-\frac{i \pi s}{2}} \int_{\hat C^\alpha_\eta} \frac{1}{1-e^{-u}} u^{s-1}du
 = e^{-\frac{i \pi s}{2}}\int_{\hat C^\alpha_\eta} \frac{e^{-u}}{1-e^{-u}}u^{s-1}du +  e^{-\frac{i \pi s}{2}}\int_{\hat C^\alpha_\eta} u^{s-1}du 
	\\ \label{2.14}
& = e^{-\frac{i\pi s}{2}} \int_{\hat C^{\alpha}_{\eta}} \frac{e^{-u}}{1-e^{-u}} u^{s-1}du + \frac{1}{s} (\alpha^s-\eta^s).
\end{align}
The sum of $I_2$ and the first term in the rhs of equation (\ref{2.14}) can be evaluated using Cauchy's theorem:
$$I_1+I_2=(2\pi)^s \sum^{\left[ \frac{\eta}{2 \pi}\right]}_{m=1}  m^{s-1}+ \frac{1}{s} (\alpha^s-\eta^s).$$
Hence equation (\ref{2.1}) follows.

Multiplying equation (\ref{2.1}) by $\chi(1-s)$ and using the identity
\begin{align}\label{2.15}
  \chi(s)\chi(1-s) = 1, \quad s \in {\C},
\end{align}  
together with the functional equation
\begin{align}\label{2.16}
\zeta(1-s) = \chi(1-s)\zeta(s), \quad s\in {\C},
\end{align}
we find the equation (\ref{2.3}). 
\end{proof}

\chapter{The Asymptotics of the Riemann Zeta Function for $t \leq \eta < \infty$}\label{sec3}
We are mainly interested in the behavior of $\zeta(s)$ in the critical strip $0 \leq \re s \leq 1$, thus we will henceforth assume that $\sigma \in [0,1]$. It is possible to extend the results of this chapter to the case when $\sigma$ belongs to any bounded interval, but the proofs are more complicated (e.g. equation (\ref{DNterm}) involves powers of $\sigma$ that cannot be neglected if $\sigma > 1$).

\begin{theorem}[{\bf The asymptotics to all orders for the case $(1+\epsilon)t < \eta$}]\label{th3.1}
Let $\zeta(s)$, $s = \sigma + it$, $\sigma, t \in \R$, denote the Riemann zeta function. Then,
\begin{align}\label{zetaformula2.1}
\zeta(s) = &\; \sum_{n =1}^{[\frac{\eta}{2\pi}]} n^{-s} - \frac{1}{1-s} \left(\frac{\eta}{2\pi}\right)^{1-s} 
	\\ \nonumber
&+ \frac{e^{-\frac{i\pi (1-s)}{2}}}{(2\pi)^{1-s}}\sum_{n=1}^\infty  \sum_{j=0}^{N-1} e^{-nz - it\ln{z}} \left(\frac{1}{n + \frac{it}{z}}\frac{d}{dz}\right)^j\frac{z^{-\sigma}}{n + \frac{it}{z}} \Biggr|_{z = i\eta}
	\\\nonumber
& + \frac{e^{\frac{i\pi (1-s)}{2}}}{(2\pi)^{1-s}}\sum_{n=1}^\infty
\sum_{j=0}^{N-1} e^{-nz - it\ln{z}} \left(\frac{1}{n + \frac{it}{z}}\frac{d}{dz}\right)^j\frac{z^{-\sigma}}{n + \frac{it}{z}} \Biggr|_{z = -i\eta}
	\\ \nonumber
& + O\biggl((2N +1)!! N \Bigl(\frac{1+\epsilon}{\epsilon}\Bigr)^{2(N+1)} \eta^{-\sigma - N}  \biggr),
	\\ \nonumber
& \hspace{2cm} (1+\epsilon)t < \eta < \infty, \quad \epsilon > 0, \quad 0 \leq \sigma \leq 1, \quad N \geq 2, \quad t \to \infty,	
\end{align}
where the error term is uniform for all $\eta, \epsilon, \sigma, N$ in the above ranges and $(2N +1)!!$ is defined by
$$(2N +1)!! = 1\cdot 3 \cdot 5 \cdot \cdots \cdot (2N-1)(2N+1).$$

For $N=3$ equation (\ref{zetaformula2.1}) simplifies to
\begin{align}\label{zetaformula2.1b}
\zeta(s) = &\; \sum_{n =1}^{[\frac{\eta}{2\pi}]} n^{-s} - \frac{1}{1-s} \left(\frac{\eta}{2\pi}\right)^{1-s} 
	\\ \nonumber
& + \frac{2i\eta^{-s}}{(2\pi)^{1-s}} \biggl\{
-i\arg(1 - e^{i \eta})
+ \frac{t - i \sigma}{\eta} \re \Li_{2}(e^{i\eta}) 
	\\ \nonumber
&\hspace{2.1cm} +\frac{1}{\eta^2} \bigl[ i t^2 +  (2\sigma +1) t - i\sigma(\sigma+1)\bigr]\im \Li_{3}(e^{i\eta}) \biggr\}
	\\\nonumber
&    + O\biggl(\frac{1}{\eta^{3 + \sigma}}\biggl(t^3 + \Bigl(\frac{1+\epsilon}{\epsilon}\Bigr)^{8}\biggr)\biggr),
    	\\ \nonumber
& \hspace{2cm} (1+\epsilon)t < \eta < \infty, \quad \epsilon > 0, \quad 0 \leq \sigma \leq 1, \quad t \to \infty,	
\end{align}
where the error term is uniform for all $\eta, \epsilon, \sigma$ in the above ranges and the polylogarithm $\Li_m(z)$ is defined by
\begin{align}\label{polylogdef}
  \Li_m(z) = \sum_{k=1}^\infty \frac{z^k}{k^m}, \qquad m \geq 1.
\end{align}  
Similarly, it is straightforward for any $N \geq 4$ to derive an asymptotic formula for $\zeta(s)$ analogous to (\ref{zetaformula2.1b}) with an error term of order 
$$O\biggl(\frac{1}{\eta^{N + \sigma}}\biggl(t^N + \Bigl(\frac{1+\epsilon}{\epsilon}\Bigr)^{2(N+1)}\biggr)\biggr).$$
\end{theorem}
\begin{proof}
Suppose first that $0 < t \leq \eta < \infty$, $\epsilon > 0$, $0 \leq \sigma \leq 1$, $N \geq 2$. All error terms of the form $O(\cdot)$ will be uniform with respect to $\eta, \epsilon, \sigma, N$ unless otherwise specified.
The proof is based on equation (\ref{2.3}), i.e.,
\begin{align}\label{2.22}
\zeta(1-s) = \sum_{n =1}^{[\frac{\eta}{2\pi}]} n^{s-1} - \frac{1}{s} \left(\frac{\eta}{2\pi}\right)^s + G_L(t,\sigma;\eta) + G_U(t, \sigma; \eta), \qquad 0 < \eta < \infty,
\end{align}
where $G_L$ and $G_U$ are defined by
\begin{align}\label{GLdef}
  G_L(t, \sigma; \eta) = \frac{e^{\frac{i\pi s}{2}}}{(2\pi)^s} \int_{-i\eta}^{\infty e^{i\phi_1}} \frac{e^{-z}z^{s-1} dz}{1 - e^{-z}}, \qquad -\frac{\pi}{2} < \phi_1 < \frac{\pi}{2}, \quad 0 < \eta < \infty,
\end{align}
and
\begin{align}\label{GUdef}
  G_U(t, \sigma; \eta) = \frac{e^{-\frac{i\pi s}{2}}}{(2\pi)^s} \int_{i\eta}^{\infty e^{i\phi_2}} \frac{e^{-z}z^{s-1} dz}{1 - e^{-z}}, \qquad -\frac{\pi}{2} < \phi_2 < \frac{\pi}{2}, \quad 0 < \eta < \infty.
\end{align}
The assumption $\eta \notin 2\pi \Z$ implies that $G_L$ and $G_U$ are well-defined.

\medskip
\noindent
{\bf The asymptotics of $G_L$} \nopagebreak
\medskip

\noindent
Using the expansion
\begin{align}\label{geometricsum}
  \frac{1}{1 - e^{-z}} = \sum_{n=0}^\infty e^{-nz}, \qquad \re{z} > 0,
\end{align}
in the definition (\ref{GLdef}) of $G_L$, we find
\begin{align}\label{2.24}
G_L(t, \sigma; \eta) = \frac{e^{\frac{i\pi s}{2}}}{(2\pi)^s}\sum_{n=1}^\infty \int_{-i\eta}^{\infty e^{i\phi_1}} e^{-nz} z^{s-1} dz, \qquad -\frac{\pi}{2} < \phi_1 < \frac{\pi}{2}, \quad 0 < \eta < \infty.
\end{align}
The interchange of the integration and the summation in (\ref{2.24}) is allowed because the sum on the rhs of (\ref{geometricsum}) converges absolutely and uniformly in $\re{z} \geq \delta$ for any $\delta > 0$.
Integration by parts shows that the integral in (\ref{2.24}) can be written as
\begin{align} \nonumber
\int_{-i\eta}^{\infty e^{i\phi_1}} e^{-nz} z^{s-1} dz
= &\; \sum_{j=0}^{N-1} e^{-nz + it\ln{z}} \left(\frac{1}{n - \frac{it}{z}}\frac{d}{dz}\right)^j\frac{z^{\sigma-1}}{n - \frac{it}{z}} \Biggr|_{z = -i\eta}
	\\ \label{GLIBP}
& + \int_{-i\eta}^{\infty e^{i\phi_1}} e^{-nz + it\ln{z}} D_N dz,
\end{align}
where $D_N = D_N(z, n, \sigma, t)$ is short-hand notation for
\begin{align}\label{DNdef}
  D_N := \left(\frac{d}{dz} \frac{1}{n - \frac{it}{z}}\right)^N z^{\sigma-1}.
\end{align}  
We claim that
\begin{align}\label{GLerrorestimate}
\frac{e^{\frac{i\pi s}{2}}}{(2\pi)^s} \sum_{n=1}^\infty
\int_{-i\eta}^{\infty e^{i\phi_1}} e^{-nz + it\ln{z}} D_N dz  = O\biggl(  (2N -1)!!  (N+1)^2 \eta^{\sigma -1 - N} \biggr).
\end{align}
In order to derive equation (\ref{GLerrorestimate}) we first note that we can write $D_N$ in the form
\begin{align}\label{DNterm}
  D_N = \sum_{b=0}^N \sum_{c=0}^N A_{bc}^{(N)} z^{\sigma-1} \frac{t^b (nz)^{N-b} \sigma^c}{(nz - it)^{2N}}
\end{align}
where $A_{bc}^{(N)}$, $b=0, \dots, N$, $c=0, \dots, N$, are integers which satisfy
$$|A_{bc}^{(N)}| \leq (2N-1)!! = 1 \cdot 3 \cdot 5 \cdot \cdots \cdot (2N-1), \qquad b=0, \dots, N, \quad c=0, \dots, N.$$
Indeed, if $D_N$ can be written as in (\ref{DNterm}), then the computation
\begin{align*}
&\frac{d}{dz}\biggl(\frac{1}{n - \frac{it}{z}} A_{bc}^{(N)} z^{\sigma-1} \frac{t^b (nz)^{N-b} \sigma^c}{(nz - it)^{2N}}
\biggr)
	\\
& \hspace{2cm} =A_{bc}^{(N)} z^{\sigma-1} \frac{t^b(nz)^{N-b} \sigma^c}{(n z - it)^{2(N+1)}}
[it(b - N) - it\sigma - nz(1 + b + N) +  nz\sigma)]
\end{align*}
shows that $D_{N+1}$ also can be written as in (\ref{DNterm}).

In order to estimate the lhs of (\ref{GLerrorestimate}), we choose $\phi_1 = 0$. Then, since $t \leq |z|$ and $|n - \frac{it}{z}| \geq n$ for all $z$ on the contour, we can estimate
\begin{align*}
D_N & = O\biggl((2N -1)!! |z|^{\sigma -1}\sum_{b=0}^N \sum_{c=0}^N \biggl| \left(\frac{t}{nz}\right)^b \frac{(nz)^{N} \sigma^c}{z^{2N}(n - \frac{it}{z})^{2N}}\biggr|\biggr)
	\\ \nonumber
& = O\biggl((2N -1)!! |z|^{\sigma -1 - N} (N+1)^2 \frac{1}{n^N} \biggr), \qquad
z = u - i\eta, \quad u \in [0, \infty).
\end{align*}
Thus, the lhs of (\ref{GLerrorestimate}) is
\begin{align}\nonumber
& O\biggl(  e^{-\frac{\pi t}{2}} \sum_{n=1}^\infty\int_{-i\eta}^{\infty}
|e^{-nz}e^{it\ln{z}} | (2N -1)!!  (N+1)^2 |z|^{\sigma -1 - N} \frac{1}{n^N} |dz| \biggr)
	\\ \nonumber
& = O\biggl(  (2N -1)!!  (N+1)^2  \sum_{n=1}^\infty\frac{1}{n^N}
\int_{-i\eta}^{\infty} e^{-n\re{z}} |z|^{\sigma -1 - N}  |dz| \biggr)
	\\ \nonumber
& = O\biggl(  (2N -1)!!  (N+1)^2 \eta^{\sigma -1 - N} \sum_{n=1}^\infty\frac{1}{n^N}
\int_{0}^{\infty} e^{-nu} du \biggr)
	\\ \nonumber
& = O\biggl(  (2N -1)!!  (N+1)^2 \eta^{\sigma -1 - N} \sum_{n=1}^\infty\frac{1}{n^{N+1}} \biggr)
	\\ \nonumber
& = O\biggl(  (2N -1)!!  (N+1)^2 \eta^{\sigma -1 - N} \biggr).
\end{align}
This proves (\ref{GLerrorestimate}).

Equations (\ref{2.24}), (\ref{GLIBP}), and (\ref{GLerrorestimate}) imply that $G_L$ satisfies
\begin{align}\label{GLfinal}
G_L(t, \sigma; \eta) = & \; \frac{e^{\frac{i\pi s}{2}}}{(2\pi)^s}\sum_{n=1}^\infty  \sum_{j=0}^{N-1} e^{-nz + it\ln{z}} \left(\frac{1}{n - \frac{it}{z}}\frac{d}{dz}\right)^j\frac{z^{\sigma-1}}{n - \frac{it}{z}} \biggr|_{z = -i\eta}
	\\ \nonumber
&+ O\biggl(  (2N -1)!!  (N+1)^2 \eta^{\sigma -1 - N} \biggr), \qquad t \leq \eta < \infty, \quad t \to \infty.
\end{align}

\medskip
\noindent
{\bf The asymptotics of $G_U$}\nopagebreak
\medskip

\noindent
We now let $\epsilon > 0$ and suppose that $0 < (1+\epsilon) t < \eta < \infty$.
In analogy with equations (\ref{2.24}) and (\ref{GLIBP}) we now find
\begin{align}\label{GUtsigmaeta}
  G_U(t, \sigma; \eta) = \frac{e^{-\frac{i\pi s}{2}}}{(2\pi)^s} \sum_{n=1}^\infty \int_{i\eta}^{\infty e^{i\phi_2}} e^{-nz} z^{s-1} dz, \qquad -\frac{\pi}{2} < \phi_2 < \frac{\pi}{2}, \quad 0 < \eta < \infty,
\end{align}
and
\begin{align}\nonumber
\int_{i\eta}^{\infty e^{i\phi_2}} e^{-nz} z^{s-1} dz
= &\; \sum_{j=0}^{N-1} e^{-nz + it\ln{z}} \left(\frac{1}{n - \frac{it}{z}}\frac{d}{dz}\right)^j\frac{z^{\sigma-1}}{n - \frac{it}{z}} \biggr|_{z = i\eta}
	\\  \label{GUIBP}
& + \int_{i\eta}^{\infty e^{i\phi_2}} e^{-nz + it\ln{z}} D_N dz,
\end{align}
where $D_N$ is defined by (\ref{DNdef}).
We claim that
\begin{align}\label{GU2errorestimate}
\frac{e^{-\frac{i\pi s}{2}}}{(2\pi)^s}
 \sum_{n=1}^\infty &\int_{i\eta}^{\infty e^{i\phi_2}} e^{-nz + it\ln{z}} D_N dz   
= O\biggl((2N -1)!!  N \Bigl(\frac{1+\epsilon}{\epsilon}\Bigr)^{2N} \eta^{\sigma - N}  \biggr).
\end{align}
Indeed, since $D_N$ can be written in the form given by (\ref{DNterm}),
Jordan's lemma implies that we may choose $\phi_2 = \pi/2$. Then, using the estimate
\begin{align*}
D_N & = O\biggl((2N -1)!! |z|^{\sigma -1}\sum_{b=0}^N \sum_{c=0}^N \biggl| \left(\frac{t}{nz}\right)^b \frac{(nz)^{N} \sigma^c}{z^{2N}(n - \frac{it}{z})^{2N}}\biggr|\biggr)
	\\ \nonumber
& = O\biggl((2N -1)!! |z|^{\sigma -1 - N} (N+1)^2 \frac{n^N}{(n-\frac{t}{\eta})^{2N}} \biggr), \qquad  z = i\lambda, \quad \lambda \in [\eta, \infty),
\end{align*}
we find that the lhs of (\ref{GU2errorestimate}) is
\begin{align*}
O\biggl(e^{\frac{\pi t}{2}} &\sum_{n=1}^\infty\int_{i\eta}^{i\infty} |e^{-nz + it\ln{z}} D_N dz| \biggr)
	\\
& = 
O\biggl(e^{\frac{\pi t}{2}} \sum_{n=1}^\infty\int_{i\eta}^{i \infty} e^{- \frac{\pi t}{2}} 
(2N -1)!! |z|^{\sigma -1 - N} (N+1)^2 \frac{n^N}{(n-\frac{t}{\eta})^{2N}}  |dz| \biggr)
	\\
& = 
O\biggl((2N -1)!!  (N+1)^2 \int_{\eta}^{\infty} 
\lambda^{\sigma -1 - N} d\lambda  \sum_{n=1}^\infty \frac{n^N}{(n- \frac{1}{1+\epsilon})^{2N}}   \biggr)
	\\
& = 
O\biggl((2N -1)!!  (N+1)^2 
\frac{\eta^{\sigma - N}}{N - \sigma} \Bigl(\frac{1+\epsilon}{\epsilon}\Bigr)^{2N} \biggr)
	\\
& = 
O\biggl((2N -1)!!  N \Bigl(\frac{1+\epsilon}{\epsilon}\Bigr)^{2N} \eta^{\sigma - N}  \biggr).
	\\
\end{align*}
In the third equality above we have used the fact that
\begin{align}\label{sumestimate}
\sum_{n=1}^\infty \frac{n^N}{(n-\frac{1}{1+\epsilon})^{2N}} 
& = \sum_{n=1}^\infty \frac{1}{n^N} \frac{1}{(1-\frac{1}{(1+\epsilon)n})^{2N}} 
 	\\ \nonumber
& \leq \frac{1}{(1-\frac{1}{1+\epsilon})^{2N}}  \sum_{n=1}^\infty \frac{1}{n^N} 
= O\biggl(\Bigl(\frac{1+\epsilon}{\epsilon}\Bigr)^{2N}\biggr), \qquad N \geq 2.
\end{align}
This proves (\ref{GU2errorestimate}).

We next claim that (\ref{GU2errorestimate}) can be modified as follows:
\begin{align}\label{GU2errorestimate2}
\frac{e^{-\frac{i\pi s}{2}}}{(2\pi)^s}
 \sum_{n=1}^\infty&\int_{i\eta}^{\infty e^{i\phi_2}} e^{-nz + it\ln{z}} D_N dz   
 	\\ \nonumber
&  = O\biggl((2N +1)!!  N \Bigl(\frac{1+\epsilon}{\epsilon}\Bigr)^{2(N+1)} \eta^{\sigma - N - 1}  \biggr).
\end{align}
Indeed, integration by parts shows that the lhs of (\ref{GU2errorestimate}) equals
\begin{align*}
&\frac{e^{-\frac{i\pi s}{2}}}{(2\pi)^s}
 \sum_{n=1}^\infty
 e^{-nz + it\ln{z}}\frac{1}{n - \frac{it}{z}} D_N \Bigg|_{z = i\eta}
  +\frac{e^{-\frac{i\pi s}{2}}}{(2\pi)^s}
 \sum_{n=1}^\infty
 \int_{i\eta}^{\infty e^{i\phi_2}} e^{-nz + it\ln{z}} D_{N+1} dz.
\end{align*}
By (\ref{GU2errorestimate}), the second term is 
$$O\biggl((2N +1)!!  (N+1)  \Bigl(\frac{1+\epsilon}{\epsilon}\Bigr)^{2(N+1)} \eta^{\sigma - N - 1}  \biggr),$$
while, in view of (\ref{DNterm}) and (\ref{sumestimate}), the first term satisfies
\begin{align*}
 \frac{e^{-\frac{i\pi s}{2}}}{(2\pi)^s}&
 \sum_{n=1}^\infty
 e^{-ni\eta + it\ln{i\eta}}\frac{1}{n - \frac{t}{\eta}}
\sum_{b=0}^N \sum_{c=0}^N A_{bc}^{(N)} (i\eta)^{\sigma-1} \frac{t^b (ni\eta)^{N-b} \sigma^c}{(ni\eta - it)^{2N}}
	\\
&= O\biggl(\eta^{\sigma-1-N}
 \sum_{b=0}^N \sum_{c=0}^N |A_{bc}^{(N)}| 
 \sum_{n=1}^\infty
  \frac{t^b}{\eta^b}\frac{ n^{N-b}}{(n - \frac{t}{\eta})^{2N + 1}}
\biggr)
	\\
&= O\biggl( \eta^{\sigma-1-N} (2N -1)!!  (N+1)^2  \sum_{n=1}^\infty \frac{ n^{N}}{(n - \frac{1}{1+\epsilon})^{2N + 1}}\biggr)
	\\
&= O\biggl( \eta^{\sigma-1-N} (2N -1)!!  (N+1)^2 \Bigl(\frac{1+\epsilon}{\epsilon}\Bigr)^{2N+1} \biggr).
\end{align*}	
This proves (\ref{GU2errorestimate2}).

Equations (\ref{GUtsigmaeta}), (\ref{GUIBP}), and (\ref{GU2errorestimate2}) imply that $G_U$ satisfies
\begin{align} \label{GUfinal}
G_U(t, \sigma; \eta) 
= &\; \frac{e^{-\frac{i\pi s}{2}}}{(2\pi)^s}\sum_{n=1}^\infty
\sum_{j=0}^{N-1} e^{-nz + it\ln{z}} \left(\frac{1}{n - \frac{it}{z}}\frac{d}{dz}\right)^j\frac{z^{\sigma-1}}{n - \frac{it}{z}} \biggr|_{z = i\eta}
	\\ \nonumber
& +  O\biggl((2N +1)!! N \Bigl(\frac{1+\epsilon}{\epsilon}\Bigr)^{2(N+1)} \eta^{\sigma - N - 1}  \biggr), \qquad (1+\epsilon)t < \eta < \infty.
\end{align}

\bigskip
\noindent
{\bf Proof of (\ref{zetaformula2.1})}\nopagebreak
\medskip

\noindent
Substituting the expressions (\ref{GLfinal}) and (\ref{GUfinal}) for $G_L$ and $G_U$ into (\ref{2.22}), we find
\begin{align*}
\zeta(1-s) = &\; \sum_{n =1}^{[\frac{\eta}{2\pi}]} n^{s-1} - \frac{1}{s} \left(\frac{\eta}{2\pi}\right)^s 
	\\ \nonumber
&+ \frac{e^{\frac{i\pi s}{2}}}{(2\pi)^s}\sum_{n=1}^\infty  \sum_{j=0}^{N-1} e^{-nz + it\ln{z}} \left(\frac{1}{n - \frac{it}{z}}\frac{d}{dz}\right)^j\frac{z^{\sigma-1}}{n - \frac{it}{z}} \biggr|_{z = -i\eta}
	\\
& + \frac{e^{-\frac{i\pi s}{2}}}{(2\pi)^s}\sum_{n=1}^\infty
\sum_{j=0}^{N-1} e^{-nz + it\ln{z}} \left(\frac{1}{n - \frac{it}{z}}\frac{d}{dz}\right)^j\frac{z^{\sigma-1}}{n - \frac{it}{z}} \biggr|_{z = i\eta}
	\\ \nonumber
& + O\biggl((2N +1)!!  N \Bigl(\frac{1+\epsilon}{\epsilon}\Bigr)^{2(N+1)} \eta^{\sigma - N - 1}  \biggr),
	\\ 
& \hspace{5cm}  (1+\epsilon) t < \eta < \infty, \quad t \to \infty.
\end{align*}
Replacing $\sigma$ by $1 - \sigma$ and taking the complex conjugate of the resulting equation, we find
(\ref{zetaformula2.1}).

\bigskip
\noindent
{\bf Proof of (\ref{zetaformula2.1b})}\nopagebreak
\medskip

\noindent
Letting $N = 3$ in (\ref{zetaformula2.1}), we find
\begin{align*}
\zeta(s) = &\; \sum_{n =1}^{[\frac{\eta}{2\pi}]} n^{-s} - \frac{1}{1-s} \left(\frac{\eta}{2\pi}\right)^{1-s} 
	\\ \nonumber
&+ \frac{e^{-\frac{i\pi (1-s)}{2}}}{(2\pi)^{1-s}}(i\eta)^{-s}\sum_{n=1}^\infty 
e^{-i\eta n}\biggl\{
\frac{1}{n}\frac{1}{1 + q}
+ \frac{1}{\eta n^2} \frac{i (q (\sigma -1)+\sigma )}{(1+q)^3}
	\\
&\hspace{2cm} +\frac{1}{\eta^2 n^3} \frac{-q^2 (\sigma -1)^2+q \left(-2 \sigma ^2+\sigma +2\right)-\sigma  (\sigma
   +1)}{(1+q)^5} \biggr\}
   	\\ \nonumber
&+ \frac{e^{\frac{i\pi (1-s)}{2}}}{(2\pi)^{1-s}}(-i\eta)^{-s}\sum_{n=1}^\infty 
e^{i\eta n}\biggl\{
\frac{1}{n}\frac{1}{1 - q}
+ \frac{1}{\eta n^2}  \frac{i (q (\sigma -1)-\sigma )}{(1-q)^3}
	\\
&\hspace{2cm}  -\frac{1}{\eta^2 n^3} \frac{q^2 (\sigma -1)^2+q \left(-2 \sigma ^2+\sigma +2\right)+\sigma  (\sigma
   +1)}{(1 - q)^5} \biggr\}
   	\\
& + O\biggl(\eta^{-\sigma - 3} \Bigl(\frac{1+\epsilon}{\epsilon}\Bigr)^{8}\biggr), \qquad q = \frac{t}{n\eta}.
\end{align*}
The above sums can be expanded as series in $q$ with coefficients expressible in terms of the polylogarithms $\Li_m(z)$ defined by (\ref{polylogdef}).
 For example, 
 \begin{align*}
\sum_{n=1}^\infty \frac{e^{i\eta n}}{n(1 - q)}
& = \sum_{n=1}^\infty \frac{e^{i\eta n}}{n}\sum_{k=0}^\infty \biggl(\frac{t}{n\eta}\biggr)^k
	\\
& = \sum_{k=0}^\infty \biggl(\frac{t}{\eta}\biggr)^k \sum_{n=1}^\infty \frac{e^{i\eta n}}{n^{k+1}}
	\\
& = \sum_{k=0}^{\infty} \left(\frac{t}{\eta}\right)^k \Li_{k+1}\left(e^{i \eta}\right),
\end{align*}
where the interchange of the two sums can be justified by separating the terms where $k = 0$ and noticing that the remaining double sum is absolutely convergent. More precisely,
 \begin{align*}
\lim_{N\to \infty} \lim_{M \to \infty} \sum_{n=1}^N \frac{e^{i\eta n}}{n}\sum_{k=0}^M \biggl(\frac{t}{n\eta}\biggr)^k
&=
\lim_{N\to \infty} \left(\sum_{n=1}^N \frac{e^{i\eta n}}{n} +  \lim_{M \to \infty}  \sum_{n=1}^N \frac{e^{i\eta n}}{n}\sum_{k=1}^M \biggl(\frac{t}{n\eta }\biggr)^k\right)
	\\
& = \Li_1(e^{i\eta})  + \sum_{n=1}^\infty \frac{e^{i\eta n}}{n}\sum_{k=1}^\infty \biggl(\frac{t}{n\eta}\biggr)^k
	\\
& = \Li_1(e^{i\eta}) + \sum_{k=1}^\infty \sum_{n=1}^\infty  \biggl(\frac{t}{\eta}\biggr)^k \frac{e^{i\eta n}}{n^{k+1}}
	\\
& = \sum_{k=0}^\infty \biggl(\frac{t}{\eta}\biggr)^k \Li_{k+1}\left(e^{i \eta}\right).
\end{align*}
Similarly, for any $j\geq0$ and $l,m\geq 1$,
\begin{align*}
\sum_{n=1}^\infty \frac{e^{i\eta n}}{n^l} \frac{q^j}{(1 - q)^m}
&= \sum_{n=1}^\infty \frac{e^{i\eta n}}{n^l}\sum_{k=0}^\infty {{k+m-1}\choose{k}} \biggl(\frac{t}{n\eta}\biggr)^{k+j}
	\\
& =  \sum_{k=0}^\infty {{k+m-1}\choose{k}} \biggl(\frac{t}{\eta}\biggr)^{k+j}
\Li_{k+l+j}(e^{i\eta}).
\end{align*}
Letting $\eta  \to -\eta$ in this equation, we find
\begin{align*}
\sum_{n=1}^\infty \frac{e^{-i\eta n}}{n^l}\frac{(-q)^j}{(1 + q)^m}
 =  \sum_{k=0}^\infty {{k+m-1}\choose{k}} \biggl(-\frac{t}{\eta}\biggr)^{k+j}
\Li_{k+l+j}(e^{-i\eta}).
\end{align*}
This gives the following expression for the Riemann zeta function:
\begin{align*}
& \zeta(s)  =  \sum_{n =1}^{[\frac{\eta}{2\pi}]} n^{-s} - \frac{1}{1-s} \left(\frac{\eta}{2\pi}\right)^{1-s} 
	\\ \nonumber
&+ \frac{e^{-\frac{i\pi (1-s)}{2}}}{(2\pi)^{1-s}}(i\eta)^{-s}\biggl\{
\sum_{k=0}^\infty \biggl(-\frac{t}{\eta}\biggr)^k \Li_{k+1}(e^{-i\eta})
+ \frac{i\sigma}{\eta} \sum_{k=0}^\infty {{k+2}\choose{k}}\biggl(-\frac{t}{\eta}\biggr)^k \Li_{k+2}(e^{-i\eta})
	\\
&\hspace{3.3cm} - \frac{i(\sigma-1)}{\eta} \sum_{k=0}^\infty {{k+2}\choose{k}} \biggl(-\frac{t}{\eta}\biggr)^{k+1} \Li_{k+3}(e^{-i\eta})
	\\
&\hspace{3.3cm} - \frac{\sigma(\sigma+1)}{\eta^2} \sum_{k=0}^\infty{{k+4}\choose{k}} \biggl(-\frac{t}{\eta}\biggr)^k \Li_{k+3}(e^{-i\eta})
	\\
&\hspace{3.3cm} +\frac{2\sigma^2 - \sigma -2}{\eta^2} \sum_{k=0}^\infty {{k+4}\choose{k}}\biggl(-\frac{t}{\eta}\biggr)^{k+1} \Li_{k+4}(e^{-i\eta})
	\\
&\hspace{3.3cm} - \frac{(\sigma -1)^2}{\eta^2} \sum_{k=0}^\infty {{k+4}\choose{k}} \biggl(-\frac{t}{\eta}\biggr)^{k+2} \Li_{k+5}(e^{-i\eta})
 \biggr\}
   	\\ \nonumber
&+ \frac{e^{\frac{i\pi (1-s)}{2}}}{(2\pi)^{1-s}}(-i\eta)^{-s}\biggl\{
\sum_{k=0}^\infty \biggl(\frac{t}{\eta}\biggr)^k \Li_{k+1}(e^{i\eta})
- \frac{i\sigma}{\eta}\sum_{k=0}^\infty {{k+2}\choose{k}} \biggl(\frac{t}{\eta}\biggr)^k \Li_{k+2}(e^{i\eta})
	\\
&\hspace{3.3cm} + \frac{i(\sigma-1)}{\eta}\sum_{k=0}^\infty {{k+2}\choose{k}} \biggl(\frac{t}{\eta}\biggr)^{k+1} \Li_{k+3}(e^{i\eta})
	\\
&\hspace{3.3cm} - \frac{\sigma(\sigma+1)}{\eta^2}\sum_{k=0}^\infty {{k+4}\choose{k}} \biggl(\frac{t}{\eta}\biggr)^{k} \Li_{k+3}(e^{i\eta})
	\\
&\hspace{3.3cm} + \frac{2\sigma^2 - \sigma -1}{\eta^2}\sum_{k=0}^\infty {{k+4}\choose{k}} \biggl(\frac{t}{\eta}\biggr)^{k+1} \Li_{k+4}(e^{i\eta})
	\\
&\hspace{3.3cm} - \frac{(\sigma-1)^2}{\eta^2}\sum_{k=0}^\infty {{k+4}\choose{k}} \biggl(\frac{t}{\eta}\biggr)^{k+2} \Li_{k+5}(e^{i\eta})
 \biggr\}
   	\\
& + O\biggl(\eta^{-\sigma - 3} \Bigl(\frac{1+\epsilon}{\epsilon}\Bigr)^{8}\biggr).
\end{align*}
By including only the first few terms in the sums over $k$, we find an expansion for $\zeta(s)$ in powers of $\frac{t}{\eta}$. For example, including only the terms of order larger than $O(\frac{t^3}{\eta^{3 + \sigma}})$ we have
\begin{align*}
\zeta(s) = &\; \sum_{n =1}^{[\frac{\eta}{2\pi}]} n^{-s} - \frac{1}{1-s} \left(\frac{\eta}{2\pi}\right)^{1-s} 
	\\ \nonumber
&- \frac{i\eta^{-s}}{(2\pi)^{1-s}}\biggl\{
 \Li_{1}(e^{-i\eta}) -\frac{t}{\eta} \Li_{2}(e^{-i\eta})
+\biggl(\frac{t}{\eta}\biggr)^2 \Li_{3}(e^{-i\eta})
	\\
&+ \frac{i\sigma}{\eta}  \Li_{2}(e^{-i\eta})
- \frac{3i\sigma t}{\eta^2}  \Li_{3}(e^{-i\eta})
	\\
& + \frac{i(\sigma-1)t}{\eta^2} \Li_{3}(e^{-i\eta})
- \frac{\sigma(\sigma+1)}{\eta^2} \Li_{3}(e^{-i\eta})
 \biggr\}
   	\\ \nonumber
&+ \frac{i\eta^{-s}}{(2\pi)^{1-s}} \biggl\{
\Li_{1}(e^{i\eta})
+ \frac{t}{\eta} \Li_{2}(e^{i\eta})
+ \biggl(\frac{t}{\eta}\biggr)^2 \Li_{3}(e^{i\eta})
	\\
&- \frac{i\sigma}{\eta}  \Li_{2}(e^{i\eta})
- \frac{3 i\sigma t}{\eta^2} \Li_{3}(e^{i\eta})
	\\
& + \frac{i(\sigma-1)t}{\eta^2} \Li_{3}(e^{i\eta})
- \frac{\sigma(\sigma+1)}{\eta^2}  \Li_{3}(e^{i\eta})
 \biggr\}
   	\\
& + O\biggl(\frac{1}{\eta^{3 + \sigma}}\biggl(t^3 + \Bigl(\frac{1+\epsilon}{\epsilon}\Bigr)^{8}\biggr)\biggr).
\end{align*}
After simplification we find
\begin{align*}
\zeta(s) = &\; \sum_{n =1}^{[\frac{\eta}{2\pi}]} n^{-s} - \frac{1}{1-s} \left(\frac{\eta}{2\pi}\right)^{1-s} 
	\\ \nonumber
& + \frac{2i\eta^{-s}}{(2\pi)^{1-s}} \biggl\{
i \im \Li_{1}(e^{i\eta})
+  \frac{t}{\eta} \re \Li_{2}(e^{i\eta}) 
+ i \biggl(\frac{t}{\eta}\biggr)^2 \im \Li_{3}(e^{i\eta})
	\\
& \hspace{2cm} - \frac{i\sigma}{\eta} \re \Li_{2}(e^{i\eta})
+ \frac{3\sigma t}{\eta^2} \im \Li_{3}(e^{i\eta})
	\\
&  \hspace{2cm} - \frac{(\sigma-1) t}{\eta^2} \im \Li_{3}(e^{i\eta}) 
- \frac{i\sigma(\sigma+1)}{\eta^2} \im \Li_{3}(e^{i\eta})
 \biggr\}
 	\\
&  + O\biggl(\frac{1}{\eta^{3 + \sigma}}\biggl(t^3 + \Bigl(\frac{1+\epsilon}{\epsilon}\Bigr)^{8}\biggr)\biggr),
\end{align*}
The identity $\Li_{1}\left(e^{i \eta}\right) = -\ln(1 - e^{i \eta})$ now yields (\ref{zetaformula2.1b}).
\end{proof}   

\begin{theorem}[The asymptotic expansion to all orders for the case $t = \eta$] \label{th3.2}
Let $\zeta(s)$, $s = \sigma + it$, $\sigma, t \in \R$, denote the Riemann zeta function. Then,
\begin{align}\nonumber
\zeta(s) = & \sum_{n =1}^{[\frac{t}{2\pi}]} n^{-s} - \frac{1}{1-s} \left(\frac{t}{2\pi}\right)^{1-s} 
	\\ \nonumber
& + \frac{e^{-\frac{i\pi (1-s)}{2}}}{(2\pi)^{1-s}}\sum_{n=1}^\infty  \sum_{j=0}^{N-1} e^{-nz - it\ln{z}} \left(\frac{1}{n + \frac{it}{z}}\frac{d}{dz}\right)^j\frac{z^{-\sigma}}{n + \frac{it}{z}} \biggr|_{z = it}
	\\ \nonumber
& + \frac{e^{\frac{i\pi (1-s)}{2}}}{(2\pi)^{1-s}}\sum_{n=2}^\infty
\sum_{j=0}^{N-1} e^{-nz - it\ln{z}} \left(\frac{1}{n + \frac{it}{z}}\frac{d}{dz}\right)^j\frac{z^{-\sigma}}{n + \frac{it}{z}} \biggr|_{z = -it}
	\\ \nonumber
& + \left(\frac{t}{2\pi}\right)^{1-s} e^{it}
\sum_{k=0}^{2N}   \frac{\overline{c_k(1-\sigma)}\Gamma(\frac{k+1}{2})}{t^{\frac{k+1}{2}}}
 +  O\biggl((2N +1)!! N 2^{2N} t^{-\sigma - N}  \biggr),
	\\ \label{zetaformula2.2}
& \hspace{2cm}  0 \leq \sigma \leq 1, \quad N \geq 2, \quad t \to \infty,	
\end{align}
where the error term is uniform for all $\sigma, N$ in the above ranges and the coefficients $c_k(\sigma)$ are given by equation (\ref{rhoquotientexpansion}) below. The first few of the $c_k$'s are given by
\begin{align} \nonumber
& c_0(\sigma) = \frac{1-i}{2},
	\\\nonumber
& c_1(\sigma) = \frac{i}{3}-i \sigma,
	\\\nonumber
& c_2(\sigma) = -\frac{1+i}{12}\left(6 \sigma ^2-6 \sigma +1\right),
	\\\label{firstfewcks}
& c_3(\sigma) = \frac{1}{135} \left(-45 \sigma ^3+90 \sigma ^2-45 \sigma
   +4\right),
	\\ \nonumber
& c_4(\sigma) = \frac{i-1}{432} \left(36 \sigma ^4-120 \sigma
   ^3+120 \sigma ^2-36 \sigma +1\right),
	\\\nonumber
& c_5(\sigma) = \frac{i}{5670} \left(189 \sigma ^5-945 \sigma ^4+1575
   \sigma ^3-987 \sigma ^2+168 \sigma
   +8\right),
	\\ \nonumber
& c_6(\sigma) = \frac{1+i}{194400} \left(1080 \sigma
   ^6-7560 \sigma ^5+18900 \sigma ^4-20160 \sigma ^3+8190 \sigma ^2-450 \sigma
   -139\right).
\end{align}

For $N=3$ equation (\ref{zetaformula2.2}) simplifies to  
\begin{align} \nonumber
\zeta(s) = & \sum_{n =1}^{[\frac{t}{2\pi}]} n^{-s} - \frac{1}{1-s} \left(\frac{t}{2\pi}\right)^{1-s} 
	\\ \nonumber
&- \frac{2it^{-s}e^{it}}{(2\pi)^{1-s}} \biggl\{\frac{e^{-it}}{2t^2}\bigl[-t^2 - it(\sigma -1) + (\sigma -1)^2\bigr]
- i \im \Li_1(e^{it})
+ \frac{i\sigma}{t} \re \Li_2(e^{it})
	\\\nonumber
& \hspace{1cm} + \frac{it + \sigma + \sigma^2}{t^2} i\im \Li_3(e^{it})
+\frac{2+3\sigma}{t^2}  \re \Li_4(e^{it})
+ \frac{3i}{t^2} \im \Li_5(e^{it})\biggr\}
\end{align}
\begin{align}\nonumber
&+ \frac{ t^{-s}e^{it}}{(2\pi)^{1-s}}\biggl\{ \frac{1+i}{2} \sqrt{\pi t} -\frac{1}{3} i (3 \sigma -2)
	\\\nonumber
&\hspace{1cm}
+ \frac{i-1}{24 \sqrt{t}} \sqrt{\pi } \left(6\sigma ^2-6 \sigma +1\right)
+ \frac{1}{135 t} \left(45 \sigma^3-45 \sigma ^2+4\right)
	\\\nonumber
&\hspace{1cm}
-\frac{1+i}{576 t^{\frac{3}{2}}}\sqrt{\pi } \left(36 \sigma^4-24 \sigma ^3-24 \sigma ^2+12 \sigma +1\right)
	\\\nonumber
&\hspace{1cm}+\frac{i}{2835 t^2} \bigl(189 \sigma ^5-315 \sigma ^3+42 \sigma ^2+84 \sigma -8\bigr)
	\\\nonumber
& \hspace{1cm} 
+ \frac{1 - i}{103680 t^{\frac{5}{2}}} \sqrt{\pi } \bigl(1080 \sigma^6+1080 \sigma ^5-2700 \sigma ^4-1440 \sigma ^3
	\\\label{zetaformula2.2b}
& \hspace{1cm} +1710 \sigma ^2+270 \sigma -139\bigr) \biggr\}  + O\bigl(t^{-\sigma - 3}\bigr),
\qquad 0 \leq \sigma \leq 1, \quad t \to \infty,	
\end{align}
where the error term is uniform for all $\sigma$ in the above range and the polylogarithm $\Li_m(z)$, $m \geq 1$, is defined by (\ref{polylogdef}).
Similarly, it is straightforward for any $N \geq 4$ to derive an asymptotic formula for $\zeta(s)$ analogous to (\ref{zetaformula2.2b}) with an error term of order $O\bigl(t^{-\sigma-N}\bigr)$.
\end{theorem}
\begin{proof}
Setting $\eta = t$ in (\ref{2.22}), we find
\begin{align}\label{zetaformulaaist}
\zeta(1-s) = \sum_{n =1}^{[\frac{t}{2\pi}]} n^{s-1} - \frac{1}{s} \left(\frac{t}{2\pi}\right)^s + G_L(t,\sigma;t) + G_U(t, \sigma; t), \qquad 0 < t < \infty.
\end{align}
Equation (\ref{GLfinal}) is valid also when $\eta =t$ and gives the asymptotics of $G_L(t,\sigma; t)$.
On the other hand, using the expansion (\ref{geometricsum}), we write
\begin{align*}
G_U(t, \sigma; t) = G_U^{(1)}(t, \sigma; t) + G_U^{(2)}(t, \sigma; t),
\end{align*}
where
\begin{align}\label{GU1def}
& G_U^{(1)}(t, \sigma; t) = \frac{e^{-\frac{i\pi s}{2}}}{(2\pi)^s}\int_{it}^{\infty e^{i\phi_2}} 
e^{-z} z^{s-1} dz, \qquad -\frac{\pi}{2} < \phi_2 < \frac{\pi}{2}, \quad 0 < t < \infty,
\end{align}
and
\begin{align}\label{GU2def}
& G_U^{(2)}(t, \sigma; t) = \frac{e^{-\frac{i\pi s}{2}}}{(2\pi)^s}\sum_{n=2}^\infty \int_{it}^{\infty e^{i\phi_2}} 
e^{-nz} z^{s-1} dz, \qquad -\frac{\pi}{2} < \phi_2 < \frac{\pi}{2}, \quad 0 < t < \infty.
\end{align}
The asymptotics of $G_U^{(2)}$ can be found using integration by parts, whereas the asymptotics of $G_U^{(1)}$ will be computed by considering the critical point at $z = it$.

\bigskip
\noindent
{\bf The asymptotics of $G_U^{(2)}$}\nopagebreak
\medskip
 
\noindent
Repeating the steps that led to (\ref{GUfinal}) but with the sum over $n$ only going from $2$ to $\infty$ and with $t=\eta$, we find the following analog of (\ref{GUfinal}):
\begin{align} \nonumber
G_U^{(2)}(t, \sigma; t) 
= &\; \frac{e^{-\frac{i\pi s}{2}}}{(2\pi)^s}\sum_{n=2}^\infty
\sum_{j=0}^{N-1} e^{-nz + it\ln{z}} \left(\frac{1}{n - \frac{it}{z}}\frac{d}{dz}\right)^j\frac{z^{\sigma-1}}{n - \frac{it}{z}} \biggr|_{z = it}
	\\ \label{GU2final}
& +  O\biggl((2N +1)!! N 2^{2N} t^{\sigma - N - 1}  \biggr).
\end{align}

\medskip
\noindent
{\bf The asymptotics of $G_U^{(1)}$}\nopagebreak
\medskip

\noindent
Letting in the definition (\ref{GU1def}) of $G_U^{(1)}$, $\phi_2 = 0$ and
$$z = it + t \rho, \qquad \rho \in [0, \infty),$$
we find
\begin{align*}
& G_U^{(1)} = \frac{e^{-\frac{i\pi s}{2}}}{(2\pi)^s} 
(it)^{s - 1} te^{-it} \int_{0}^{\infty}  e^{-t\left(\rho - i\ln(1 - i \rho) \right)} (1 - i\rho)^{\sigma-1} d\rho.
\end{align*}
Letting 
\begin{align}\label{vrho}
  v = \rho - i\ln(1 - i \rho),
\end{align}  
we find
\begin{align*}
& G_U^{(1)} =  \left(\frac{t}{2\pi}\right)^{s} e^{-it}
\int_\gamma  e^{-t v}  \frac{(1 - i\rho(v))^{\sigma}}{\rho(v)}dv,
\end{align*}
where $\rho(v)$ is defined by inverting (\ref{vrho}) and  $\gamma$  denotes the image of the contour $[0, \infty)$ under (\ref{vrho}). In order to ascertain that the value of $\rho(v)$ is well-defined, we consider in detail the map $\phi$ defined by
\begin{align*}
\phi: \; &\C\setminus [-i, -i\infty) \to \C,
	\\
& \rho \mapsto v = \rho - i\ln(1 - i \rho),
\end{align*}
where $[-i, -i\infty)$ is a branch cut and the principal branch is chosen for the logarithm, i.e.,
\begin{align}\label{vrholnarg}
  v = \rho - i \ln{|1-i\rho|} + \arg(1-i\rho), \qquad \arg(1 - i\rho) \in (-\pi, \pi).
\end{align}  
The function $\phi$ satisfies
\begin{align*}
& \phi(\rho) = 0 \quad \text{iff} \quad \rho = 0;
	\\
& \phi'(\rho) = \frac{\rho}{i + \rho} = 0 \quad \text{iff} \quad \rho = 0;	
	\\ 
& v = \phi(\rho) = -\frac{i\rho^2}{2} + O(\rho^3), \qquad \rho \to 0.
\end{align*}
Moreover, we claim that $\phi$ maps the first quadrant of the complex $\rho$-plane bijectively onto the region delimited by the contour $\gamma = \phi([0, \infty))$ and the positive imaginary axis in the complex $v$-plane, see figure \ref{mapphi.pdf}. Indeed, it is clear that $\phi$ maps $[0, i\infty)$ bijectively onto $[0, i\infty)$. Writing $\rho = \rho_1 + i\rho_2$, we have
\begin{align*}
  \frac{\partial \re \phi}{\partial \rho_1} = \frac{|\rho|^2 + \rho_2}{\rho_1^2 + (\rho_2 +1)^2}, \qquad
    \frac{\partial \re \phi}{\partial \rho_2} = \frac{\rho_1}{\rho_1^2 + (\rho_2 +1)^2}, \qquad
\end{align*}
This shows that each level curve of the function  $\re \phi$ in the first quadrant intersects the positive $\rho$-axis in a unique point and that the vector
\begin{align}\label{tangenttolevelcurve}
  \bigl(-\rho_1, |\rho|^2 + \rho_2\bigr)
\end{align}
is tangent to the level curve passing through $\rho$. Since
\begin{align*}
  \frac{\partial \im \phi}{\partial \rho_1} = \frac{-\rho_1}{\rho_1^2 + (\rho_2 +1)^2}, \qquad
    \frac{\partial \im \phi}{\partial \rho_2} = \frac{|\rho|^2 + \rho_2}{\rho_1^2 + (\rho_2 +1)^2},
\end{align*}
we infer that the derivative of $\im \phi$ in the direction of (\ref{tangenttolevelcurve}) equals $|\rho|^2$, showing that $\im \phi$ is strictly increasing along each level curve. This proves the claim.  

Thus, if $v$ belongs to the region delimited by the contour $\gamma = \phi([0, \infty))$ and the positive imaginary axis, we may define $\rho(v)$ as the unique inverse image of $v$ under $\phi$ which belongs to the first quadrant.  
Using analyticity to deform the contour $\gamma$ to the positive real axis, we obtain
\begin{align*}
& G_U^{(1)} =  \left(\frac{t}{2\pi}\right)^{s} e^{-it}
\int_0^\infty  e^{-t v}  \frac{(1 - i\rho(v))^{\sigma}}{\rho(v)}dv.
\end{align*}

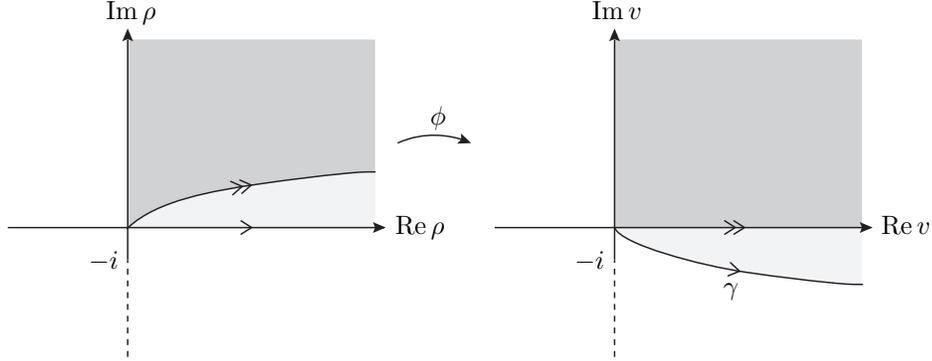
\begin{figure}
\bigskip \begin{overpic}[width=.92\textwidth]{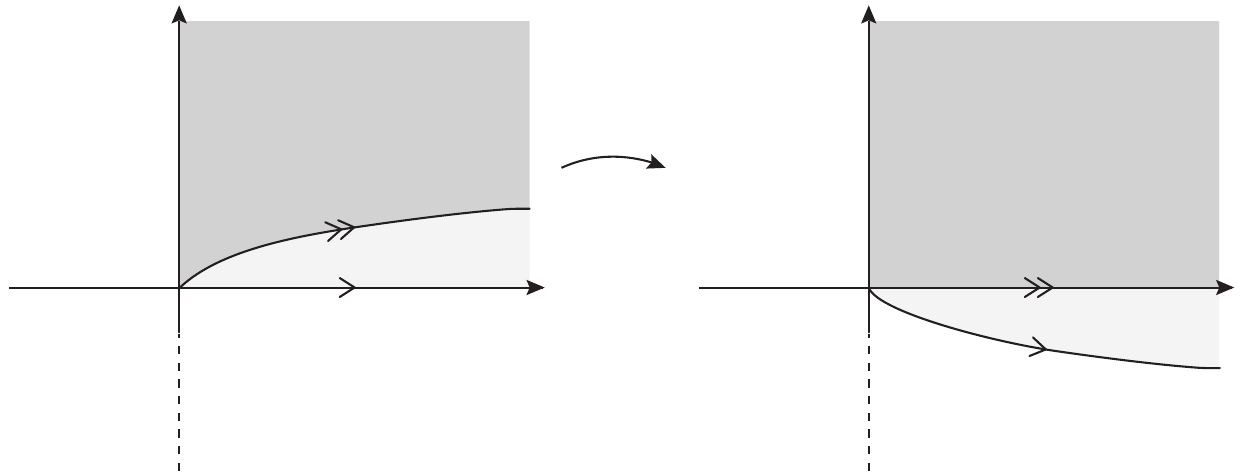}
 \put(45,14.5){$\re \rho$}
 \put(12,39){$\im \rho$}
 \put(100.5,14.5){$\re v$}
 \put(67.5,39){$\im v$}
 \put(49,27){$\phi$}
 \put(82.5,7.5){$\gamma$}
 \put(10,10.5){$-i$}
 \put(65.5,10.5){$-i$}
   \end{overpic}
   \caption{The function $\phi$ maps the indicated shaded areas in the complex $\rho$-plane bijectively onto the corresponding shaded areas in the $v$-plane.}\label{mapphi.pdf}
\end{figure}

We claim that there exist coefficients $\{c_k(\sigma)\}_0^{N-1}$ such that
\begin{align}\label{rhoquotientexpansion}
 \frac{(1 - i\rho(v))^{\sigma}}{\rho(v)}
 = \frac{1}{\sqrt{v}} \biggl(\sum_{k=0}^{N-1} c_k(\sigma) v^{k/2} + r_N(\sqrt{v})\biggr), \qquad v \geq 0,
\end{align}
where $r_N(\sqrt{v})$ satisfies
\begin{align}\label{rNintegral}
\int_{0}^{\infty}  e^{-t v} \frac{r_N(\sqrt{v})}{\sqrt{v}} dv = O\biggl(\frac{\Gamma(\frac{N+1}{2})}{t^{\frac{N+1}{2}}}\biggr).
\end{align}
Indeed, since $\phi$ has a double zero at $\rho =0$, there exists a neighborhood $V$ of $\rho = 0$ in the complex $\rho$-plane and a Riemann surface $\Sigma$ which is a two-sheeted cover of $V$ with local parameter $\lambda:=\sqrt{v}$, such that the map $\rho \mapsto \sqrt{v}$ is bijective and holomorphic $V \to \Sigma$. 
The function $\phi$ is analytic away from $[-i, -i\infty)$ and $\phi'(\rho) \neq 0$ for all $\rho \neq 0$, thus the Riemann surface $\Sigma$ and the map
$\sqrt{v} \mapsto \rho$ can be analytically extended as long as  $\phi^{-1}(v) \cap [-i, -i\infty) = \emptyset$. The function
$\phi$ maps the segments $[-i + 0, -i\infty+0)$ and $[-i-0, -i\infty-0)$ onto the lines $\re{v} = -\pi$ and $\re{v} = \pi$ respectively, hence $\sqrt{v} \mapsto \rho$ is well-defined and analytic for all $|v| < \pi$.
Thus, the map
$$f:\lambda \mapsto  \lambda \frac{(1 - i\rho)^{\sigma}}{\rho} = \lambda\frac{e^{\sigma\ln(1 - i\rho)}}{\rho}$$ 
is holomorphic from the open disk $\{0 < |\lambda| < \sqrt{\pi}\} \subset \Sigma$ to $\C$.
Let
$$c_k(\sigma) = \frac{f^{(k)}(0)}{k!}, \qquad k = 0, \dots, N-1,$$
and define $r_N(\lambda)$ so that (\ref{rhoquotientexpansion}) holds. 
Then $r_N(\lambda)$ satisfies
$$r_N(\lambda) = \sum_{k=N}^{\infty} \frac{f^{(k)}(0)}{k!} \lambda^{k}$$
for $0 \leq |\lambda| < \sqrt{\pi}$. 
In order to prove (\ref{rNintegral}), it is sufficient to show that there exists a constant $C$ independent of $N$  such that
\begin{align}\label{rNCvN2}
|r_N(\sqrt{v})| \leq C v^{\frac{N}{2}}, \qquad v \geq 0, \quad N \geq 2, \quad 0 \leq \sigma \leq 1.
\end{align}
Cauchy's estimates imply that the numbers $c_k(\sigma) = \frac{f^{(k)}(0)}{k!}$ are uniformly bounded for all $k \geq 0$ and $\sigma \in [0, 1]$. Thus, an application of Cauchy's estimates in a disk of radius $|\lambda | = r < \sqrt{\pi}$ yields
$$|r_N(\lambda)| \leq M_r \sum_{k=N}^{\infty} \biggl(\frac{|\lambda|}{r}\biggr)^k
\leq M_r \biggl(\frac{|\lambda|}{r}\biggr)^N  \frac{1}{1 - \frac{|\lambda|}{r}}, \qquad |\lambda | < r,$$
where $M_r =  \text{max}_{|\lambda| =r} |f(\lambda)|$. Choosing any $r \in (1, \sqrt{\pi})$ we see that such a $C$ exists for $0 \leq v \leq 1$. On the other hand, equations (\ref{vrholnarg}) and (\ref{rhoquotientexpansion}) imply that such a $C$ exists also for $v \geq 1$. This proves (\ref{rNCvN2}) and hence also (\ref{rhoquotientexpansion}) and (\ref{rNintegral}).

Equation (\ref{rhoquotientexpansion}) together with the identity
$$\int_0^\infty e^{-tv} v^\alpha dv = \frac{\Gamma(\alpha+1)}{t^{\alpha+1}}, \qquad \alpha > -1,$$
imply the following asymptotic expansion of $G_U^{(1)}$:
\begin{align}\label{GU1asymptotics}
& G_U^{(1)}(t, \sigma; t) =  \left(\frac{t}{2\pi}\right)^{s} e^{-it}
\sum_{k=0}^{N-1}   \frac{c_k\Gamma(\frac{k+1}{2})}{t^{\frac{k+1}{2}}}
+ O\biggl(\frac{\Gamma(\frac{N+1}{2})}{t^{\frac{N+1}{2} - \sigma}}\biggr).
\end{align}

\medskip
\noindent
{\bf Proof of (\ref{zetaformula2.2})}\nopagebreak
\medskip

\noindent
Equations (\ref{GLfinal}), (\ref{GU2final}), and (\ref{GU1asymptotics}) give the asymptotics of $G_L$ and $G_U$. Substitution into (\ref{zetaformulaaist}) yields
\begin{align*}
\zeta(1-s) = & \sum_{n =1}^{[\frac{t}{2\pi}]} n^{s-1} - \frac{1}{s} \left(\frac{t}{2\pi}\right)^s 
	\\
& + \frac{e^{\frac{i\pi s}{2}}}{(2\pi)^s}\sum_{n=1}^\infty  \sum_{j=0}^{N-1} e^{-nz + it\ln{z}} \left(\frac{1}{n - \frac{it}{z}}\frac{d}{dz}\right)^j\frac{z^{\sigma-1}}{n - \frac{it}{z}} \biggr|_{z = -it}
	\\ \nonumber
& + \frac{e^{-\frac{i\pi s}{2}}}{(2\pi)^s}\sum_{n=2}^\infty
\sum_{j=0}^{N-1} e^{-nz + it\ln{z}} \left(\frac{1}{n - \frac{it}{z}}\frac{d}{dz}\right)^j\frac{z^{\sigma-1}}{n - \frac{it}{z}} \biggr|_{z = it}
	\\
& + \left(\frac{t}{2\pi}\right)^{s} e^{-it}
\sum_{k=0}^{2N}   \frac{c_k(\sigma)\Gamma(\frac{k+1}{2})}{t^{\frac{k+1}{2}}}
 +  O\biggl((2N +1)!! N 2^{2N} t^{\sigma - N - 1}  \biggr),
\end{align*}
where we have replaced $N$ by $2N+1$ in (\ref{GU1asymptotics}) and used that $\Gamma(N+1) = N! \leq (2N +1)!!$ to eliminate one of the error terms.
Replacing $\sigma$ by $1 - \sigma$ and taking the complex conjugate of the resulting equation, we find
(\ref{zetaformula2.2}).

\bigskip
\noindent
{\bf Proof of (\ref{zetaformula2.2b})}\nopagebreak
\medskip

\noindent
Letting $N = 3$ in (\ref{zetaformula2.2}) and using the expressions in (\ref{firstfewcks}) for $\{c_k\}_0^6$, we find
that the term on the rhs of (\ref{zetaformula2.2}) involving the $c_k$'s yields the term involving the second curly bracket on the rhs of (\ref{zetaformula2.2b}). On the other hand, 
using that
$$\frac{1}{n + \frac{it}{z}} \frac{d}{dz} \frac{z^{-\sigma}}{n + \frac{it}{z}}
= -\frac{z^{1-\sigma } (n \sigma  z+i (\sigma -1) t)}{(n z+i t)^3}$$
and
$$\biggl(\frac{1}{n + \frac{it}{z}} \frac{d}{dz}\biggr)^2 \frac{z^{-\sigma}}{n + \frac{it}{z}}
= \frac{z^{1-\sigma } \left(n^2 \sigma  (\sigma +1) z^2+i n \left(2 \sigma ^2-\sigma
   -2\right) t z-(\sigma -1)^2 t^2\right)}{(n z+i t)^5},$$
long but straightforward computations show that the two terms involving the double sums on the rhs of (\ref{zetaformula2.2}) yield the term involving the first curly bracket on the rhs of (\ref{zetaformula2.2b}). This proves
 (\ref{zetaformula2.2b}).
 \end{proof}

\begin{remark}\upshape
Using the basic identity (\ref{1.2}), it is possible to analyze the case of $0 < \eta < t$ in a way that is similar to the case of $\eta > t$. Actually, this approach is implemented in \cite{FF2018} where the asymptotics of the Hurwitz zeta function to all orders is computed. However, in the next section we find it convenient to analyze the case of $0 < \eta < t$ by implementing the approach introduced by Siegel. We expect that Siegel's approach can be used also to obtain asymptotic formulas for the Hurwitz and Lerch zeta functions. In this regard, we note that an analog of the approximate functional equation (\ref{1.6}) for the Lerch zeta function is derived in \cite{GLS2003} and that a Riemann-Siegel type integral representation for the Lerch zeta function is presented in \cite{BS2012}. The large $t$ asymptotics of the second moment of the Lerch function is derived in \cite{GLS2003b}. Another analog of (\ref{1.6}) for the Lerch function as well as an alternative proof of the large t asymptotics of its second moment can be found in \cite{M2017}.
\end{remark}

\chapter{The Asymptotics of the Riemann Zeta Function for $0 < \eta < t$}\label{sec4}
In this chapter, we consider the asymptotics of $\zeta(s)$ as $t \to \infty$ with $0 < \eta < t$. Theorem \ref{zetath2} and its corollary treat the cases $\epsilon < \eta < \sqrt{t}$ and $2\pi \sqrt{t} < \eta < \frac{2\pi}{\epsilon} t$ under the assumption that $\dist(\eta, 2\pi \Z) > \epsilon$ for some $\epsilon > 0$. The case when $\eta$ is of the same order as $\sqrt{t}$ and the case when $\dist(\eta, 2\pi \Z) \to 0$ are covered by theorem \ref{ZETATH} and its corollary.
Throughout this chapter, we assume that the branch cut for the logarithm runs along the positive real axis. Recall the standing assumption (\ref{basicassumption}) that $\eta, t \notin 2 \pi \Z$.

\begin{theorem}[{\bf The asymptotics to all orders for the case $\epsilon < \eta < \sqrt{t}$}]\label{zetath2}
For every $\epsilon > 0$, there exists a constant $A > 0$ such that
\begin{align}\label{zetaformula2}
& \zeta(s) 
= \sum_{n=1}^{[\frac{t}{\eta}]} \frac{1}{n^s} 
+ \chi(s) \sum_{n=1}^{[\frac{\eta}{2\pi}]} \frac{1}{n^{1-s}}
 	\\ \nonumber
& - \frac{e^{-i\pi s}\Gamma(1-s)}{2\pi i}e^{-([\frac{t}{\eta}]+1)i\eta} e^{\frac{i \pi}{2}(s-1)} \eta^{s-1} e^{\frac{i\pi}{4}} \sum_{k=0}^{[\frac{N-1}{2}]} \frac{\varphi^{(2k)}(0)}{(2k)!} i^k 
 \biggl(\frac{2\eta^2}{t}\biggr)^{k + \frac{1}{2}} \Gamma\biggl(k + \frac{1}{2}\biggr)
	\\ \nonumber
& + 
e^{-i\pi s}\Gamma(1-s) e^{-\frac{\pi t}{2}}\eta^{\sigma-1}
	\\\nonumber
& \quad\times \begin{cases}	
O\biggl(\left(\frac{2N}{t}\right)^{\frac{N}{6}} \frac{\eta}{\sqrt{t}}\biggr), & \epsilon < \eta < t^{\frac{1}{3}} < \infty, \quad 1 \leq N < \frac{At}{\eta^3}, \\
  O\biggl(N e^{- \frac{At}{\eta^2}} + \left(\frac{N\eta^2}{t}\right)^{\frac{N + 1}{2}}\biggr), & t^{\frac{1}{3}} < \eta < \sqrt{t} < \infty, \quad 1 \leq N < \frac{At}{\eta^2}, \end{cases}  
 	\\ \nonumber
& \hspace{4cm} \dist(\eta, 2\pi \Z) > \epsilon, \quad 0 \leq \sigma \leq 1, \quad  t \to \infty,
\end{align}
where the error terms are uniform for all $\eta, \sigma, N$ within the above ranges, and the function $\varphi(z)$ is defined by
$$\varphi(z) = \frac{e^{(s-1)\ln(1 + \frac{z}{i\eta}) - \frac{it}{2\eta^2}z^2 - [\frac{t}{\eta}]z}}{e^z - e^{-i\eta}}.$$

For $N=3$ equation (\ref{zetaformula2}) simplifies to
\begin{align}\nonumber
\zeta(s) 
= &\; \sum_{n=1}^{[\frac{t}{\eta}]} \frac{1}{n^s} 
+ \chi(s) \sum_{n=1}^{[\frac{\eta}{2\pi}]} \frac{1}{n^{1-s}}
 	\\ \nonumber
& - \frac{e^{-i\pi s}\Gamma(1-s)}{2\pi i}e^{-([\frac{t}{\eta}]+1)i\eta} e^{\frac{i \pi}{2}(s-1)} \eta^{s-1} e^{\frac{i\pi}{4}} 
\biggl\{
\frac{\sqrt{2\pi}}{1 - e^{-i\eta}} \frac{\eta}{\sqrt{t}} 
	\\ \nonumber
&+
\frac{1}{(1-e^{-i \eta})^3}\biggl[\frac{(1-e^{-i \eta})^2 \bigl((-\eta [\frac{t}{\eta}] -i (\sigma -1)+t)^2+\sigma -1\bigr)}{\eta^2}
	\\ \nonumber
&\hspace{2.5cm} -\frac{2 (1-e^{-i \eta}) (-\eta  [\frac{t}{\eta}] -i (\sigma -1)+t)}{\eta }+e^{-i \eta }+1\biggr]
   \frac{i\sqrt{\pi}\eta^3}{\sqrt{2}t^{\frac{3}{2}}}
\biggr\}
	\\ \nonumber
& + 
e^{-i\pi s}\Gamma(1-s) e^{-\frac{\pi t}{2}}\eta^{\sigma-1}
	\\\nonumber
& \quad\times \begin{cases}	
O\bigl(\frac{\eta}{t}\bigr), & \epsilon < \eta < t^{\frac{1}{3}} < \infty, \quad 3\eta^3 < At, \\
  O\bigl(e^{- \frac{At}{\eta^2}} + \frac{\eta^4}{t^2}\bigr), & t^{\frac{1}{3}} < \eta < \sqrt{t} < \infty, \quad 3\eta^2 < At, \end{cases}  
 	\\ \label{zetaformula2b}
& \hspace{4cm}  \dist(\eta, 2\pi \Z) > \epsilon, \quad 0 \leq \sigma \leq 1, \quad  t \to \infty,
\end{align}
where the error term is uniform for all $\eta, \sigma$  in the above ranges. 
\end{theorem}
\begin{remark}\upshape
 It will be shown below (see equations (\ref{bnestimate}) and (\ref{bnestimatecase2})) that
$$\frac{\varphi^{(n)}(0)}{n!} = O\biggl(\frac{t^{\frac{n}{3}}}{\eta^n}\biggr) + O(1), \qquad n \geq 0,  \quad
\text{(not uniformly in $n$)}.$$
Thus, the term in (\ref{zetaformula2}) which involves $\varphi^{(n)}(0)$ is of order 
$$
\begin{cases}
O\biggl(e^{-i\pi s}\Gamma(1-s) e^{-\frac{\pi t}{2}} \frac{\eta^\sigma}{\sqrt{t}}\frac{1}{t^{\frac{n}{6}}}\biggr), &
\epsilon < \eta < t^{\frac{1}{3}}, \\
O\biggl(e^{-i\pi s}\Gamma(1-s) e^{-\frac{\pi t}{2}} \eta^{\sigma-1}\left(\frac{2\eta^2}{t}\right)^{\frac{n+1}{2}}\biggr), &
t^{\frac{1}{3}} < \eta < \sqrt{t}, \end{cases} \qquad
\text{(not uniformly in $n$)}.$$
It follows that, as expected, the first term left out of the sum is smaller than the error term.
\end{remark}

\begin{proof}[Proof of theorem \ref{zetath2}.]
We start with the expression for $\zeta(s)$ given on page 82 of Titchmarsh \cite{T1986}:
\begin{align}\label{zetaCjs}
\zeta(s) = \sum_{n=1}^{[\frac{t}{\eta}]} \frac{1}{n^s} 
+ \chi(s) \sum_{n=1}^{[\frac{\eta}{2\pi}]} \frac{1}{n^{1-s}}
+ \frac{e^{-i\pi s}\Gamma(1-s)}{2\pi i} \sum_{j=1}^4 I_j, \qquad 0 < \eta < t,
\end{align}
where 
$$I_j =  \int_{C_j} \frac{w^{s-1} e^{-[\frac{t}{\eta}]w}}{e^w - 1} dw,$$
the contours $C_1, C_2, C_3, C_4$ are the straight lines joining $\infty, c\eta + i\eta(1+c), -c\eta + i\eta(1-c), -c\eta - (2[\frac{\eta}{2\pi}] + 1)\pi i, \infty$, where $c$ is an absolute constant, $0 < c \leq 1/2$, see figure \ref{Cjs.pdf}. We henceforth let $c = 2^{-3/2}$ so that the total length of the contour  $C_2$ is $\eta$. 
The complex powers in (\ref{zetaCjs}) are defined using a branch cut that runs along the {\it positive} real axis, i.e., $w^{s-1} = e^{(s-1)(\ln |w| + i \arg w)}$ with $\arg w \in (0, 2\pi)$. We will adopt this choice of branch cut along the positive real axis in the remainder of this chapter.
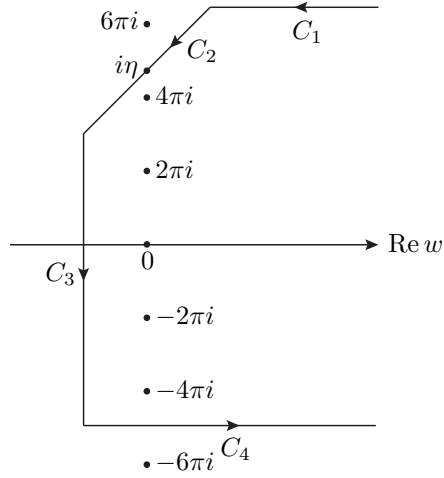
\begin{figure}
\bigskip
 \begin{overpic}[width=.4\textwidth]{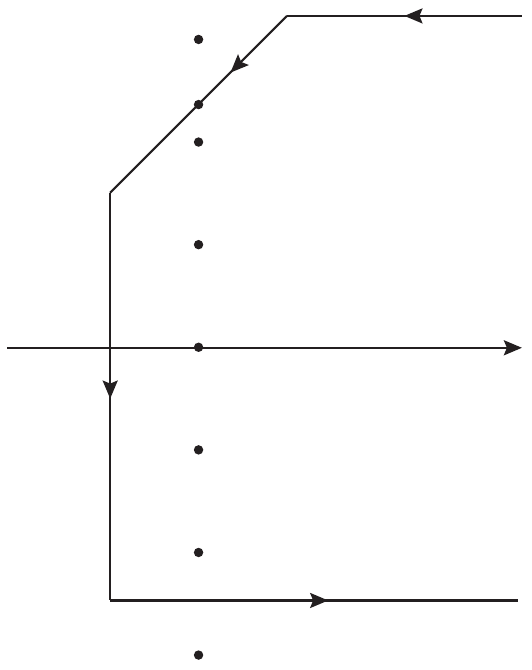}
 \put(80,46.5){$\re w$}
 \put(60, 91){$C_1$}
 \put(38, 87){$C_2$}
 \put(8.5, 40.5){$C_3$}
 \put(45, 4){$C_4$}
 \put(23, 84){$i\eta$}
 \put(31.5, 1){$-6\pi i$}
 \put(31.5, 16){$-4\pi i$}
 \put(31.5, 31.5){$-2\pi i$}
 \put(28.5, 43){$0$}
 \put(31.5, 62){$2\pi i$}
 \put(31.5, 77.5){$4\pi i$}
 \put(20, 93){$6\pi i$}
   \end{overpic}
   \caption{The critical point $i\eta$ and the integration contours $\{C_j\}_1^4$ in the complex $w$-plane.}\label{Cjs.pdf}
\end{figure}
The contributions to $\zeta(s)$ from $I_1$, $I_3$, and $I_4$ are exponentially small by the arguments on pages 82-83 of \cite{T1986} (in these arguments, it is assumed that $x = t/\eta \leq y = \eta/2\pi$, however this assumption is not used for the derivations of these estimates and therefore we do not need to make it here.)

Let $\epsilon > 0$ be given and suppose that $4\epsilon < \eta < \sqrt{t}$, $0 \leq \sigma \leq 1$, $m = [t/\eta]$, $\dist(\eta, 2\pi \Z) > 4\epsilon$, and $N \geq 1$. All error terms of the form $O(\cdot)$ will be uniform with respect to $\eta, \sigma, N$ (but not with respect to $\epsilon$). In order to analyze the integral $I_2$, we write
$$I_2 = e^{-(m+1)i\eta} (i\eta)^{s-1} \int_{C_2} e^{\frac{it}{2\eta^2}(w - i\eta)^2} \varphi(w - i\eta) dw,$$
where
$$\varphi(z) = \frac{e^{(s-1)\ln(1 + \frac{z}{i\eta}) - \frac{t}{\eta} z - \frac{it}{2\eta^2}z^2 + (\frac{t}{\eta} - [\frac{t}{\eta}])z}}{e^z - e^{-i\eta}}.$$
Defining $\phi(z)$ by
\begin{align}\label{phidef}
\phi(z) = e^{(s-1)\ln\bigl(1 + \frac{z}{i\eta}\bigr) - \frac{t}{\eta} z - \frac{it}{2\eta^2}z^2},
\end{align}
we have
$$\varphi(z) = \frac{\phi(z) e^{(\frac{t}{\eta} - [\frac{t}{\eta}])z}}{e^z - e^{-i\eta}}.$$
We split the contour $C_2$ as follows:
$$C_2 = C_2^\epsilon \cup C_2^r,$$
where $C_2^\epsilon$ denotes the segment of $C_2$ of length $2\epsilon$ which consists of the points within a distance $\epsilon$ from $i\eta$. We write
$$I_2 = I_2^\epsilon + I_2^r,$$
where
\begin{align*}
& I_2^\epsilon = e^{-(m+1)i\eta} (i\eta)^{s-1} \int_{C_2^\epsilon} e^{\frac{it}{2\eta^2}(w - i\eta)^2} \varphi(w - i\eta) dw
\end{align*}
and
\begin{align*}
& I_2^r = e^{-(m+1)i\eta} (i\eta)^{s-1} \int_{C_2^r} e^{\frac{it}{2\eta^2}(w - i\eta)^2} \varphi(w - i\eta) dw.
\end{align*}

We claim that $I_2^r$ can be estimated as follows:
\begin{align}\label{I2restimate}
  I_2^r = O\biggl(e^{-\frac{\pi t}{2}} \frac{\eta^{\sigma}}{\sqrt{t}} e^{-\frac{t \epsilon^2}{12 \eta^2}}\biggr).
\end{align}
Indeed, the change of variables $w = i\eta + \lambda e^{\frac{i\pi}{4}}$ gives
$$I_2^r = O\biggl(e^{-\frac{\pi t}{2}} \eta^{\sigma -1} \biggl(\int_{-\frac{\eta}{2}}^{-\epsilon} + \int_\epsilon^{\frac{\eta}{2}}
\biggr) e^{-\frac{t}{2\eta^2} \lambda^2} |\varphi(\lambda e^{\frac{i\pi}{4}})| d\lambda\biggr).$$
We have
\begin{align}\label{etetaAestimate}
  \biggl|\frac{e^{(\frac{t}{\eta} - [\frac{t}{\eta}])z}}{e^z - e^{-i\eta}}\biggr| = O(1), \qquad
z = \lambda e^{\frac{i\pi}{4}}, \quad \epsilon < |\lambda | < \frac{\eta}{2}.
\end{align}
Moreover, the definition of $\phi(z)$ implies
\begin{align}\nonumber
  \ln{\phi(z)} & = (s-1)\ln\left(1 + \frac{z}{i\eta}\right) - \frac{t}{\eta}z - \frac{it}{2\eta^2}z^2
  	\\ \label{lnphizestimate}
& = (\sigma - 1)\ln\left(1 + \frac{z}{i\eta}\right)	- \frac{it}{\eta^2}z^2 \sum_{k=1}^\infty \frac{(-1)^{k-1}}{k+2} \left(\frac{z}{i\eta}\right)^k.
\end{align}
Hence,
\begin{align}\label{relogphiw}
  \re{\ln{\phi(z)}} \leq |\sigma -1| \ln\frac{3}{2} + \frac{t}{\eta^2}|z|^2\frac{1}{3} \frac{|z|}{\eta} \frac{1}{1- 1/2} = |\sigma -1| \ln\frac{3}{2} + \frac{2t|z|^3}{3\eta^3}, \qquad |z| \leq \frac{\eta}{2},
\end{align}  
and so
\begin{align}\label{philambdaeipi4}
|\phi(\lambda e^{\frac{i\pi}{4}})| < e^{|\sigma -1| \ln{\frac{3}{2}} + \frac{2t }{3\eta^3}|\lambda|^3}, \qquad
\epsilon < |\lambda | < \frac{\eta}{2}.
\end{align}
Equations (\ref{etetaAestimate}) and (\ref{philambdaeipi4}) imply that
$$|\varphi(\lambda e^{\frac{i\pi}{4}})| = O\Bigl(e^{\frac{2t}{3\eta^3}|\lambda|^3}\Bigr), \qquad \epsilon < |\lambda | < \frac{\eta}{2}.$$
This yields
\begin{align*}
I_2^r & = O\biggl(e^{-\frac{\pi t}{2}} \eta^{\sigma -1}
 \int_\epsilon^{\frac{\eta}{2}}
e^{-\frac{t}{2\eta^2}\lambda^2 + \frac{2t}{3\eta^3}|\lambda|^3} d\lambda\biggr)
	\\
& =
O\biggl(e^{-\frac{\pi t}{2}} \eta^{\sigma -1}
 \int_\epsilon^{\frac{\eta}{2}}
e^{-\frac{t}{2\eta^2}\lambda^2 + \frac{t}{3\eta^2} \lambda^2} d\lambda\biggr)
= O\biggl(e^{-\frac{\pi t}{2}} \eta^{\sigma -1}
 \int_\epsilon^{\frac{\eta}{2}}
e^{-\frac{t}{6\eta^2}\lambda^2 } d\lambda\biggr).
\end{align*}
Splitting the integrand as
$$e^{-\frac{t}{6\eta^2}\lambda^2} 
= e^{-\frac{t}{12\eta^2}\lambda^2} \times e^{-\frac{t}{12\eta^2}\lambda^2}$$
and noting that 
$$\int_\epsilon^{\frac{\eta}{2}} e^{-\frac{t}{12\eta^2}\lambda^2} d\lambda
\leq \int_0^\infty e^{-\frac{t}{12\eta^2}\lambda^2} d\lambda
= \sqrt{3\pi}\frac{ \eta}{\sqrt{t}},$$
we find
\begin{align*}
I_2^r & =
O\biggl(e^{-\frac{\pi t}{2}} \eta^{\sigma -1} e^{-\frac{t}{12 \eta^2}\epsilon^2} \frac{\eta}{\sqrt{t}} \biggr).
\end{align*}
This proves (\ref{I2restimate}).

We now consider $I_2^\epsilon$. 
In view of the assumption $\dist(\eta, 2\pi \Z) > 4\epsilon$, we have
\begin{align}\label{expquotientepsilon}
\biggl|\frac{e^{(\frac{t}{\eta} - [\frac{t}{\eta}])z}}{e^z - e^{-i\eta}}\biggr| = O(1), \qquad
|z| < 2\epsilon.
\end{align}
In particular, $\varphi(z)$ is analytic for $|z| < 2\epsilon$. Thus we can write
$$\varphi(z) = \sum_{n=0}^{\infty} b_n z^n = \sum_{n=0}^{N-1} b_n z^n + s_N(z), \qquad |z| < 2 \epsilon,$$
where
\begin{align}\label{sNexpression}
s_N(z) = \frac{z^N}{2\pi i}\int_{\Gamma} \frac{\varphi(w) dw}{w^N(w-z)},
\end{align}
and $\Gamma$ is a counterclockwise contour which encircles $0$ and $z$ but none of the poles of $\varphi$. 
We claim that\footnote{Here and below $A > 0$ denotes a generic constant independent of $t, \eta, \sigma, N$, which may change within a computation.}
\begin{align} \nonumber
I_2^\epsilon = &\; e^{-(m+1)i\eta} (i\eta)^{s-1} \int_{C_2^\epsilon} e^{\frac{it}{2\eta^2}(w - i\eta)^2} \sum_{n=0}^{N-1} b_n (w - i\eta)^n dw
	\\ \label{I2epsilonmidstep}
&+  \begin{cases}	
O\biggl( \eta^{\sigma-1}e^{-\frac{\pi t}{2}} \left(\frac{2N}{t}\right)^{\frac{N}{6}} \frac{\eta}{\sqrt{t}}\biggr), & 1 \leq N < \frac{At}{\eta^3}, \\
O\biggl(e^{-\frac{\pi t}{2}}\eta^{\sigma-1} \left(\frac{N\eta^2}{t}\right)^{\frac{N + 1}{2}}\biggr), &   t^{\frac{1}{3}} < \eta < \sqrt{t}.\end{cases}  
\end{align}
In order to establish (\ref{I2epsilonmidstep}), we need to estimate the error term
$$e^{-(m+1)i\eta} (i\eta)^{s-1} \int_{C_2^\epsilon} e^{\frac{it}{2\eta^2}(w - i\eta)^2} s_N(w - i\eta) dw.$$

Let us first consider the case $1 \leq N < \frac{At}{\eta^3}$. 
Let $|z| < \frac{40}{21}\epsilon$ (so that $\frac{21}{20}|z| < 2\epsilon$) and let $\Gamma$ in (\ref{sNexpression}) be a circle with center $w =0$ and radius $\rho_N$, where
$$\frac{21}{20} |z| \leq \rho_N < 2\epsilon.$$
Then, since $\varphi(z)$  is analytic in the disk $|z| < 2\epsilon$, equations (\ref{sNexpression}), (\ref{relogphiw}), and (\ref{expquotientepsilon}) yield
$$s_N(z) = O\left(\frac{2\pi \rho_N}{\rho_N - |z|}|z|^N \rho_N^{-N}e^{\frac{2t\rho_N^3}{3\eta^3}}\right)
= O\left(|z|^N \rho_N^{-N}e^{\frac{2t\rho_N^3}{3\eta^3}}\right), \qquad |z| < \frac{40}{21}\epsilon.$$
The function $\rho^{-N}e^{\frac{2t\rho^3}{3\eta^3}}$ has the minimum $(\frac{2et}{N\eta^3})^{N/3}$ for $\rho = (\frac{N}{2t})^{1/3}\eta$; $\rho_N$ can have this value if
$$\frac{21}{20}|z| \leq \biggl(\frac{N}{2t}\biggr)^{1/3}\eta < 2\epsilon.$$
The assumption $1 \leq N < \frac{At}{\eta^3}$ implies that for $A$ sufficiently small, we have $(\frac{N}{2t})^{1/3}\eta < 2\epsilon$.
Hence, letting $\rho_N = (\frac{N}{2t})^{1/3}\eta$, we find
\begin{align}\label{sNestimate}
  s_N(z) = 
  O\left(|z|^N \biggl(\frac{2et}{N\eta^3}\biggr)^{\frac{N}{3}}\right), 
  \qquad N < \frac{At}{\eta^3}, \quad |z| \leq \frac{20}{21}\biggl(\frac{N}{2t}\biggr)^{\frac{1}{3}} \eta < \frac{40}{21}\epsilon,
\end{align}
For $|z| < \frac{40}{21}\epsilon$ we can also take $\rho_N = \frac{21}{20}|z|$, which yields
\begin{align}\label{sNestimate2}
s_N(z) 
= O\left(\Bigl(\frac{20}{21}\Bigr)^N e^{\frac{2}{3}\frac{t}{\eta^3} (\frac{21}{20}|z|)^3}\right) 
= O\left(e^{\frac{4}{5}\frac{t}{\eta^3} |z|^3}\right) 
= O\left(e^{\frac{2}{5}\frac{t}{\eta^2} |z|^2}\right), \qquad |z| < \frac{40}{21}\epsilon,
\end{align}
where we have used that $\frac{|z|}{\eta} < \frac{|z|}{4\epsilon} < \frac{1}{2}$ in the last step.
Using (\ref{sNestimate}) and (\ref{sNestimate2}), we estimate
\begin{align}\nonumber
& e^{-(m+1)i\eta} (i\eta)^{s-1} \int_{C_2^\epsilon} e^{\frac{it}{2\eta^2}(w - i\eta)^2} s_N(w - i\eta) dw
	\\ \nonumber
& = O\biggl(\eta^{\sigma -1}e^{-\frac{\pi t}{2}}\biggl\{ \int_0^{A(\frac{N}{t})^{\frac{1}{3}} \eta} e^{-\frac{t}{2\eta^2}\lambda^2} \lambda^N \left(\frac{2et}{N\eta^3}\right)^{\frac{N}{3}} d\lambda
	\\\nonumber
& \hspace{1cm} + \int_{A(\frac{N}{t})^{\frac{1}{3}} \eta}^{\epsilon} e^{-\frac{t}{2\eta^2}\lambda^2 + \frac{2}{5}\frac{t}{\eta^2}\lambda^2} d\lambda \biggr\}
\biggr)
	\\\nonumber
& = O\biggl(\eta^{\sigma -1} e^{-\frac{\pi t}{2}} \biggl\{\left(\frac{2et}{N\eta^3}\right)^{\frac{N}{3}} 2^{\frac{N -1}{2}} \left(\frac{\eta^2}{t}\right)^{\frac{N+1}{2}} \Gamma\left(\frac{N +1}{2}\right)
+ \int_{A(\frac{N}{t})^{\frac{1}{3}}\eta}^{\epsilon} e^{-\frac{t}{10\eta^2}\lambda^2} d\lambda \biggr\} \biggr)
	\\\label{C2epsestimate}
& = O\biggl(\eta^{\sigma -1} e^{-\frac{\pi t}{2}} \left(\frac{2et}{N\eta^3}\right)^{\frac{N}{3}} 2^{\frac{N -1}{2}} \left(\frac{\eta^2}{t}\right)^{\frac{N+1}{2}} \Gamma\left(\frac{N +1}{2}\right) \biggr), \qquad 1\leq N < \frac{At}{\eta^3},
\end{align}
where we have used the following estimate to find the last equality:
\begin{align*}
 \int_{A(\frac{N}{t})^{\frac{1}{3}}\eta}^{\epsilon} e^{-\frac{t}{10\eta^2}\lambda^2} d\lambda
& = O\biggl(e^{-\frac{t}{20\eta^2}(A(\frac{N}{t})^{\frac{1}{3}}\eta)^2}
  \int_{A(\frac{N}{t})^{\frac{1}{3}}\eta}^{\infty} e^{-\frac{t}{20\eta^2}\lambda^2} d\lambda\biggr)
 	\\ \nonumber
&  = O\biggl(e^{-\frac{A^2 t^{\frac{1}{3}}}{20}N^{\frac{2}{3}}} \frac{\eta}{\sqrt{t}}\biggr)
= O(e^{-\frac{A^2}{20} t^{\frac{1}{3}}}).
\end{align*}
Using the expansion
\begin{align}\label{Gammaexpansion}
\Gamma(x) = e^{x(\ln(x) - 1)}\biggl(\sqrt{\frac{2\pi}{x}} + O\biggl(\frac{1}{x^{3/2}}\biggr)\biggr), \qquad x \to \infty,
\end{align}
and the inequality $\frac{4}{e} < 2$,
it follows that the rhs of (\ref{C2epsestimate}) is
\begin{align}\label{Oetasigmaminusone}
O\biggl(\eta^{\sigma -1} e^{-\frac{\pi t}{2}} \left(\frac{2N}{t}\right)^{\frac{N}{6}} \frac{\eta}{\sqrt{t}}\biggr).
\end{align}
This proves (\ref{I2epsilonmidstep}) in the case when $1\leq N < \frac{At}{\eta^3}$.

We now consider the case when $t^{\frac{1}{3}} < \eta < \sqrt{t}$.
Using (\ref{expquotientepsilon}) and (\ref{relogphiw}) in the representation (\ref{sNexpression}) with $\Gamma$ a circle with center $w=0$ and radius $2^{1/3}\epsilon$, we find that 
\begin{align}\label{sNestimate3}
s_N(z) = O\Big(|z|^N e^{\frac{4t\epsilon^3}{3\eta^3}}\Big), \qquad |z| < \epsilon.
\end{align}
Now $t^{\frac{1}{3}} < \eta$ implies that $e^{\frac{4t\epsilon^3}{3\eta^3}} = O(1)$, so we can estimate
\begin{align}\nonumber
e^{-(m+1)i\eta} & (i\eta)^{s-1} \int_{C_2^\epsilon} e^{\frac{it}{2\eta^2}(w - i\eta)^2} s_N(w - i\eta) dw
	\\ \nonumber
& = O\biggl(\eta^{\sigma -1}e^{-\frac{\pi t}{2}} \int_0^{\epsilon} e^{-\frac{t}{2\eta^2}\lambda^2} \lambda^Nd\lambda \biggr)
	\\\nonumber
& = O\biggl(\eta^{\sigma -1} e^{-\frac{\pi t}{2}} 2^{\frac{N -1}{2}} \left(\frac{\eta^2}{t}\right)^{\frac{N+1}{2}} \Gamma\left(\frac{N +1}{2}\right) \biggr)
	\\ \nonumber
& = O\biggl(\eta^{\sigma -1} e^{-\frac{\pi t}{2}} \left(\frac{\eta^2}{t}\right)^{\frac{N+1}{2}} N^{\frac{N}{2}} \biggr), \qquad
t^{\frac{1}{3}} < \eta < \sqrt{t}.
\end{align}
This completes the proof of (\ref{I2epsilonmidstep}).

We next claim that, up to a small error term, the contour $C_2^\epsilon$ in the integral in (\ref{I2epsilonmidstep}) can be replaced by the infinite line $C_2'$, where $C_2'$ denotes the infinite straight line of which $C_2$ is a part. More precisely, we claim that there exists an $A >0$ such that
\begin{align}\label{C2epsilonC2prime}
e^{-(m+1)i\eta} (i\eta)^{s-1} \int_{C_2'\backslash C_2^\epsilon}& e^{\frac{it}{2\eta^2}(w - i\eta)^2} \sum_{n=0}^{N-1} b_n (w - i\eta)^n dw
	\\ \nonumber
& = \begin{cases}
  O\Bigl( e^{-\frac{\pi t}{2} - \frac{At}{\eta^2} }\Bigr), & 4\epsilon < \eta < t^{\frac{1}{3}}, \quad 1 \leq N < \frac{At}{\eta^3} , \\
  O\Bigl(\eta^{\sigma -1}e^{-\frac{\pi t}{2} - \frac{At}{\eta^2}} N\Bigr), & t^{\frac{1}{3}} < \eta < \sqrt{t}, \quad 1 \leq N < \frac{At}{\eta^2}.
\end{cases}  
 \end{align}
In order to prove (\ref{C2epsilonC2prime}), we note that the coefficient of $b_n$ on the lhs of (\ref{C2epsilonC2prime}) is
\begin{align}\label{bncoefficientbigOh}
O\biggl(\eta^{\sigma -1}e^{-\frac{\pi t}{2}} \int_\epsilon^\infty e^{-\frac{t}{2\eta^2} \lambda^2} \lambda^n d\lambda\biggr).
\end{align}
We write the integrand as
$$e^{-\frac{t}{4\eta^2} \lambda^2} \lambda^n \times e^{-\frac{t}{4\eta^2} \lambda^2}.$$
The first factor is steadily decreasing for $\lambda > \sqrt{\frac{2n}{t}} \eta$, and so it decreases throughout the interval of integration provided that $n < N < \frac{t\epsilon^2}{2\eta^2}$. 
The term in (\ref{bncoefficientbigOh}) is then
\begin{align}\label{bncoefficientbigOh2}
O\biggl(\eta^{\sigma -1}e^{-\frac{\pi t}{2}} 
e^{-\frac{t}{4\eta^2} \epsilon^2} \epsilon^n
\int_\epsilon^\infty e^{-\frac{t}{4\eta^2} \lambda^2} d\lambda\biggr)
= 
O\biggl(\eta^{\sigma -1}e^{-\frac{\pi t}{2}} 
e^{-\frac{t}{4\eta^2} \epsilon^2} \epsilon^n \frac{\eta}{\sqrt{t}}\biggr).
\end{align}

Let us assume that $4\epsilon < \eta < t^{\frac{1}{3}}$ and $1 \leq N < At/\eta^3$. In this case, choosing $|z| < \frac{20}{21}\bigl(\frac{N}{2t}\bigr)^{\frac{1}{3}} \eta$, equation (\ref{sNestimate}) yields
\begin{align}\label{bnestimate}
b_n = (s_n(z) - s_{n+1}(z))z^{-n} 
= O\left(\left(\frac{2et}{n\eta^3}\right)^{\frac{n}{3}}\right), \qquad N < \frac{At}{\eta^3},\quad 1 \leq n \leq N-1.
\end{align}
Multiplying the rhs of (\ref{bncoefficientbigOh2}) by  $b_n$ and summing from $0$ to $N-1$, we find that the total error is
$$O\left(  \eta^{\sigma -1}e^{-\frac{\pi t}{2}} e^{-\frac{t}{4\eta^2} \epsilon^2} \frac{\eta}{\sqrt{t}}
\biggl(b_0 +  \sum_{n=1}^{N-1}\epsilon^n \left(\frac{2et}{n\eta^3}\right)^{\frac{n}{3}}\biggr)\right).$$
Now the function $(\frac{t}{n\eta^3})^{\frac{n}{3}}$ increases steadily up to $n = \frac{t}{\eta^3 e}$, so that if $n < A \frac{t}{\eta^3}$, where $A < 1/e$, it is of order
$$O\Bigl(e^{\frac{1}{3}\frac{At}{\eta^3}\ln\frac{1}{A}}\Bigr).$$
Hence, if $\epsilon^3 < \frac{1}{2e}$, the total error is
\begin{align*}
& O\left(\frac{\eta^\sigma}{\sqrt{t}} e^{-\frac{\pi t}{2}} e^{-\frac{t}{4\eta^2} \epsilon^2}\left(b_0 +  e^{\frac{1}{3}\frac{At}{\eta^3}\ln\frac{1}{A}} \sum_{n=1}^{N-1} \big(\epsilon(2e)^{1/3}\big)^n\right)\right)
	\\
& = O\left(\frac{\eta^\sigma}{\sqrt{t}} e^{-\frac{\pi t}{2}} e^{-\frac{t}{\eta^2}(\frac{\epsilon^2}{4} - \frac{1}{3}\frac{A}{\eta}\ln\frac{1}{A})}\right).
\end{align*}
Since $4\epsilon < \eta$ and $A \ln\frac{1}{A}  \to 0$ as $A \to 0$, we can choose $A > 0$ such that $\frac{\epsilon^2}{4} - \frac{1}{3}\frac{A}{\eta}\ln\frac{1}{A} > 0$. 
This proves (\ref{C2epsilonC2prime}) in the case when $4\epsilon < \eta < t^{\frac{1}{3}}$ and $1 \leq N < \frac{At}{\eta^3}$.

In the case when $t^{\frac{1}{3}} < \eta < \sqrt{t}$ and $1 \leq N < \frac{At}{\eta^2}$, we instead use (\ref{sNestimate3}) to find
\begin{align}\label{bnestimatecase2} 
b_n = (s_n(z) - s_{n+1}(z))z^{-n} 
= O\Bigl(e^{\frac{4t\epsilon^3}{3\eta^3}}\Bigr) = O(1), \qquad n \geq 0.
\end{align}
Then, the estimate (\ref{bncoefficientbigOh2}) implies that the total error is
\begin{align*}
O\left(\eta^{\sigma -1}e^{-\frac{\pi t}{2} - \frac{t}{4\eta^2}\epsilon^2} N\right)
	 = O\left(\eta^{\sigma -1}e^{-\frac{\pi t}{2} - \frac{At}{\eta^2}} N\right).
\end{align*}
This completes the proof of (\ref{C2epsilonC2prime}).

We finally analyze the sum
\begin{align*}
& e^{-(m+1)i\eta} (i\eta)^{s-1} \int_{C_2'} e^{\frac{it}{2\eta^2}(w - i\eta)^2} \sum_{n=0}^{N-1} b_n (w - i\eta)^n dw
	\\
& = -e^{-(m+1)i\eta} (i\eta)^{s-1}  \sum_{n=0}^{N-1} b_n e^{\frac{\pi i}{4}(n+1)}
\int_{-\infty}^\infty e^{-\frac{t}{2\eta^2} \lambda^2} \lambda^n d\lambda
	\\
& = -e^{-(m+1)i\eta} (i\eta)^{s-1}  \sum_{n=0}^{N-1} b_n e^{\frac{\pi i}{4}(n+1)}
2^{\frac{n-1}{2}} (1 + (-1)^n) \biggl(\frac{\eta^2}{t}\biggr)^{\frac{n+1}{2}} \Gamma\biggl(\frac{n+1}{2}\biggr)
	\\
& = -e^{-(m+1)i\eta} (i\eta)^{s-1} e^{\frac{i\pi}{4}} \sum_{k=0}^{[\frac{N-1}{2}]} b_{2k} i^k 
 \biggl(\frac{2\eta^2}{t}\biggr)^{k + \frac{1}{2}} \Gamma\biggl(k + \frac{1}{2}\biggr).
\end{align*}
Together with equations (\ref{zetaCjs}), (\ref{I2restimate}), (\ref{I2epsilonmidstep}), and (\ref{C2epsilonC2prime}), this yields equation (\ref{zetaformula2}) with $\epsilon$ replaced with $4 \epsilon$. Since $\epsilon > 0$ was arbitrary, the proof is complete.
\end{proof}

\begin{corollary}[{\bf The asymptotics to all orders for the case $2\pi \sqrt{t} < \eta < \frac{2\pi}{\epsilon} t$}]\label{zetacor2}
For every $\epsilon > 0$, there exists a constant $A > 0$ such that
\begin{align} \nonumber
\zeta(&s) 
= \sum_{n=1}^{[\frac{t}{\eta}]} \frac{1}{n^s} 
+ \chi(s) \sum_{n=1}^{[\frac{\eta}{2\pi}]} \frac{1}{n^{1-s}}
 	\\ \nonumber
& + \chi(s) \frac{e^{i\pi(1-s)}\Gamma(s)}{2\pi i}e^{([\frac{\eta}{2\pi}]+1)i\frac{2\pi t}{\eta}} e^{\frac{i \pi s}{2}} \biggl(\frac{2\pi t}{\eta}\biggr)^{-s} e^{-\frac{i\pi}{4}}
	\\ \nonumber
&\times  \sum_{k=0}^{[\frac{N-1}{2}]} \frac{\psi^{(2k)}(0)}{(2k)!} (-i)^k 
 \biggl(\frac{8\pi^2 t}{\eta^2}\biggr)^{k + \frac{1}{2}} \Gamma\biggl(k + \frac{1}{2}\biggr)
	\\ \nonumber
& + 
\chi(s) e^{i\pi(1-s)}\Gamma(s) e^{-\frac{\pi t}{2}}\biggl(\frac{\eta}{2\pi t}\biggr)^{\sigma}
	\\\nonumber
& \quad\times \begin{cases}	
O\biggl(\left(\frac{2N}{t}\right)^{\frac{N}{6}} \frac{2\pi t}{\eta\sqrt{t}}\biggr), & 2\pi t^{\frac{2}{3}} < \eta < \frac{2\pi}{\epsilon} t < \infty, \quad 1 \leq N < \frac{A\eta^3}{t^2}, \\
  O\biggl(N e^{- \frac{A\eta^2}{4\pi^2 t}} + \left(\frac{N 4\pi^2 t}{\eta^2}\right)^{\frac{N + 1}{2}}\biggr), 
  & 2\pi \sqrt{t} < \eta < 2\pi t^{\frac{2}{3}} < \infty, \quad 1 \leq N < \frac{A\eta^2}{t}, \end{cases}  
 	\\ \label{zetacor2formula}
& \hspace{4cm} \dist\biggl(\frac{2\pi t}{\eta}, 2\pi \Z\biggr) > \epsilon, \quad 0 \leq \sigma \leq 1, \quad  t \to \infty, 
\end{align}
where the error terms are uniform for all $\eta, \sigma, N$ within the above ranges, and the function $\psi(z)$ is defined by
$$\psi(z) = \frac{e^{-s\ln(1 + \frac{i\eta z}{2\pi t}) + \frac{i\eta^2}{8\pi^2 t}z^2 - [\frac{\eta}{2\pi}]z}}{e^z - e^{\frac{2\pi i t}{\eta}}}.$$
\end{corollary}
\begin{proof}
We replace $\sigma$ by $1- \sigma$ in (\ref{zetaformula2}) and take the complex conjugate of both sides. We then multiply the resulting equation by $\chi(s)$ and use the identities
\begin{align}\nonumber
 & \overline{\zeta(\bar{s})} = \zeta(s), \qquad
  \overline{\Gamma(\bar{s})} = \Gamma(s), \qquad
  \overline{\chi(\bar{s})} = \chi(s), 
  	\\ \label{baridentities}
& \chi(s)\chi(1-s) = 1, \qquad \chi(s)\zeta(1-s) = \zeta(s).	
\end{align}
This yields the following equation
\begin{align}\nonumber
\zeta(s) 
=&\;  \chi(s) \sum_{n=1}^{[\frac{t}{\eta}]} \frac{1}{n^{1-s}} 
+ \sum_{n=1}^{[\frac{\eta}{2\pi}]} \frac{1}{n^{s}}
 	\\ \nonumber
& + \chi(s)\frac{e^{i\pi(1-s)}\Gamma(s)}{2\pi i}e^{([\frac{t}{\eta}]+1)i\eta} e^{\frac{i \pi s}{2}} \eta^{-s} e^{-\frac{i\pi}{4}}
	\\ \nonumber
& \quad \times \sum_{k=0}^{[\frac{N-1}{2}]} \biggl(\frac{d^{2k}}{dz^{2k}}\bigg|_{z=0} \frac{e^{-s\ln(1 - \frac{z}{i\eta}) + \frac{it}{2\eta^2}z^2 - [\frac{t}{\eta}]z}}{e^z - e^{i\eta}}\biggr) \frac{1}{(2k)!} (-i)^k 
 \biggl(\frac{2\eta^2}{t}\biggr)^{k + \frac{1}{2}} \Gamma\biggl(k + \frac{1}{2}\biggr)
	\\ \nonumber
& + \chi(s)
e^{i\pi(1-s)}\Gamma(s) e^{-\frac{\pi t}{2}}\eta^{-\sigma}
	\\\nonumber
& \quad\times \begin{cases}	
O\biggl(\left(\frac{2N}{t}\right)^{\frac{N}{6}} \frac{\eta}{\sqrt{t}}\biggr), & \epsilon < \eta < t^{\frac{1}{3}} < \infty, \quad 1 \leq N < \frac{At}{\eta^3}, \\
  O\biggl(N e^{- \frac{At}{\eta^2}} + \left(\frac{N\eta^2}{t}\right)^{\frac{N + 1}{2}}\biggr), & t^{\frac{1}{3}} < \eta < \sqrt{t} < \infty, \quad 1 \leq N < \frac{At}{\eta^2}, \end{cases}  
 	\\ \nonumber
& \hspace{4cm} \dist(\eta, 2\pi \Z) > \epsilon, \quad 0 \leq \sigma \leq 1, \quad  t \to \infty,
\end{align}
Replacing $\eta$ by $\frac{2\pi t}{\eta}$, we find (\ref{zetacor2formula}).
\end{proof}

\begin{theorem}[{\bf The asymptotics to all orders for the case $\epsilon \sqrt{t} < \eta < t$}]\label{ZETATH}
For every $\epsilon > 0$, there exists a constant $A>0$ such that
\begin{align}\label{zetaformula}
\zeta(s) 
=& \sum_{n=1}^{[\frac{t}{\eta}]} \frac{1}{n^s} 
+ \chi(s) \sum_{n=1}^{[\frac{\eta}{2\pi}]} \frac{1}{n^{1-s}}
 	\\ \nonumber
& + e^{-i\pi s}\Gamma(1-s)\biggl\{ e^{\frac{i\pi(s-1)}{2}}\eta^{s-1}  e^{\frac{2t}{\eta}[\frac{\eta}{2\pi}]\pi i - it 
- \frac{it}{2\eta^2}(2[\frac{\eta}{2\pi}]\pi - \eta)^2} S_N(s, \eta)
	\\ \nonumber
 & \hspace{4cm} + O\biggl( e^{-\frac{\pi t}{2}} \left(\frac{3N}{t}\right)^{\frac{N}{6}} \frac{\eta^\sigma}{\sqrt{t}}\biggr)\biggr\},
 	\\ \nonumber
& \qquad \epsilon \sqrt{t} < \eta < t , \quad 0 \leq \sigma \leq 1, \quad 1 \leq N < At, \quad  t \to \infty,
\end{align}
where the error term is uniform for all $\eta, \sigma, N$ in the above ranges and the function $S_N(s, \eta)$ is defined by
\begin{align}\label{SNdef}
S_N(s, \eta)= \sum_{n=0}^{N-1} a_n \sum_{k = 0}^{n} {n \choose{k}}& \biggl(2\Bigl[\frac{\eta}{2\pi}\Bigr]\pi i - i\eta\biggr)^k
	\\ \nonumber
&\times  \partial_2^{n-k} \Phi \left(-\frac{2\pi t}{\eta^2}, \frac{2t}{\eta} - \frac{2\pi t}{\eta^2}\Bigl[\frac{\eta}{2\pi}\Bigr] - \Bigl[\frac{t}{\eta}\Bigr] - \frac{1}{2}\right)
\end{align}
with $\partial_2$ denoting differentiation with respect to the second argument. The coefficients $a_n$ are determined by the recurrence formula
\begin{align}\label{anrecurrence}
i\eta(n+1)a_{n+1} = (\sigma - n - 1)a_n - \frac{it}{\eta^2}a_{n-2}, \qquad n = 0,1,2, \dots,
\end{align}
together with the initial conditions $a_{-2} = a_{-1} = 0$ and $a_0 = 1$, and the function $\Phi(\tau, u)$ is defined by
\begin{align}\label{Phidef}
\Phi(\tau, u) = \int_{0 \nwarrow 1} \frac{e^{\pi i \tau x^2 + 2\pi i u x}}{e^{\pi i x} - e^{-\pi i x}} dx, \qquad \tau < 0, \quad u \in \C,
\end{align}
with the contour $0 \nwarrow 1$ denoting a straight line parallel to $e^{3\pi i/4}$ which crosses the real axis between $0$ and $1$.\footnote{In the particular case of $\eta = \sqrt{2\pi t}$, an analog of (\ref{zetaformula}) is stated without proof in the recent paper \cite[Eq. (1.4)]{M2017}.}

In the particular case when $p,q > 0$ are integers, we have
\begin{align}\label{Phiformula}
\Phi\left(-\frac{p}{q}, u\right) = &\; \frac{1}{1 - (-1)^q e^{-\pi i q p - 2\pi i q u}}\Bigg\{\sum_{n=0}^{q-1} (-1)^n e^{-\pi i n^2 \frac{p}{q} - 2\pi i n u} 
	\\\nonumber
& + \frac{(-1)^q e^{-\pi i q p - 2\pi i q u}}{\sqrt{p/q}} e^{\frac{3\pi i}{4}} \sum_{n=0}^{p-1} e^{\frac{\pi i q}{p}(u + n + \frac{1}{2})^2}\Bigg\}.
\end{align}

For $N=3$ equation (\ref{zetaformula}) simplifies to
\begin{align}\label{zetaformulab}
\zeta(s) 
=& \sum_{n=1}^{[\frac{t}{\eta}]} \frac{1}{n^s} 
+ \chi(s) \sum_{n=1}^{[\frac{\eta}{2\pi}]} \frac{1}{n^{1-s}}
 	\\ \nonumber
& + e^{-i\pi s}\Gamma(1-s)\biggl\{ e^{\frac{i\pi(s-1)}{2}}\eta^{s-1}  e^{\frac{2t}{\eta}[\frac{\eta}{2\pi}]\pi i - it 
- \frac{it}{2\eta^2}(2[\frac{\eta}{2\pi}]\pi - \eta)^2} 
	\\ \nonumber
& \times \biggl[\Phi + \frac{\sigma -1}{i\eta} \biggl(
  \partial_2 \Phi  
  + 
\biggl(2\Bigl[\frac{\eta}{2\pi}\Bigr]\pi i - i\eta\biggr)   \Phi  \biggr)
  	\\ \nonumber
&  -\frac{(\sigma -2) (\sigma -1)}{2 \eta ^2} \biggl(
  \partial_2^{2} \Phi  
 + 2 \biggl(2\Bigl[\frac{\eta}{2\pi}\Bigr]\pi i - i\eta\biggr)  \partial_2 \Phi  
  + \biggl(2\Bigl[\frac{\eta}{2\pi}\Bigr]\pi i - i\eta\biggr)^2 \Phi     \biggr) \biggr]
  	\\ \nonumber
& + O\biggl( e^{-\frac{\pi t}{2}} \frac{\eta^\sigma}{t} \biggr)\biggr\},
\qquad 
 \epsilon \sqrt{t} < \eta < t , \quad 0 \leq \sigma \leq 1, \quad  t \to \infty,
\end{align}
where the error term is uniform for all $\eta, \sigma$ in the above ranges and $\Phi$ and its partial derivatives are evaluated at the point (\ref{Phievaluationpoint}).
\end{theorem}
\begin{remark} \upshape
1. If $\eta = \sqrt{\frac{2\pi t}{b}}$ where $b = p/q > 0$ is a rational number, then equations (\ref{zetaformula})-(\ref{Phiformula}) provide an explicit asymptotic expansion of $\zeta(s)$ to all orders. The particular case $b = 1$ is the famous case analyzed by Siegel in \cite{S1932} using ideas from Riemann's unpublished notes. We note that, building on the work of \cite{K1993}, an elementary (but formal) derivation of the Riemann-Siegel formula is presented in \cite{B1995}.

2. The results of theorem \ref{ZETATH} are convenient for studying the higher-order asymptotics of $\zeta(s)$ in the case when $\eta = \text{constant} \times \sqrt{t}$. However, if the order of $\eta$ is strictly smaller or larger than the order of $\sqrt{t}$, the results of theorem \ref{ZETATH} are less convenient, because in this case the asymptotics of $\zeta(s)$ involves the asymptotics of the sums in (\ref{Phiformula}) as $p$ and/or $q$ tend to infinity.
In this case, the alternative representation of the asymptotics of theorem \ref{zetath2}, which avoids the appearance of the above sums, is usually more convenient. 

3. The error term in (\ref{zetaformula}) is uniform with respect to $\eta, \sigma, N$. More precisely, this means that equation (\ref{zetaformula}) is equivalent to the following statement: For every $\epsilon>0$, there exist constants $A>0$, $C>0$, $K>0$ such that the inequality
\begin{align*}
\Biggl|& \zeta(s) 
- \sum_{n=1}^{[\frac{t}{\eta}]} \frac{1}{n^s} 
- \chi(s) \sum_{n=1}^{[\frac{\eta}{2\pi}]} \frac{1}{n^{1-s}}
 	\\ \nonumber
& \qquad- e^{-i\pi s}\Gamma(1-s) e^{\frac{i\pi(s-1)}{2}}\eta^{s-1}  e^{\frac{2t}{\eta}[\frac{\eta}{2\pi}]\pi i - it - \frac{it}{2\eta^2}(2[\frac{\eta}{2\pi}]\pi -\eta)^2} S_N(s, \eta) \Biggr| 
	\\
& \qquad\qquad \leq C |e^{-i\pi s}\Gamma(1-s)| e^{-\frac{\pi t}{2}} \left(\frac{3N}{t}\right)^{\frac{N}{6}} \frac{\eta^{\sigma}}{\sqrt{t}},
\end{align*}
holds for all $t > K$ and all $\sigma, N, \eta$ such that $0 \leq \sigma \leq 1$, $1 \leq N < At$, and $\epsilon \sqrt{t} < \eta < t$.

4. The main difference between the proofs of theorems \ref{zetath2} and \ref{ZETATH} can be explained as follows. Consider the representation (\ref{zetaCjs}) of $\zeta(s)$.
If $\eta$ is of strictly smaller order than $t^{1/2}$, the main contribution to the asymptotics of $\zeta(s)$ comes from the part of the contour $C_2$ that lies within a distance $\epsilon$ of the critical point $i\eta$; the remaining part of $C_2$ gives a contribution which is suppressed by a factor of the form $e^{-\text{const} \times \frac{t}{\eta^2}}$ cf. Eq. (\ref{I2restimate}). 
Thus, the proof of theorem \ref{zetath2} relies on a direct study of the integral near the critical point $i\eta$ using the method of steepest descent. 
On the other hand, in the case when $\eta \sim \sqrt{t}$, the asymptotic expansion of $\zeta(s)$ depends on the integral along all of the contour $C_2$. The contour $C_2$ is a line segment of total length $\eta$ centered on the critical point $i\eta$. In the proof of theorem \ref{ZETATH}, we handle this integral by relating it to the function $\Phi(\tau, u)$ defined in (\ref{Phidef}). 

5. The coefficients $a_n$ satisfy the estimate
$$a_n = O\biggl(\frac{t^{[\frac{n}{3}]}}{\eta^n}\biggr), \qquad \text{(not uniformly in $n$);}$$
indeed, assuming this estimate up to $n$, we find
$$a_{n+1} = O\biggl(\frac{t^{[\frac{n}{3}]}}{\eta^{n+1}}\biggr)
+ O\biggl(\frac{t}{\eta^3} \frac{t^{[\frac{n-2}{3}]}}{\eta^{n-2}}\biggr)
= O\biggl(\frac{t^{[\frac{n+1}{3}]}}{\eta^{n+1}}\biggr), \qquad \text{(not uniformly in $n$)}.$$
\end{remark}

\noindent
\begin{proof}[Proof of Theorem \ref{ZETATH}.]
Let $t > 0$, $0 \leq \sigma \leq 1$, and $m = [t/\eta]$. As in the proof of theorem \ref{zetath2}, $\zeta(s)$ is given by (\ref{zetaCjs}) where the contributions from $I_1, I_3$, and $I_4$ are exponentially small, and the error terms of the form $O(\cdot)$ are uniform with respect to the variables $\eta, \sigma, N$ as $t \to \infty$.
It remains to analyze the integral $I_2$.
In the neighborhood of $w = i\eta$ we have
\begin{align*}
(s-1)\ln\frac{w}{i\eta}
& = (s-1)\left(\frac{w-i\eta}{i\eta} - \frac{1}{2}\left(\frac{w-i\eta}{i\eta}\right)^2 + \cdots \right)
	\\
& = \frac{t}{\eta}(w - i\eta) + \frac{it}{2\eta^2}(w - i\eta)^2 + \cdots.
\end{align*}
Hence we write
$$e^{(s-1)\ln{\frac{w}{i\eta}} }
= e^{\frac{t}{\eta}(w - i\eta) + \frac{it}{2\eta^2} (w - i\eta)^2} \phi(w - i\eta),$$
where $\phi(z)$ is defined by (\ref{phidef}). Define $\{a_n\}_0^\infty$ by 
\begin{align*}
\phi(z) = \sum_{n=0}^\infty a_n z^n, \qquad |z| < \eta,
\end{align*}
The identity
$$\frac{d\phi}{dz}
= \left(\frac{s-1}{i\eta + z} - \frac{t}{\eta} - \frac{it}{\eta^2} z \right) \phi(z)$$
implies
$$(i\eta + z)\sum_{n=1}^\infty n a_n z^{n-1}
= \left[s-1 - (i\eta + z)\left(\frac{t}{\eta} + \frac{itz}{\eta^2}\right)\right]\sum_{n=0}^\infty a_n z^n.$$
Hence the coefficients $a_n$ are determined in succession by the recurrence formula (\ref{anrecurrence}) supplemented with the conditions $a_{-2} = a_{-1} = 0$ and $a_0 = 1$. The first few coefficients are given by
\begin{align*}
& a_0 = 1, \qquad a_1 = \frac{\sigma -1}{i\eta}, \qquad a_2 = -\frac{(\sigma -2) (\sigma -1)}{2 \eta ^2},
	\\
& a_3 = \frac{-2 t+i (\sigma -3) (\sigma -2) (\sigma -1)}{6 \eta ^3}, \qquad
	\\
& a_4 = \frac{(\sigma -4) (\sigma -3) (\sigma -2) (\sigma -1)+2 i (4 \sigma -7) t}{24 \eta ^4}.
\end{align*}	

Representing $\phi(z)$ in the form
$$\phi(z) = \sum_{n=0}^{N-1} a_n z^n + r_N(z),$$
we find
$$r_N(z) = \frac{z^N}{2\pi i} \int_\Gamma \frac{\phi(w)}{w^N(w-z)}dw, \qquad |z| < \eta,$$
where $\Gamma$  is a contour contained in the disk of radius $\eta$ centered at the origin which encircles the points $0$ and $z$ once. 

Let $|z| < \frac{4}{7}\eta$ (so that $\frac{21}{20}|z| < \frac{3}{5}\eta$) and let $\Gamma$ be a circle with center $w =0$ and radius $\rho_N$, where
$$\frac{21}{20} |z| \leq \rho_N \leq \frac{3}{5} \eta.$$
Equation (\ref{lnphizestimate}) implies the following analog of (\ref{relogphiw}):
\begin{align*}
  \re{\ln{\phi(z)}} \leq |\sigma -1| \ln\frac{8}{5} + \frac{5t|z|^3}{6\eta^3}, \qquad |z| \leq \frac{3}{5}\eta.
\end{align*}  
Hence
$$r_N(z) = O\left(\frac{2\pi \rho_N}{\rho_N - |z|}|z|^N \rho_N^{-N}e^{\frac{5t\rho_N^3}{6\eta^3}}\right)
= O\left(|z|^N \rho_N^{-N}e^{\frac{5t\rho_N^3}{6\eta^3}}\right), \qquad |z| < \frac{4}{7}\eta.$$
The function $\rho^{-N}e^{\frac{5t\rho^3}{6\eta^3}}$ has the minimum $(\frac{5et}{2N\eta^3})^{N/3}$ for $\rho = (\frac{2N}{5t})^{1/3}\eta$; $\rho_N$ can have this value if
$$\frac{21}{20}|z| \leq \biggl(\frac{2N}{5t}\biggr)^{\frac{1}{3}}\eta \leq \frac{3}{5}\eta.$$
Hence, 
\begin{align} \label{rNestimate}
 r_N(z) & = O\left(|z|^N \biggl(\frac{5et}{2N\eta^3}\biggr)^{\frac{N}{3}}\right),
  \qquad N < At, \quad |z| \leq \frac{20}{21}\biggl(\frac{2N}{5t}\biggr)^{\frac{1}{3}} \eta.
\end{align}
For $|z| < \frac{4}{7}\eta$ we can also take $\rho_N = \frac{21}{20}|z|$, which yields
\begin{align}\label{rNestimate2}
r_N(z) = O\left(\Bigl(\frac{20}{21}\Bigr)^N e^{\frac{5}{6}\frac{t}{\eta^3} (\frac{21}{20}|z|)^3}\right)
= O\left(e^{\frac{14}{29}\frac{t}{\eta^2} |z|^2}\right), \qquad |z| < \frac{\eta}{2}.
\end{align}
where we have used the inequalities $\frac{|z|}{\eta} < \frac{1}{2}$ and $\frac{5}{12}(\frac{21}{20})^3 < \frac{14}{29}$ in the last step.

We write $I_2$ in the form
\begin{align}\label{I2expression}
  I_2 = I_2^S + I_2^R,
\end{align} 
where $I_2^S$ and $I_2^R$ denote the integrals
\begin{align}
I_2^S  & =  \int_{C_2} (i\eta)^{s-1}\frac{e^{\frac{t}{\eta}(w - i\eta) + \frac{it}{2\eta^2} (w - i\eta)^2 - mw}}{e^w - 1} \sum_{n=0}^{N-1} a_n (w - i\eta)^n dw
\end{align}
and 
\begin{align}
I_2^R = \int_{C_2} (i\eta)^{s-1}\frac{e^{\frac{t}{\eta}(w - i\eta) + \frac{it}{2\eta^2} (w - i\eta)^2 - mw}}{e^w - 1} r_N(w - i\eta) dw.
\end{align}

We claim that
\begin{align}\label{I2Restimate}
I_2^R = O\biggl( e^{-\frac{\pi t}{2}} \left(\frac{3N}{t}\right)^{\frac{N}{6}} \frac{\eta^\sigma}{\sqrt{t}}\biggr).
\end{align}
Indeed, suppose first that $|e^w - 1| > A$ on $C_2$. Then we have for $w\in C_2$ with $|w - i\eta| < \epsilon$, the estimate
$$\biggl|\frac{e^{(\frac{t}{\eta}-m)w}}{e^w -1}\biggr| < \frac{e^\epsilon}{A},$$
whereas for $w\in C_2$ with $|w - i\eta| > \epsilon$, we have the estimates
\begin{align}\label{expquotientestimate}
\frac{e^{(\frac{t}{\eta}-m)w}}{e^w -1} = \begin{cases} O\biggl(\frac{e^{(\frac{t}{\eta}-m-1)\re{w}}}{1 - e^{-\re{w}}}\biggr), & \re{w} \geq 0, \\
O\biggl(\frac{e^{(\frac{t}{\eta}-m)\re{w}}}{ e^{\re{w}} - 1}\biggr), & \re{w} < 0. \end{cases}
\end{align}
It follows from these estimates that the factor $\frac{e^{(\frac{t}{\eta}-m)w}}{e^w -1}$ in the integrand of $I_2^R$ is of $O(1)$ on $C_2$.
We now make the change of variables $w = i\eta + \lambda e^{\frac{i\pi}{4}}$ and split the integral into two integrals; in the first integral $|\lambda | \leq A(\frac{N}{t})^{\frac{1}{3}} \eta$, while in the second integral $A(\frac{N}{t})^{\frac{1}{3}} \eta \leq |\lambda| \leq \eta/2$.
In the first integral, we use the estimate (\ref{rNestimate}) of $r_N$, whereas in the second integral, we use the estimate (\ref{rNestimate2}). This yields
\begin{align*}
I_2^R = &\; O\biggl(\eta^{\sigma -1}e^{-\frac{\pi t}{2}}\biggl\{ \int_0^{A(\frac{N}{t})^{\frac{1}{3}} \eta} e^{-\frac{t}{2\eta^2}\lambda^2} \lambda^N \left(\frac{5et}{2N\eta^3}\right)^{\frac{N}{3}} d\lambda
	\\
& + \int_{A(\frac{N}{t})^{\frac{1}{3}} \eta}^{\frac{\eta}{2}} e^{-\frac{t}{2\eta^2}\lambda^2 + \frac{14}{29}\frac{t}{\eta^2}\lambda^2} d\lambda \biggr\}\biggr).
\end{align*}
Steps almost identical to those leading from (\ref{C2epsestimate}) to (\ref{Oetasigmaminusone}) now show that 
\begin{align*}
I_2^R = O\biggl(\eta^{\sigma -1} e^{-\frac{\pi t}{2}} \left(\frac{3N}{t}\right)^{\frac{N}{6}} \frac{\eta}{\sqrt{t}}\biggr).	
\end{align*}
The case where the contour goes near a pole gives a similar result.
Indeed, if the contour passes near the pole at $w = 2w\pi i$, say within a distance $1/2$ of it, we take it around an arc of the circle $|w - i\eta| = 1$. Letting $w = i\eta + e^{i\theta}$ and using (\ref{rNestimate}), the total contribution to $I_2^R$ from this arc is given by 
\begin{align*}
& O\biggl(\eta^{\sigma -1}e^{-\frac{\pi t}{2}} \int_0^{2\pi} \biggl| e^{(\frac{t}{\eta} - [\frac{t}{\eta}])e^{i\theta} + \frac{it}{2\eta^2}e^{2i\theta}} r_N(e^{i\theta}) \biggr| d\theta\biggr)
	\\
& = O\biggl(\eta^{\sigma -1}e^{-\frac{\pi t}{2}} 
\biggl(\frac{5et}{2N\eta^3}\biggr)^{\frac{N}{3}}
\int_0^{2\pi} e^{-\frac{t}{2\eta^2}\sin(2\theta)} d\theta\biggr)
	\\
& = O\biggl(\eta^{\sigma -1} e^{-\frac{\pi t}{2}} \biggl(\frac{5et}{2N\eta^3}\biggr)^{\frac{N}{3}}\biggr)
= O\biggl(\eta^{\sigma -1} e^{-\frac{\pi t}{2}} \left(\frac{3N}{t}\right)^{\frac{N}{6}} \frac{\eta}{\sqrt{t}}\biggr),
\end{align*}
where we used the assumption that $\eta > \epsilon \sqrt{t}$.
This proves (\ref{I2Restimate}).

We next consider $I_2^S$. We claim that there exists a constant $A>0$ such that
\begin{align}\label{I2Sestimate}
  I_2^S =  \int_{C_2'} (i\eta)^{s-1}\frac{e^{\frac{t}{\eta}(w - i\eta) + \frac{i}{2}\frac{t}{\eta^2} (w - i\eta)^2 - mw}}{e^w - 1} \sum_{n=0}^{N-1} a_n (w - i\eta)^n dw + O(e^{-\frac{\pi t}{2} - At}),
\end{align}
where $C_2'$ denotes the infinite straight line of which $C_2$ is a part.
Indeed, in view of (\ref{expquotientestimate}), if we replace $C_2$ by $C_2'$, the integral multiplying $a_n$ in the expression for $I_2^S$ changes by
\begin{align}\label{ancoefficientbigOh}
O\left(\eta^{\sigma -1} e^{-\frac{\pi t}{2}} \int_{\frac{\eta}{2}}^\infty e^{-\frac{t}{2\eta^2}\lambda^2} \lambda^n  d\lambda\right).
\end{align}
We write the integrand as
$$e^{-\frac{t}{4\eta^2} \lambda^2} \lambda^n \times e^{-\frac{t}{4\eta^2} \lambda^2}.$$
The first factor is steadily decreasing for $\lambda > \sqrt{\frac{2n}{t}} \eta$, and so it decreases throughout the interval of integration provided that $n < N < At$, with $A$ sufficiently small. The term in (\ref{ancoefficientbigOh}) is then
\begin{align}\label{wholetermO}
O\left(\eta^{\sigma -1} e^{-\frac{\pi t}{2}} e^{-\frac{t}{4\eta^2} \frac{\eta^2}{4}} \left(\frac{\eta}{2}\right)^n \int_{\frac{\eta}{2}}^\infty e^{-\frac{t}{4\eta^2}\lambda^2} d\lambda\right)
= O\left(\eta^{\sigma -1} e^{-\frac{\pi t}{2}} e^{-\frac{t}{4\eta^2} \frac{\eta^2}{4}} \left(\frac{\eta}{2}\right)^n \frac{\eta}{\sqrt{t}} \right).
\end{align}
Also, by (\ref{rNestimate}), choosing $z$ such that $|z| \leq \frac{20}{21} (\frac{2}{5t})^{\frac{1}{3}}\eta$, we find
\begin{align}\nonumber
a_n & = (r_n(z) - r_{n+1}(z))z^{-n} 
	\\ \nonumber
& = O\left(\left(|z|^n \left(\frac{5et}{2n\eta^3}\right)^{\frac{n}{3}} - |z|^{n+1} \left(\frac{5et}{2(n+1)\eta^3}\right)^{\frac{n+1}{3}}\right)|z|^{-n}\right)
	\\\label{anestimate}
&=  O\left(\left(\frac{5et}{2n\eta^3}\right)^{\frac{n}{3}}\right), \qquad 1 \leq n \leq N-1.
\end{align}
Multiplying (\ref{wholetermO}) by $a_n$ and summing from  $0$ to $N-1$, we find that the total error is
\begin{align*}
 & O\left(\eta^{\sigma -1} e^{-\frac{\pi t}{2} - \frac{t}{16}}\frac{\eta}{\sqrt{t}}
\biggl(a_0+ \sum_{n=1}^{N-1} \left(\frac{\eta}{2}\right)^n \left(\frac{5 e t}{2 n \eta^3}\right)^{\frac{n}{3}} \biggr)\right)
  	\\
& = 
O\left(\eta^{\sigma -1} e^{-\frac{\pi t}{2} - \frac{t}{16}}\frac{\eta}{\sqrt{t}}
\biggl(1 + \sum_{n=1}^{N-1} \left(\frac{5 e t}{16 n }\right)^{\frac{n}{3}} \biggr)\right).
\end{align*}
Now the factor $(t/n)^{\frac{n}{3}}$ increases steadily up to $n = t/e$, and so if $n < A t$, where $A < 1/e$, it is
$$O(e^{\frac{1}{3}t A\ln\frac{1}{A}}).$$
Hence if $N < At$, with $A$ sufficiently small, the total error is
$$O(e^{-\frac{\pi t}{2} - At}).$$
This proves (\ref{I2Sestimate}).

Finally, we analyze the sum 
\begin{align}\label{finalsum}
 (i\eta)^{s-1} \sum_{n=0}^{N-1} a_n\int_{C_2'} \frac{e^{\frac{t}{\eta}(w - i\eta) + \frac{it}{2\eta^2}(w - i\eta)^2 - m w}}{e^w - 1} (w - i\eta)^n dw.
 \end{align}
The integral in (\ref{finalsum}) may be expressed as
\begin{align*}
  -\int_L e^{\frac{t}{\eta}(w + 2[\frac{\eta}{2\pi}]\pi i - i\eta) + \frac{it}{2\eta^2}(w + 2[\frac{\eta}{2\pi}]\pi i - i\eta)^2 - [\frac{t}{\eta}] w}
  \frac{(w + 2[\frac{\eta}{2\pi}]\pi i - i\eta)^n}{e^w - 1} dw,
\end{align*}
where $L$ is a line in the direction $\arg{w} = \pi/4$, passing between $0$ and $2\pi i$.
This is $n!$ times the coefficient of $\xi^n$ in the following expression:
\begin{align}\label{intwithxi}
 & -\int_L e^{\frac{t}{\eta}(w + 2[\frac{\eta}{2\pi}]\pi i - i\eta) + \frac{it}{2\eta^2}(w + 2[\frac{\eta}{2\pi}]\pi i - i\eta)^2 - [\frac{t}{\eta}] w + \xi(w + 2[\frac{\eta}{2\pi}]\pi i - i\eta)} \frac{dw}{e^w -1}
  	\\ \nonumber
& = - e^{\frac{t}{\eta}2[\frac{\eta}{2\pi}]\pi i - it + \frac{it}{2\eta^2}(2[\frac{\eta}{2\pi}]\pi i - i\eta)^2 + \xi(2[\frac{\eta}{2\pi}]\pi i - i\eta)}
\int_{L} e^{\frac{it}{2\eta^2} w^2 + (\frac{2t}{\eta} - \frac{2\pi t}{\eta^2}[\frac{\eta}{2\pi}] - [\frac{t}{\eta}] + \xi)w}\frac{dw}{e^w - 1}
	\\ \nonumber
& = e^{\frac{t}{\eta}2[\frac{\eta}{2\pi}]\pi i - it + \frac{it}{2\eta^2}(2[\frac{\eta}{2\pi}]\pi i - i\eta)^2}2\pi i \Phi\left(-\frac{2\pi t}{\eta^2}, \frac{2t}{\eta} - \frac{2\pi t}{\eta^2}\Bigl[\frac{\eta}{2\pi}\Bigr] - \Bigl[\frac{t}{\eta}\Bigr] + \xi - \frac{1}{2}\right)
	\\\nonumber
&\quad \times e^{\xi(2[\frac{\eta}{2\pi}]\pi i - i\eta)}, 
\end{align}
where the function $\Phi(\tau, u)$ is defined by (\ref{Phidef}), i.e.,
  $$\Phi(\tau, u) = \int_{0 \nwarrow 1} \frac{e^{\pi i \tau x^2 + 2\pi i u x}}{e^{\pi i x} - e^{-\pi i x}} dx
  = \frac{1}{2\pi i}\int_{-L} \frac{e^{-\frac{i\tau}{4\pi}w^2 + (u + \frac{1}{2})w}}{e^w -1} dw, \qquad \tau < 0, \quad u \in \C.$$
We rewrite the expression on the rhs of (\ref{intwithxi}) as
\begin{align*}
e^{\frac{t}{\eta}2[\frac{\eta}{2\pi}]\pi i - it + \frac{it}{2\eta^2}(2[\frac{\eta}{2\pi}]\pi i - i\eta)^2}2\pi i 
& \sum_{l=0}^\infty \partial_2^l\Phi \left(-\frac{2\pi t}{\eta^2}, \frac{2t}{\eta} - \frac{2\pi t}{\eta^2}\Bigl[\frac{\eta}{2\pi}\Bigr] - \Bigl[\frac{t}{\eta}\Bigr] - \frac{1}{2}\right)\frac{\xi^l}{l!}
	\\
&\times \sum_{k=0}^\infty \frac{\xi^k(2[\frac{\eta}{2\pi}]\pi i - i\eta)^k}{k!}.
\end{align*}
It follows that the expression in (\ref{finalsum}) is given by
\begin{align}\label{finalsumexpression}
 & 2\pi i (i\eta)^{s-1}    e^{\frac{t}{\eta}2[\frac{\eta}{2\pi}]\pi i - it + \frac{it}{2\eta^2}(2[\frac{\eta}{2\pi}]\pi i - i\eta)^2} S_N(s, \eta)
\end{align}
where $S_N(s, \eta)$ is defined in (\ref{SNdef}).
Equations (\ref{I2Restimate}), (\ref{I2Sestimate}), and (\ref{finalsumexpression}) show that $I_2$ satisfies
\begin{align*}
I_2 = &\; 2\pi i (i\eta)^{s-1}    e^{\frac{t}{\eta}2[\frac{\eta}{2\pi}]\pi i - it + \frac{it}{2\eta^2}(2[\frac{\eta}{2\pi}]\pi i - i\eta)^2} S_N(s, \eta) 
	\\
& +O\biggl( e^{-\frac{\pi t}{2}} \left(\frac{3N}{t}\right)^{\frac{N}{6}} \frac{\eta^\sigma}{\sqrt{t}}\biggr)
+ O(e^{-\frac{\pi t}{2} - At}).
\end{align*}
Equation (\ref{zetaformula}) follows by substituting this expression into (\ref{zetaCjs}).

In order to prove (\ref{Phiformula}) we suppose that $u \in \C$ and $\tau < 0$. We claim that $\Phi$ satisfies the two recursion relations
\begin{align}\label{identity1}
  \Phi(\tau, u) = \Phi(\tau, u + 1) - \frac{e^{\frac{3\pi i}{4}}}{\sqrt{|\tau|}} e^{-\frac{\pi i}{\tau}(u + \frac{1}{2})^2}
\end{align}
and
\begin{align}\label{identity2}
  \Phi(\tau, u) = 1 - e^{\pi i \tau - 2\pi i u}\Phi(\tau, u - \tau).
\end{align}
Repetitive use of these two equations yields the following identities:
\begin{align}\label{identity1N}
\Phi(\tau, u) = \Phi(\tau, u + N) - \frac{e^{\frac{3\pi i}{4}}}{\sqrt{|\tau|}} \sum_{n=0}^{N-1} e^{-\frac{\pi i}{\tau}(u + n + \frac{1}{2})^2}
\end{align} 
and
\begin{align}\label{identity2N}
  \Phi(\tau, u) = \sum_{n=0}^{N-1} (-1)^n e^{\pi i n^2 \tau - 2\pi i n u} + (-1)^N e^{\pi i N^2 \tau - 2\pi i N u} \Phi(\tau, u - N \tau)
\end{align}
valid for all $N \geq 0$. If $\tau = -p/q$, we apply (\ref{identity1N}) with $N = p$ and (\ref{identity2N}) with $N = q$ to find
\begin{align}\label{identity1pq}
\Phi(\tau, u) = \Phi(\tau, u + p) - \frac{e^{\frac{3\pi i}{4}}}{\sqrt{|\tau|}} \sum_{n=0}^{p-1} e^{-\frac{\pi i}{\tau}(u + n + \frac{1}{2})^2}
\end{align} 
and
\begin{align}\label{identity2pq}
  \Phi(\tau, u) = \sum_{n=0}^{q-1} (-1)^n e^{\pi i n^2 \tau - 2\pi i n u} + (-1)^q e^{\pi i q^2 \tau - 2\pi i q u} \Phi(\tau, u + p).
\end{align}
Eliminating $\Phi(\tau, u + p)$ from these two equations, we find (\ref{Phiformula}).

It remains to prove (\ref{identity1}) and (\ref{identity2}). In order to prove (\ref{identity1}), we note that
\begin{align*}
\Phi(\tau, u+1) - \Phi(\tau,u) & = \int_{0 \nwarrow 1} e^{\pi i \tau x^2} \frac{e^{2 \pi i (u+1)x} - e^{2\pi i u x}}{e^{\pi i x} - e^{-\pi i x}} dx
	\\
& = \int_{0 \nwarrow 1} e^{\pi i \tau x^2 + 2\pi i (u + \frac{1}{2}) x} dx
	\\
& = e^{-\frac{\pi i}{\tau} (u + \frac{1}{2})^2} \int_{0 \nwarrow 1} e^{\pi i \tau (x + \frac{u + \frac{1}{2}}{\tau})^2} dx
	\\
& = e^{-\frac{\pi i}{\tau} (u + \frac{1}{2})^2} \int_{0 \nwarrow 1} e^{\pi i \tau x^2} dx.
\end{align*}
The identities
$$\int_{0 \nwarrow 1} e^{\pi i \tau x^2} dx = \int_{0 \nwarrow 1} e^{-\pi i y^2} \frac{dy}{\sqrt{|\tau|}}
= \frac{e^{\frac{3\pi i}{4}}}{\sqrt{|\tau|}}$$
imply equation (\ref{identity1}).

The identity (\ref{identity2}) is established as follows:
\begin{align*}
\Phi(\tau, u) & = \int_{0 \nwarrow 1} \frac{e^{\pi i \tau x^2 + 2\pi i u x}}{e^{\pi i x} - e^{-\pi i x}} dx
	\\
& = 1 + \int_{-1 \nwarrow 0} \frac{e^{\pi i \tau x^2 + 2\pi i u x}}{e^{\pi i x} - e^{-\pi i x}} dx
	\\
& = 1 + \int_{0 \nwarrow 1} \frac{e^{\pi i \tau (x-1)^2 + 2\pi i u (x-1)}}{e^{\pi i (x-1)} - e^{-\pi i (x-1)}} dx
	\\
& = 1 - e^{\pi i \tau - 2\pi i u} \int_{0 \nwarrow 1} \frac{e^{\pi i \tau x^2 + 2\pi i (u-\tau) x}}{e^{\pi i x} - e^{-\pi ix}} dx
	\\
& = 1 - e^{\pi i \tau - 2\pi i u} \Phi(\tau, u - \tau),
\end{align*}
where the second equality follows from the fact that the residue of the integrand at $x =0$ is $\frac{1}{2\pi i}$.
\end{proof}

Theorem \ref{ZETATH} is valid for $\epsilon \sqrt{t} < \eta < t$. As a corollary, we can find an analogous result valid for $2\pi < \eta < \frac{\sqrt{t}}{\epsilon}$.

\begin{corollary}[{\bf The asymptotics to all orders for the case $2\pi < \eta < \frac{\sqrt{t}}{\epsilon}$}]\label{zetacor}
For every $\epsilon > 0$, there exists a constant $A>0$ such that
\begin{align}\label{zetaformula3}
\zeta(s) 
=& \sum_{n=1}^{[\frac{t}{\eta}]} \frac{1}{n^{s}}
+ \chi(s) \sum_{n=1}^{[\frac{\eta}{2\pi}]} \frac{1}{n^{1-s}} 
 	\\ \nonumber
& + \chi(s) e^{i\pi(1-s)}\Gamma(s)\biggl\{ e^{\frac{i\pi s}{2}}\biggl(\frac{\eta}{2\pi t}\biggr)^{s}  
e^{-i\eta[\frac{t}{\eta}]  + it 
+ \frac{i\eta^2}{2t}([\frac{t}{\eta}] - \frac{t}{\eta})^2} \overline{S_N\biggl(1-\bar{s}, \frac{2\pi t}{\eta}\biggr)}
	\\ \nonumber
 & \hspace{4cm} + O\biggl( e^{-\frac{\pi t}{2}} \left(\frac{3N}{t}\right)^{\frac{N}{6}} \frac{(\frac{2\pi t}{\eta})^{1-\sigma}}{\sqrt{t}}\biggr)\biggr\},
 	\\ \nonumber
& \qquad 
2\pi < \eta < \frac{\sqrt{t}}{\epsilon}, \quad 0 \leq \sigma \leq 1, \quad 1 \leq N < At, \quad  t \to \infty,
\end{align}
where the error term is uniform for all $\eta, \sigma, N$ in the above ranges and $S_N$ is given by (\ref{SNdef}).
\end{corollary}
\begin{proof}
We replace $\sigma$ by $1- \sigma$ in (\ref{zetaformula}) and take the complex conjugate of both sides. We then multiply the resulting equation by $\chi(s)$ and use the identities (\ref{baridentities}).
This yields the following equation
\begin{align}\nonumber
\zeta(s) =&\; \chi(s) \sum_{n=1}^{[\frac{t}{\eta}]} \frac{1}{n^{1-s}} + \sum_{n=1}^{[\frac{\eta}{2\pi}]} \frac{1}{n^{s}}
 	\\ \nonumber
& + \chi(s) e^{i\pi(1-s)}\Gamma(s)\biggl\{ (e^{\frac{i\pi s}{2}}\eta^{-s}  e^{-\frac{2t}{\eta}[\frac{\eta}{2\pi}]\pi i + it 
+ \frac{it}{2\eta^2}(2[\frac{\eta}{2\pi}]\pi - \eta)^2} \overline{S_N(1-\bar{s}, \eta)}
	\\ \nonumber
 & \hspace{4cm} + O\biggl(e^{-\frac{\pi t}{2}} \left(\frac{3N}{t}\right)^{\frac{N}{6}} \frac{\eta^{1-\sigma}}{\sqrt{t}}\biggr)\biggr\},
 	\\ \nonumber & \qquad 
 \epsilon \sqrt{t} < \eta < t , \quad 0 \leq \sigma \leq 1, \quad 1 \leq N < At, \quad  t \to \infty.
\end{align}
Replacing $\eta$ by $\frac{2\pi t}{\eta}$, we find (\ref{zetaformula3}).
\end{proof}

\chapter{Consequences of the Asymptotic Formulae}\label{sec5}
Using theorems \ref{th3.1} and \ref{th3.2} we can compute several interesting sums.

\begin{theorem}\label{th5.1}
Define the polylogarithm $\Li_m(z)$ by (\ref{polylogdef}). The following relation holds:
\begin{align}\label{sumzetaformula2.1}
\sum_{n = [\frac{\eta_1}{2\pi}] + 1}^{[\frac{\eta_2}{2\pi}]} n^{-s} = &\;
\frac{1}{1-s}\biggl[ \left(\frac{\eta_2}{2\pi}\right)^{1-s} - \left(\frac{\eta_1}{2\pi}\right)^{1-s} \biggr]
	\\ \nonumber
&+ \frac{e^{-\frac{i\pi (1-s)}{2}}}{(2\pi)^{1-s}}\sum_{n=1}^\infty  \sum_{j=0}^{N-1} e^{-nz - it\ln{z}} \left(\frac{1}{n + \frac{it}{z}}\frac{d}{dz}\right)^j\frac{z^{-\sigma}}{n + \frac{it}{z}} \Biggr|_{z = i\eta_2}^{i\eta_1}
	\\\nonumber
& + \frac{e^{\frac{i\pi (1-s)}{2}}}{(2\pi)^{1-s}}\sum_{n=1}^\infty
\sum_{j=0}^{N-1} e^{-nz - it\ln{z}} \left(\frac{1}{n + \frac{it}{z}}\frac{d}{dz}\right)^j\frac{z^{-\sigma}}{n + \frac{it}{z}} \Biggr|_{z = -i\eta_2}^{-i\eta_1}
	\\ \nonumber
& + O\biggl((2N +1)!! N \Bigl(\frac{1+\epsilon}{\epsilon}\Bigr)^{2(N+1)} \eta_1^{-\sigma - N}  \biggr),
	\\ \nonumber
& (1+\epsilon)t < \eta_1 < \eta_2 < \infty, \quad \epsilon > 0, \quad 0 \leq \sigma \leq 1, \quad N \geq 2, \quad t \to \infty,	
\end{align}
where the error term is uniform for all $\eta_1, \eta_2, \epsilon, \sigma,N$ in the above ranges.

For $N=3$ equation (\ref{sumzetaformula2.1}) simplifies to  
\begin{align}\label{sumzetaformula2.1b}
\sum_{n = [\frac{\eta_1}{2\pi}] + 1}^{[\frac{\eta_2}{2\pi}]} n^{-s}  = &\; 
\frac{1}{1-s}\biggl[ \left(\frac{\eta_2}{2\pi}\right)^{1-s} - \left(\frac{\eta_1}{2\pi}\right)^{1-s} \biggr]
	\\ \nonumber
& + \frac{2i\eta_1^{-s}}{(2\pi)^{1-s}} \biggl\{
-i\arg(1 - e^{i \eta_1})
+ \frac{t - i \sigma}{\eta_1} \re \Li_{2}(e^{i\eta_1}) 
	\\ \nonumber
&\hspace{2.1cm} +\frac{1}{\eta_1^2} \bigl[ i t^2 -  3i \sigma t - (\sigma-1) t - i\sigma(\sigma+1)\bigr]\im \Li_{3}(e^{i\eta_1}) \biggr\}
	\\ \nonumber
& - \frac{2i\eta_2^{-s}}{(2\pi)^{1-s}} \biggl\{
-i\arg(1 - e^{i \eta_2})
+ \frac{t - i \sigma}{\eta_2} \re \Li_{2}(e^{i\eta_2}) 
	\\ \nonumber
&\hspace{2.1cm} +\frac{1}{\eta_2^2} \bigl[ i t^2 -  3i \sigma t - (\sigma-1) t - i\sigma(\sigma+1)\bigr]\im \Li_{3}(e^{i\eta_2}) \biggr\}	
	\\\nonumber
&    + O\biggl(\frac{1}{\eta_1^{3 + \sigma}}\biggl(t^3 + \Bigl(\frac{1+\epsilon}{\epsilon}\Bigr)^{8}\biggr)\biggr),
    	\\ \nonumber
& \hspace{1.5cm} (1+\epsilon)t < \eta_1 < \eta_2 < \infty, \quad \epsilon > 0, \quad 0 \leq \sigma \leq 1, \quad t \to \infty,
\end{align}
where the error term is uniform for all $\eta_1, \eta_2, \epsilon, \sigma$ in the above ranges.
\end{theorem}
\begin{proof}
Replacing in equation (\ref{zetaformula2.1}) $\eta$ with $\eta_1$ and subtracting the resulting equation from the equation obtained from (\ref{zetaformula2.1}) by replacing $\eta$ with $\eta_2$, we find (\ref{sumzetaformula2.1}). Equation (\ref{sumzetaformula2.1b}) follows in a similar way from (\ref{zetaformula2.1b}).
\end{proof}

\begin{theorem}\label{th5.2} 
Define the polylogarithm $\Li_m(z)$ by (\ref{polylogdef}). 
The following relation holds:
\begin{align}\label{th52formula}
\sum_{n =[\frac{t}{2\pi}] +1}^{[\frac{\eta}{2\pi}]} n^{-s} 
= &\; \frac{1}{1-s} \biggl[\left(\frac{\eta}{2\pi}\right)^{1-s} - \left(\frac{t}{2\pi}\right)^{1-s}\biggr] 
	\\ \nonumber
&+ \frac{e^{-\frac{i\pi (1-s)}{2}}}{(2\pi)^{1-s}}\sum_{n=1}^\infty  \sum_{j=0}^{N-1} e^{-nz - it\ln{z}} \left(\frac{1}{n + \frac{it}{z}}\frac{d}{dz}\right)^j\frac{z^{-\sigma}}{n + \frac{it}{z}} \Biggr|_{z = i\eta}^{it}
	\\\nonumber
&- \frac{e^{\frac{i\pi (1-s)}{2}}}{(2\pi)^{1-s}}\sum_{n=1}^\infty
\sum_{j=0}^{N-1} e^{-nz - it\ln{z}} \left(\frac{1}{n + \frac{it}{z}}\frac{d}{dz}\right)^j\frac{z^{-\sigma}}{n + \frac{it}{z}} \Biggr|_{z = -i\eta}
	\\ \nonumber
& + \frac{e^{\frac{i\pi (1-s)}{2}}}{(2\pi)^{1-s}}\sum_{n=2}^\infty
\sum_{j=0}^{N-1} e^{-nz - it\ln{z}} \left(\frac{1}{n + \frac{it}{z}}\frac{d}{dz}\right)^j\frac{z^{-\sigma}}{n + \frac{it}{z}} \biggr|_{z = -it}
	\\ \nonumber
& + \left(\frac{t}{2\pi}\right)^{1-s} e^{it}
\sum_{k=0}^{2N}   \frac{\overline{c_k(1-\sigma)}\Gamma(\frac{k+1}{2})}{t^{\frac{k+1}{2}}}
	\\ \nonumber
 & +  O\biggl( \frac{(2N +1)!! N2^{2N}}{t^{\sigma + N}}  
 + \frac{(2N +1)!! N (\frac{1+\epsilon}{\epsilon})^{2(N+1)}}{\eta^{\sigma + N}}  \biggr),
	\\ \nonumber
& \hspace{1cm} (1+\epsilon)t < \eta < \infty, \quad \epsilon > 0, \quad 0 \leq \sigma \leq 1, \quad N \geq 2, \quad t \to \infty,	
\end{align}
where the error term is uniform for all $\eta, \epsilon, \sigma, N$ in the above ranges and the coefficients $c_k(\sigma)$ are defined in (\ref{rhoquotientexpansion}).

For $N=3$ equation (\ref{th52formula}) simplifies to  
\begin{align}\label{th52formulab}
\sum_{n =[\frac{t}{2\pi}] +1}^{[\frac{\eta}{2\pi}]} n^{-s} 
= \;&  \frac{1}{1-s} \biggl[\left(\frac{\eta}{2\pi}\right)^{1-s} - \left(\frac{t}{2\pi}\right)^{1-s}\biggr] 
+ \Bigl\{ \cdots \Bigr\} 
	\\ \nonumber
& + O\biggl(\frac{1}{\eta^{3 + \sigma}}\biggl(t^3 + \Bigl(\frac{1+\epsilon}{\epsilon}\Bigr)^{8}\biggr) + \frac{1}{t^{3 + \sigma}}\biggr),
	\\ \nonumber
& \hspace{2cm} (1+\epsilon)t < \eta < \infty, \quad \epsilon > 0, \quad 0 \leq \sigma \leq 1, \quad t \to \infty,	
\end{align}
where the error term is uniform for all $\eta, \epsilon, \sigma$  in the above ranges and $\{ \cdots\}$ denotes an explicit expression involving $\sigma$, $t$, and the polylogarithms $\Li_m(e^{it})$, $m = 1, \dots, 5$.
\end{theorem}
\begin{proof}
Equation (\ref{th52formula}) follows by subtracting equation (\ref{zetaformula2.2}) from equation (\ref{zetaformula2.1}). Similarly, equation (\ref{th52formulab}) follows by subtracting equation (\ref{zetaformula2.2b}) from equation (\ref{zetaformula2.1b}).

\end{proof}

The following equation is given on page 78 of Titchmarsh \cite{T1986}:
\begin{align}\label{Titchmarshpage78}
  \sum_{x < n \leq N} \frac{1}{n^s} 
  \sim \chi(s) \sum_{\frac{t}{2\pi N} < n \leq \frac{t}{2\pi x}} \frac{1}{n^{1-s}}.
\end{align}
Theorem \ref{ZETATH} yields a precise version of the relation (\ref{Titchmarshpage78}) with an explicit error estimate.

\begin{theorem}\label{th5.3}
For every $\epsilon > 0$, there exists a constant $A>0$ such that
\begin{align}\label{sumzetaformula}
\sum_{n=[\frac{t}{\eta_2}]+1}^{[\frac{t}{\eta_1}]} \frac{1}{n^s}  
= \;& \chi(s) \sum_{n=[\frac{\eta_1}{2\pi}]+1}^{[\frac{\eta_2}{2\pi}]} \frac{1}{n^{1-s}}
 	\\ \nonumber
& + e^{-i\pi s}\Gamma(1-s)\biggl\{ e^{\frac{i\pi(s-1)}{2}}\eta^{s-1}  e^{\frac{2t}{\eta}[\frac{\eta}{2\pi}]\pi i - it 
- \frac{it}{2\eta^2}(2[\frac{\eta}{2\pi}]\pi - \eta)^2} S_N(s, \eta)
	\\ \nonumber
 & \hspace{4cm} + O\biggl( e^{-\frac{\pi t}{2}} \left(\frac{3N}{t}\right)^{\frac{N}{6}} \frac{\eta^\sigma}{\sqrt{t}}\biggr)\biggr\}\Biggr|_{\eta = \eta_1}^{\eta_2},
 	\\ \nonumber
& \qquad \epsilon \sqrt{t} < \eta_1 < \eta_2 < t , \quad 0 \leq \sigma \leq 1, \quad 1 \leq N < At, \quad  t \to \infty,
\end{align}
where the error term is uniform for all $\eta_1, \eta_2, \sigma, N$ in the above ranges and the function $S_N(s, \eta)$ is defined by (\ref{SNdef}).

For $N=3$ equation (\ref{sumzetaformula}) simplifies to  
\begin{align}\label{sumzetaformulab}
\sum_{n=[\frac{t}{\eta_2}]+1}^{[\frac{t}{\eta_1}]} \frac{1}{n^s}  = \;& \chi(s) \sum_{n=[\frac{\eta_1}{2\pi}]+1}^{[\frac{\eta_2}{2\pi}]} \frac{1}{n^{1-s}}
 	\\ \nonumber
& + e^{-i\pi s}\Gamma(1-s)\biggl\{ e^{\frac{i\pi(s-1)}{2}}\eta^{s-1}  e^{\frac{t}{\eta}2[\frac{\eta}{2\pi}]\pi i - it 
- \frac{it}{2\eta^2}(2[\frac{\eta}{2\pi}]\pi - \eta)^2} 
	\\ \nonumber
&\qquad \times \biggl[\Phi + \frac{\sigma -1}{i\eta} \biggl(
  \partial_2 \Phi  
  + 
\biggl(2\Bigl[\frac{\eta}{2\pi}\Bigr]\pi i - i\eta\biggr)   \Phi  \biggr)
  	\\ \nonumber
&  \qquad\qquad -\frac{(\sigma -2) (\sigma -1)}{2 \eta ^2} \biggl(
  \partial_2^{2} \Phi  
 + 2 \biggl(2\Bigl[\frac{\eta}{2\pi}\Bigr]\pi i - i\eta\biggr)  \partial_2 \Phi  
	\\ \nonumber
&\qquad\qquad  + \biggl(2\Bigl[\frac{\eta}{2\pi}\Bigr]\pi i - i\eta\biggr)^2 \Phi     \biggr) \biggr]
+ O\biggl( e^{-\frac{\pi t}{2}} \frac{\eta^\sigma}{t} \biggr)\biggr\}\Biggr|_{\eta = \eta_1}^{\eta_2},
 	\\ \nonumber
& \qquad 
 \epsilon \sqrt{t} < \eta_1 < \eta_2 < t , \quad 0 \leq \sigma \leq 1, \quad t \to \infty,
\end{align}
where the error term is uniform for all $\eta_1, \eta_2, \sigma$ in the above ranges, $\Phi$ is defined by (\ref{Phidef}), and $\Phi$ and its partial derivatives are evaluated at the point (\ref{Phievaluationpoint}).
\end{theorem}
\begin{proof}
Replacing in equation (\ref{zetaformula}) $\eta$ with $\eta_1$ and subtracting the resulting equation from the equation obtained from (\ref{zetaformula}) by replacing $\eta$ with $\eta_2$, we find (\ref{sumzetaformula}). Equation (\ref{sumzetaformulab}) follows in a similar way from (\ref{zetaformulab}).
\end{proof}

\begin{remark}\upshape
The relation (\ref{Titchmarshpage78}) is a particular case of (\ref{sumzetaformula}) (let $x = t/\eta_2$ and $N = t/\eta_1$).
\end{remark}

\part{Asymptotics to all Orders of a Two-Parameter Generalization of the Riemann Zeta Function}

\chapter{An Exact Representation for $\Phi(u,v,\beta)$}\label{sec6}
The two-parameter generalization $\Phi(u,v,\beta)$ of $\zeta(s)$ was defined in (\ref{X1.1}). In Theorem \ref{theorem2.1}, we derived an exact representation for $\zeta(s)$. In this chapter, we prove Theorem \ref{exactPhith} which provides an analogous representation for $\Phi(u,v,\beta)$. 

Recall that we make the basic assumption (\ref{basicassumption}); in particular $\eta \notin 2 \pi \Z$. The branch cut for the logarithm is assumed to run along the negative real axis.

\begin{theorem}[An exact representation for $\Phi$]\label{exactPhith}
  Let $\Phi(u,v,\beta)$ be defined by (\ref{X1.1}). 
Then
\begin{align}\nonumber
  \Phi(u&, v,\beta) = (e^{-i\pi u} - 1)e^{-i\pi v}(2\pi e^{\frac{i\pi}{2}})^{u+v-1} \sum_{m=1}^{[\frac{\eta}{2\pi}]} m^{u-1} (m+\beta)^{v-1}
  	\\ \nonumber
& + e^{-i\pi v} \bigg((1 - e^{i\pi u})\int_{-i\eta}^{\infty e^{i\phi_1}} + (e^{-i\pi u} - 1)\int_{i\eta}^{\infty e^{i\phi_2}}\bigg)z^{u-1} (z + 2i\pi \beta)^{v-1} \frac{dz}{e^z - 1}
	\\ \nonumber
& + \int_{L^3} \frac{z^{u-1}}{e^{-z} -1} \big[(z-2i\pi \beta)^{v-1} + e^{-i\pi v} e^{-z} (z + 2i\pi \beta)^{v-1}\big] dz
	\\\nonumber
& + (e^{-i \pi u} - 1) \int_{\hat{C}_\eta^\alpha} \frac{z^{u-1}}{e^{-z} -1} \big[(z-2i\pi \beta)^{v-1} + e^{-i\pi v}e^{-z} (z + 2i\pi \beta)^{v-1}\big]dz,
	\\ \label{X2.12}
&\hspace{1cm} u, v \in \C, \quad \beta > 0, \quad 0<\eta<2\pi \beta, \quad - \frac{\pi}{2} < \phi_j < \frac{\pi}{2}, \ \ j = 1,2,
\end{align}
where $0 < \alpha < 2\pi \min(1,\beta)$ and the contours $L_3$ and $\hat{C}_\eta^\alpha$ with the orientations shown in figures \ref{L1L2L3.pdf} and \ref{Calphaa.pdf} are defined in (\ref{2.7}) and (\ref{2.12c}), respectively.
\end{theorem}
\begin{proof}
We decompose the contour $H_\alpha$ in the definition (\ref{X1.1}) of $\Phi$ into the union of the three contours $\{L_j\}_1^3$ defined in (\ref{2.7}) with the orientation shown in figure \ref{L1L2L3.pdf}:
$$H_\alpha = L_1 \cup L_2 \cup L_3.$$

Proceeding in analogy with the proof of Theorem \ref{theorem2.1}, we write the integral along $L_1$ as follows:
\begin{align}\nonumber
& \int_{L_1} z^{u-1}(z-2i\pi \beta)^{v-1}\frac{dz}{e^{-z} -1} 
 = \int_{i\alpha}^{\infty e^{i\pi}} z^{u-1} (z - 2i\pi \beta)^{v-1} \frac{dz}{e^{-z} - 1}
 	\\ \nonumber
& = e^{i\pi u} \int_{-i\alpha}^\infty \zeta^{u-1}(-\zeta - 2i\pi \beta)^{v-1} \frac{d\zeta}{e^\zeta -1} 
	\\ \label{X2.15}
& = e^{i\pi u} e^{-i\pi(v-1)} \int_{-i\alpha}^\infty \zeta^{u-1} (\zeta + 2i\pi \beta)^{v-1} \frac{d\zeta}{e^\zeta -1},
\end{align}
where the first equality is a consequence of Cauchy's theorem, the second equality is a consequence of the substitution $z = e^{i\pi} \zeta$, and the third equality uses the fact that $-\zeta - 2i\pi \beta$ is in the third quadrant. 

Similarly, the integral along $L_2$ can be written in the form
\begin{align}\label{X2.16}
  \int_{L^2} z^{u-1}(z - 2i\pi \beta)^{v-1} \frac{dz}{e^{-z} -1} 
  = e^{-i\pi u} e^{-i\pi(v-1)} \int_{\infty}^{i\alpha} \zeta^{u-1} (\zeta + 2i\pi \beta)^{v-1} \frac{d\zeta}{e^\zeta -1}.
\end{align}
The starting point for deriving this identity is the application of Cauchy's theorem in the domain enclosed by $L_2$ and by the ray from $\infty e^{-i\pi}$ to $-i\alpha$, i.e., in the shaded domain $2$ of figure \ref{domains123.pdf}.

The integral along $L_3$ can be written as follows:
\begin{align} \nonumber
  \int_{L_3} z^{u-1} &(z - 2i\pi \beta)^{v-1} \frac{dz}{e^{-z} - 1}
  = -e^{-i\pi(v-1)} \int_{L_3} \zeta^{u-1}(\zeta + 2i\pi \beta)^{v-1} \frac{d\zeta}{e^\zeta -1}
  	\\ \nonumber
&  + \int_{L_3} \zeta^{u-1} \bigg[\frac{(\zeta - 2i\pi \beta)^{v-1}}{e^{-\zeta} -1} + e^{-i\pi(v-1)} \frac{(\zeta + 2i\pi \beta)^{v-1}}{e^\zeta -1}\bigg] d\zeta
  	\\ \nonumber
 = & -e^{-i\pi(v-1)} \bigg(\int_{-i\alpha}^\infty + \int_\infty^{i\alpha}\bigg) \zeta^{u-1}(\zeta + 2i\pi \beta)^{v-1} \frac{d\zeta}{e^\zeta -1}
	\\\label{X2.17}
& + \int_{L_3} \frac{\zeta^{u-1}}{e^{-\zeta} -1} \big[(\zeta - 2i\pi \beta)^{v-1} + e^{-i\pi v}e^{-\zeta}(\zeta + 2i\pi \beta)^{v-1}\big] d\zeta,
\end{align}
where the first equality is an identity and the second equality follows from an application of Cauchy's theorem in the domain enclosed by $L_3$ and the two rays $(-i\alpha,\infty)$ and $(\infty, i\alpha)$, i.e., in the shaded domain $3$ of figure \ref{domains123.pdf}. 

Adding equations (\ref{X2.15})-(\ref{X2.17}) we obtain
\begin{align}\nonumber
\Phi(u,v,\beta) = & -2ie^{-i\pi v} \sin\Big(\frac{\pi u}{2}\Big) \bigg(e^{\frac{i\pi u}{2}} \int_{-i\alpha}^\infty + e^{-\frac{i\pi u}{2}} \int_{i\alpha}^\infty\bigg) z^{u-1}(z+2i\pi \beta)^{v-1} \frac{dz}{e^z -1}
	\\ \label{X2.18}
& + \int_{L_3} \frac{z^{u-1}}{e^{-z} -1}\big[(z - 2i\pi \beta)^{v-1} + e^{-i\pi v}e^{-z}(z + 2i\pi \beta)^{v-1}\big] dz.
\end{align}
\begin{figure}
 \begin{overpic}[width=.5\textwidth]{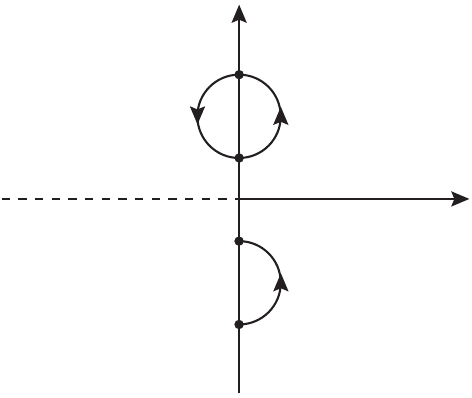}
 \put(44,70){$i\eta$}
 \put(43,47){$i\alpha$}
 \put(38.5,32){$-i\alpha$}
 \put(39,14){$-i\eta$}
 \put(62.5,58){$C_\alpha^\eta$}
 \put(31.5,58.5){$\hat{C}_\eta^\alpha$}
 \put(62.5,22.5){$C_{-\eta}^{-\alpha}$}
 \put(101.5,40.5){$\re z$} 
   \end{overpic}
   \caption{The contours $C^\eta_\alpha$, $C^{-\alpha}_{-\eta}$ and $\hat C^\alpha_\eta$.}\label{XCalphaa.pdf}
\end{figure}
Defining the contours $C_\alpha^\eta, C_{-\eta}^{-\alpha}, \hat{C}_\eta^\alpha$ as in (\ref{Cetaalphadef}) with the orientations shown in figure \ref{XCalphaa.pdf}, the second integral in the rhs of (\ref{X2.18}) can be rewritten as the sum of an integral along the contour $C_\alpha^\eta$ plus an integral along the ray $(i\eta, \infty e^{i\phi_1})$. Similarly, the first integral in the rhs of (\ref{X2.18}) can be rewritten as the sum of an integral along the curve $-C_{-\eta}^{-\alpha}$ plus an integral along the ray $(-i\eta, \infty e^{i\phi_2})$. 
Hence the first two terms in the rhs of (\ref{X2.18}) yield the second line in the rhs of (\ref{X2.12}), as well as the additional term $I$ defined by
\begin{align}\label{X2.20}
  I = -2ie^{-i\pi v} \sin\Big(\frac{\pi u}{2}\Big) \bigg(-e^{\frac{i\pi u}{2}}\int_{C_{-\eta}^{-\alpha}} + e^{-\frac{i\pi u}{2}} \int_{C_\alpha^\eta}\bigg) z^{u-1}(z+ 2i\pi \beta)^{v-1} \frac{dz}{e^z -1}.
\end{align}
Letting $z = \zeta e^{-i\pi}$ in the integral involving $C_{-\eta}^{-\alpha}$ and using the fact that $-\zeta + 2i\pi \beta$ is in the first quadrant we find
\begin{align}\label{X2.21}
-e^{\frac{i\pi u}{2}}\int_{C_{-\eta}^{-\alpha}} z^{u-1}(z+2i\pi \beta)^{v-1} \frac{dz}{e^z -1}
= e^{-\frac{i\pi u}{2}} e^{i\pi v} \int_{\hat{C}_\eta^\alpha} z^{u-1}(z - 2i\pi \beta)^{v-1} \frac{dz}{e^{-z} -1}.
\end{align}
Using the above equation in the rhs of (\ref{X2.20}) and then adding and subtracting in the resulting equation the term
$$-2i\sin\Big(\frac{\pi u}{2}\Big) e^{-\frac{i\pi u}{2}} e^{-i\pi v} \int_{\hat{C}_\eta^\alpha} z^{u-1} (z + 2i\pi \beta)^{v-1} \frac{dz}{e^z -1},$$
we find that $I$ is given by
\begin{align*}
  I = & -2i\sin\Big(\frac{\pi u}{2}\Big) e^{-\frac{i\pi u}{2}} \bigg\{e^{-i\pi v} \int_{C_\alpha^\eta \cup \hat{C}_\eta^\alpha} z^{u-1}(z + 2i\pi \beta)^{v-1} \frac{dz}{e^z -1}
  	\\
& + \int_{\hat{C}_\eta^\alpha} z^{u-1}\bigg(\frac{(z-2i\pi \beta)^{v-1}}{e^{-z} -1} - \frac{e^{-i\pi v}(z + 2i\pi \beta)^{v-1}}{e^z -1}\bigg)\bigg\}.
\end{align*}
That is,
\begin{align}\nonumber
  I = &\; (e^{-i\pi u} -1)\bigg\{e^{-i\pi v}\int_{C_\alpha^\eta \cup \hat{C}_\eta^\alpha} z^{u-1}(z + 2i \pi \beta)^{v-1} \frac{dz}{e^z -1}
  	\\ \label{X2.22}
& + \int_{\hat{C}_\eta^\alpha} \frac{z^{u-1}}{e^{-z} -1} \big[(z-2i\pi \beta)^{v-1} + e^{-i\pi v} e^{-z} (z + 2i\pi \beta)^{v-1}\big]\bigg\}.
\end{align}
Cauchy's theorem implies that the first integral in the curly bracket in (\ref{X2.22}) equals
$$(2\pi e^{\frac{i\pi}{2}})^{u+v-1} \sum_{m=1}^{[\frac{\eta}{2\pi}]} m^{u-1} (m+\beta)^{v-1}.$$
Thus (\ref{X2.18}) with the aid of equation (\ref{X2.22}) becomes equation (\ref{X2.12}).
\end{proof}

\chapter{The Asymptotics of $\Phi(u,v,\beta)$}\label{sec7}
In Theorem \ref{zetath2}, we established an asymptotic formula for $\zeta(s)$ for $\epsilon < \eta < \sqrt{t}$. In this chapter, we establish an analogous asymptotic formula for $\Phi$.

Let $u = \sigma_1 + it$ and $v = \sigma_2 - it$. The function $\Phi(u,v,\beta)$ was defined in (\ref{X1.1}) by
\begin{align}\label{XPhidefrepeat}
\Phi(u,v,\beta) = \int_{H_\alpha} z^{u-1}(z-2i\pi \beta)^{v-1} \frac{dz}{e^{-z} -1}, \qquad u,v \in \C,  \quad \beta \in \R \setminus \{0\},
\end{align}
where $H_\alpha$, $0 < \alpha < 2\pi\min(1, |\beta|)$, denotes the Hankel contour (\ref{X1.2}) surrounding the negative real axis in the counterclockwise direction, and the complex powers are defined by $w^a = e^{a(\ln|w| + i\arg w)}$ with $\arg w \in (-\pi, \pi]$, i.e., the branch cut runs along the negative real axis. In this and the following chapter, we prefer to have the branch cut along the positive real axis. Therefore we change variables $z = -w$ in (\ref{XPhidefrepeat}) to get
\begin{align}\label{PhihatHalphaexpression}
\Phi(u,v,\beta) = -e^{-i\pi(u+v)} \int_{\hat{H}_\alpha} w^{u-1}(w + 2i\pi \beta)^{v-1} \frac{dw}{e^{w} -1}, \quad \; \; u,v \in \C,  \;\; \beta \in \R \setminus \{0\},
\end{align}
where $\hat{H}_\alpha$ denotes the Hankel contour surrounding the positive real axis in the counterclockwise direction, i.e., 
\begin{align}\label{hatHalphadef}
\hat{H}_\alpha = \left\{r + i0 \, | \,  \alpha < r <\infty\right\} \cup  \left\{\alpha e^{i\theta} \, | \, 0 <\theta < 2\pi\right\} \cup \left\{r-i0 \, | \,  \alpha<r<\infty\right\},
\end{align}
and the complex powers are defined by $w^a = e^{a(\ln|w| + i\arg w)}$ with $\arg w \in [0,2\pi)$, i.e., the branch cut now runs along the positive real axis.

Using the identity
$$\frac{1}{e^w -1} = \sum_{n=1}^m e^{-nw} + \frac{e^{-mw} }{e^w -1} , \qquad m = 1,2, \dots,$$
we may write 
\begin{align}\nonumber
\Phi(u,v,\beta) = -e^{-i\pi(u+v)} \bigg\{&\sum_{n=1}^m \int_{\hat{H}_\alpha} w^{u-1}(w + 2i\pi \beta)^{v-1} e^{-nw} dw
	\\ \label{XPhisumC}
& + \int_{\hat{H}_\alpha}  \frac{w^{u-1}(w + 2i\pi \beta)^{v-1} e^{-mw} dw}{e^{w} -1}\bigg\}.
\end{align}
We are interested in the asymptotic behavior of $\Phi(u,v,\beta)$ as $t \to \infty$.
Thus we note that the second integral on the rhs of (\ref{XPhisumC}) has critical points at the solutions of
$$\frac{d}{dw}\ln\big(w^{u-1}(w + 2i\pi \beta)^{v-1} e^{-mw}\big) = 0,$$
that is, at
$$w = \frac{-2 i \pi  \beta m+u+v-2 \pm \sqrt{(-2 i \pi  \beta m+u+v-2)^2+8 i \pi  \beta
   m (u-1)}}{2 m}.$$
Anticipating that $m$ will grow like $t^\alpha$, $\alpha > 0$, we find that for large $t$ the critical points are approximately given by
$$i \beta \pi \bigg(-1\pm \sqrt{1 + \frac{2t}{\beta \pi m}}\bigg).$$
By performing a steepest descent analysis we can find the asymptotics of $\Phi$ to all orders. The idea is that given some $\eta > 0$, we choose $m$ so that the critical point for large $t$ lies at $i\eta$.
Solving the equation 
$$\eta = \beta \pi \bigg(-1 + \sqrt{1 + \frac{2t}{\beta \pi m}}\bigg)$$
for $m$ we find that $m$ should be given by
$$t\bigg(\frac{1}{\eta} - \frac{1}{\eta + 2\pi \beta}\bigg).$$
But since $m$ has to be an integer, we instead use the approximate definition
$$m = [x] \quad \text{where} \quad x = t\bigg(\frac{1}{\eta} - \frac{1}{\eta + 2\pi \beta}\bigg).$$

\begin{theorem}[{\bf The asymptotics of $\Phi(u,v,\beta)$ to all orders for  $\epsilon < \eta < \sqrt{t}$}]\label{Xzetath2}
For every $\epsilon > 0$, there exists a constant $A > 0$ such that
\begin{align} \nonumber
& -e^{i\pi(\sigma_1 + \sigma_2)} \Phi(u,v,\beta) 
 = \sum_{n=1}^m \int_{\hat{H}_\alpha} w^{u-1}(w + 2i\pi \beta)^{v-1} e^{-nw} dw
	\\  \nonumber
& + i (2\pi)^{\sigma_1 + \sigma_2 - 1} e^{\frac{\pi i}{2}(\sigma_1 + \sigma_2)} \sum_{n=1}^{[\frac{\eta}{2\pi}]}   n^{u-1} (n+\beta)^{v-1} 
 	\\ \nonumber
&-e^{-(m+1)i\eta} (i\eta)^{u-1} (i\eta + 2\pi i \beta)^{v-1}  e^{\frac{i\pi}{4}} 
	\\ \nonumber
& \hspace{3cm} \times \sum_{k=0}^{[\frac{N-1}{2}]} \frac{\varphi^{(2k)}(0)}{(2k)!} i^k 
\biggl(\frac{t}{2}\bigg(\frac{1}{\eta^2} - \frac{1}{(\eta + 2\pi \beta)^2}\bigg)\biggr)^{-k - \frac{1}{2}} \Gamma\biggl(k + \frac{1}{2}\biggr)
	\\ \nonumber
& + \eta^{\sigma_1 -1} (\eta + 2\pi \beta)^{\sigma_2 -1}
	\\ \nonumber
& \hspace{.5cm} \times \begin{cases}	
O\biggl(\left(\frac{22N}{t}\right)^{\frac{N}{6}} \frac{\eta}{\sqrt{t}}\biggr), & \epsilon < \eta < t^{\frac{1}{3}} < \infty, \quad 1 \leq N < \frac{At}{\eta^3}, \\
  O\biggl(N e^{- \frac{At}{\eta^2}} + \big(\frac{3t}{4N\eta^2}\big)^{-\frac{N+1}{2}} \biggr), & t^{\frac{1}{3}} < \eta < \sqrt{t} < \infty, \quad 1 \leq N < \frac{At}{\eta^2}, \end{cases}   \qquad t \to \infty,
 	\\ \label{Xzetaformula2}	
& \beta > 0, \quad \epsilon < \eta < 2\pi \beta - \epsilon < \infty, \quad \dist(\eta, 2\pi \Z) > \epsilon, \quad \sigma_1 \in [0,1], \quad \sigma_2 \in [0,1], 
\end{align}
where $\hat{H}_\alpha$ denotes the Hankel contour defined in (\ref{hatHalphadef}),
\begin{align*}
& m = \bigg[t\bigg(\frac{1}{\eta} - \frac{1}{\eta + 2\pi \beta}\bigg)\bigg],
	\\
& \varphi(z) = \frac{e^{(u-1)\ln(1 + \frac{z}{i\eta}) + (v-1)\ln(1 + \frac{z}{i\eta + 2\pi i\beta}) - x z - \frac{it}{2}(\frac{1}{\eta^2} - \frac{1}{(\eta + 2\pi \beta)^2})z^2 + (x-m)z}}{e^z - e^{-i\eta}},
\end{align*}
and the error terms are uniform with respect to $\eta, \sigma_1, \sigma_2, \beta, N$ in the given ranges.
\end{theorem}
\begin{proof}
Let $\epsilon > 0$ be given and suppose that $8\epsilon < \eta < \min(\sqrt{t}, 2\pi \beta - 8 \epsilon)$, $\sigma_1 \in [0,1]$, $\sigma_2 \in [0,1]$, $\dist(\eta, 2\pi \Z) > 8\epsilon$, and $N \geq 1$ (since $\epsilon > 0$ is arbitrary, we can replace $8\epsilon$ with $\epsilon$ at the end). All error terms of the form $O(\cdot)$ will be uniform with respect to $\eta, \sigma_1, \sigma_2, \beta, N$ in the given ranges, but not with respect to $\epsilon$. We let $A>0$ denote a generic constant which can change within a computation.

Let
$$x = t\bigg(\frac{1}{\eta} - \frac{1}{\eta + 2\pi \beta}\bigg), \qquad m = [x].$$
Since $\eta < \sqrt{t}$ and $\eta < 2\pi \beta$, we have
$$x = \frac{t}{\eta} \frac{1}{\frac{\eta}{2\pi \beta} +1} \geq \frac{\sqrt{t}}{2}. $$
In particular, $m \to \infty$ as $t \to \infty$.
We deform the Hankel contour  $\hat{H}_\alpha$ into the straight lines $C_j$, $j = 1, \dots, 4$ joining $\infty$, $c\eta + i\eta(1+c)$, $-c\eta + i\eta(1-c)$, $-c\eta - i\eta$, $\infty$, where $0 < c \leq 1/2$ is an absolute constant, see figure \ref{Cjs.pdf}.
Then
\begin{align}\label{XintCww}
\int_{\hat{H}_\alpha}  \frac{w^{u-1}(w + 2i\pi \beta)^{v-1} e^{-mw}}{e^{w} -1} dw
= &-2\pi i \sum_{\substack{n = - [\frac{\eta}{2\pi}]\\n \neq 0}}^{[\frac{\eta}{2\pi}]} (2\pi i n)^{u-1} (2\pi i(n+\beta))^{v-1} + \sum_{j=1}^4 I_j,
\end{align}
where
$$I_j = \int_{C_j} \frac{w^{u-1}(w + 2i\pi \beta)^{v-1} e^{-mw}}{e^{w} -1} dw, \qquad j = 1, \dots, 4.$$
Note that
\begin{align*}
-2\pi i \sum_{\substack{n = - [\frac{\eta}{2\pi}]\\n \neq 0}}^{[\frac{\eta}{2\pi}]} & (2\pi i n)^{u-1} (2\pi i(n+\beta))^{v-1}
= - (2\pi)^{\sigma_1 + \sigma_2 - 1} e^{\frac{i\pi}{2}(u+v -1)} \sum_{n=1}^{[\frac{\eta}{2\pi}]}  n^{u-1} (n+\beta)^{v-1}
	\\
& -(2\pi)^{\sigma_1 + \sigma_2 - 1} e^{\frac{\pi i}{2}v}  e^{\frac{3 \pi i}{2}(u-1)} \sum_{n = 1}^{[\frac{\eta}{2\pi}]} n^{u-1} (-n+\beta)^{v-1}
	\\
 = &\;  i (2\pi)^{\sigma_1 + \sigma_2 - 1} e^{\frac{\pi i}{2}(\sigma_1 + \sigma_2)} \sum_{n=1}^{[\frac{\eta}{2\pi}]}   n^{u-1} \big[(n+\beta)^{v-1} - e^{\pi i u} (-n+\beta)^{v-1}\big]
	\\
 = &\;  i (2\pi)^{\sigma_1 + \sigma_2 - 1} e^{\frac{\pi i}{2}(\sigma_1 + \sigma_2)} \sum_{n=1}^{[\frac{\eta}{2\pi}]}   n^{u-1} (n+\beta)^{v-1} + O(e^{-At}).	
\end{align*}
In view of (\ref{XPhisumC}) and (\ref{XintCww}) this gives the first two terms on the rhs of (\ref{Xzetaformula2}). It remains to analyze the contributions from the $I_j$'s. 

We first prove that $I_1, I_3, I_4$ are exponentially small as $t \to \infty$.
Let $w = \rho e^{i\phi}$, $0 < \phi < 2\pi$, and $w + 2i\pi \beta = \tilde{\rho} e^{i\tilde{\phi}}$, $0 < \tilde{\phi} < 2\pi$. Then
$$|w^{u-1}| = \rho^{\sigma_1 - 1} e^{-t\phi}, \qquad
|(w + 2i\pi \beta)^{v-1}| = \tilde{\rho}^{\sigma_2 - 1} e^{t\tilde{\phi}}.$$
Also,
$$x\eta = \frac{tQ}{1 + Q}, \quad \text{where} \quad Q := \frac{2\pi \beta}{\eta} > 1.$$

\bigskip
\noindent
{\bf Analysis of the integral $I_4$}\nopagebreak
\medskip

\noindent  
For $w \in C_4$ we have 
$$|e^w -1|>A, \qquad \rho \geq \eta, \qquad \tilde{\rho} \geq 2\pi \beta - \eta, \quad
\phi \geq \pi + \arctan\frac{\eta}{c\eta},$$
and
$$\tilde{\phi} \leq \pi - \arctan\frac{2\pi \beta - \eta}{c\eta} = \pi - \arctan\frac{Q - 1}{c}.$$
Hence
\begin{align}\nonumber
|w^{u-1}(w + 2i\pi \beta)^{v-1}|
& = \rho^{\sigma_1 -1} \tilde{\rho}^{\sigma_2 -1} e^{-t(\phi - \tilde{\phi})}
	\\ \nonumber
& \leq \eta^{\sigma_1 -1} (2\pi \beta -\eta)^{\sigma_2 -1}e^{-t(\arctan\frac{1}{c} + \arctan\frac{Q - 1}{c})}
	\\ \label{XC4est}
&= O(e^{-t(\arctan\frac{1}{c} + \arctan\frac{Q - 1}{c})}).
\end{align}
Thus,
\begin{align*}
I_4 & = \int_{C_4} \frac{w^{u-1}(w + 2i\pi \beta)^{v-1} e^{-mw} }{e^{w} -1} dw 
	\\
& = O\bigg(e^{-t(\arctan\frac{1}{c} + \arctan\frac{Q - 1}{c})} \int_{-c\eta}^\infty e^{-mw_1}dw_1\bigg)
	\\
& = O(e^{-t(\arctan\frac{1}{c} + \arctan\frac{Q - 1}{c}) + mc\eta})
	\\
& = O(e^{-t(\arctan\frac{1}{c} + \arctan\frac{Q - 1}{c}) + xc\eta})  = O(e^{-tF(Q)}),
\end{align*}
where 
$$F(Q) = \arctan\frac{1}{c} + \arctan\frac{Q - 1}{c} - \frac{cQ}{1 + Q}.$$
Since $\frac{c^2}{4} < Q$, we have
$$F'(Q) = \frac{4 c Q-c^3}{(Q+1)^2\left(c^2+(Q-1)^2\right)} > 0, \qquad Q > 1.$$
Hence 
$$\inf_{Q > 1} F(Q) \geq F(1) = \arctan\bigg(\frac{1}{c}\bigg) - \frac{c}{2} > A > 0.$$
It follows that $I_4 = O(e^{-At})$.

\bigskip
\noindent
{\bf Analysis of the integral $I_3$}\nopagebreak
\medskip

\noindent 
For $w = w_1 + iw_2 \in C_3$ we have 
$$\phi =  \pi - \arctan \frac{w_2}{c\eta}, \qquad \tilde{\phi} = \pi - \arctan \frac{w_2 + 2\pi \beta}{c\eta}.$$
The function $\phi - \tilde{\phi}$ assumes its maximum at $w_2 = -\pi \beta$ and decreases symmetrically as $w_2$ moves away from this point. Its minimum on $C_3$ is therefore assumed at the upper endpoint where $w_2 = (1-c)\eta$.
For $w_2 = (1-c)\eta$ we have
\begin{align}\label{XphiminustildephionC3}
\phi - \tilde{\phi} = \arctan \frac{w_2 + 2\pi \beta}{c\eta} - \arctan \frac{w_2}{c\eta}
= \arctan \bigg(\frac{1-c + Q}{c}\bigg) - \arctan \frac{1-c}{c}.
\end{align}
Hence
\begin{align*}
w^{u-1}(w + 2i\pi \beta)^{v-1} e^{-mw}
& = O\Big((c\eta)^{\sigma_1 -1} (c\eta)^{\sigma_2 -1} e^{-t(\arctan (\frac{1-c + Q}{c}) - \arctan \frac{1-c}{c})} e^{mc\eta}\Big)
	\\
& = O\Big(e^{-t(\arctan (\frac{1-c + Q}{c}) - \arctan \frac{1-c}{c})} e^{xc\eta}\Big)
= O\big(e^{-tF(Q)}\big),
\end{align*}
where 
$$F(Q) = \arctan \bigg(\frac{1-c + Q}{c}\bigg) - \arctan\bigg(\frac{1-c}{c}\bigg) - \frac{cQ}{1 + Q}.$$
Since $c < 1+Q$, we have
$$F'(Q) = -\frac{2 c^2 (c-Q-1)}{(Q+1)^2 (c^2 + (1+Q-c)^2)} > 0, \qquad Q > 1.$$
Hence 
$$\inf_{Q > 1} F(Q) \geq F(1) = \arctan\left(\frac{2-c}{c}\right) - \arctan\left(\frac{1 - c}{c}\right) -\frac{c}{2} > A > 0.$$
It follows that $I_3= O(e^{-At})$.

\bigskip
\noindent
{\bf Analysis of the integral $I_1$}\nopagebreak
\medskip

\noindent  
For $w \in C_1$ we have $w = w_1 + i(1+c)\eta$ and $w_1 \geq c\eta$. Also,
\begin{align}\nonumber
& \rho \geq \eta, \qquad \tilde{\rho} \geq \eta + 2\pi \beta,
	\\ \nonumber
& \phi = \arctan\frac{(1+c)\eta}{w_1} = \arctan\frac{1+c}{P}, 
	\\ \label{XonC1}
& \tilde{\phi} = \arctan\frac{(1+c)\eta + 2\pi \beta}{w_1} = \arctan\frac{1+c + Q}{P},
\end{align}
where $P := w_1/\eta \geq c$.
There exists $A > 0$ such that
$$1 \leq (1-A)e^{c\epsilon} \leq (1-A)e^{w_1},$$
and so
$$|e^w -1| \geq e^{w_1} -1 \geq A e^{w_1}.$$
Thus
\begin{align*}
\frac{w^{u-1}(w + 2i\pi \beta)^{v-1} e^{-mw}}{e^{w} -1}  
& = O\big(\eta^{\sigma_1 -1}(\eta+2\pi \beta)^{\sigma_2 -1} e^{-t(\phi - \tilde{\phi}) - (m+1)w_1}\big)
	\\
& = O\big(\eta^{\sigma_1 -1}(\eta+2\pi \beta)^{\sigma_2 -1} e^{-t(\phi - \tilde{\phi}) - xw_1}\big) 
	\\
& = O\big(\eta^{\sigma_1 -1}(\eta+2\pi \beta)^{\sigma_2 -1} e^{-tF(P,Q)}\big),
\end{align*}
where
$$F(P,Q) = \arctan\frac{1+c}{P} - \arctan\frac{1+c + Q}{P} + P - \frac{P}{1 + Q}.$$
Now
$$\frac{\partial F}{\partial Q} = \frac{P}{(Q+1)^2}-\frac{1}{P
   \left(\frac{(c+Q+1)^2}{P^2}+1\right)} > 0, \qquad P\geq c, \quad Q \geq 1,$$
so $F$ assumes its minimum over the domain $\{P\geq c, Q \geq 1\}$ on the boundary where $Q =1$.
For definiteness, let us henceforth set
$$c = \frac{1}{8}.$$ 
Then
$$\frac{\partial [F(P,1) - \frac{P}{16}]}{\partial P} > 0 \quad \text{for} \quad P \geq c, \quad \text{and} \quad F(c,1) - \frac{c}{16} > 0,$$
so that $F(P,1) > P/16$ for $P \geq c$.
Hence
\begin{align*}
I_1 & = \int_{C_1} \frac{w^{u-1}(w + 2i\pi \beta)^{v-1} e^{-mw} }{e^{w} -1} dw
	\\
& = O\bigg( \eta^{\sigma_1 -1}(\eta+2\pi \beta)^{\sigma_2 -1} \int_{c \eta}^{\infty}  e^{-tP/16} dw_1\bigg)
	\\
& = O\bigg(\eta^{\sigma_1}(\eta+2\pi \beta)^{\sigma_2 -1}  \int_{c}^{\infty}  e^{-tP/16} dP\bigg)
	\\
& = O\bigg(\eta^{\sigma_1}(\eta+2\pi \beta)^{\sigma_2 -1} t^{-1} e^{-c t/16}\bigg).
\end{align*}
It follows that $I_1= O(e^{-At})$.

\bigskip
\noindent
{\bf Analysis of the integral $I_2$}\nopagebreak
\medskip

\noindent
It remains to analyze the integral $I_2$. We write
$$I_2 = e^{-(m+1)i\eta} (i\eta)^{u-1} (i\eta + 2\pi i \beta)^{v-1} \int_{C_2} e^{\frac{it}{2}(\frac{1}{\eta^2} - \frac{1}{(\eta + 2\pi \beta)^2})(w - i\eta)^2} \varphi(w - i\eta) dw,$$
where
$$\varphi(z) = \frac{e^{(u-1)\ln(1 + \frac{z}{i\eta}) + (v-1)\ln(1 + \frac{z}{i\eta + 2\pi i\beta}) - x z - \frac{it}{2}(\frac{1}{\eta^2} - \frac{1}{(\eta + 2\pi \beta)^2})z^2 + (x-m)z}}{e^z - e^{-i\eta}}.$$
Defining $\phi(z)$ by
\begin{align}\label{Xphidef}
\phi(z) = e^{(u-1)\ln\bigl(1 + \frac{z}{i\eta}\bigr) + (v-1)\ln(1 + \frac{z}{i\eta + 2\pi i\beta}) - x z - \frac{it}{2}(\frac{1}{\eta^2} - \frac{1}{(\eta + 2\pi \beta)^2})z^2},
\end{align}
we have
$$\varphi(z) = \frac{\phi(z) e^{(x-m)z}}{e^z - e^{-i\eta}}.$$
We split the contour $C_2$ as follows:
$$C_2 = C_2^\epsilon \cup C_2^r,$$
where $C_2^\epsilon$ denotes the segment of $C_2$ of length $2\epsilon$ which consists of the points within a distance $\epsilon$ from $i\eta$. We write
$$I_2 = I_2^\epsilon + I_2^r,$$
where
\begin{align*}
& I_2^\epsilon = e^{-(m+1)i\eta} (i\eta)^{u-1} (i\eta + 2\pi i \beta)^{v-1} \int_{C_2^\epsilon} e^{\frac{it}{2}(\frac{1}{\eta^2} - \frac{1}{(\eta + 2\pi \beta)^2})(w - i\eta)^2} \varphi(w - i\eta) dw
\end{align*}
and
\begin{align*}
& I_2^r = e^{-(m+1)i\eta} (i\eta)^{u-1} (i\eta + 2\pi i \beta)^{v-1} \int_{C_2^r} e^{\frac{it}{2}(\frac{1}{\eta^2} - \frac{1}{(\eta + 2\pi \beta)^2})(w - i\eta)^2} \varphi(w - i\eta) dw.
\end{align*}

We claim that $I_2^r$ can be estimated as follows:
\begin{align}\label{XI2restimate}
  I_2^r = O\biggl(\frac{\eta^{\sigma_1}(\eta + 2\pi \beta)^{\sigma_2 -1}}{\sqrt{t}} e^{-\frac{t \epsilon^2}{32 \eta^2}}\biggr).
\end{align}
Indeed, the change of variables $w = i\eta + \lambda e^{\frac{i\pi}{4}}$ gives
$$I_2^r = O\biggl(\eta^{\sigma_1 -1} (\eta + 2\pi \beta)^{\sigma_2 - 1} \biggl(\int_{-\sqrt{2} c \eta}^{-\epsilon} + \int_\epsilon^{\sqrt{2} c \eta}
\biggr) e^{-\frac{t}{2}(\frac{1}{\eta^2} - \frac{1}{(\eta + 2\pi \beta)^2}) \lambda^2} |\varphi(\lambda e^{\frac{i\pi}{4}})| d\lambda\biggr).$$
Now there exists an $A > 0$ such that
\begin{align}\label{XetetaAestimate}
  \biggl|\frac{e^{(x-m)z}}{e^z - e^{-i\eta}}\biggr| < A, \qquad
z = \lambda e^{\frac{i\pi}{4}}, \quad \epsilon < |\lambda | < \sqrt{2}c \eta.
\end{align}
Moreover, the definition of $\phi(z)$ implies
\begin{align*}
  \ln{\phi(z)} = &\; (u-1)\ln\left(1 + \frac{z}{i\eta}\right) + (v-1)\ln\left(1 + \frac{z}{i\eta + 2\pi i \beta}\right)
  	\\
&   - tz\bigg(\frac{1}{\eta} - \frac{1}{\eta + 2\pi \beta}\bigg) - \frac{itz^2}{2}\bigg(\frac{1}{\eta^2} - \frac{1}{(\eta + 2\pi \beta)^2}\bigg)
  	\\
 = &\; (\sigma_1 - 1)\ln\left(1 + \frac{z}{i\eta}\right)
+(\sigma_2 - 1)\ln\left(1 + \frac{z}{i\eta + 2\pi i \beta}\right)
	\\
& - it z^2 \sum_{j=1}^\infty \frac{(-1)^{j-1}}{j+2} \bigg[\frac{1}{\eta^2} \left(\frac{z}{i\eta}\right)^j - \frac{1}{(\eta + 2\pi \beta)^2} \left(\frac{z}{i\eta + 2\pi i \beta}\right)^j\bigg].
\end{align*}
Let $\alpha \in (\sqrt{2}c, 1)$ be a constant. Then, for $|z| \leq \alpha \eta$, we have
\begin{align}\nonumber
  \re{\ln{\phi(z)}} & \leq |\sigma_1 -1| \ln(1+\alpha) + |\sigma_2 -1|\ln\bigg(1 + \frac{\alpha}{1 + \frac{2\pi \beta}{\eta}}\bigg)
  + \frac{2t}{\eta^2}|z|^2\frac{1}{3} \frac{|z|}{\eta} \sum_{j=0}^\infty \alpha^j
  	\\ \label{Xrelogphiw}
& \leq |\sigma_1 -1| \ln(1+\alpha) + |\sigma_2 -1| \ln\bigg(1 + \frac{\alpha}{1 + \frac{2\pi \beta}{\eta}}\bigg)
  + \frac{2t|z|^3}{3(1-\alpha)\eta^3},
\end{align}  
and so
\begin{align}\label{Xphilambdaeipi4}
|\phi(\lambda e^{\frac{i\pi}{4}})| < e^{|\sigma_1 -1| \ln(1+\alpha) + |\sigma_2 -1|\ln(1 + \frac{\alpha}{1 + Q})
 + \frac{2t}{3(1-\alpha)\eta^3}\lambda^3}, \qquad
\epsilon < |\lambda | < \alpha \eta.
\end{align}
Equations (\ref{XetetaAestimate}) and (\ref{Xphilambdaeipi4}) imply that
$$|\varphi(\lambda e^{\frac{i\pi}{4}})| = O\Bigl(e^{\frac{2t}{3(1-\alpha)\eta^3}|\lambda|^3}\Bigr), \qquad \epsilon < |\lambda | < \sqrt{2} c \eta.$$
This yields
\begin{align*}
I_2^r & = O\biggl(\eta^{\sigma_1 -1} (\eta + 2\pi \beta)^{\sigma_2 -1}
 \int_\epsilon^{\sqrt{2}c\eta}
e^{-\frac{t}{2}(\frac{1}{\eta^2} - \frac{1}{(\eta + 2\pi \beta)^2})\lambda^2 + \frac{2t}{3(1-\alpha)\eta^3}|\lambda|^3} d\lambda\biggr)
	\\
& =
O\biggl(\eta^{\sigma_1 -1} (\eta + 2\pi \beta)^{\sigma_2 -1}
 \int_\epsilon^{\sqrt{2} c \eta} e^{-\frac{t}{2}(\frac{1}{\eta^2} - \frac{1}{(\eta + 2\pi \beta)^2})\lambda^2 + \frac{2t \sqrt{2} c}{3(1-\alpha)\eta^2} \lambda^2} d\lambda\biggr).
\end{align*}
Now $\frac{1}{\eta + 2\pi \beta} < \frac{1}{2\eta}$, so 
\begin{align}\label{X1overeta2}
\frac{1}{\eta^2} - \frac{1}{(\eta + 2\pi \beta)^2} > \frac{3}{4\eta^2}.
\end{align}
Hence
\begin{align*}
I_2^r = O\biggl(\eta^{\sigma_1 -1} (\eta + 2\pi \beta)^{\sigma_2 -1}
 \int_\epsilon^{\sqrt{2} c \eta}
e^{-\frac{t\lambda^2}{\eta^2}(\frac{3}{8} - \frac{2\sqrt{2}c}{3(1-\alpha)})} d\lambda\biggr).
\end{align*}
For definiteness, let us henceforth set 
$$\alpha = \frac{3}{5}.$$ 
Then $\frac{3}{8} - \frac{2\sqrt{2}c}{3(1-\alpha)} > \frac{1}{16}$. Hence
\begin{align*}
I_2^r =  O\biggl(\eta^{\sigma_1 -1} (\eta + 2\pi \beta)^{\sigma_2 -1}
 \int_\epsilon^{\sqrt{2} c \eta}
e^{-\frac{t\lambda^2}{16\eta^2}} d\lambda\biggr).
\end{align*}

Splitting the integrand as
$$e^{-\frac{t}{16\eta^2}\lambda^2} 
= e^{-\frac{t}{32\eta^2}\lambda^2} \times e^{-\frac{t}{32\eta^2}\lambda^2}$$
and noting that 
$$\int_\epsilon^{\sqrt{2} c \eta} e^{-\frac{t}{32\eta^2}\lambda^2} d\lambda
\leq \int_0^{\infty} e^{-\frac{t}{32\eta^2}\lambda^2} d\lambda
= \sqrt{8\pi} \frac{ \eta}{\sqrt{t}},$$
we find
\begin{align*}
I_2^r & =
O\biggl(\eta^{\sigma_1 -1} (\eta + 2\pi \beta)^{\sigma_2 -1} e^{-\frac{t}{32 \eta^2}\epsilon^2} \frac{\eta}{\sqrt{t}} \biggr).
\end{align*}
This proves (\ref{XI2restimate}).

We now consider $I_2^\epsilon$. In view of the assumption $\dist(\eta, 2\pi \Z) > 8\epsilon$, we have
\begin{align}\label{Xexpquotientepsilon}
\biggl|\frac{e^{(x - m)w}}{e^z - e^{-i\eta}}\biggr| = O(1), \qquad
|z| < 2\epsilon.
\end{align}
In particular, $\varphi(z)$ is analytic for $|z| < 2\epsilon$. Thus we can write
$$\varphi(z) = \sum_{n=0}^{\infty} b_n z^n = \sum_{n=0}^{N-1} b_n z^n + s_N(z), \qquad |z| < 2 \epsilon,$$
where
\begin{align}\label{XsNexpression}
s_N(z) = \frac{z^N}{2\pi i}\int_{\Gamma} \frac{\varphi(w) dw}{w^N(w-z)}
\end{align}
and $\Gamma$ is a counterclockwise contour which encircles $0$ and $z$ but none of the poles of $\varphi$. 
We claim that
\begin{align}\nonumber
I_2^\epsilon = &\; e^{-(m+1)i\eta} (i\eta)^{u-1} (i\eta + 2\pi i \beta)^{v-1} \int_{C_2^\epsilon} e^{\frac{it}{2}(\frac{1}{\eta^2} - \frac{1}{(\eta + 2\pi \beta)^2}) (w - i\eta)^2} \sum_{n=0}^{N-1} b_n (w - i\eta)^n dw
	\\  \label{XI2epsilonmidstep}
&+  \begin{cases}	
O\biggl(\eta^{\sigma_1 -1}(\eta + 2\pi \beta)^{\sigma_2 -1} \left(\frac{22N}{t}\right)^{\frac{N}{6}} \frac{\eta}{\sqrt{t}}\biggr), & 1 \leq N < \frac{At}{\eta^3}, \\
O\biggl(\eta^{\sigma_1 -1}(\eta + 2\pi \beta)^{\sigma_2 -1}  \big(\frac{3t}{4N\eta^2}\big)^{-\frac{N+1}{2}} \biggr), \quad &   t^{\frac{1}{3}} < \eta < \sqrt{t}.\end{cases}  
\end{align}
In order to establish (\ref{XI2epsilonmidstep}), we need to estimate the error term
$$e^{-(m+1)i\eta} (i\eta)^{u-1} (i\eta + 2\pi i \beta)^{v-1} \int_{C_2^\epsilon} e^{\frac{it}{2}(\frac{1}{\eta^2} - \frac{1}{(\eta + 2\pi \beta)^2})(w - i\eta)^2} s_N(w - i\eta) dw.$$

Let us first consider the case $1 \leq N < \frac{At}{\eta^3}$. 
Let $|z| < \frac{40}{21}\epsilon$ (so that $\frac{21}{20}|z| < 2\epsilon$) and let $\Gamma$ in (\ref{XsNexpression}) be a circle with center $w =0$ and radius $\rho_N$, where
$$\frac{21}{20} |z| \leq \rho_N < 2\epsilon.$$
Since $\varphi(z)$  is analytic in the disk $|z| < 2\epsilon$, equations (\ref{XsNexpression}), (\ref{Xrelogphiw}), and (\ref{Xexpquotientepsilon}) yield
$$s_N(z) = O\left(\frac{2\pi \rho_N}{\rho_N - |z|}|z|^N \rho_N^{-N}e^{\frac{2t\rho_N^3}{3(1-\alpha)\eta^3}}\right)
= O\left(|z|^N \rho_N^{-N}e^{\frac{2t\rho_N^3}{3(1-\alpha)\eta^3}}\right), \qquad |z| < \frac{40}{21}\epsilon.$$
The function $\rho^{-N}e^{\frac{2t\rho^3}{3(1-\alpha)\eta^3}}$ has the minimum $(\frac{2et}{N(1-\alpha)\eta^3})^{N/3}$ for $\rho = (\frac{(1-\alpha)N}{2t})^{1/3}\eta$; $\rho_N$ can have this value if
$$\frac{21}{20}|z| \leq \biggl(\frac{(1-\alpha)N}{2t}\biggr)^{1/3}\eta < 2\epsilon.$$
The assumption $1 \leq N < \frac{At}{\eta^3}$ implies that for $A$ sufficiently small, we have 
$$\bigg(\frac{(1-\alpha)N}{2t}\bigg)^{1/3}\eta < 2\epsilon.$$
Hence, letting $\rho_N = (\frac{(1-\alpha)N}{2t})^{1/3}\eta$, we find
\begin{align}\label{XsNestimate}
  s_N(z) = 
  O\left(|z|^N \biggl(\frac{2et}{N(1-\alpha)\eta^3}\biggr)^{\frac{N}{3}}\right), 
  \qquad 1 \leq N < \frac{At}{\eta^3}, \;\; |z| \leq \frac{20}{21}\biggl(\frac{(1-\alpha)N}{2t}\biggr)^{\frac{1}{3}} \eta, 
\end{align}
For $|z| < \frac{40}{21}\epsilon$ we can also take $\rho_N = \frac{21}{20}|z|$, which yields 
\begin{align}\nonumber
s_N(z) & = O\left(\Bigl(\frac{20}{21}\Bigr)^N e^{\frac{2}{3(1-\alpha)}\frac{t}{\eta^3} (\frac{21}{20}|z|)^3}\right) 
= O\left(e^{\frac{2}{3}(\frac{21}{20})^3 \frac{1}{8(1-\alpha)} \frac{t}{\eta^2} |z|^2}\right)
	\\ \label{XsNestimate2}
& = O\left(e^{\frac{1}{10(1-\alpha)} \frac{t}{\eta^2} |z|^2}\right), 
\qquad |z| \leq \epsilon,
\end{align}
where we have used that $\frac{|z|}{\eta} < \frac{|z|}{8\epsilon} < \frac{1}{8}$ in the second step.
Using (\ref{XsNestimate}) and (\ref{XsNestimate2}), we estimate
\begin{align}\nonumber
& e^{-(m+1)i\eta} (i\eta)^{u-1}(i\eta + 2\pi i\beta)^{v-1} \int_{C_2^\epsilon} e^{\frac{it}{2}(\frac{1}{\eta^2} - \frac{1}{(\eta + 2\pi \beta)^2})(w - i\eta)^2} s_N(w - i\eta) dw
	\\ \nonumber
 = &\; O\biggl(\eta^{\sigma_1 -1}(\eta + 2\pi \beta)^{\sigma_2 -1} 
\biggl\{ \int_0^{A(\frac{N}{t})^{\frac{1}{3}} \eta} e^{-\frac{t}{2}(\frac{1}{\eta^2} - \frac{1}{(\eta + 2\pi \beta)^2})\lambda^2} \lambda^N \left(\frac{2et}{N(1-\alpha)\eta^3}\right)^{\frac{N}{3}} d\lambda
	\\\nonumber
& + \int_{A(\frac{N}{t})^{\frac{1}{3}} \eta}^{\epsilon} e^{-\frac{t}{2}(\frac{1}{\eta^2} - \frac{1}{(\eta + 2\pi \beta)^2})\lambda^2 + \frac{1}{10(1-\alpha)}\frac{t}{\eta^2}\lambda^2} d\lambda \biggr\}
\biggr)
	\\\nonumber
 = &\; O\biggl(\eta^{\sigma_1 -1}(\eta + 2\pi \beta)^{\sigma_2 -1}  \biggl\{\left(\frac{2et}{N(1-\alpha)\eta^3}\right)^{\frac{N}{3}} 2^{\frac{N -1}{2}} 
	\\ \label{Xem1ieta}
&\times  \left(t \bigg(\frac{1}{\eta^2} - \frac{1}{(\eta + 2\pi \beta)^2}\bigg)\right)^{-\frac{N+1}{2}}\Gamma\left(\frac{N +1}{2}\right) 
+ \int_{A(\frac{N}{t})^{\frac{1}{3}}\eta}^{\epsilon} e^{-\frac{t}{8\eta^2}\lambda^2} d\lambda \biggr\} \biggr),
\end{align}
where we have used the following estimate which is a consequence of (\ref{X1overeta2}) to find the last equality:
$$\frac{1}{2}\bigg(\frac{1}{\eta^2} - \frac{1}{(\eta + 2\pi \beta)^2}\bigg) - \frac{1}{10(1-\alpha)}\frac{1}{\eta^2} > \frac{1}{\eta^2} \bigg(\frac{3}{8} - \frac{1}{10(1-\alpha)}\bigg) = \frac{1}{8\eta^2}.$$
Since 
\begin{align*}
 \int_{A(\frac{N}{t})^{\frac{1}{3}}\eta}^{\epsilon} e^{-\frac{t}{8\eta^2}\lambda^2} d\lambda
& = O\biggl(e^{-\frac{t}{16\eta^2}(A(\frac{N}{t})^{\frac{1}{3}}\eta)^2}
  \int_{A(\frac{N}{t})^{\frac{1}{3}}\eta}^{\infty} e^{-\frac{t}{16\eta^2}\lambda^2} d\lambda\biggr)
 	\\ \nonumber
&  = O\biggl(e^{-\frac{A^2 t^{\frac{1}{3}} N^{\frac{2}{3}}}{16}} \frac{\eta}{\sqrt{t}}\biggr)
= O\Big(e^{-\frac{A^2}{16} t^{\frac{1}{3}}}\Big),
\end{align*}
the rhs of (\ref{Xem1ieta}) is
\begin{align}\nonumber
O\biggl(&\eta^{\sigma_1 -1}(\eta + 2\pi \beta)^{\sigma_2 -1} \left(\frac{2et}{N(1-\alpha)\eta^3}\right)^{\frac{N}{3}}2^{\frac{N -1}{2}}
	\\ \label{XC2epsestimate} 
& \times  \left(t \bigg(\frac{1}{\eta^2} - \frac{1}{(\eta + 2\pi \beta)^2}\bigg)\right)^{-\frac{N+1}{2}} \Gamma\left(\frac{N +1}{2}\right) \biggr), \qquad 1\leq N < \frac{At}{\eta^3}.
\end{align}
Using (\ref{X1overeta2}) and the asymptotic expression (\ref{Gammaexpansion}) for the Gamma function, it follows that the expression in (\ref{XC2epsestimate}) is
\begin{align}\label{XOetasigmaminusone}
O\biggl(\eta^{\sigma_1 -1}(\eta + 2\pi \beta)^{\sigma_2 -1} \left(\bigg(\frac{2e}{1-\alpha}\bigg)^2\frac{N}{t}\right)^{\frac{N}{6}} \bigg(\frac{4}{3}\bigg)^{\frac{N}{2}} e^{-\frac{N}{2}} \frac{\eta}{\sqrt{t}}\biggr).
\end{align}
Since $(\frac{2e}{1-\alpha})^2 (\frac{4}{3})^3 e^{-3}  < 22$, this proves (\ref{XI2epsilonmidstep}) in the case when $1\leq N < \frac{At}{\eta^3}$.

We now consider the case when $t^{\frac{1}{3}} < \eta < \sqrt{t}$.
Using (\ref{Xrelogphiw}) and (\ref{Xexpquotientepsilon}) in the representation (\ref{XsNexpression}) with $\Gamma$ a circle with center $w=0$ and radius $2^{1/3}\epsilon$, we find that 
\begin{align}\label{XsNestimate3}
s_N(z) = O\Big(|z|^N e^{\frac{4t\epsilon^3}{3(1-\alpha)\eta^3}}\Big), \qquad |z| < \epsilon.
\end{align}
Now $t^{\frac{1}{3}} < \eta$ implies that $e^{\frac{4t\epsilon^3}{3(1-\alpha)\eta^3}} = O(1)$, so we can estimate
\begin{align}\nonumber
&e^{-(m+1)i\eta} (i\eta)^{u-1} (i\eta + 2\pi i \beta)^{v -1}  \int_{C_2^\epsilon} e^{\frac{it}{2}(\frac{1}{\eta^2} - \frac{1}{(\eta + 2\pi \beta)^2}) (w - i\eta)^2} s_N(w - i\eta) dw
	\\ \nonumber
& = O\biggl(\eta^{\sigma_1 -1}(\eta + 2\pi \beta)^{\sigma_2 -1} \int_0^{\epsilon} e^{-\frac{t}{2} (\frac{1}{\eta^2} - \frac{1}{(\eta + 2\pi \beta)^2}) \lambda^2} \lambda^Nd\lambda \biggr)
	\\\nonumber
& = O\biggl(\eta^{\sigma_1 -1}(\eta + 2\pi \beta)^{\sigma_2 -1} 2^{\frac{N -1}{2}} 
\left(t\bigg(\frac{1}{\eta^2} - \frac{1}{(\eta + 2\pi \beta)^2}\bigg)\right)^{-\frac{N+1}{2}} \Gamma\left(\frac{N +1}{2}\right) \biggr)
	\\ \nonumber
& = O\biggl(\eta^{\sigma_1 -1}(\eta + 2\pi \beta)^{\sigma_2 -1}  \left(t\bigg(\frac{1}{\eta^2} - \frac{1}{(\eta + 2\pi \beta)^2}\bigg)\right)^{-\frac{N+1}{2}} N^{\frac{N}{2}} \biggr), \qquad
t^{\frac{1}{3}} < \eta < \sqrt{t}.
\end{align}
In view of (\ref{X1overeta2}), this completes the proof of (\ref{XI2epsilonmidstep}).

We next claim that, up to a small error term, the contour $C_2^\epsilon$ in the integral in (\ref{XI2epsilonmidstep}) can be replaced by the infinite line $C_2'$, where $C_2'$ denotes the infinite straight line of which $C_2$ is a part. More precisely, we claim that there exists an $A >0$ such that
\begin{align} \nonumber
& e^{-(m+1)i\eta} (i\eta)^{u-1}(i\eta + 2\pi i\beta)^{v-1} \int_{C_2'\backslash C_2^\epsilon}e^{\frac{it}{2}(\frac{1}{\eta^2} - \frac{1}{(\eta + 2\pi \beta)^2})(w - i\eta)^2} \sum_{n=0}^{N-1} b_n (w - i\eta)^n dw
	\\ \label{XC2epsilonC2prime}
& = \begin{cases}
  O\Bigl( \eta^{\sigma_1 -1}(\eta + 2\pi \beta)^{\sigma_2 -1}e^{- \frac{At}{\eta^2} }\Bigr), & 8\epsilon < \eta < t^{\frac{1}{3}}, \quad 1 \leq N < \frac{At}{\eta^3} , \\
  O\Bigl(\eta^{\sigma_1 -1}(\eta + 2\pi \beta)^{\sigma_2 -1} e^{- \frac{At}{\eta^2}} N\Bigr), & t^{\frac{1}{3}} < \eta < \sqrt{t}, \quad 1 \leq N < \frac{At}{\eta^2}.
\end{cases}  
 \end{align}
In order to prove (\ref{XC2epsilonC2prime}), we note that the coefficient of $b_n$ on the lhs of (\ref{XC2epsilonC2prime}) is
\begin{align}\label{XbncoefficientbigOh}
O\biggl(\eta^{\sigma_1 -1}(\eta + 2\pi \beta)^{\sigma_2 -1} \int_\epsilon^\infty e^{-\frac{3t}{8\eta^2} \lambda^2} \lambda^n d\lambda\biggr).
\end{align}
We write the integrand as
$$e^{-\frac{3t}{16\eta^2} \lambda^2} \lambda^n \times e^{-\frac{3t}{16\eta^2} \lambda^2}.$$
The first factor is steadily decreasing for $\lambda > \sqrt{\frac{8n}{3t}} \eta$, and so it decreases throughout the interval of integration provided that $n < N < \frac{3t\epsilon^2}{8\eta^2}$. 
The term in (\ref{XbncoefficientbigOh}) is then
\begin{align}\nonumber
& O\biggl(\eta^{\sigma_1 -1} (\eta + 2\pi \beta)^{\sigma_2 -1}e^{-\frac{3t}{16\eta^2} \epsilon^2} \epsilon^n
\int_\epsilon^\infty e^{-\frac{3t}{16\eta^2} \lambda^2} d\lambda\biggr)
	\\
& \label{XbncoefficientbigOh2}
= O\biggl(\eta^{\sigma_1 -1}  (\eta + 2\pi \beta)^{\sigma_2 -1} 
e^{-\frac{3t}{16\eta^2} \epsilon^2} \epsilon^n \frac{\eta}{\sqrt{t}}\biggr).
\end{align}

Let us assume that $8\epsilon < \eta < t^{\frac{1}{3}}$ and $1 \leq N < At/\eta^3$. In this case, choosing $|z| < \frac{20}{21}\big(\frac{(1-\alpha)N}{2t}\big)^{\frac{1}{3}} \eta$, equation (\ref{XsNestimate}) yields
\begin{align}\label{Xbnestimate}
b_n = (s_n(z) - s_{n+1}(z))z^{-n} 
= O\left(\left(\frac{2et}{n(1-\alpha)\eta^3}\right)^{\frac{n}{3}}\right), \qquad N < \frac{At}{\eta^3},\quad n \geq 1.
\end{align}
Multiplying the rhs of (\ref{XbncoefficientbigOh2}) by  $b_n$ and summing from $0$ to $N-1$, we find that the total error is
$$O\left(  \eta^{\sigma_1 -1}(\eta + 2\pi \beta)^{\sigma_2 -1} e^{-\frac{3t}{16\eta^2} \epsilon^2} \frac{\eta}{\sqrt{t}}
\biggl(b_0 +  \sum_{n=1}^{N-1}\left(\frac{2e\epsilon^3 t}{(1-\alpha) n\eta^3}\right)^{\frac{n}{3}}\biggr)\right).$$
Now the function $(\frac{t}{n\eta^3})^{\frac{n}{3}}$ increases steadily up to $n = \frac{t}{\eta^3 e}$, so that if $n < A \frac{t}{\eta^3}$, where $A < 1/e$, it is of order
$$O\Bigl(e^{\frac{1}{3}\frac{At}{\eta^3}\ln\frac{1}{A}}\Bigr).$$
Hence, choosing $\epsilon >0$ and $A > 0$ so small that $\epsilon^3 < \frac{1-\alpha}{2e}$ and $\frac{1}{3\eta}A\ln\frac{1}{A} < \frac{3\epsilon^2}{16}$, the total error is
$$O\Bigl( \eta^{\sigma_1 -1}(\eta + 2\pi \beta)^{\sigma_2 -1} e^{- \frac{At}{\eta^2}}\Bigr).$$
This proves (\ref{XC2epsilonC2prime}) in the case when $8\epsilon < \eta < t^{\frac{1}{3}}$ and $1 \leq N < At/\eta^3$.

In the case when $t^{\frac{1}{3}} < \eta < \sqrt{t}$ and $1 \leq N < \frac{At}{\eta^2}$, we instead use (\ref{XsNestimate3}) to find
\begin{align}\label{Xbnestimatecase2} 
b_n = (s_n(z) - s_{n+1}(z))z^{-n} 
= O\Bigl(e^{\frac{4t\epsilon^3}{3(1-\alpha)\eta^3}}\Bigr) = O(1), \qquad n \geq 0.
\end{align}
Then, the total error is
\begin{align*}
& O\left(\eta^{\sigma_1 -1}(\eta + 2\pi \beta)^{\sigma_2 -1} e^{-\frac{3t}{16\eta^2}\epsilon^2}  N\right), \qquad  t^{\frac{1}{3}} < \eta < \sqrt{t},
\end{align*}
which completes the proof of (\ref{XC2epsilonC2prime}).

We finally analyze the sum
\begin{align*}
& e^{-(m+1)i\eta} (i\eta)^{u-1} (i\eta + 2\pi i\beta)^{v-1} \int_{C_2'} e^{\frac{it}{2}(\frac{1}{\eta^2} - \frac{1}{(\eta + 2\pi \beta)^2}) (w - i\eta)^2} \sum_{n=0}^{N-1} b_n (w - i\eta)^n dw
	\\
 = &-e^{-(m+1)i\eta} (i\eta)^{u-1} (i\eta + 2\pi i \beta)^{v-1} \sum_{n=0}^{N-1} b_n e^{\frac{\pi i}{4}(n+1)}
\int_{-\infty}^\infty e^{-\frac{t}{2}(\frac{1}{\eta^2} - \frac{1}{(\eta + 2\pi \beta)^2}) \lambda^2} \lambda^n d\lambda
	\\
= & -e^{-(m+1)i\eta} (i\eta)^{u-1} (i\eta + 2\pi i \beta)^{v-1}
	\\
& \times \sum_{n=0}^{N-1} b_n e^{\frac{\pi i}{4}(n+1)}
2^{\frac{n-1}{2}} (1 + (-1)^n) \biggl(t\bigg(\frac{1}{\eta^2} - \frac{1}{(\eta + 2\pi \beta)^2}\bigg)\biggr)^{-\frac{n+1}{2}} \Gamma\biggl(\frac{n+1}{2}\biggr)
	\\
= & -e^{-(m+1)i\eta} (i\eta)^{u-1} (i\eta + 2\pi i \beta)^{v-1}  e^{\frac{i\pi}{4}} 
	\\
& \times \sum_{k=0}^{[\frac{N-1}{2}]} b_{2k} i^k 
\biggl(\frac{t}{2}\bigg(\frac{1}{\eta^2} - \frac{1}{(\eta + 2\pi \beta)^2}\bigg)\biggr)^{-k - \frac{1}{2}} \Gamma\biggl(k + \frac{1}{2}\biggr).
\end{align*}
Together with equations (\ref{XI2restimate}), (\ref{XI2epsilonmidstep}), and (\ref{XC2epsilonC2prime}), this yields (\ref{Xzetaformula2}) with $\epsilon$ replaced with $8 \epsilon$. Since $\epsilon > 0$ was arbitrary, the proof is complete.
\end{proof}

In the special case when $N = 6$  and $\eta = \pi$, theorem \ref{Xzetath2} reduces to the following result which will be used in chapter \ref{hurwitzsec}.

\begin{corollary}\label{Phiasymptoticscorollary}
Let $u = \sigma_1 + it$ and $v = \sigma_2 - it$. 
Then
\begin{align} \nonumber
& -e^{i\pi(\sigma_1 + \sigma_2)} \Phi(u,v,\beta) 
 = \sum_{n=1}^m \int_{\hat{H}_\alpha} w^{u-1}(w + 2i\pi \beta)^{v-1} e^{-nw} dw
 	\\ \nonumber
&+ \frac{(-1)^{m+1} e^{\frac{\pi i}{4}}e^{\frac{\pi i}{2}(\sigma_1 + \sigma_2)} \pi^{\sigma_1 + \sigma_2-\frac{1}{2}} ( 1 + 2\beta)^v}{2\sqrt{2\beta (1+ \beta) t}}
 + O(\beta^{\sigma_2 -1} t^{-3/2}),    \qquad t \to \infty, 
\end{align}
uniformly for $\beta \in [\frac{3}{4}, \infty)$, $\sigma_2 \in [0,1]$, and for $\sigma_1$ in compact subsets of $(0,1]$, where
\begin{align*}
& m = \bigg[\frac{t}{\pi(1 + (2 \beta)^{-1})}\bigg].
\end{align*}
\end{corollary}
\begin{proof}
The expansion follows by setting $N = 6$  and $\eta = \pi$ in (\ref{Xzetaformula2}). Indeed, a direct computation shows that $\varphi''(0) = O(1)$  and $\varphi^{(4)}(0) = O(t)$ uniformly for $\sigma_1, \sigma_2, \beta$ in the given ranges; the terms involving $\varphi''(0)$ and $\varphi^{(4)}(0)$ in (\ref{Xzetaformula2}) are therefore $O(\beta^{\sigma_2 -1} t^{-3/2})$.
Since $\varphi(0) = 1/2$, the corollary follows.
\end{proof}

\begin{remark}
For $\sigma_1 > 0$, the Hankel contour integrals occurring in theorem \ref{Xzetath2} and corollary  \ref{Phiasymptoticscorollary} can be alternatively written as
\begin{align*}
\int_{\hat{H}_\alpha} w^{u-1}(w + 2i\pi \beta)^{v-1} e^{-nw} dw
= (e^{2\pi i u} - 1)\int_0^\infty w^{u-1}(w + 2i\pi \beta)^{v-1} e^{-nw} dw
\end{align*}
where the contribution from the factor $e^{2\pi i u} = O(e^{-2\pi t})$ is exponentially small. This follows by collapsing the contour $\hat{H}_\alpha$ onto the positive real axis.
\end{remark}

Theorem \ref{Xzetath2} and Corollary \ref{Phiasymptoticscorollary} provide asymptotic formulas for $\Phi(u,v,\beta)$ as $t \to \infty$. The following theorem, which we will need in chapter \ref{hurwitzsec}, considers the large $t$ behavior of $\Phi(v,u,\beta)$.

\begin{theorem}[The asymptotics of $\Phi(v,u,\beta)$]
Let $u = \sigma_1 + it$ and $v = \sigma_2 - it$. Then, for every $N \geq 1$, 
\begin{align}\label{Phivuestimate}
\Phi(v, u, \beta) = O(\beta^{\sigma_1 -1} e^{2\pi t} t^{-N}), \qquad t \to \infty,
\end{align}
uniformly for $\sigma_1,  \sigma_2 \in [0,1]$ and $\beta \geq \frac{1}{2}$.
\end{theorem}
\begin{proof}
Suppose $\beta \geq 1/2$ and $\sigma_1,  \sigma_2 \in [0,1]$. By (\ref{PhihatHalphaexpression}), we have
$$\Phi(v,u,\beta) = -e^{-i \pi(u+v)} \int_{\hat{H}_1} w^{v-1} (w + 2\pi i \beta)^{u-1} \frac{dw}{e^w -1},$$
where
\begin{align}\label{wv1wu1estimate}
|w^{v-1} (w + 2\pi i \beta)^{u-1}| = |w|^{\sigma_2-1}|w + 2\pi i \beta|^{\sigma_1 -1} e^{t \arg w - t\arg (w + 2\pi i \beta)},
\end{align}
with $\arg w, \arg (w + 2\pi i \beta) \in [0, 2\pi)$. 
There exists a constant $c > 0$ independent of $\beta \geq 1/2$ such that
$$\arg w - \arg (w + 2\pi i \beta) \leq 2\pi - c,$$
for $w$ on the upper side of $[1, \infty)$ and on the unit circle. Hence
\begin{align}\label{Phivubeta}
\Phi(v,u,\beta) = -e^{-i \pi(\sigma_1 + \sigma_2)} e^{2 i\pi v} \int_1^\infty w^{v-1} (w + 2\pi i \beta)^{u-1} \frac{dw}{e^w -1} + O(\beta^{\sigma_1 -1} e^{(2\pi-c)t}).
\end{align}
Using (\ref{wv1wu1estimate}) and the fact that $(e^w -1)^{-1} \leq 2e^{-w}$ for $w \geq 1$, we find
\begin{align*}
\bigg|\int_1^\infty w^{v-1} (w + 2\pi i \beta)^{u-1} \frac{dw}{e^w -1}\bigg|
& \leq 2\int_1^\infty w^{\sigma_2-1}|w + 2\pi i \beta|^{\sigma_1 -1} e^{- t\arg (w + 2\pi i \beta)} \frac{dw}{e^w}
	\\
& \leq 2\beta^{\sigma_1 -1}\int_1^\infty e^{- t \arctan\frac{2\pi \beta}{w}} e^{-w} dw.
\end{align*}
The inequality
$$\arctan\frac{2\pi \beta}{w} \geq \frac{1}{w}, \qquad w \geq 1, \quad \beta \geq \frac{1}{2},$$
then implies
\begin{align}\label{int1inftywv1}
\bigg|\int_1^\infty w^{v-1} (w + 2\pi i \beta)^{u-1} \frac{dw}{e^w -1}\bigg|
\leq 2\beta^{\sigma_1 -1}\int_1^\infty e^{- \frac{t}{w}} e^{-w} dw.
\end{align}
Iterated integration by parts using
$$ \frac{\partial }{\partial w}\bigg(\frac{1}{t} e^{- \frac{t}{w}}\bigg) = \frac{e^{- \frac{t}{w}}}{w^2},$$
gives
\begin{align}\nonumber
\int_1^\infty e^{- \frac{t}{w}} e^{-w} dw
= & \; \frac{e^{-t-1}}{-t} 
+ \frac{1}{-t} \int_1^\infty \frac{e^{- \frac{t}{w}}}{w^2} w^2 \frac{d}{dw}(w^2 e^{-w}) dw
	\\\nonumber
=& \; \frac{e^{-t-1}}{-t} + \frac{e^{-t-1}}{(-t)^2} 
+ \frac{1}{(-t)^2}\int_1^\infty \frac{e^{- \frac{t}{w}}}{w^2} \bigg(w^2\frac{d}{dw} \bigg)^2(w^2 e^{-w}) dw
	\\\nonumber
=& \; \cdots =  \sum_{k = 1}^{N} \frac{e^{-t}}{(-t)^k} \bigg(w^2 \frac{d}{dw} \bigg)^{k-1}(w^2 e^{-w})\bigg|_{w = 1}
	\\ \label{inttwIBP}
& + \frac{1}{(-t)^N}\int_1^\infty \frac{e^{- \frac{t}{w}}}{w^2} \bigg(w^2 \frac{d}{dw} \bigg)^N (w^2 e^{-w}) dw
= O(t^{-N}).
\end{align}
The theorem follows from (\ref{Phivubeta}), (\ref{int1inftywv1}), and (\ref{inttwIBP}).
\end{proof}

\chapter{More Explicit Asymptotics of $\Phi(u,v, \beta)$}\label{sec8}
The asymptotic formula for $\Phi(u,v, \beta)$ derived in Theorem \ref{Xzetath2} involves the sum
$$\sum_{n=1}^m \int_{\hat{H}_\alpha} w^{u-1}(w + 2i\pi \beta)^{v-1} e^{-nw} dw.$$
Our goal in this chapter is to obtain a more explicit asymptotic formula for $\Phi$ (see corollary \ref{Phicor}) by computing the asymptotic behavior of this sum in the range where $\beta > t^{1+\epsilon}$.
We recall that $u = \sigma_1 + it$ and $v = \sigma_2 - it$.

\begin{theorem}[{\bf The asymptotics to all orders of the integral \\ $\int_{\hat{H}_\alpha} w^{u-1}(w + 2i\pi \beta)^{v-1} e^{-nw} dw$}]\label{sec8th}
For every $\epsilon > 0$, there exists a constant $A > 0$ such that
\begin{align} \nonumber
\int_{\hat{H}_\alpha} & w^{u-1}(w + 2i\pi \beta)^{v-1} e^{-nw} dw
 = -(i\eta)^{u-1} (i\eta + 2\pi i \beta)^{v-1} e^{-in\eta}  e^{\frac{i\pi}{4}}
 	\\\nonumber
& \times  \sum_{k=0}^{[\frac{N-1}{2}]}  \frac{\phi^{(2k)}(0)}{(2k)!}  i^k 
\biggl(\frac{t}{2}\bigg(\frac{1}{\eta^2} - \frac{1}{(\eta + 2\pi \beta)^2}\bigg)\biggr)^{-k - \frac{1}{2}} \Gamma\biggl(k + \frac{1}{2}\biggr)
	\\\nonumber
& + O\bigg(\eta^{\sigma_1} (\eta + 2\pi \beta)^{\sigma_2 -1} \bigg(\frac{7N}{t}\bigg)^{\frac{N}{6}} \frac{1}{\sqrt{t}}\bigg), \qquad t \to \infty,
 	\\ \label{Xzetaformula}	
&  t^{1+\epsilon} < \beta < \infty, \quad n = 1, 2, \dots, [t], \quad \sigma_1 \in [0,1], \quad \sigma_2 \in [0,1], \quad 1 \leq N < At,
\end{align}
where $\hat{H}_\alpha$ denotes the Hankel contour defined in (\ref{hatHalphadef}), $\eta$ and $\phi(z)$ are given by
\begin{align}\label{Xetadefbeta}
\eta = \beta \pi \bigg(-1 + \sqrt{1 + \frac{2t}{\beta \pi n}}\bigg),
\end{align}
\begin{align}\label{Xphizdef}
\phi(z) = e^{(u-1)\ln(1 + \frac{z}{i\eta}) + (v-1)\ln(1 + \frac{z}{i\eta + 2\pi i\beta}) - nz - \frac{it}{2}(\frac{1}{\eta^2} - \frac{1}{(\eta + 2\pi \beta)^2})z^2},
\end{align}
and the error terms are uniform with respect to $n,\sigma_1, \sigma_2, \beta, N$ in the given ranges.
\end{theorem}
\begin{proof}
Let $\epsilon > 0$ be given and suppose that $\beta > t^{1+\epsilon}$, $n = 1, \dots, [t]$, $\sigma_1 \in [0,1]$, $\sigma_2 \in [0,1]$, and $N \geq 1$. All error terms of the form $O(\cdot)$ will be uniform with respect to $n, \sigma_1, \sigma_2, \beta, N$ (but not with respect to $\epsilon$).
Let $\eta$ be given by (\ref{Xetadefbeta}).
The function $x (\sqrt{1+1/x} -1)$ increases from $0$ to $1/2$ as $x$ goes from $0$ to $\infty$. Hence, using that
$$\eta 
= \frac{2t}{n} \frac{\beta \pi n}{2t}\bigg(\sqrt{1 + \frac{2t}{\beta \pi n}} - 1\bigg)$$ 
we infer that, given any $\delta > 0$, we have
\begin{align}\label{Xetatn}
(1- \delta) \frac{t}{n} < \eta < \frac{t}{n}
\end{align}
for all sufficiently large $t$. 

As in the proof of Theorem \ref{Xzetath2}, we deform the Hankel contour  $\hat{H}_\alpha$ into the straight lines $C_j$, $j = 1, \dots, 4$ joining $\infty$, $i\eta + (1+i)c\eta$, $i\eta - (1+i)c\eta$, $-c\eta - i\eta$, $\infty$, where $0 < c \leq 1/2$ is an absolute constant.
This implies
\begin{align*}
\int_{\hat{H}_\alpha} w^{u-1}(w + 2i\pi \beta)^{v-1} e^{-nw} dw =  \sum_{j=1}^4 I_j,
\end{align*}
where
$$I_j = \int_{C_j} w^{u-1}(w + 2i\pi \beta)^{v-1} e^{-nw} dw.$$

We first prove that $I_1, I_3, I_4$ are exponentially small as $t \to \infty$.
We let $Q = 2\pi \beta/\eta$. Then, by (\ref{Xetatn}), $Q > \frac{2\pi \beta n}{t} > A n t^{\epsilon} \to \infty$ as $t \to \infty$.

Using (\ref{XC4est}) and (\ref{Xetatn}), we find that $I_4$ is exponentially small:
\begin{align*}
I_4 & = \int_{C_4} w^{u-1}(w + 2i\pi \beta)^{v-1} e^{-nw} dw
 = O\bigg(e^{-t(\arctan\frac{1}{c} + \arctan\frac{Q - 1}{c})} \int_{-c\eta}^\infty e^{-nw_1}dw_1\bigg)
	\\
& = O(e^{-t(\arctan\frac{1}{c} + \arctan\frac{Q - 1}{c}) + nc\eta}) 
=
O(e^{-t(\arctan\frac{1}{c} + \arctan\frac{Q - 1}{c} - c)})
	\\
& =
O(e^{-t(\arctan\frac{1}{c} - c)}) = O(e^{-At}).
\end{align*}

We next consider $I_3$. Letting $w = \rho e^{i\phi}$ and $w + 2i\pi \beta = \tilde{\rho} e^{i\tilde{\phi}}$, equations (\ref{XphiminustildephionC3}) and (\ref{Xetatn}) yield
\begin{align*}
w^{u-1}(w + 2i\pi \beta)^{v-1} &e^{-nw}
 = O((c\eta)^{\sigma_1 -1} (c\eta)^{\sigma_2 -1} e^{-t(\arctan (\frac{1-c + Q}{c}) - \arctan \frac{1-c}{c})} e^{nc\eta})
	\\
&= O(e^{-t(\arctan (\frac{1-c + Q}{c}) - \arctan \frac{1-c}{c} - c)})
= O(e^{-tF(Q)}),
\end{align*}
where 
$$F(Q) = \arctan \bigg(\frac{1-c + Q}{c}\bigg) - \arctan\bigg(\frac{1-c}{c}\bigg) - c.$$
Since $Q \to \infty$ and $F(Q) > A > 0$ for all sufficiently large $Q$, it follows that $I_3= O(e^{-At})$.

For $w \in C_1$ we have $w = w_1 + i(1+c)\eta$ with $w_1 \geq c\eta$. Equations (\ref{XonC1}) and (\ref{Xetatn}) imply
\begin{align*}
w^{u-1}(w + 2i\pi \beta)^{v-1} e^{-nw} 
& = O\big(\eta^{\sigma_1 -1}(\eta+2\pi \beta)^{\sigma_2 -1} e^{-t(\phi - \tilde{\phi}) - nw_1}\big)
	\\
& = O\big(e^{-t(\phi - \tilde{\phi} + \frac{nw_1}{t})}\big)
= O\big(e^{-t(\phi - \tilde{\phi} + (1-\delta) P)}\big)
= O\big(e^{-tF(P,Q)}\big)
\end{align*}
where $P = w_1/\eta \geq c$ and
$$F(P,Q) = \arctan\frac{1+c}{P} - \arctan\frac{1+c + Q}{P} + (1-\delta)P.$$
Now the rhs of the inequality
$$F(P,Q) - A P > \arctan\frac{1+c}{P} - \frac{\pi}{2} + (1-\delta)P - AP$$
is an increasing function of  $P \geq c$ and $\arctan\frac{1+c}{c} - \frac{\pi}{2} + (1-\delta)c - Ac > 0$ for $\delta$ and $A$ small enough. Hence $F(P,Q) > AP$.
It follows that
\begin{align*}
I_1 & = \int_{C_1} w^{u-1}(w + 2i\pi \beta)^{v-1} e^{-nw} dw 
	\\
& = O\bigg(\int_{c \eta}^{\infty}  e^{-tAP} dw_1\bigg) = O\bigg(\eta \int_{c}^{\infty}  e^{-tAP} dP\bigg) 
= O\big(\eta t^{-1} e^{-c A t}\big) = O\big(e^{-c A t}\big).
\end{align*}
 
It remains to analyze the integral $I_2$. Using that $n = t(\frac{1}{\eta} - \frac{1}{\eta + 2\pi \beta})$, we write 
$$I_2 = (i\eta)^{u-1} (i\eta + 2\pi i \beta)^{v-1} e^{-it\eta(\frac{1}{\eta} - \frac{1}{\eta + 2\pi \beta})} \int_{C_2} e^{\frac{it}{2}(\frac{1}{\eta^2} - \frac{1}{(\eta + 2\pi \beta)^2})(w - i\eta)^2} \phi(w - i\eta) dw,$$
where $\phi(z)$ is defined in (\ref{Xphizdef}). 

Define $\{a_j\}_0^\infty$ by 
\begin{align*}
\phi(z) = \sum_{j=0}^\infty a_j z^j, \qquad |z| < \eta.
\end{align*}
Then
$$\phi(z) = \sum_{j=0}^{N-1} a_j z^j + r_N(z),  \qquad |z| < \eta,$$
where
$$r_N(z) = \frac{z^N}{2\pi i} \int_\Gamma \frac{\phi(w)}{w^N(w-z)}dw,$$
and $\Gamma$  is a counterclockwise contour contained in the disk of radius $\eta$ centered at the origin which encircles the points $0$ and $z$ once. 

Let $\alpha \in (\sqrt{2}c, 1)$ be a constant. Let $|z| < \frac{20\alpha}{21}\eta$ (so that $\frac{21}{20}|z| < \alpha\eta$) and let $\Gamma$ be a circle with center $w =0$ and radius $\rho_N$, where
$$\frac{21}{20} |z| \leq \rho_N \leq \alpha \eta.$$
By (\ref{Xrelogphiw}),
$$r_N(z) = O\left(\frac{2\pi \rho_N |z|^N}{\rho_N^N(\rho_N - |z|)} e^{\frac{2t\rho_N^3}{3(1-\alpha)\eta^3}}\right)
= O\left(|z|^N \rho_N^{-N}e^{\frac{2t\rho_N^3}{3(1-\alpha)\eta^3}}\right), \qquad |z| \leq \alpha \eta.$$
The function $\rho^{-N}e^{\frac{2t\rho^3}{3(1-\alpha)\eta^3}}$ has the minimum $(\frac{2et}{(1-\alpha)N\eta^3})^{N/3}$ for $\rho = (\frac{(1-\alpha)N}{2t})^{1/3}\eta$; $\rho_N$ can have this value if
$$\frac{21}{20}|z| \leq \bigg(\frac{(1-\alpha)N}{2t}\bigg)^{\frac{1}{3}}\eta \leq \alpha\eta.$$
Hence, 
\begin{align} \label{XrNestimate}
 r_N(z) & = O\left(|z|^N \biggl(\frac{2et}{(1-\alpha)N\eta^3}\biggr)^{\frac{N}{3}}\right),
  \qquad N \leq \frac{2\alpha^3}{1-\alpha} t, \quad |z| \leq \frac{20}{21}\biggl(\frac{(1-\alpha)N}{2t}\biggr)^{\frac{1}{3}} \eta.
\end{align}
For $|z| < \frac{20\alpha}{21}\eta$ we can also take $\rho_N = \frac{21}{20}|z|$, which yields
\begin{align}\label{XrNestimate2pre}
r_N(z) = O\left(\Bigl(\frac{20}{21}\Bigr)^N e^{\frac{2 t}{3(1-\alpha)\eta^3} (\frac{21}{20}|z|)^3}\right)
= O\left(e^{\frac{2 \alpha}{3(1-\alpha)}(\frac{21}{20})^2 \frac{t}{\eta^2} |z|^2}\right), \qquad |z| < \frac{20\alpha}{21}\eta.
\end{align}
For definiteness, we henceforth set
$$c = \frac{1}{8}, \qquad \alpha = \frac{1}{4}.$$
Then (\ref{XrNestimate2pre}) yields
\begin{align}\label{XrNestimate2}
r_N(z) = O\left(e^{\frac{t}{4\eta^2} |z|^2}\right), \qquad |z| < \frac{20\alpha}{21}\eta.
\end{align}

We write $I_2$ in the form
\begin{align}\label{XI2expression}
  I_2 = I_2^S + I_2^R,
\end{align} 
where $I_2^S$ and $I_2^R$ denote the integrals
\begin{align*}
I_2^S   = & \; (i\eta)^{u-1} (i\eta + 2\pi i \beta)^{v-1} e^{-it\eta(\frac{1}{\eta} - \frac{1}{\eta + 2\pi \beta})}  
	\\
& \times \int_{C_2} e^{ \frac{it}{2}(\frac{1}{\eta^2} - \frac{1}{(\eta + 2\pi \beta)^2})(w - i\eta)^2} \sum_{j=0}^{N-1} a_j (w - i\eta)^j dw
\end{align*}
and 
\begin{align}\nonumber
I_2^R = &\; (i\eta)^{u-1} (i\eta + 2\pi i \beta)^{v-1} e^{-it\eta(\frac{1}{\eta} - \frac{1}{\eta + 2\pi \beta})} 
	\\\label{XI2Rdef}
& \times \int_{C_2} e^{ \frac{it}{2}(\frac{1}{\eta^2} - \frac{1}{(\eta + 2\pi \beta)^2})(w - i\eta)^2} r_N(w - i\eta) dw.
\end{align}
We claim that
\begin{align}\label{XI2Restimate}
I_2^R = O\biggl( \frac{\eta^{\sigma_1}(\eta + 2\pi \beta)^{\sigma_2 -1}}{\sqrt{t}} \bigg(\frac{7N}{t}\bigg)^{\frac{N}{6}}\biggr).
\end{align}
To prove (\ref{XI2Restimate}) we make the change of variables $w = i\eta + \lambda e^{\frac{i\pi}{4}}$ in (\ref{XI2Rdef}) and split the integral into two integrals; in the first integral $|\lambda | \leq A(\frac{N}{t})^{\frac{1}{3}} \eta$, while in the second integral $A(\frac{N}{t})^{\frac{1}{3}} \eta \leq |\lambda| \leq \sqrt{2} c \eta$.
In the first integral, we use the estimate (\ref{XrNestimate}) of $r_N$, whereas in the second integral, we use the estimate (\ref{XrNestimate2}). This yields
\begin{align*}
I_2^R = &\; O\biggl(\eta^{\sigma_1 -1} (\eta + 2\pi \beta)^{\sigma_2 -1} \biggl\{ \int_0^{A(\frac{N}{t})^{\frac{1}{3}} \eta} e^{-\frac{t}{2}(\frac{1}{\eta^2} - \frac{1}{(\eta + 2\pi \beta)^2})\lambda^2} \lambda^N \left(\frac{2et}{(1-\alpha)N\eta^3}\right)^{\frac{N}{3}} d\lambda
	\\
&+ \int_{A(\frac{N}{t})^{\frac{1}{3}} \eta}^{\sqrt{2} c \eta} e^{-\frac{t}{2}(\frac{1}{\eta^2} - \frac{1}{(\eta + 2\pi \beta)^2})\lambda^2 + \frac{t}{4\eta^2}\lambda^2} d\lambda \biggr\}\biggr).
\end{align*}
Since $\eta < 2\pi \beta$, the estimate (\ref{X1overeta2}) holds; thus steps almost identical to those leading from (\ref{Xem1ieta}) to (\ref{XOetasigmaminusone}) now show that 
\begin{align*}
I_2^R = O\biggl(\eta^{\sigma_1 -1} (\eta + 2\pi \beta)^{\sigma_2 -1} \left(\bigg(\frac{2e}{1-\alpha}\bigg)^2\frac{N}{t}\right)^{\frac{N}{6}} \bigg(\frac{4}{3}\bigg)^{\frac{N}{2}} e^{-\frac{N}{2}}  \frac{\eta}{\sqrt{t}}\biggr).	
\end{align*}
Since $(\frac{2e}{1-\alpha})^2  (\frac{4}{3})^3 e^{-3}  < 7$, this proves (\ref{XI2Restimate}).

We next consider $I_2^S$. We claim that there exists a constant $A>0$ such that
\begin{align}\nonumber
  I_2^S = & \; (i\eta)^{u-1} (i\eta + 2\pi i \beta)^{v-1} e^{-it\eta(\frac{1}{\eta} - \frac{1}{\eta + 2\pi \beta})}
 	\\ \label{XI2Sestimate}
& \times  \int_{C_2'}  e^{\frac{it}{2}(\frac{1}{\eta^2} - \frac{1}{(\eta + 2\pi \beta)^2})(w - i\eta)^2}  \sum_{j=0}^{N-1} a_j (w - i\eta)^j dw + O(e^{- At}),
\end{align}
where $C_2'$ denotes the infinite straight line of which $C_2$ is a part.
Indeed, if we replace $C_2$ by $C_2'$, the coefficient multiplying $a_j$ in the expression for $I_2^S$ changes by
\begin{align}\label{XancoefficientbigOh}
O\left(\eta^{\sigma_1 -1}(\eta + 2\pi \beta)^{\sigma_2 -1}  \int_{\sqrt{2} c \eta}^\infty e^{-\frac{t}{2}(\frac{1}{\eta^2} - \frac{1}{(\eta + 2\pi \beta)^2})\lambda^2} \lambda^j  d\lambda\right)
= O\left( \int_{\sqrt{2} c \eta}^\infty e^{-\frac{3t}{8\eta^2}\lambda^2} \lambda^j  d\lambda\right).
\end{align}
We can write the integrand as
$$e^{-\frac{3t}{16\eta^2}\lambda^2} \lambda^j \times e^{-\frac{3t}{16\eta^2}\lambda^2},$$
and the first factor is steadily decreasing for $\lambda > \sqrt{\frac{8j}{3t}} \eta$, and so it decreases throughout the interval of integration provided that $j < N < At$ with $0 < A \leq 3c^2/4$. The rhs of (\ref{XancoefficientbigOh}) is then
\begin{align}\label{XwholetermO}
O\left(e^{-\frac{3t}{16\eta^2} 2c^2 \eta^2} (\sqrt{2} c \eta)^j \int_{\sqrt{2} c \eta}^\infty e^{-\frac{3t}{16\eta^2}\lambda^2} d\lambda\right)
= O\left(e^{-\frac{3c^2 t}{8}} (\sqrt{2} c \eta)^j \frac{\eta}{\sqrt{t}} \right).
\end{align}
Also, by (\ref{XrNestimate}), choosing $z$ with $|z|$ small enough, we find
\begin{align}\nonumber
a_j & = (r_j(z) - r_{j+1}(z))z^{-j} 
	\\ \nonumber
& = O\left(\left(|z|^j \left(\frac{2et}{(1-\alpha)j\eta^3}\right)^{\frac{j}{3}} - |z|^{j+1} \left(\frac{2et}{(1-\alpha)(j+1)\eta^3}\right)^{\frac{j+1}{3}}\right)|z|^{-j}\right)
	\\\label{Xanestimate}
&=  O\left(\left(\frac{2et}{(1-\alpha)j\eta^3}\right)^{\frac{j}{3}}\right), \qquad N < At, \quad 1 \leq j \leq N - 1.
\end{align}
Multiplying (\ref{XwholetermO}) by $a_j$ and summing from  $0$ to $N-1$, we find that the total error is
\begin{align*}
 & O\left(e^{- \frac{3c^2 t}{8}}\frac{\eta}{\sqrt{t}}
\biggl(a_0+ \sum_{j=1}^{N-1} (\sqrt{2} c \eta)^j \left(\frac{2 e t}{(1-\alpha) j \eta^3}\right)^{\frac{j}{3}} \biggr)\right)
  	\\
& = 
O\left(e^{- \frac{3c^2 t}{8}}\frac{\eta}{\sqrt{t}}
\biggl(1 + \sum_{j=1}^{N-1} \left(\frac{2^{5/2} c^3 e t}{(1-\alpha) j }\right)^{\frac{j}{3}} \biggr)\right).
\end{align*}
The function $(t/j)^{\frac{j}{3}}$ increases steadily up to $j = t/e$, and so if $j < A t$, where $A < 1/e$, it is
$$O(e^{\frac{1}{3}t A\ln\frac{1}{A}}).$$
Moreover, $\frac{2^{5/2} c^3 e }{(1-\alpha) } < 1$. Hence if $N < At$, with $A$ sufficiently small, the total error is $O(e^{- At})$.
This proves (\ref{XI2Sestimate}).

We finally analyze the integral
\begin{align*}
&   \int_{C_2'}  e^{\frac{it}{2}(\frac{1}{\eta^2} - \frac{1}{(\eta + 2\pi \beta)^2})(w - i\eta)^2}  \sum_{j=0}^{N-1} a_j (w - i\eta)^j dw
 	\\
& = - \sum_{j=0}^{N-1} a_j e^{\frac{\pi i}{4}(j+1)}
\int_{-\infty}^\infty e^{-\frac{t}{2}(\frac{1}{\eta^2} - \frac{1}{(\eta + 2\pi \beta)^2}) \lambda^2} \lambda^j d\lambda
	\\
& = - \sum_{j=0}^{N-1} a_j e^{\frac{\pi i}{4}(j+1)}
2^{\frac{j-1}{2}} (1 + (-1)^j) \biggl(t\bigg(\frac{1}{\eta^2} - \frac{1}{(\eta + 2\pi \beta)^2}\bigg)\biggr)^{-\frac{j+1}{2}} \Gamma\biggl(\frac{j+1}{2}\biggr)
	\\
& = - e^{\frac{i\pi}{4}} \sum_{k=0}^{[\frac{N-1}{2}]} a_{2k} i^k 
\biggl(\frac{t}{2}\bigg(\frac{1}{\eta^2} - \frac{1}{(\eta + 2\pi \beta)^2}\bigg)\biggr)^{-k - \frac{1}{2}} \Gamma\biggl(k + \frac{1}{2}\biggr).
\end{align*}
Together with equations (\ref{XI2Restimate}) and (\ref{XI2Sestimate}), this yields (\ref{Xzetaformula}).

\end{proof}

Summing the asymptotic formula of theorem \ref{sec8th} from  $n=1$ to $n= m$, we obtain the following corollary. 

\begin{corollary}[{\bf The asymptotics of $ \sum_{n=1}^m \int_{\hat{H}_\alpha} w^{u-1}(w + 2i\pi \beta)^{v-1} e^{-nw} dw$}]\label{Xsumcor}
For every $\epsilon > 0$, there exists a constant $A > 0$ such that
\begin{align} \nonumber
\sum_{n=1}^m \int_{\hat{H}_\alpha} & w^{u-1}(w + 2i\pi \beta)^{v-1} e^{-nw} dw
 = J + O(t^{1/6} \beta^{\sigma_2 -1}), \qquad t \to \infty,
 	\\ \label{Xsumasymptotics}
&  m = 1, \dots, [t], \quad t^{1+\epsilon} < \beta < \infty, \quad \sigma_1 \in [0,1], \quad \sigma_2 \in [0,1], 
\end{align}
where $J$ is defined by
\begin{align}\label{XJdef}
J = e^{\frac{i\pi}{2}(\sigma_1 + \sigma_2)} e^{\frac{i\pi}{4}} \sqrt{\frac{\pi}{2t}} (\pi \beta)^{\sigma_1 + \sigma_2 - 1} \sum_{n=1}^m (\sqrt{R} -1)^{\sigma_1} (\sqrt{R} + 1)^{\sigma_2} \bigg(\frac{2t}{\pi n \beta}\bigg)^{it} F,
\end{align}
with
\begin{align}\label{XRFdef}
R := 1 + \frac{2t}{\pi n \beta}, \qquad
F := \frac{(1 + \sqrt{R})^{-2it}}{R^{1/4}} e^{i\pi n \beta ( 1- \sqrt{R})},
\end{align}
and the error term is uniform with respect to $m, \beta, \sigma_1, \sigma_2$ in the given ranges.
\end{corollary}
\begin{proof}
Summing equation (\ref{Xzetaformula}) with $N =2$ from $n = 1$ to $n = m$, we find
\begin{align*}
\sum_{n=1}^m \int_{\hat{H}_\alpha} & w^{u-1}(w + 2i\pi \beta)^{v-1} e^{-nw} dw
 = J + O\biggl(\frac{1}{t^{5/6}} \sum_{n=1}^m \eta^{\sigma_1} (\eta + 2\pi \beta)^{\sigma_2 -1} \biggr),
\end{align*}
where
\begin{align*}
& J = - \sum_{n=1}^m (i\eta)^{u-1} (i\eta + 2\pi i \beta)^{v-1} e^{-in\eta}  e^{\frac{i\pi}{4}} 
\biggl(\frac{t}{2}\bigg(\frac{1}{\eta^2} - \frac{1}{(\eta + 2\pi \beta)^2}\bigg)\biggr)^{- \frac{1}{2}} \sqrt{\pi},
	\\
& \eta = \beta \pi \bigg(-1 + \sqrt{1 + \frac{2t}{\beta \pi n}}\bigg).
\end{align*}
Straightforward algebra shows that $J$ can be written as in (\ref{XJdef}).
Equation (\ref{Xsumasymptotics}) follows because, by (\ref{Xetatn}),
$$\sum_{n=1}^m \eta^{\sigma_1} (\eta + 2\pi \beta)^{\sigma_2 -1}
= O\bigg(\sum_{n=1}^m \Big(\frac{t}{n}\Big)^{\sigma_1} \beta^{\sigma_2 -1}\bigg)
= O\big(t^{\sigma_1} \beta^{\sigma_2 -1} t^{1-\sigma_1}\big)
= O\big(t \beta^{\sigma_2 -1}\big).$$
\end{proof}

Substituting the result of corollary \ref{Xsumcor} into the asymptotic expansion of theorem \ref{Xzetath2}, we arrive at the following asymptotic formula for $\Phi(u,v,\beta)$. For brevity of presentation, we state the result only to leading order.

\begin{corollary}[{\bf The asymptotics of $\Phi(u,v,\beta)$}]\label{Phicor}
For every $\epsilon > 0$, there exists a constant $A > 0$ such that
\begin{align} \nonumber
\Phi&(u,v,\beta) 
 = -e^{-\frac{i\pi}{2}(\sigma_1 + \sigma_2)} e^{\frac{i\pi}{4}} \sqrt{\frac{\pi}{2t}} (\pi \beta)^{\sigma_1 + \sigma_2 -1} 
 	\\ \nonumber
& \times \sum_{n=1}^{[t(\frac{1}{\eta} - \frac{1}{\eta + 2\pi \beta})]} (\sqrt{R} -1)^{\sigma_1} (\sqrt{R} + 1)^{\sigma_2} \bigg(\frac{2t}{\pi n \beta}\bigg)^{it} F
	\\ \nonumber
& - i (2\pi)^{\sigma_1 + \sigma_2 - 1} e^{-\frac{\pi i}{2}(\sigma_1 + \sigma_2)} \sum_{n=1}^{[\frac{\eta}{2\pi}]}   n^{u-1} (n+\beta)^{v-1} + O(t^{1/6} \beta^{\sigma_2 -1}), \qquad  t \to \infty,
	\\ \label{XPhiasymptoticssumsum}
&1 \leq \eta < t^{\frac{1}{3} - \epsilon} < \infty, \quad
 t^{1+\epsilon} < \beta < \infty, \quad \dist(\eta, 2\pi \Z) > \epsilon, \quad  \sigma_1 \in [0,1], \quad \sigma_2 \in [0,1], 
\end{align}
where $R$ and $F$ are defined in (\ref{XRFdef}) and the error term is uniform with respect to $\eta, \beta, \sigma_1, \sigma_2$ in the given ranges.
\end{corollary}
\begin{proof}
Letting $N =2$ in Theorem \ref{Xzetath2} and using Corollary \ref{Xsumcor} we find equation (\ref{XPhiasymptoticssumsum}) but with the error term
\begin{align*} \nonumber
 & O(t^{1/6} \beta^{\sigma_2 -1}) + \eta^{\sigma_1 -1} (\eta + 2\pi \beta)^{\sigma_2 -1} 
O\biggl(t^{-\frac{1}{3}} \frac{\eta}{\sqrt{t}}\biggr)
	\\
&+ e^{-\big([t(\frac{1}{\eta} - \frac{1}{\eta + 2\pi \beta})]+1\big)i\eta}  \frac{(i\eta)^{u-1} (i\eta + 2\pi i \beta)^{v-1}e^{\frac{i\pi}{4}}}{1 - e^{-i\eta}} \biggl(\frac{t}{2}\bigg(\frac{1}{\eta^2} - \frac{1}{(\eta + 2\pi \beta)^2}\bigg)\biggr)^{- \frac{1}{2}} \sqrt{\pi}.
\end{align*} 
Since this error term is
\begin{align*} \nonumber
O\big(t^{1/6} \beta^{\sigma_2 -1}\big)
+ O\big(\eta^{\sigma_1} \beta^{\sigma_2 -1}  t^{-\frac{5}{6}} \big)
+ O\bigg(\eta^{\sigma_1-1} \beta^{\sigma_2-1} \frac{\eta}{\sqrt{t}}\bigg)
=O\big(t^{1/6} \beta^{\sigma_2 -1}\big)
\end{align*} 
the result follows.
\end{proof}

\chapter[Fourier coefficients of the product]{Fourier coefficients of the product \\ of two Hurwitz zeta functions}\label{hurwitzsec}
The erratic behavior of the Riemann zeta function in the critical strip makes it difficult to obtain good asymptotic information on $\zeta(s)$  for $0 \leq \re s \leq 1$. 
Much of the literature on the asymptotics of $\zeta(s)$ therefore instead focuses on the asymptotics of various related quantities in which the erratic behavior has been smoothed out by averaging. The prime example is the study of the large $T$  behavior of the moments (see \cite[chapter VII]{T1986} and e.g. \cite{CFKRS2005, S2009})
$$\int_1^T |\zeta(\sigma + it)|^{2k} dt.$$
A second example involves the small $\delta$ behavior of the integral
$$\int_0^\infty |\zeta(\sigma + it)|^{2k} e^{-\delta t} dt,$$
where $\delta > 0$, see \cite[Section 7.12]{T1986}. 
Another class of examples, which is relevant for the present chapter, considers the large $t$ behavior of various quantities involving the Hurwitz zeta function $\zeta(s,\alpha)$ averaged with respect to the extra parameter $\alpha$. 

Recall that the Hurwitz zeta function $\zeta(s,\alpha)$ is defined as the analytic continuation of the sum
$$\zeta(s,\alpha) = \sum_{n=0}^\infty \frac{1}{(n + \alpha)^s}, \qquad \re s > 1, \quad \re \alpha > 0.$$
For each $\alpha$ with $\re \alpha > 0$, $\zeta(s,\alpha)$ is a meromorphic function of $s \in \C$ with a simple pole at $s = 1$ and no other poles.
For each $s \in \C \setminus \{1\}$, $\zeta(s, \alpha)$ is analytic in the half-plane $\re \alpha > 0$.
These properties of the Hurwitz function follow easily from the representation
\begin{align}\label{hurwitzHankel}
\zeta(s, \alpha) = \frac{\Gamma(1-s)}{2\pi i} \int_{H_1} \frac{e^{\alpha z} z^{s-1}}{1 - e^z} dz, \qquad s \in \C \setminus \{1\}, \quad \re \alpha > 0,
\end{align}
where $H_1$ is the Hankel contour surrounding the negative real axis defined in (\ref{X1.2}).

For many purposes it is more convenient to work with the {\it modified} Hurwitz function $\zeta_1(s,\alpha)$ defined by
\begin{align}\label{modifiedhurwitzshift}
\zeta_1(s, \alpha) = \zeta(s, \alpha) - \alpha^{-s} = \zeta(s, \alpha +1),
\end{align}
which is regular at $\alpha = 0$ (in fact, analytic for $\re \alpha > -1$). 
For each $s\in \C \setminus\{1\}$, $\zeta_1(s,\alpha)$ is a smooth function of $\alpha \in [0,1]$. Thus, for each choice of $u,v \in \C \setminus\{1\}$, the product $\zeta_1(u,\alpha)\zeta_1(v,\alpha)$ can be represented by its Fourier series
\begin{align}\label{zeta1zeta1Fourierseries}
\zeta_1(u,\alpha)\zeta_1(v,\alpha) = \sum_{n \in \Z} q_n(u,v) e^{2\pi i n \alpha}, \qquad 0 < \alpha < 1,
\end{align}
where the Fourier coefficients $q_n(u,v)$ are defined by
$$q_n(u,v) = \int_0^1 \zeta_1(u,\alpha) \zeta_1(v,\alpha) e^{-2\pi i n \alpha}  d\alpha, \qquad n \in \Z.$$
By standard Fourier analysis, the series in (\ref{zeta1zeta1Fourierseries}) converges in $L^2([0,1])$ and converges pointwise for each $\alpha \in (0,1)$. 

In this chapter, we will show the following two theorems on the Fourier coefficients $q_n(u,v)$. The proofs are presented in sections  \ref{qnidentitysec} and \ref{qnasymptoticssec}, respectively.

\begin{theorem}[Expression for $q_n(u,v)$]\label{QNTH}
For each $n \in \Z$, the $n$th Fourier coefficient $q_n(u,v)$ can be expressed as follows for $u,v \in \C \setminus \{1\}$ with $\re u, \re v < 2$:
\begin{align}\label{qnidentity}
q_n(u,v) = &\;  b_n(u+v) + R_n(u,v) + R_n(v,u) - T_n(u,v) - T_n(v,u),
\end{align}
where, for any $n \in \Z$, $b_n(s)$, $R_n(u,v)$, and $T_n(u,v)$ are the meromorphic continuations of the functions
\begin{subequations}\label{bnRnTndef}
\begin{align}\label{bndef}
& b_n(s) =  \int_1^\infty \alpha^{-s}  e^{-2\pi i n \alpha}  d\alpha, \qquad \re s > 1,
	\\\label{Rndef}
& R_n(u,v) = \int_0^\infty \alpha^{-v} \zeta_1(u,\alpha) e^{-2i\pi n \alpha} d\alpha, \qquad \re (u+v) > 2, \quad \re v < 1,
	\\\label{Tndef}
& T_n(u,v) = \int_0^1 \alpha^{-v}  \zeta_1(u, \alpha) e^{-2\pi i n\alpha} d\alpha, \qquad \re v < 1.
\end{align}
\end{subequations}
These meromorphic continuations are given by
\begin{subequations}\label{bnRnTncontinuation}
\begin{align}\label{bncontinuation}
 b_n(s) = & \begin{cases} \Gamma(1-s, 2\pi i n)(2\pi i n)^{s-1}, & n \neq 0, 
	\\
\frac{1}{s-1}, \quad & n = 0, 
\end{cases}
	\\ \label{Rncontinuation}
 R_n(u,v) = &\; \frac{\Gamma(1-v)}{\Gamma(u)}  \times \begin{cases}  
    \frac{ie^{\pi i v}\Phi(u,v, n)}{2\sin(\pi u)},  &  n \geq 1, \\
 \frac{ie^{-\pi i v}\Phi(u,v, n)}{2\sin(\pi u)}, & n \leq -1, \\
\Gamma(u+v-1)\zeta(u+v-1), & n = 0,
\end{cases}
	\\\label{Tncontinuation}
T_n(u,v) = &\; \frac{\zeta(u) - 1}{1 - v} + \frac{u}{1-v} \int_0^1 \alpha^{1-v}  \zeta_1(u+1, \alpha) e^{-2i\pi n \alpha} d\alpha
	\\\nonumber
& + \frac{2i\pi n}{1-v} \int_0^1 \alpha^{1-v}  \zeta_1(u,\alpha) e^{-2i\pi n \alpha} d\alpha, \qquad \re v < 2, \quad n \in \Z,
\end{align}
\end{subequations}
where $\Phi(u,v,n)$ is the function defined in (\ref{X1.1}) and $\Gamma(s, z)$ denotes the incomplete Gamma function:
$$\Gamma(s,z) = \int_z^\infty r^{s-1} e^{-r} dr, \qquad \Gamma(s,0) = \Gamma(s) = \int_0^\infty r^{s-1} e^{-r} dr.$$
\end{theorem}

\begin{theorem}[Large $t$ asymptotics of $q_n(u,v)$]\label{QNASYMPTOTICSTH}
Let $u = \sigma_1 + it$ and $v = \sigma_2 - it$. For each integer $n \geq 1$, the $n$th Fourier coefficient $q_n(u,v)$ satisfies the following asymptotic formula as $t \to \infty$:
\begin{align}\label{qnasymptoticsPhi}
q_n(u,v) = &\;  b_n(u+v) + \frac{\Gamma(1-v)}{\Gamma(u)} \frac{ie^{\pi i v}\Phi(u,v, n)}{2\sin(\pi u)}  
- \frac{\zeta(u)}{1-v} -  \frac{\zeta(v)}{1-u}
	\\ \nonumber
& 
- \frac{\zeta(u+1) u}{(1-v)(2-v)}  - \frac{\zeta(v+1) v}{(1-u)(2-u)}  + O(t^{-1}), \qquad \sigma_1, \sigma_2 \in [0,1],
\end{align}
where the error term is uniform for $\sigma_1, \sigma_2$ in the given range.

For $n = 0$, the following formulas hold as $t \to \infty$:
\begin{itemize}
\item For $\sigma_1, \sigma_2 \in [0,1]$ with $\sigma_1 + \sigma_2 \neq 1$,
\begin{align}\nonumber
q_0(u, v) = 
&\;  \frac{1}{\sigma_1 + \sigma_2 - 1} + \Gamma(\sigma_1 + \sigma_2 -1)\zeta(\sigma_1 + \sigma_2 -1)\bigg(\frac{\Gamma(1-u)}{\Gamma(v)} + \frac{\Gamma(1-v)}{\Gamma(u)}\bigg) 
 	\\ \label{q0asymptoticsneq1}
& - \frac{\zeta(u)}{1- v} - \frac{\zeta(v)}{1- u} - \frac{\zeta(u+1) u}{(1-v)(2-v)}  - \frac{\zeta(v+1)v}{(1-u)(2-u)} 
+ O(t^{-1}),  
\end{align}
where the error term is uniform for all $\sigma_1, \sigma_2 \in [0,1]$ such that $\sigma_1 + \sigma_2 \neq 1$.
 
\item For $\sigma_1, \sigma_2 \in [0,1]$ with $\sigma_1 + \sigma_2 = 1$,
\begin{align} \label{q0asymptoticseq1}
q_0(u,v) = &\; \ln \frac{t}{2\pi}  + \gamma 
- \frac{\zeta(u)}{1- v} - \frac{\zeta(v)}{1- u} 
	\\\nonumber
& - \frac{\zeta(u+1) u}{(1-v)(2-v)}  - \frac{\zeta(v+1) v}{(1-u)(2-u)} 
+ O(t^{-1}), 
\end{align}
where the error term is uniform for all $\sigma_1, \sigma_2 \in [0,1]$ such that $\sigma_1 + \sigma_2 = 1$, and $\gamma$ denotes the Euler constant.
\end{itemize}
\end{theorem}


\begin{remark}\upshape
The results of theorem \ref{QNTH} and \ref{QNASYMPTOTICSTH} are well-known in the case when $n = 0$. More precisely, the formulas (\ref{q0asymptoticseq1}) and (\ref{q0asymptoticsneq1}) can be found in Theorem 1 and Theorem 2 of \cite{Z1994}, respectively.
They are included here for the sake of comparison; the new content lies in the case of nonzero $n$. 
Note however that the terms $\frac{\zeta(u+1)  u}{(1-v)(2-v)}$ and $\frac{\zeta(v+1) v}{(1-u)(2-u)}$ in (\ref{q0asymptoticsneq1}) and (\ref{q0asymptoticseq1}) are missing in \cite{Z1994}. These terms must be included in order to make the asymptotic formulas uniform for $\sigma_1$ and/or $\sigma_2$ near $0$, because it is known that $\zeta(1+it)$ is unbounded as $t \to \infty$; in fact, $\zeta(1+it) = \Omega(\ln \ln t)$, see Titchmarsh \cite[Theorem 8.5]{T1986}.  
\end{remark}

\begin{remark}\upshape
The integral in (\ref{Rndef}) converges whenever $\re (u+v) > 2$ and $\re v < 1$, because for each $s \in \C \setminus \{1\}$, $\zeta_1(s,\alpha)$ satisfies the large $\alpha$ asymptotics
\begin{align}\label{Hurwitzalphaasymptotics}
\zeta_1(s, \alpha) = \frac{\alpha^{1-s}}{s-1} - \frac{\alpha^{-s}}{2} + O(\alpha^{-s-1}), 
\end{align}
uniformly as $\alpha \to \infty$ in the half-plane $\re \alpha > 0$, see e.g. \cite[Eq. (25.11.43)]{NIST}.
\end{remark}

\section{Corollaries}
By substituting the asymptotic formula for $\Phi(u,v,n)$ found in corollary \ref{Phiasymptoticscorollary} into (\ref{qnasymptoticsPhi}), we obtain the following corollary of theorem \ref{QNASYMPTOTICSTH}.

\begin{corollary}
Let $u = \sigma_1 + it$ and $v = \sigma_2 - it$. For each integer $n \geq 1$, 
\begin{align}\nonumber
q_n(u,v) =  &\; b_n(\sigma_1 + \sigma_2)	
 +  \frac{\Gamma(1-v)}{\Gamma(u)} \frac{ie^{\pi i v}}{2\sin(\pi u)} 
 	\\\label{qnasymptotics}
&\times \bigg\{-e^{-\pi i(\sigma_1 + \sigma_2)}\sum_{j=1}^m  \int_{\hat{H}_\alpha} w^{u-1}(w + 2i\pi n)^{v-1} e^{-jw} dw
 	\\ \nonumber
&+\frac{(-1)^{m} e^{\frac{\pi i}{4}}e^{-\frac{\pi i}{2}(\sigma_1 + \sigma_2)} \pi^{\sigma_1 + \sigma_2-\frac{1}{2}} ( 1 + 2n)^v}{2\sqrt{2n (1+ n) t}} 
  + O(t^{-3/2})
 \bigg\}
	\\\nonumber
& - \frac{\zeta(u)}{1- v} - \frac{\zeta(v)}{1- u} - \frac{\zeta(u+1) u}{(1-v)(2-v)}  - \frac{\zeta(v+1)v}{(1-u)(2-u)} + O(t^{-1}),
	\\\nonumber
& \hspace{5cm} \sigma_1,\sigma_2 \in [0,1], \quad t \to \infty,
\end{align}
where $m = [\frac{t}{\pi(1 + (2 n)^{-1})}]$ and the error term is uniform for $\sigma_1, \sigma_2 \in [0,1]$.
\end{corollary}

In the special case of $\sigma_ 1 = \sigma_2$, we obtain asymptotic estimates for the Fourier coefficients 
$$Q_n(u) =  \int_0^1 |\zeta_1(u,\alpha)|^2 e^{-2\pi i n \alpha} d\alpha$$ 
of $|\zeta_1(s, \alpha)|^2$ as corollaries.

\begin{corollary}[Large $t$ asymptotics of $Q_0(u)$]\label{Q0cor}
For $\sigma \in [0,1]$ with $\sigma \neq 1/2$,
\begin{align}\nonumber
Q_0(\sigma + it) = 
&\;  \frac{1}{2 \sigma - 1} + 2\Gamma(2 \sigma -1)\zeta(2 \sigma -1) \sin(\pi \sigma) t^{1-2\sigma}
 	\\ \label{Q0asymptoticsneq1}
& - 2 \re \frac{\zeta(\sigma + it)}{1- \sigma + it}  - 2\frac{\im \zeta(\sigma+1+it)}{t}
+ O(t^{-1}), \qquad t\to \infty,
\end{align}
where the error term is uniform with respect to $\sigma$ in the given range. For $\sigma = 1/2$,
\begin{align} \label{Q0asymptoticseq1}
Q_0\bigg(\frac{1}{2} + it\bigg) = &\; \ln \frac{t}{2\pi}  + \gamma 
- 2 \re \frac{\zeta(\frac{1}{2} + it)}{\frac{1}{2} + it} + O(t^{-1}), \qquad  t\to \infty.
\end{align}
\end{corollary}
\begin{proof}
Let $u = \sigma + it$ and $v = \sigma - it$. Then (see (\ref{Gammaasymptotics}))
$$\frac{\Gamma(1-u)}{\Gamma(v)} + \frac{\Gamma(1-v)}{\Gamma(u)}
= 2\sin(\pi \sigma)t^{1-2\sigma}  + O(t^{-2\sigma-1})$$
uniformly for $\sigma \in [0,1]$. Hence (\ref{Q0asymptoticsneq1}) and (\ref{Q0asymptoticseq1}) follow from equations (\ref{q0asymptoticsneq1}) and (\ref{q0asymptoticseq1}), respectively. 
\end{proof}

\begin{corollary}\label{Q0boundcor}
As $t \to \infty$,
\begin{align}\label{Q0bound}
Q_0(\sigma + it) = \begin{cases} O(t^{1-2\sigma}), & \sigma \in [0, 1/2), \\
O(\ln t), & \sigma = 1/2, \\
O(1), & \sigma \in (1/2, 1].
\end{cases}
\end{align}
If the error terms $O(t^{1-2\sigma})$ and $O(1)$ on the rhs are replaced with $O(t^{1-2\sigma} \ln t)$ and $O(\ln t)$, respectively, then (\ref{Q0bound}) holds uniformly for $\sigma \in [0,1]$. 
\end{corollary}
\begin{proof}
Equation (\ref{Q0bound}) follows immediately from corollary \ref{Q0cor}. From this corollary it also follows that the error terms in (\ref{Q0bound}) are uniform with respect to $\sigma \in [0,\frac{1}{2} - \delta] \cup \{\frac{1}{2}\} \cup [\frac{1}{2} + \delta, 1]$ for any fixed $\delta > 0$. We will complete the proof by showing that 
\begin{align}\label{Q0lnttbound}
Q_0(\sigma + it) = O(t^{1-2\sigma} \ln t), \qquad t \to \infty,
\end{align}
uniformly for $\sigma \in [0, 1/2)$ and
\begin{align}\label{Q0lnttbound2}
Q_0(\sigma + it) = O(\ln t), \qquad t \to \infty,
\end{align}
uniformly for $\sigma \in (1/2,1]$.
By (\ref{Q0asymptoticsneq1}), we have
\begin{align}\nonumber
Q_0(\sigma + it) = \frac{1}{2 \sigma - 1} + 2\Gamma(2 \sigma -1)\zeta(2 \sigma -1) \sin(\pi \sigma) t^{1-2\sigma}
+ O(1), \qquad t \to \infty,
\end{align}
uniformly for $\sigma \in [0, 1/2) \cup (1/2,1]$.
Since
$$t^{1-2\sigma} = 1 - 2 \int_{1/2}^\sigma t^{1-2\sigma'} (\ln t) d\sigma' = 
\begin{cases}
1 + O((1 - 2\sigma) t^{1-2\sigma} \ln t), & \sigma \in  [0, 1/2), \\
1 + O((2\sigma - 1) \ln t), & \sigma \in  (1/2,1],
\end{cases}$$
uniformly for $\sigma$ in the given ranges, and since
$$\frac{1}{2 \sigma - 1} + 2\Gamma(2 \sigma -1)\zeta(2 \sigma -1) \sin(\pi \sigma), \quad
2\Gamma(2 \sigma -1)\zeta(2 \sigma -1) \sin(\pi \sigma)|2\sigma -1|,$$
are bounded functions of $\sigma \in [0,1]$, we find (\ref{Q0lnttbound}) and (\ref{Q0lnttbound2}). 
\end{proof}

\begin{corollary}[Large $t$ asymptotics of $Q_n(u)$]\label{Qncor}
Let $u = \sigma + it$. For each integer $n \geq 1$, 
\begin{align}\nonumber
Q_n(u) = &\;\overline{Q_{-n}(u)} =  b_n(2\sigma) -  
 t^{1 - 2\sigma}  i e^{-\pi i \sigma}  \sum_{j=1}^m  \int_{\hat{H}_\alpha}  w^{u-1}(w + 2i\pi n)^{\bar{u}-1} e^{-jw} dw
 	\\ \label{Qnasymptotics}
&+ t^{\frac{1}{2} - 2\sigma}  \frac{i(-1)^{m} e^{\frac{\pi i}{4}}\pi^{2 \sigma- \frac{1}{2}} ( 1 + 2n)^{\bar{u}}}{2\sqrt{2n (1+ n)}}
- 2 \re \frac{\zeta(\sigma + it)}{1- \sigma + it}  - 2\frac{\im \zeta(\sigma+1+it)}{t}
	\\ \nonumber
&  + O(t^{-1}) + O(t^{- \frac{1}{2} - 2\sigma}), \qquad t\to \infty, \quad \sigma \in [0,1],
\end{align}
where $m = [\frac{t}{\pi(1 + (2 n)^{-1})}]$ and the error term is uniform for $\sigma$ in the given range.
\end{corollary}
\begin{proof}
Let $u = \sigma + it$ and $v = \sigma - it$. By (\ref{Gammaasymptotics}), the coefficient of $\Phi(u,v,n)$ in (\ref{qnasymptoticsPhi}) satisfies
$$\frac{\Gamma(1-v)}{\Gamma(u)} \frac{ie^{\pi i v}}{2\sin(\pi u)} = t^{1-2\sigma} e^{\frac{\pi i}{2}(2\sigma + 1)}\big(1 + O(t^{-2})\big), \qquad t \to \infty,$$
uniformly for $\sigma \in [0,1]$. Hence there exist constants $T \geq 1$ and $c > 0$ such that
$$\bigg|\frac{\Gamma(1-v)}{\Gamma(u)} \frac{ie^{\pi i v}}{2\sin(\pi u)}\bigg| \geq c t^{1-2\sigma},$$
for all $t \geq T$ and all $\sigma \in [0,1]$. Since the terms not involving $\Phi$ on the rhs of (\ref{qnasymptoticsPhi}) are $O(1)$, and clearly $|Q_n(u)| \leq Q_0(u)$ for each integer $n$, we conclude from (\ref{qnasymptoticsPhi}) and corollary \ref{Q0boundcor} that
\begin{align}\label{t1sigmaPhiuvn}
t^{1 - 2\sigma} \Phi(u,v,n) = O( Q_0(u)) + O(1) = \begin{cases} O(t^{1-2\sigma} \ln t), & \sigma \in [0, 1/2), \\
O(\ln t), & \sigma \in [1/2, 1].
\end{cases}
\end{align}
uniformly for $\sigma \in [0,1]$.
Equation (\ref{Qnasymptotics}) now follows from (\ref{qnasymptotics}).
\end{proof}

From (\ref{t1sigmaPhiuvn}) and corollary \ref{Q0boundcor}, we also obtain an estimate for $\Phi$. 

\begin{corollary}\label{Phiboundcor}
For each integer $n \geq 1$,
\begin{align}\label{Philargetbound}
\Phi(\sigma + it, \sigma - it,n) = \begin{cases} O(1), & \sigma \in [0, 1/2), \\
O(\ln t), & \sigma = 1/2, \\
O(t^{2\sigma - 1}), & \sigma \in (1/2, 1].
\end{cases}
\end{align}
If the error terms $O(1)$ and $O(t^{2\sigma - 1})$ on the rhs are replaced with $O(\ln t)$ and $O(t^{2\sigma - 1} \ln t)$, respectively, then (\ref{Philargetbound}) holds uniformly for $\sigma \in [0,1]$. 
\end{corollary}

\section{Asymptotics of the zeroth Fourier coefficient}\label{zerothsec}
In order to put theorems \ref{QNTH} and \ref{QNASYMPTOTICSTH} in context, let us first consider the case $n = 0$ of the zeroth Fourier coefficient.
The study of the large $t$ asymptotics of the zeroth Fourier coefficient
$$Q_0(s) = \int_0^1 |\zeta_1(s, \alpha)|^2 d\alpha, \qquad s = \sigma + it,$$
i.e., of the mean square of $\zeta_1(s,\alpha)$ as a function  of $\alpha \in [0,1]$, has a long history. Koksma and Lekkerkerker showed in the 1952 paper \cite{KL1952} that $Q_0(\frac{1}{2} + it) = O(\ln t)$ and that
$$Q_0(\sigma + it) = \frac{1}{2\sigma -1} + O\bigg(t^{1-2\sigma}\bigg(\frac{1}{2\sigma -1} + \ln t\bigg)\bigg) \quad \text{uniformly for} \quad  \frac{1}{2} < \sigma \leq 1.$$
In the case when $\sigma = 1/2$, Balasubramanian established \cite{B1979} the asymptotic formula
$$Q_0\bigg(\frac{1}{2} + it\bigg) = \ln t + O(\ln \ln t)$$
and in the 1980s the error term in this formula was improved to $O(1)$ by Rane \cite{R1983}. Sitaramachandrarao was awarded a prize (announced in Hardy-Ramanujan Journal {\bf 10} (1987), p. 28) for proving the more precise result
$$Q_0\bigg(\frac{1}{2} + it\bigg) = \ln \frac{t}{2\pi}  + \gamma + O(t^{-3/16}(\ln t)^{3/8}).$$
Zhang independently obtained \cite{Z1990} a similar result but with the slightly larger error term $O(t^{-3/16}(\ln t)^{11/8})$, which was later improved to $O(t^{-7/36}(\ln t)^{25/18})$ in \cite{Z1991}. Finally, in \cite{A1992} and \cite{Z1994}, Andersson and Zhang independently arrived at the asymptotic expansion stated in (\ref{Q0asymptoticseq1}).

Around the same time as Andersson and Zhang obtained the expansion (\ref{Q0asymptoticseq1}), Katsurada and Matsumoto started using Atkinson's dissection argument applied to the product $\zeta_1(u,\alpha) \zeta_1(v,\alpha)$ with two independent complex variables $u$ and $v$ in order to analyze $Q_0$, see \cite{KM1993}. In \cite{KM1996}, they refined their approach and found an alternative proof of (\ref{Q0asymptoticseq1}) together with a number of other related results.
The derivation of \cite{KM1996} is based on the following expression for the zeroth Fourier coefficient (see \cite[Eq. (2.1)]{KM1996}):
\begin{align}\label{q0identity}
q_0(u,v) = &\; \frac{1}{u+v-1} + R_0(u,v) + R_0(v,u) 
 - T_0(u,v) - T_0(v,u),
	\\\nonumber
& \hspace{3cm} -N+1 < \re u, \re v < N+1,
\end{align}
where
\begin{align}\label{R0def}
& R_0(u,v) = \Gamma(u+v-1)\zeta(u+v-1) \frac{\Gamma(1-v)}{\Gamma(u)}
\end{align}
and, for any integer $N \geq 0$,
\begin{align}\nonumber
 T_0(u,v) = &\; \sum_{k=0}^{N-1} \frac{(u)_k}{(1-v)_{k+1}}(\zeta(u+k) - 1)
	\\ \label{T0def}
& + \frac{(u)_N}{(1-v)_N} \sum_{l=1}^\infty l^{1-u-v} \int_l^\infty \beta^{u+v - 2} (1+ \beta)^{-u-N} d\beta,
\end{align}
with the Pochhammer symbol $(a)_k$ being defined by
$$(a)_k = \frac{\Gamma(a+k)}{\Gamma(a)} = a(a+1)(a+2) \cdots (a+k-1).$$
It is not hard to see that $T_0$ is independent of the choice of $N$ (see Lemma \ref{R0T0lemma}). 

It appears that the identity (\ref{q0identity}) provides the most powerful approach available for determining the asymptotics of $Q_0(s)$. For example, the expansion of corollary \ref{Q0cor} for $Q_0(\sigma + it)$ with $\sigma \neq 1/2$ is easily obtained from (\ref{q0identity}) by taking $u = \bar{v} = \sigma + it$ and $N = 1$, and noting that the sum over $l$ in (\ref{T0def}) is $O(t^{-1})$; the expansion for $\sigma = 1/2$ is obtained in a similar way by taking $N =1$ and $\sigma = (1 + \delta)/2$ and letting $\delta \to 0$, see \cite{KM1996} or the proof of theorem \ref{QNASYMPTOTICSTH}.

The expression (\ref{qnidentity}) for $q_n$ given in theorem \ref{QNTH} is the natural generalization to any integer $n \in \Z$ of the expression (\ref{q0identity}) for $q_0$. Indeed, since $b_0(s) = \frac{1}{s-1}$ clearly is the analytic continuation of the integral $\int_1^\infty \alpha^{-s} d\alpha$, this follows from (\ref{bnRnTndef}) and the following lemma.

\begin{lemma}\label{R0T0lemma}
The functions $R_0$ and $T_0$ defined in (\ref{R0def}) and (\ref{T0def}) are the analytic continuations of the integrals
\begin{align*}
& R_0(u,v) = \int_0^\infty \alpha^{-v} \zeta_1(u, \alpha) d\alpha, \qquad \re (u+v) > 2, \quad \re v < 1,
	\\ 
& T_0(u,v) = \int_0^1 \alpha^{-v}  \zeta_1(u, \alpha)  d\alpha, \qquad  \re v < 1.
\end{align*}
\end{lemma}
\begin{proof}
Employing the integral representation 
\begin{align}\label{X1.4}
  \zeta_1(u,\alpha) = \frac{1}{\Gamma(u)}\int_0^\infty \frac{e^{-\alpha r} r^{u-1}}{e^r -1} dr, \qquad \re u > 1, \quad \re \alpha > -1,
\end{align}
and using Fubini's theorem to change the order of integration, we find, for $\re (u+v) > 2$ and $\re v < 1$,
\begin{align*}
\int_0^\infty \alpha^{-v} \zeta_1(u,\alpha) d\alpha
= \frac{1}{\Gamma(u)} \int_0^\infty \frac{ r^{u-1}}{e^r -1} \int_0^\infty \alpha^{-v} e^{-\alpha r} d\alpha dr.
\end{align*}
The change of variables $\beta = \alpha r$ shows that the rhs can be rewritten as
$$\frac{1}{\Gamma(u)} \int_0^\infty \frac{r^{u + v-2}}{e^r -1} \int_0^\infty \beta^{-v} e^{-\beta} d\beta dr.$$
Since the $\beta$-integral equals $\Gamma(1-v)$, we arrive at
\begin{align*}
& \int_0^\infty \alpha^{-v} \zeta_1(u,\alpha) d\alpha
 = \frac{\Gamma(1-v)}{\Gamma(u)} \int_0^\infty  \frac{r^{u + v- 2}}{e^r - 1}dr
	\\
& =  \frac{\Gamma(1-v)}{\Gamma(u)}\Gamma(u+v-1)\zeta(u+v-1), \qquad  \re (u+v) > 2, \quad \re v < 1,
\end{align*}
which proves the statement for $R_0$.

Integrating by parts repeatedly using the facts that $\zeta_1(s,1) = \zeta(s) -1$ and
\begin{align}\label{partialalphazeta1}
\partial_\alpha \zeta_1(s, \alpha) = -s\zeta_1(s+1, \alpha),
\end{align}
we obtain, for any integer $N \geq 1$ and $\re v < 1$,
\begin{align}\nonumber
& \int_0^1 \alpha^{-v}  \zeta_1(u, \alpha)  d\alpha = \frac{\zeta(u) - 1}{1-v} + \frac{u}{1-v} \int_0^1 \alpha^{1-v}  \zeta_1(u+1, \alpha)  d\alpha
	\\ \label{int01alphaminusv}
& = \cdots = \sum_{k = 0}^{N-1} \frac{(u)_k (\zeta(u + k) - 1)}{(1-v)_{k+1}} + \frac{(u)_N}{(1-v)_N}\int_0^1 \alpha^{N-v}   \zeta_1(u+N, \alpha) d\alpha.
\end{align}
The statement for $T_0$ will follow once we prove that
\begin{align}\label{T0midstep}
\int_0^1 \alpha^{N-v}\zeta_1(u+N, \alpha)   d\alpha = \sum_{l=1}^\infty l^{1-u-v} \int_l^\infty \beta^{u+v  - 2} (1+ \beta)^{-u-N} d\beta
\end{align}
for $\re v < N+1$ and $\re u > -N +1$. But for $\re u > -N +1$ we can use the representation 
\begin{align}\label{zeta1sum}
\zeta_1(s,\alpha) = \sum_{l=1}^\infty \frac{1}{(l + \alpha)^s}, \qquad \re s > 1, \quad \re \alpha > -1,
\end{align}
to write the left-hand side of (\ref{T0midstep}) as
$$\sum_{l=1}^\infty \int_0^1 \alpha^{N-v} (l + \alpha)^{-u-N}  d\alpha.$$
The change of variables $\alpha = l/\beta$ completes the proof of (\ref{T0midstep}) and hence also of the lemma.
\end{proof}

\begin{remark}
The formula for $T_n$ given in (\ref{Tncontinuation}) is obtained by integrating by parts once in the integral (\ref{Tndef}). 
If we instead integrate by parts $N$ times, we find
\begin{align}\nonumber
T_n(u,v) = & \int_0^1  \alpha^{-v} \zeta_1(u, \alpha)  e^{-2\pi i n\alpha}d\alpha 
 = \sum_{k=0}^{N-1} \frac{(-1)^k}{(1-v)_{k+1}} \frac{\partial^k}{\partial \alpha^k}\bigg|_{\alpha = 1}( \zeta_1(u, \alpha) e^{-2\pi i n\alpha}) 
	\\ \label{preTnIBP}
& + \frac{(-1)^N}{(1-v)_N} \int_0^1 \alpha^{N-v}   \frac{\partial^N}{\partial \alpha^N}(\zeta_1(u, \alpha) e^{-2\pi i n\alpha}) d\alpha 
	\\\nonumber
= &\; \sum_{k=0}^{N-1} \sum_{j = 0}^k \begin{pmatrix} k \\ j \end{pmatrix} \frac{(u)_j(2\pi i n)^{k-j}}{(1-v)_{k+1}}  (\zeta(u + j) - 1)
	\\\nonumber
& + \sum_{j = 0}^N 
\begin{pmatrix} N \\ j \end{pmatrix} 
\frac{(u)_j(2\pi i n)^{N-j}}{(1-v)_N} 
\int_0^1 \alpha^{N-v} \zeta_1(u + j, \alpha)  e^{-2\pi i n \alpha} d\alpha.
\end{align}
Using that (cf. (\ref{T0midstep})), for $2 \leq j \leq N$ with $\re v < N+1$ and $\re u > -j +1$, 
\begin{align*}
& \int_0^1 \alpha^{N-v}\zeta_1(u+j, \alpha) e^{-2\pi i n \alpha}  d\alpha 
	\\
& = \sum_{l=1}^\infty l^{1-u-v+N-j} \int_l^\infty \beta^{u+v  - 2-N+j} (1+ \beta)^{-u-j} e^{-2\pi i n l/\beta} d\beta,
\end{align*}
this leads to the following expression for $T_n$ when $\re u > -1$ and $\re v < N+1$, which is the natural extension to any integer $n$ of the expression (\ref{T0def}) for $T_0$:
\begin{align}\nonumber
 T_n(u,v) = & \sum_{k=0}^{N-1} \sum_{j = 0}^k \begin{pmatrix} k \\ j \end{pmatrix} \frac{(u)_j(2\pi i n)^{k-j}}{(1-v)_{k+1}}  (\zeta(u + j) - 1)
+ \sum_{j = 2}^N 
\begin{pmatrix} N \\ j \end{pmatrix} 
\frac{(u)_j(2\pi i n)^{N-j}}{(1-v)_N} 
	\\\label{TnIBP}
& 
\times \sum_{l=1}^\infty l^{1-u-v+N-j} \int_l^\infty \beta^{u+v  - 2-N+j} (1+ \beta)^{-u-j} e^{-2\pi i n l/\beta} d\beta
	\\ \nonumber
& + \frac{(2\pi i n)^{N}}{(1-v)_N} 
\int_0^1\alpha^{N-v}  \zeta_1(u, \alpha) e^{-2\pi i n \alpha} d\alpha
	\\\nonumber
& + \frac{u(2\pi i n)^{N-1}}{(1-v)_N} 
\int_0^1\alpha^{N-v}  \zeta_1(u+1, \alpha) e^{-2\pi i n \alpha} d\alpha.
\end{align}
\end{remark}

\section{Proof of Theorem \ref{QNTH}}\label{qnidentitysec}

\begin{lemma}[Fourier series of $\zeta_1(u, \alpha) \zeta_1(v, \alpha)$]
Suppose $u,v \in \C \setminus \{1\}$ and $\delta \in (0, \pi/2)$.  Then the Fourier series representation of $\zeta_1(u, \cdot) \zeta_1(v, \cdot) \in C^\infty([0,1])$ is given by
\begin{align}\label{zeta1zeta1fourier}
\zeta_1(u, \alpha) \zeta_1(v, \alpha) = \sum_{n \in \Z} q_n(u,v) e^{2\pi i n \alpha}, 
\end{align}
where
\begin{align}\label{qnexpression}
& q_n(u,v) = \int_1^{e^{-(\sgn n) i \delta} \infty}  \big(\alpha^{-u-v} + \alpha^{-v} \zeta_1(u, \alpha) + \alpha^{-u} \zeta_1(v,\alpha) \big) e^{-2\pi i n \alpha} d\alpha, 
	\\\nonumber
& \hspace{9cm} n \in \Z \setminus \{0\},
\end{align}
and
$$q_0(u,v) = \int_0^1  \zeta_1(u, \alpha) \zeta_1(v,\alpha) d\alpha.$$
The series in (\ref{zeta1zeta1fourier}) converges in $L^2([0,1])$ and converges pointwise for each $\alpha \in (0,1)$. 
If $\re(u+v) > 2$, the Fourier coefficients can be expressed as
\begin{align}\label{qnexpression2}
& q_n(u,v) = \int_1^\infty \big(\alpha^{-u-v} + \alpha^{-v} \zeta_1(u, \alpha) +  \alpha^{-u} \zeta_1(v,\alpha)\big) e^{-2\pi i n \alpha}  d\alpha, 
	\\\nonumber
&\hspace{4cm} n \in \Z, \quad u,v \in \C \setminus \{1\}, \quad \re(u+v) > 2.
\end{align}

\end{lemma}
\begin{proof}
Equation (\ref{qnexpression}) follows from (\ref{modifiedhurwitzshift}) and easy algebra. Indeed, for $u,v \in \C \setminus \{1\}$ and $n \geq 1$, we compute
\begin{align*}
q_n(u,v) 
 = &\; \int_0^1 \zeta_1(u,\alpha) \zeta_1(v,\alpha) e^{-2\pi i n \alpha}  d\alpha
 	\\
 = &\; \bigg(\int_0^{e^{-i\delta}\infty}  - \int_1^{e^{-i\delta}\infty} \bigg) \zeta_1(u,\alpha) \zeta_1(v,\alpha) e^{-2\pi i n \alpha}  d\alpha
 	\\
= &\; \int_0^{e^{-i\delta}\infty} \zeta(u, \alpha +1) \zeta(v, \alpha +1) e^{-2\pi i n \alpha} d\alpha
	\\
&- \int_1^{e^{-i\delta}\infty} (\zeta(u,\alpha) - \alpha^{-u})(\zeta(v,\alpha) - \alpha^{-v}) e^{-2\pi i n \alpha}  d\alpha.
\end{align*}
Changing variables $\beta = \alpha +1$ in the integral from $0$ to $e^{-i\delta} \infty$, we infer that
\begin{align*}
q_n(u,v) = &\; \int_1^{e^{-i\delta}\infty}  \zeta(u,\beta) \zeta(v,\beta) e^{-2\pi i n \beta} d\beta
- \int_1^{e^{-i\delta}\infty} \zeta(u,\alpha)\zeta(v,\alpha)  e^{-2\pi i n \alpha}  d\alpha
	\\
& + \int_1^{e^{-i\delta}\infty} \big( \alpha^{-v}  \zeta(u,\alpha)
+ \alpha^{-u}  \zeta(v,\alpha)  - \alpha^{-u-v} \big)  e^{-2\pi i n \alpha}  d\alpha
	\\
= &\; \int_1^{e^{-i\delta}\infty} \big(\alpha^{-v}  \zeta_1(u,\alpha) 
+ \alpha^{-u}  \zeta_1(v,\alpha)  + \alpha^{-u-v} \big) e^{-2\pi i n \alpha} d\alpha,
\end{align*}
which proves (\ref{qnexpression}) for $n \geq 1$; the same computation applies also when $n \leq -1$ if $\delta$ is replaced by $-\delta$.
The expression in (\ref{qnexpression2}) follows because, by (\ref{Hurwitzalphaasymptotics}), $\zeta(s, \alpha) = \frac{\alpha^{1-s}}{s-1} + O(\alpha^{-s})$ as $\alpha \to \infty$, so if $\re(u+v) > 3$ and $u,v \neq 1$, the above proof can be repeated for each $n \in \Z$ with $\delta = 0$; the resulting formula can then be extended to $\re(u+v) > 2$ by analytic continuation.
\end{proof}

\begin{proof}[Proof of Theorem \ref{QNTH}]
Equations (\ref{qnidentity}) and (\ref{bnRnTndef}) are a direct consequence of (\ref{qnexpression2}).
It remains to establish the meromorphic continuations in (\ref{bnRnTncontinuation}).
For $n = 0$, this was carried out already in section \ref{zerothsec} (see Lemma \ref{R0T0lemma} and equation (\ref{int01alphaminusv}) with $N = 1$).

Suppose $n$ is a nonzero integer. Making the change of variables $r = 2\pi i n \alpha$ in the integral in (\ref{bndef}) and deforming the contour back to the real axis, we obtain, for $\re s > 1$,
$$b_n(s) = \int_{2\pi i n}^{\infty} r^{-s} (2\pi i n)^{s-1} e^{-r} dr = \Gamma(1-s, 2\pi i n)(2\pi i n)^{s-1}.$$
This proves (\ref{bncontinuation}). 

To determine the analytic continuation of $R_n$, we substitute the integral representation (\ref{X1.4}) for $\zeta_1$ into (\ref{Rndef}) and change the order of integration. This gives, for $\re(u + v) > 2$ and $\re v < 1$,
\begin{align*}
\int_0^\infty \alpha^{-v} \zeta_1(u,\alpha) e^{-2\pi i n\alpha} d\alpha
= \frac{1}{\Gamma(u)} \int_0^\infty \frac{ r^{u-1}}{e^r -1} \int_0^\infty \alpha^{-v} e^{-\alpha (r + 2\pi i n)} d\alpha dr.
\end{align*}
The change of variables $\beta = \alpha(r + 2\pi i n)$ shows that the rhs can be rewritten as
$$\frac{1}{\Gamma(u)} \int_0^\infty \frac{r^{u-1}}{e^r -1} (r+2\pi i n)^{v-1} \int_0^\infty \beta^{-v} e^{-\beta} d\beta dr,$$
where we have used the decay of $e^{-\beta}$ to deform the contour back to the positive real axis. 
Since the $\beta$-integral equals $\Gamma(1-v)$, we arrive at the formula
\begin{align}\label{int0inftyalphaminusv}
& \int_0^\infty \alpha^{-v} \zeta_1(u,\alpha) e^{-2\pi i n \alpha} d\alpha
= \frac{\Gamma(1-v)}{\Gamma(u)} \int_0^\infty  \frac{r^{u-1}}{e^r - 1} (r + 2\pi i n)^{v-1} dr, 
	\\ \nonumber
&\hspace{6cm} \re(u + v) > 2, \quad \re v < 1.
\end{align}
The integral on the rhs can be expressed in terms of $\Phi$: 
\begin{align}\label{X1.6}
  \int_0^\infty \frac{r^{u-1}}{e^r - 1}(r + 2\pi i n)^{v-1} dr = \Phi(u,v, n) \times \begin{cases} 
  \frac{1}{e^{-i\pi(u+v)} - e^{i\pi(u-v)}}, \quad &  n > 0, \\
  \frac{1}{e^{-i\pi(u-v)} - e^{i\pi(u+v)}}, & n < 0,
  \end{cases} 
\end{align}
for $\re u >1$ and $v \in \C$.
Indeed, for $\re u > 1$ the contribution from the circular part of the Hankel contour $H_\alpha$ in (\ref{X1.1}) vanishes in the limit $\alpha \to 0$; hence (\ref{X1.1}) can be written as
\begin{align*}
  \Phi(u,v,n) = (e^{i\pi u} - e^{-i\pi u}) \int_0^\infty \frac{r^{u-1}}{e^r -1}(-r-2\pi i n)^{v-1}dr.
\end{align*}
For $n > 0$ (resp. $n < 0$), $-r-2\pi i n$ lies in the third (resp. second) quadrant, and so $(-r-2\pi i n)^{v-1} = e^{-(\sgn n)i\pi(v-1)} (r+2\pi i n)^{v-1}$, which proves (\ref{X1.6}).
The expression (\ref{Rncontinuation}) for $R_n$ follows from (\ref{int0inftyalphaminusv}) and (\ref{X1.6}).

The expression (\ref{Tncontinuation}) for $T_n$ is the special case $N =1$ of (\ref{preTnIBP}).
\end{proof}

\section{Proof of Theorem \ref{QNASYMPTOTICSTH}}\label{qnasymptoticssec}

\begin{lemma}\label{Rnlemma}
Let $u = \sigma_1 + it$ and $v = \sigma_2 - it$. Then, for every $N \geq 1$, $R_n(v, u) = O(t^{-N})$ uniformly for $n \geq 1$ and $\sigma_1, \sigma_2 \in [0,1]$ as $t \to \infty$. 
\end{lemma}
\begin{proof}
This is an immediate consequence of the expression (\ref{Rncontinuation}) for $R_n$, the estimate (\ref{Phivuestimate}) of $\Phi(v,u,\beta)$, and the fact that (see (\ref{Gamma1minuss}))
$$\frac{\Gamma(1-u)}{\Gamma(v)}\frac{ie^{\pi i u}}{2\sin(\pi v)} = O(t^{1 - \sigma_1 - \sigma_2}e^{-2\pi t}), \qquad t \to \infty,$$
uniformly for $\sigma_1, \sigma_2 \in [0,1]$. 
\end{proof}

\begin{lemma}\label{Tnlemma}
Let $u = \sigma_1 + it$, $v = \sigma_2 - it$, and $n \in \Z$. Then, as $t \to \infty$,
\begin{subequations}\label{Tnasymptotics}
\begin{align}\label{Tnasymptoticsa}
& T_n(u,v) = \frac{\zeta(u) - 1}{1-v} + \frac{(\zeta(u+1) - 1)u}{(1-v)(2-v)}  + O(t^{-1}), 
	\\\label{Tnasymptoticsb}
& T_n(v, u) = \frac{\zeta(v) - 1}{1-u} + \frac{(\zeta(v+1) - 1)v}{(1-u)(2-u)} + O(t^{-1}), 
\end{align}
\end{subequations}
uniformly for $\sigma_1, \sigma_2 \in [0, 1]$.
\end{lemma}
\begin{proof}
Consider the expression (\ref{TnIBP}) for $T_n$ for some integer $N \geq 3$. 
An integration by parts shows that the integral with respect to  $d \beta$ in (\ref{TnIBP}) can be written as follows for $\re v < N+1$:
\begin{align}\label{intE1E2E3}
& \int_l^\infty \beta^{u+v-2-N+j} (1+ \beta)^{-u-j} e^{-2\pi i n l/\beta} d\beta
= \frac{(1+ l)^{1-u-j}}{u+j-1} l^{u+v-2-N+j} 
	\\ \nonumber
& - \int_l^\infty \frac{(1+ \beta)^{1-u-j}}{1-u-j} \frac{\partial}{\partial \beta}\big(\beta^{u+v-2-N+j} e^{-2\pi i n l/\beta}\big)d\beta
 = E_1 + E_2 + E_3,
\end{align}
where we have used the short-hand notation
\begin{align*}
& E_1 := \frac{(1+ l)^{1-u-j}}{u+j-1} l^{u+v-2-N+j},
	\\
& E_2 := (u+v-2-N+j) \int_l^\infty \frac{(1+ \beta)^{1-u-j}}{u+j-1} \beta^{u+v-3-N+j} e^{-2\pi i n l/\beta}d\beta	,
	\\
& E_3 := 2\pi i n l \int_l^\infty \frac{(1+ \beta)^{1-u-j}}{u+j-1} \beta^{u+v-4-N+j} e^{-2\pi i n l/\beta}d\beta.
\end{align*}
Substituting (\ref{intE1E2E3}) into (\ref{TnIBP}), we see that $T_n(u,v)$ can be expressed for $\re u > -1$ and $\re v < N+1$ as 
\begin{align}\label{TnIBP2}
T_n(u,v) = &\; \sum_{k=0}^{N-1} \sum_{j = 0}^k \begin{pmatrix} k \\ j \end{pmatrix} \frac{(u)_j(2\pi i n)^{k-j}}{(1-v)_{k+1}}  (\zeta(u + j) - 1)
 + \sum_{m=1}^5 F_m,
\end{align}
where we have used the short-hand notation
\begin{align*}
& F_m := \sum_{j = 2}^N 
\begin{pmatrix} N \\ j \end{pmatrix} 
\frac{(u)_j(2\pi i n)^{N-j}}{(1-v)_N} 
\sum_{l=1}^\infty l^{1-u-v+N-j} E_m, \qquad m = 1,2,3,
	\\
& F_4 := \frac{(2\pi i n)^{N}}{(1-v)_N} 
\int_0^1\alpha^{N-v}  \zeta_1(u, \alpha) e^{-2\pi i n \alpha} d\alpha,
	\\
& F_5 := \frac{u(2\pi i n)^{N-1}}{(1-v)_N} 
\int_0^1\alpha^{N-v}  \zeta_1(u+1, \alpha) e^{-2\pi i n \alpha} d\alpha.
\end{align*}
The first two terms on the rhs of (\ref{Tnasymptoticsa}) come from the terms with $k = j =0$ and $k = j = 1$ in the double sum in (\ref{TnIBP2}). All other terms in this double sum are $O(t^{-1})$ by standard estimates of $\zeta(s)$. Thus the asymptotic formula (\ref{Tnasymptoticsa}) for $T_n(u,v)$ will follow if we can show that $F_m = O(t^{-1})$ uniformly for $\sigma_1, \sigma_2 \in [0,1]$ for each $m = 1,\dots, 5$. 

The terms $F_1$, $F_2$, and $F_3$ are easily estimated:
\begin{align*}
|F_1|& \leq \frac{C}{t} \sum_{j = 2}^N \sum_{l=1}^\infty l^{-1} (1+l)^{1-\sigma_1-j}
= O(t^{-1}),
	\\
 |F_2| & \leq C \sum_{j=2}^N \sum_{l=1}^\infty l^{1-\sigma_1 - \sigma_2 + N - j} \int_l^\infty \frac{(1+ \beta)^{1- \sigma_1 -j} }{|u+j-1|} \beta^{\sigma_1 + \sigma_2 -3-N+j} d\beta
	\\
&\leq 
\frac{C}{t} \sum_{j=2}^N \sum_{l=1}^\infty l^{-\sigma_1 - j} = O(t^{-1}),
	\\
 |F_3| &
\leq
C \sum_{j = 2}^N \sum_{l=1}^\infty l^{2-\sigma_1 - \sigma_2 + N - j} \int_l^\infty \frac{(1+ \beta)^{1- \sigma_1 -j} }{|u+j-1|} \beta^{\sigma_1 + \sigma_2 -4-N+j} d\beta
	\\
& \leq
\frac{C}{t}
\sum_{j = 2}^N\sum_{l=1}^\infty l^{-\sigma_1-j}
= O(t^{-1})	,
\end{align*}
uniformly for $\sigma_1, \sigma_2 \in [0,1]$. 

Integrating by parts twice and using (\ref{partialalphazeta1}), we find, for $\sigma_1, \sigma_2 \in [0,1]$,
\begin{align}\label{int01alphaIBP}
 \int_0^1 \alpha^{N-v}  e^{-2\pi i n \alpha} \zeta_1(u, \alpha)  d\alpha
= &\; \partial_\alpha f(v, 1) \zeta_1(u, 1) + u f(v, 1) \zeta_1(u+1, 1)
	\\\nonumber
& + u(u+1) \int_0^1 f(v, \alpha) \zeta_1(u+2, \alpha)  d\alpha,
\end{align}
where $f(v,\alpha)$ is defined by
$$f(v, \alpha) = \int_0^\alpha \int_0^{\alpha_1} {\alpha_2}^{N-v}  e^{-2\pi i n \alpha_2} d\alpha_2 d\alpha_1.$$
The functions $f(v, \alpha)$ and $\partial_\alpha f(v, \alpha)$ are uniformly bounded for $\alpha \in [0,1]$ and $\re v \in [0,1]$. Also, by a trivial estimate of the sum in (\ref{zeta1sum}), $\zeta_1(u+2, \alpha)$ is uniformly bounded for $\alpha \in [0,1]$ and $\re u \in [0,1]$. Since $\zeta_1(u, 1) = \zeta(u) - 1$ and $\zeta_1(u+1, 1) = \zeta(u+1) - 1$, we conclude that the rhs of (\ref{int01alphaIBP}) is $O(t^2)$ uniformly for $\sigma_1, \sigma_2 \in [0,1]$. It follows that $F_4 = O(t^{2-N}) = O(t^{-1})$ uniformly for $\sigma_1, \sigma_2 \in [0,1]$ as $t \to \infty$. A similar argument shows that $F_5$ also is $O(t^{2-N})$.
This completes the proof of the estimate (\ref{Tnasymptoticsa}) for $T_n(u,v)$; the estimate (\ref{Tnasymptoticsb}) for $T_n(v,u)$ follows by interchanging $u$ and  $v$ in the above arguments. 
\end{proof}

\begin{proof}[Proof of Theorem \ref{QNASYMPTOTICSTH}]
Employing lemma \ref{Rnlemma} and lemma \ref{Tnlemma} in the expression (\ref{qnidentity}) for $q_n(u,v)$ with $n \geq 1$, we find formula (\ref{qnasymptoticsPhi}).

The formula (\ref{q0asymptoticsneq1}) for the asymptotics of $q_0(u,v)$ when $\sigma_1 + \sigma_2 \neq 1$ follows immediately by substituting the expression (\ref{Rncontinuation}) for $R_0$ into equation (\ref{qnidentity}) and using Lemma \ref{Tnlemma}.
To derive the formula (\ref{q0asymptoticseq1}) for the asymptotics of $q_0(u,v)$ when $\sigma_1 + \sigma_2 = 1$, we first set $u = \sigma + \delta + it$ and $v = 1 - \sigma + \delta - it$ in equation (\ref{qnidentity}) and take the limit $\delta \to 0$; this yields 
$$q_0(u,v) = - \ln 2\pi + \gamma + \frac{\psi(u) + \psi(v)}{2} + T_0(u,v) + T_0(v,u)$$
where $u = \sigma + it$, $v = 1 - \sigma - it$, and $\psi(s) = \Gamma'(s)/\Gamma(s)$. 
The formula (\ref{q0asymptoticseq1}) then follows from Lemma \ref{Tnlemma} and the fact that $\psi(\sigma \pm it) = \ln t \pm \frac{\pi i}{2} + O(t^{-1})$ uniformly for $\sigma \in [0,1]$ as $t \to \infty$.
\end{proof}

\section{Fourier coefficients of $\zeta(s,\alpha)$ and $\zeta_1(s, \alpha)$}
In this section, we consider, for completeness, the Fourier coefficients of $\zeta(s,\alpha)$ and $\zeta_1(s,\alpha)$. 
In particular, we will see that the coefficients $b_n(s)$ defined in (\ref{bndef}) are the Fourier coefficients of $\zeta_1(s,\alpha)$, cf. \cite{R1997, AFpreprint}.

We first compute the Fourier coefficients $a_n(s)$ of $\zeta(s,\alpha)$ defined by
$$a_n(s) = \int_0^1 \zeta(s, \alpha) e^{-2\pi i n \alpha}  d\alpha, \qquad n \in \Z.$$
The coefficients $a_n(s)$ are well-defined whenever $\re s < 1$, because the relation $\zeta(s, \alpha) = \zeta(s, \alpha +1) + \alpha^{-s}$ implies that $\zeta(s, \alpha) \sim \alpha^{-s}$ as $\alpha \downarrow 0$; hence $\zeta(s, \cdot) \in L^1([0,1])$ iff $\re s < 1$. 
For $0 < \re s < 1$, the next lemma was proved in \cite{R1997} using different techniques.

\begin{lemma}[Fourier series of $\zeta(s,\alpha)$]
Suppose $\re s < 1$. Then the Fourier series representation of $\zeta(s,\cdot) \in L^1([0,1]) \cap C^\infty((0,1])$ is given by
\begin{align}\label{hurwitzfourier}
\zeta(s,\alpha) = \sum_{n \in \Z} a_n(s) e^{2\pi i n \alpha}, \qquad a_n(s) = \begin{cases} \Gamma(1-s) (2\pi i n)^{s-1}, & n \in \Z \setminus \{0\}, \\ 
0, & n = 0.
\end{cases}
\end{align}
\end{lemma}
\begin{proof}
For $\re s < 0$, a computation using the representation (\ref{hurwitzHankel}) for $\zeta(s,\alpha)$ gives
\begin{align}\nonumber
a_n(s) & = \int_0^1 \zeta(s, \alpha) e^{-2\pi i n \alpha} d\alpha
= \frac{\Gamma(1-s)}{2\pi i} \int_0^1 e^{-2\pi i n \alpha} \int_{H_1} \frac{e^{\alpha z} z^{s-1}}{1 - e^z} dz d\alpha
	\\\nonumber
& = \frac{\Gamma(1-s)}{2\pi i} \int_{H_1}  \int_0^1 e^{\alpha (z-2\pi i n)} d\alpha \frac{ z^{s-1}}{1 - e^z} dz 
	\\\nonumber
& = \frac{\Gamma(1-s)}{2\pi i} \int_{H_1} \frac{e^{z-2\pi i n} - 1}{z -2 \pi i n} \frac{ z^{s-1}}{1 - e^z} dz 
= -\frac{\Gamma(1-s)}{2\pi i} \int_{H_1} \frac{ z^{s-1}}{z - 2 \pi i n} dz
	\\\label{ancomputation}
& 
= \begin{cases} \Gamma(1-s) (2\pi i n)^{s-1}, & n \in \Z \setminus \{0\}, \\ 
0, & n = 0,
\end{cases} \quad \re s < 0,
\end{align}
where, in the last step, we have deformed the Hankel contour $H_1$ to infinity and used the residue theorem. By writing
$$a_n(s) = \int_0^1  \zeta(s, \alpha +1) e^{-2\pi i n \alpha} d\alpha
+ \int_0^1  \alpha^{-s} e^{-2\pi i n \alpha} d\alpha, \qquad \re s < 0, \; \; n \in \Z,$$
we see that each Fourier coefficient $a_n(s)$ is analytic for $\re s < 1$.
Since the expression in (\ref{hurwitzfourier}) for $a_n(s)$ also is analytic for $\re s < 1$ for each $n$, the lemma follows from (\ref{ancomputation}) and analytic continuation.
\end{proof}

\begin{remark}
Since $\zeta(s,\cdot) \in L^1([0,1]) \cap C^\infty((0,1])$, the Fourier series in (\ref{hurwitzfourier}) converges pointwise for each $\alpha \in (0,1)$. Easy estimates show that the convergence is absolute and uniform for $s$ in compact subsets of $\re s < 0$ and $\alpha \in [0,1]$. This reflects the fact that $\zeta(s,\alpha)$ can be extended continuously to a function $\zeta(s, \cdot) \in C([0,1])$ with $\zeta(s,0) = \zeta(s, 1)$ if $\re s < 0$ (see (\ref{modifiedhurwitzshift})). 
\end{remark}

\begin{lemma}[Fourier series of $\zeta_1(s,\alpha)$]
Suppose $s \in \C \setminus \{1\}$. Then the Fourier series representation of $\zeta_1(s,\cdot) \in C^\infty([0,1])$ is given by
\begin{align}\label{modifiedhurwitzfourier}
\zeta_1(s,\alpha) = \sum_{n \in \Z} b_n(s) e^{2\pi i n \alpha}, 
\end{align}
where $b_n(s)$, $n \in \Z$, are the coefficients in (\ref{bncontinuation}).
\end{lemma}
\begin{proof}
We need to show that
\begin{align}\label{modifiedhurwitzint}
\int_0^1 \zeta_1(s, \alpha) e^{-2\pi i n \alpha}  d\alpha = \begin{cases} \Gamma(1-s, 2\pi i n)(2\pi i n)^{s-1}, & n \neq 0, \\
\frac{1}{s-1}, & n = 0, \end{cases} \quad s \in \C \setminus \{1\}.
\end{align}
Actually, since both sides of (\ref{modifiedhurwitzint}) are analytic for $s \in \C \setminus \{1\}$ for each $n \in \Z$, it is enough to show (\ref{modifiedhurwitzint}) for $\re s < 1$. 
By the definition (\ref{modifiedhurwitzshift}) of $\zeta_1$, we have
$$\int_0^1  \zeta_1(s, \alpha) e^{-2\pi i n \alpha} d\alpha
= a_n(s) - \int_0^1 \alpha^{-s} e^{-2\pi i n \alpha} d\alpha, \qquad \re s < 1.$$
The change of variables $r = 2\pi i n \alpha$ shows that
$$-\int_0^1 \alpha^{-s} e^{-2\pi i n \alpha} d\alpha
= (\Gamma(1-s, 2\pi i n) - \Gamma(1-s)) (2\pi i n)^{s-1}, \qquad \re s < 1.$$
Recalling the expression for $a_n(s)$ in (\ref{hurwitzfourier}), the lemma follows.
\end{proof}

\begin{remark}
Since $\zeta_1(s,\cdot) \in C^\infty([0,1])$, the Fourier series in (\ref{modifiedhurwitzfourier}) converges pointwise for each $\alpha \in (0,1)$. 
For fixed $s$, the incomplete Gamma function admits the asymptotic expansion (see e.g. \cite[Eq. (8.11.2)]{NIST})
$$\Gamma(s,z) = z^{s-1} e^{-z}\bigg(1 + \sum_{k=1}^N \frac{(s-1) \cdots (s-k)}{z^k} + O(z^{-N-1})\bigg), \qquad z \to \infty,$$
where the error term is uniform for $|\arg z| \leq \frac{3\pi}{2} - \delta$.
Hence
$$b_n(s) 
= (2\pi i n)^{-1}(1 + O(n^{-1})), \qquad n \to \infty,$$
showing that the Fourier series for $\zeta_1$ is not absolutely convergent for any $s \in \C \setminus \{1\}$. 
\end{remark}

\part{Representations for the Basic Sum}

\chapter{Several Representations for the Basic Sum}\label{sec9}
The purpose of this chapter is to present integral representations of the basic sum $\sum_a^b n^{s-1}$ for certain values of $a$ and $b$. 

Let the contour $C_t^\eta$, $t < \eta$, denote the semicircle from $it$ to $i\eta$ with $\re z \geq 0$. Splitting the contour of the second integral in the rhs of equation (\ref{1.2}) into the contour $C_\eta^t$ plus the ray from $it$ to $\infty \exp(i\phi_2)$, we find 
\begin{align}\nonumber
\zeta(1-s) = &\sum_{n=1}^{[\frac{\eta}{2\pi}]} n^{s-1} - \frac{\eta^s}{s(2\pi)^s} 
- \frac{e^{-\frac{i\pi s}{2}}}{(2\pi)^s} \int_{C_t^\eta} \frac{z^{s-1}}{e^z - 1} dz
	\\ \label{6.10}
& + G_L(t, \sigma; \eta) + G_U(t, \sigma; t),
\qquad 0 \leq \sigma \leq 1, \quad  0< t < \eta,
\end{align}
where $G_L$ and $G_U$ are defined in (\ref{GLdef}) and (\ref{GUdef}) respectively.
The term $G_U(t, \sigma; t)$ was computed to all orders as $t \to \infty$ in theorem \ref{th3.2}, see equations (\ref{GU2final}) and (\ref{GU1asymptotics}). The term $G_L(t, \sigma; \eta)$ with $t < \eta$ was computed to all orders as $t \to \infty$ in (\ref{GLfinal}). Thus, by comparing equation (\ref{6.10}) with the representation obtained in theorem \ref{th3.2}, it follows that we can express the fundamental sum $\sum_{n = [t/2\pi]+1}^{[\eta/2\pi]} n^{s-1}$ in terms of the integral appearing in the rhs of (\ref{6.10}), where the error is computed to all orders. For brevity of presentation, we state this result only to leading order.

\begin{lemma}[An integral representation for the basic sum]
Let $C_t^\eta$, $t < \eta$, denote the semicircle from $z = it$ to $z = i\eta$ with $\re z \geq 0$. 
For every $\epsilon > 0$,
\begin{align}\label{6.11}
  \sum_{n=[\frac{t}{2\pi}]+1}^{[\frac{\eta}{2\pi}]} n^{s-1} 
  = &\;  \frac{e^{-\frac{i\pi s}{2}}}{(2\pi)^s} \int_{C_t^\eta} \frac{z^{s-1}}{e^z - 1} dz
 + \frac{\eta^s}{s(2\pi)^s} + \frac{i\eta^{s-1}}{(2\pi)^s}\ln(1 - e^{i\eta}) 
 	\\ \nonumber
& - \frac{it^{s-1}e^{-it}}{(2\pi)^s}\ln(1 - e^{it}) + O\biggl(\frac{1}{t^{2-\sigma}} + \frac{t}{\eta^{2-\sigma}}\biggr), 	\\ \nonumber
& \hspace{2cm} (1+\epsilon)t < \eta < \infty, \quad 0 \leq \sigma \leq 1, \quad t \to \infty,
\end{align}
where the error term is uniform for all $\eta, \sigma$ in the above ranges. 
\end{lemma}
\begin{proof}
   Letting $\sigma \to 1-\sigma$ in (\ref{zetaformula2.2intro}) and taking the complex conjugate of the resulting equation, we find
\begin{align*}
\zeta(1-s) = &\; \sum_{n=1}^{[\frac{t}{2\pi}]} n^{s-1} - \frac{1}{s}\biggl(\frac{t}{2\pi}\biggr)^s 
+ \frac{it^{s-1}e^{-it}}{(2\pi)^s}(-e^{it} + 2i\im \Li_1(e^{it}))
	\\
& + \frac{t^{s-1} e^{-it} }{(2\pi)^{s}}\biggl(\frac{1-i}{2} \sqrt{\pi t} + \frac{i}{3} - i\sigma -  \frac{1 + i}{24} \frac{ \sqrt{\pi } \left(6 \sigma ^2-6 \sigma
   +1\right)}{\sqrt{t}}\biggr) + O(t^{\sigma-2}).
\end{align*}
The lemma follows by subtracting this equation from (\ref{6.10}), employing the following equations which follow from (\ref{GLfinal}), (\ref{GU1asymptotics}), and (\ref{GU2final}) respectively:
\begin{align*}
& G_L(t, \sigma; \eta) =  -\frac{i\eta^{s-1}}{(2\pi)^s} \ln(1 - e^{i\eta}) + O\biggl(\frac{t}{\eta^{2-\sigma}}\biggr), \quad (1+\epsilon)t < \eta, \; \text{(not uniformly in $\epsilon$)},
	\\
& G_U^{(1)}(t, \sigma; t) =  \frac{t^{s-1} e^{-it} }{(2\pi)^{s}}\biggl(\frac{1-i}{2} \sqrt{\pi t} + \frac{i}{3} - i\sigma -  \frac{1 + i}{24} \frac{ \sqrt{\pi } \left(6 \sigma ^2-6 \sigma
   +1\right)}{\sqrt{t}}\biggr) + O(t^{\sigma-2}),
	\\
& G_U^{(2)}(t, \sigma; t) 
= \frac{i t^{s-1} e^{-it}}{(2\pi)^{s}}\ln(1 - e^{-it}) 
+ O(t^{\sigma - 2}),
\end{align*}
and using the identity
$$\frac{i t^{s-1} e^{-it}}{(2\pi)^{s}}\ln(1 - e^{-it}) 
+ \frac{t^{s-1}e^{-it}}{(2\pi)^s}2\im \Li_1(e^{it}))
= \frac{i t^{s-1} e^{-it}}{(2\pi)^{s}}\ln(1 - e^{it}).$$
\end{proof}

It is possible to derive an alternative integral representation for the related sum $\sum_{[\eta/2\pi]+1}^{[t/2\pi]} n^{s-1}$.

\begin{lemma}[An alternative integral representation for the basic sum]\label{principalvaluelemma}
For every $\delta > 0$,
\begin{align}\label{6.12}
 &  \sum_{n=[\frac{\eta}{2\pi}]+1}^{[\frac{t}{2\pi}]} n^{s-1} 
  = \frac{2}{(2\pi)^s} \dashint_{\eta}^t \frac{\rho^{s-1}}{e^{i\rho} - 1} d\rho + O(t^{\sigma -1}), 
  	\\ \nonumber
&\hspace{1cm} 0 < \eta < t, \quad \dist(\eta, 2\pi \Z) > \delta, \quad 
\dist(t, 2\pi \Z) > \delta, \quad 0 \leq \sigma \leq 1, \quad t \to \infty, 
\end{align}
where the error term is uniform for all $\eta, \sigma$ in the above ranges and the contour in the integral denotes the principal value integral with respect to the points
$$\left\{2\pi n \; \middle| \; n \in \Z, \; \Bigl[\frac{\eta}{2\pi}\Bigr] + 1\leq n \leq \Bigl[\frac{t}{2\pi}\Bigr]\right\}.$$
\end{lemma}
\begin{proof}
Let $\hat{C}_t^{\eta}$ denote the semicircle from $it$ to $i\eta$ with $\re z \leq 0$, defined in equation (\ref{2.12c}) with $\eta$ and $\alpha$ replaced by $t$ and $\eta$ respectively, see figure \ref{Cetatcontour.pdf}.
Let $R_\eta^t(s)$ and $L_t^\eta(s)$ denote the following: 
\begin{align}\label{Retatdef}
  R_\eta^t(s) = \frac{e^{-\frac{i\pi s}{2}}}{(2\pi)^s} \int_{C_\eta^t} \frac{z^{s-1}}{e^z - 1} dz, \qquad s \in \C
\end{align}
and
\begin{align}\label{Ltetadef}
  L_t^\eta(s) = \frac{e^{-\frac{i\pi s}{2}}}{(2\pi)^s} \int_{\hat{C}_t^\eta} \frac{z^{s-1}}{e^z - 1} dz, \qquad s \in \C.
\end{align}
\begin{figure}
\bigskip
 \begin{overpic}[width=.45\textwidth]{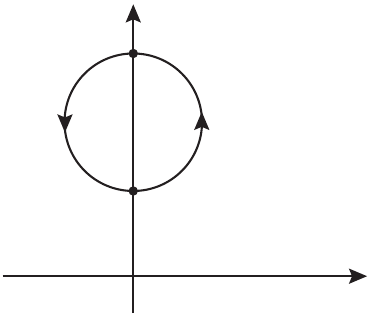}
  \put(39,73){$it$}
 \put(39,29){$i\eta$}
 \put(58,51){$C_\eta^t$}
 \put(6,51){$\hat{C}_t^{\eta}$}
 \put(102,9.5){$\re z$}
   \end{overpic}
   \caption{The contours $C_{\eta}^t$ and $\hat{C}_t^{\eta}$. }\label{Cetatcontour.pdf}
\end{figure}
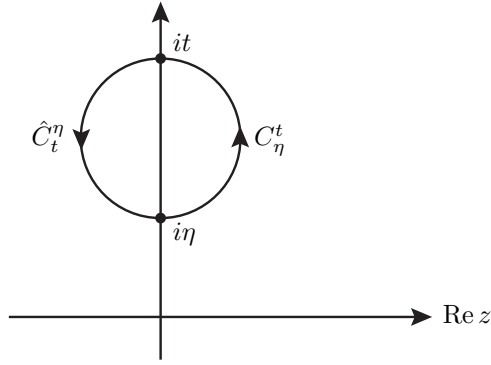
Cauchy's theorem applied in the interior of the disk whose boundary is $C_\eta^t \cup \hat{C}_t^\eta$, yields
\begin{align}\label{LplusR}
  L_t^\eta(s) + R_\eta^t(s) = \sum_{n = [\frac{\eta}{2\pi}] +1 }^{[\frac{t}{2\pi}]} n^{s-1}, \qquad s \in \C.
\end{align}
On the other hand, the Plemelj formulas imply
\begin{align}\label{RminusL}
R_\eta^t(s) - L_t^\eta(s) = \frac{2}{(2\pi)^s}\dashint_{\eta}^t \frac{\rho^{s-1}}{e^{i\rho} - 1} d\rho.
\end{align}
Using (\ref{RminusL}) to replace $R_\eta^t$ in (\ref{LplusR}), we find
\begin{align}\label{sumintLteta}
\sum_{n = [\frac{\eta}{2\pi}] +1 }^{[\frac{t}{2\pi}]} n^{s-1} = \frac{2}{(2\pi)^s}\dashint_{\eta}^t \frac{\rho^{s-1}}{e^{i\rho} - 1} d\rho + 2L_t^\eta.
\end{align}

Replacing the contour $\hat{C}_t^\eta$ in (\ref{Ltetadef}) by the union of the following three segments:
\begin{align*}
& \{ite^{i\theta} \; | \; 0 \leq \theta \leq \pi/4\}, 
\quad \{i\rho e^{\frac{\pi i}{4}} \; | \; \eta \leq \rho \leq t\},
\quad \{i\eta e^{i\theta} \; | \;  0 \leq \theta \leq \pi/4\},		
\end{align*}
we obtain
\begin{align*}
  L_t^\eta = &\; \frac{1}{(2\pi)^s}\biggl[ \int_0^{\frac{\pi}{4}} \frac{e^{i\theta(s-1)} t^{s-1}}{e^{ite^{i\theta}} -1} ite^{i\theta} d\theta
  - \int_\eta^t \frac{e^{\frac{\pi i}{4}(s-1)}\rho^{s-1}}{e^{i\rho e^{\frac{\pi i}{4}}} - 1} e^{\frac{\pi i}{4}} d\rho
  	\\
&  - \int_0^{\frac{\pi}{4}} \frac{e^{i\theta(s-1)} \eta^{s-1}}{e^{i\eta e^{i\theta}} - 1} i \eta e^{i\theta}  d\theta\biggr]
  	\\
=: &\; \frac{1}{(2\pi)^s}[J_1 + J_2 + J_3].	
\end{align*}
The assumption $\dist(t, 2\pi \Z) > \delta$ implies that 
$$|J_1| \leq \int_0^{\frac{\pi}{4}} \frac{e^{-\theta t} t^{\sigma-1}}{A} t d\theta
= \frac{t^{\sigma-1}}{A} (1 - e^{-\frac{\pi t}{4}}) = O(t^{\sigma -1}).$$
Similar computations show that $J_2$ is exponentially small and, in view of the assumption $\dist(\eta, 2\pi \Z) > \delta$, that $J_3$ is $O(\eta^\sigma t^{-1})$.
Thus,
\begin{align}\label{Ltetaestimate}
  L_t^\eta = O(t^{\sigma -1}),
\end{align}
and the lemma follows from (\ref{sumintLteta}).
\end{proof}

\begin{remark}\upshape
Lemma \ref{principalvaluelemma} implies that the leading behavior of $\zeta(s)$ is characterized by the following integral:
$$I(\eta, t, \sigma) = \dashint_\eta^t \frac{\rho^{s-1} }{e^{i\rho} - 1} d\rho.$$
\end{remark}

\begin{remark}\upshape
The estimate (\ref{Ltetaestimate}) is a consequence of the fact that the integral appearing in the definition of $L_t^\eta(s)$ does {\it not} possess any stationary points. Indeed, using
  $$\frac{1}{e^z - 1} = -\sum_{m=0}^\infty e^{mz}, \qquad \re z < 0,$$
it follows that
$$e^{-\frac{i\pi s}{2}} \int_{\hat{C}_t^\eta} \frac{z^{s-1} }{e^z - 1} dz
= \sum_{m=0}^\infty \int_\eta^t e^{im\rho + it\ln \rho} \rho^{\sigma -1} d\rho.$$
Candidates for stationary points occur at $\rho^* = -t/m$ and the inequality $\eta \leq\rho^* \leq t$ implies the non-existence of any stationary points.  
\end{remark}

The next lemma gives an alternative representation for the fundamental integral $R_\eta^t(s)$ defined in (\ref{Retatdef}).

\begin{lemma}\label{Retatlemma}
Let $C_\eta^t$ denote the semicircle from $i\eta$ to $it$ with $\re z \geq 0$ defined by equation (\ref{Cetatdef}). 
Then, for every $\delta >0$,
\begin{align}\label{Retatrepresentation}
&  R_\eta^t(s) 
= \frac{e^{-i\epsilon s}}{(2\pi)^s} \int_\eta ^t  \frac{u^{s-1}du}{e^{iue^{-i \epsilon}} - 1}
+ O\biggl(\frac{e^{\epsilon t} - 1}{t^{1-\sigma}}\biggr),
	\\ \nonumber
& 0 < \eta < t, \quad \dist(\eta, 2\pi \Z) > \delta, \quad 
\dist(t, 2\pi \Z) > \delta, \quad 0 \leq \sigma \leq 1, \quad 0 < \epsilon < 1, \quad t \to \infty, 
\end{align}
where the error term is uniform for all $\eta, \sigma, \epsilon$ in the above ranges.
\end{lemma}
\begin{proof}
Let $\epsilon \in (0, 1)$. We deform the contour $C_\eta^t$ on the lhs of (\ref{Retatrepresentation}) so that it consists of the following three pieces: 
$$\{i\eta e^{-i\theta} \; | \; 0 \leq \theta \leq \epsilon\}, \quad \{iu e^{-i\epsilon} \; | \; \eta \leq u \leq t\},
\quad \{it e^{-i\theta} \; | \; 0 \leq \theta \leq \epsilon\}.$$
This yields
\begin{align*}
e^{-\frac{i\pi s}{2}} \int_{C_\eta^t} \frac{z^{s-1}dz}{e^z - 1}
= &\; e^{-i\epsilon s} \int_\eta^t \frac{u^{s-1}du}{e^{iue^{-i\epsilon}}  - 1} 
 - i \eta^s \int_0^\epsilon \frac{e^{-i\theta s}d\theta}{e^{i \eta e^{-i\theta}} - 1} 
+ i t^s \int_0^\epsilon \frac{ e^{-i\theta s}d\theta}{e^{i t e^{-i\theta}} - 1}.
\end{align*}
The assumption that $\dist(\eta, 2\pi \Z) > \delta$ implies that there exists an $A$ such that $|e^{i\eta e^{-i\theta}} - 1| \geq A$ for all $\theta \in [0, \epsilon]$, $\epsilon \in (0,1)$, and $\eta > 0$.
Thus,
\begin{align*}
 \left|  \frac{e^{-i\theta s}}{e^{i\eta e^{-i\theta}} - 1}\right|
  \leq  \frac{e^{\theta t}}{A}, \qquad \theta \in [0, \epsilon],
\end{align*}
and so
\begin{align*}
 \left| \int_0^\epsilon \frac{e^{-i\theta s}}{e^{i \eta e^{-i\theta}}  - 1} d\theta \right|
  \leq \frac{e^{\epsilon t} -1}{A t}.
\end{align*}
Similarly, because of the assumption $\dist(t, 2\pi \Z) > \delta$,
\begin{align*}
 \left| \int_0^\epsilon \frac{e^{-i\theta s}}{e^{i t e^{-i\theta}} - 1} d\theta \right|
  \leq \frac{e^{\epsilon t} -1}{A t}.
\end{align*}
This proves the lemma. 
\end{proof}

Lemma \ref{Retatlemma} implies that
\begin{align}\label{Retatlim}  
  R_\eta^t(s) = \frac{1}{(2\pi)^s} \lim_{\substack{\epsilon \to 0 \\ \epsilon t \to 0}} \int_\eta ^t  \frac{u^{s-1}du}{e^{iue^{-i \epsilon}} - 1},
\end{align}
which in view of (\ref{LplusR}) and (\ref{Ltetaestimate}) implies the following alternative representation for the basic sum:
$$\frac{1}{(2\pi)^s} \lim_{\substack{\epsilon \to 0 \\ \epsilon t \to 0}} \int_\eta ^t  \frac{u^{s-1}du}{e^{iue^{-i \epsilon}} - 1} 
= \sum_{n = [\frac{\eta}{2\pi}] +1 }^{[\frac{t}{2\pi}]} n^{s-1} + O(t^{\sigma -1}).$$
In fact, in the remaining part of this chapter we will present arguments which suggest\footnote{The rigorous justification of some of these  arguments remains open.} that on the critical line the assumption $\epsilon t \to 0$ in (\ref{Retatlim}) can be relaxed to $\epsilon^2 t \to 0$, so that the following equation is valid:
\begin{align}\label{6.1} 
  R_\eta^t(s) = \frac{1}{(2\pi)^s} \lim_{\substack{\epsilon \to 0 \\ \epsilon^2 t \to 0}} \int_\eta^t  \frac{u^{s-1}du}{e^{iue^{-i \epsilon}} - 1},
\end{align}
where it is assumed that $\dist(\eta, 2\pi \Z) > \delta$, $\eta = O(1)$, $\dist(t, 2\pi \Z) > \delta$, and $\sigma = 1/2$. 

\begin{figure}
\bigskip
 \begin{overpic}[width=.6\textwidth]{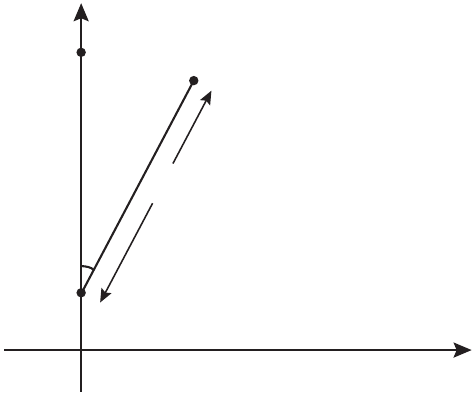}
 \put(11,71){$it$}
 \put(10.5,20.5){$i\eta$}
 \put(32,44){$t-\eta$}
 \put(18.3,29){$\epsilon$}
  \put(43,67){$z^*$}
   \end{overpic}
   \caption{The contour of integration of the integral in the rhs of equation (\ref{6.2}).}\label{zstarcontour.pdf}
\end{figure}
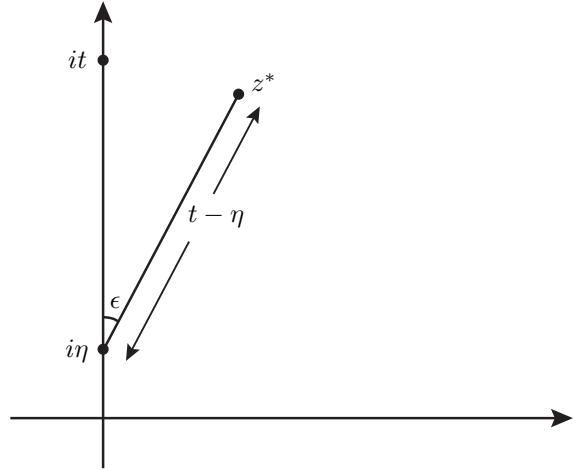

Indeed, consider the integral
$$(2\pi)^sR_\eta^t(s) = e^{-\frac{i\pi s}{2}}\int_{C_\eta^t} \frac{e^{-z}z^{s-1}dz}{1-e^{-z}}.$$
Instead of the contour $C^t_\eta$, we consider the finite ray from the point $i\eta$ to the point $z^*$, see figure \ref{zstarcontour.pdf}, where
$$z^*=i \eta +(t-\eta)e^{i(\frac{\pi}{2}-\epsilon)}=[\eta+(t-\eta)e^{-i \epsilon}]e^{\frac{i \pi}{2}}, \qquad 0 \leq \epsilon < \frac{\pi}{2}.$$
We claim that
\begin{align}\label{zstarintestimate}
& e^{-\frac{i\pi s}{2}} \int_{z^*}^{it} \frac{e^{-z}z^{s-1}dz}{1-e^{-z}}
= O(t^\sigma \epsilon e^{t\epsilon^3}).
\end{align}
Indeed, the assumption $\dist(t, 2\pi \Z) > \delta$ implies that there exists an $A>0$ independent of $\eta, \epsilon$ such that $|1-e^{-z}| > A$ on the contour from $z^*$ to $it$. 
Thus, using the parametrization $z = it e^{-i\theta}$, $0 \leq \theta \leq \epsilon$, in the lhs of (\ref{zstarintestimate}), we can prove (\ref{zstarintestimate}) as follows:
\begin{align*}
\biggl|e^{-\frac{i\pi s}{2}} \int_{z^*}^{it} \frac{e^{-z}z^{s-1}dz}{1-e^{-z}}\biggr|
& =  \biggl| i t^s \int_0^\epsilon \frac{e^{-i t e^{-i\theta}} e^{-i\theta s}d\theta}{1 - e^{-i t e^{-i\theta}}}\biggr|
\leq t^\sigma \int_0^\epsilon \frac{e^{-t \sin\theta} e^{\theta t}d\theta}{A}
	\\
& = O\biggl(t^\sigma\int_0^\epsilon e^{t\theta^3} d\theta \biggr)
= O(t^\sigma \epsilon e^{t\epsilon^3}).
\end{align*}
where we have used that
$$|\theta - \sin{\theta}| \leq \theta^3, \qquad 0 \leq \theta \leq 1.$$

It follows from (\ref{zstarintestimate}) that
\begin{align}\label{6.2}
  e^{-\frac{i\pi s}{2}}\int_{C_\eta^t} \frac{e^{-z}z^{s-1}dz}{1-e^{-z}}  = e^{-\frac{i\pi s}{2}} \lim_{\substack{\epsilon \to 0 \\ \epsilon^2 t \to 0}} \int^{z^*}_{\eta e^{\frac{i\pi}{2}}} \frac{e^{-z}z^{s-1}}{1-e^{-z}}dz.
\end{align}
We next employ the following parametrization which maps $z^*$ to $\rho=0$:
$$ z = \eta e^{\frac{i\pi}{2}} + t\left( 1 - \frac{\eta}{t} -\rho\right) e^{i\left( \frac{\pi}{2}-\epsilon\right)}=  \left[\eta + t \left( 1 - \frac{\eta}{t} -\rho\right) e^{-i\epsilon}\right] e^{\frac{i\pi}{2}},$$
$$ 0<\rho < 1 - \frac{\eta}{t}. $$
We replace in equation (\ref{6.2}) $e^{-z}$ with the expression
$$e^{-z} = e^{-ite^{-i\epsilon } + i\rho t e^{-i\epsilon}-i\eta(1- e^{-i\epsilon})} = e^{it(\rho-1) e^{-i\epsilon} + O(\eta \epsilon)}, \qquad \epsilon \to 0$$
and we also replace $z$ in the expression $z^{s-1}$ of equation (\ref{6.2}) with
\begin{align}\label{6.3}
z =t e^{\frac{i\pi}{2}} M e^{-i \mu} \left(1-\frac{\rho}{M}e^{-i(\epsilon-\mu)} \right),
\end{align}
where ($M,\mu$) are defined via the equation
\begin{align}\label{6.4}
1+ \left( \frac{t}{\eta} -1\right) e^{-i\epsilon} = \frac{t}{\eta} M e^{-i\mu}, \qquad M(t,\epsilon) >0, \quad \mu(t,\epsilon) >0.
\end{align}
Then, the rhs of equation (\ref{6.2}) is given by
\begin{subequations}
\begin{align}\label{6.5a}
t^s \lim_{\substack{\epsilon \to 0 \\ \epsilon^2 t \to 0}} e^{-i\epsilon} M^{s-1} e^{-i\mu(s-1)} \int^{1-\frac{\eta}{t}}_0
\frac{e^{it(\rho-1) e^{-i\epsilon}}}{1-e^{it(\rho-1)e^{-i\epsilon}}}
\left( 1 - \frac{\rho}{M}e^{-i\psi}\right)^{s-1}d\rho,
\end{align}
where
\begin{align}\label{6.5b}
\psi(t,\epsilon) = \epsilon - \mu.
\end{align}
\end{subequations}
The functions $M$ and $\mu$ are given by
\begin{subequations}
\begin{align}\label{6.6a}
M(t,\epsilon) = 1- \frac{\epsilon^2}{2} \left( \frac{\eta}{t} - \frac{\eta^2}{t^2}\right) + O\left( \epsilon^4\left( \frac{\eta}{t} - \frac{\eta^2}{t^2}\right)^2\right), \qquad \epsilon \to 0
\end{align}
and
\begin{align}\label{6.6b}
\mu(t,\epsilon) = \epsilon - \frac{\epsilon\eta}{t} + O(\epsilon^2), \qquad \epsilon \to 0.
\end{align}
\end{subequations}
\begin{figure}
\bigskip
 \begin{overpic}[width=.6\textwidth]{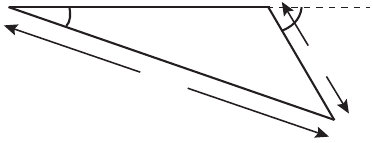}
 \put(2,39){$0$}
 \put(71,39){$1$}
 \put(21,32.3){$\mu$}
 \put(82,32){$\epsilon$}
 \put(83,21){$\frac{t}{\eta} - 1$}
 \put(40,15){$\frac{t}{\eta}M$}
   \end{overpic}
   \caption{The triangle $D_1$.}\label{triangleD1.pdf} 
\end{figure}
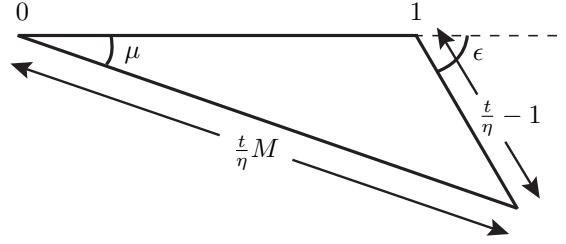
Indeed, consider the triangle $D_1$ formed by the intersection of the following finite rays, see figure \ref{triangleD1.pdf}:
$$(0,1), \quad \left( 1, 1+ \left( \frac{t}{\eta} - 1\right) e^{-i\epsilon}\right), \quad \left(0, \frac{t}{\eta} Me^{-i\mu}\right).$$
The cosine rule for the triangle $D_1$ implies the following identities:
\begin{align*}
\left( \frac{t}{\eta} M\right)^2 & = 1 + \left( \frac{t}{\eta} - 1\right)^2 - 2\left( \frac{t}{\eta} -1\right) \cos (\pi - \epsilon)
	\\ 
& = 1 + \left( \frac{t}{\eta} -1\right)^2 + 2\left( \frac{t}{\eta} -1\right) \cos \epsilon
	\\
& = \left[ 1 + \left( \frac{t}{\eta} - 1\right)\right]^2 - 2\left( \frac{t}{\eta} -1\right) (1 - \cos \epsilon).
\end{align*}
Thus,
$$ \frac{t}{\eta} M = \frac{t}{\eta} \sqrt{1 - \frac{2\left( \frac{t}{\eta}-1\right)(1-\cos\epsilon)}{\frac{t^2}{\eta^2}}},$$
or
$$ M = \sqrt{ 1 - 2\left( \frac{\eta}{t} - \frac{\eta^2}{t^2}\right)(1-\cos \epsilon)},$$
which yields (\ref{6.6a}).
The sine rule for the triangle $D_1$ implies the identity
$$ \frac{\sin\mu}{\frac{t}{\eta} -1} = \frac{\sin(\pi-\epsilon)}{\frac{t}{\eta}M}.$$
Hence,
$$ \sin \mu = \frac{1 - \frac{\eta}{t}}{M} \sin \epsilon,$$
which, using the expression (\ref{6.6a}) for $M$, yields (\ref{6.6b}).

Equation (\ref{6.6b}) implies that $\psi >0$, thus, we can define ($W,w$) via the equation
\begin{align}\label{6.7}
1 - \frac{\rho}{M}e^{-i\psi} = 1+ \frac{\rho}{M} e^{i(\pi-\psi)} = We^{iw}, \quad W(\rho,t,\epsilon)>0, \quad w(\rho,t,\epsilon)>0.
\end{align}
The functions $W$ and $w$ are given by
\begin{subequations}
\begin{align}\label{6.8a}
W(\rho,t,\epsilon) = 1 - \rho + O\left( \frac{\epsilon^2\eta}{t}\right) + O\left( \frac{\epsilon^2\eta^2}{t^2}\right), \qquad \epsilon \to 0
\end{align}
and
\begin{align}\label{6.8b}
w(\rho,t,\epsilon) =\frac{\epsilon \eta}{t} \frac{\rho}{1-\rho} \left[ 1 + O\left( \frac{\epsilon^2\eta}{t}\right) + O\left( \frac{\epsilon^2\eta^2}{t^2}\right)\right], \qquad \epsilon \to 0.
\end{align}
\end{subequations}
Indeed, consider the triangle $D_2$ formed by the intersection of the following finite rays, see figure \ref{triangleD2.pdf}:
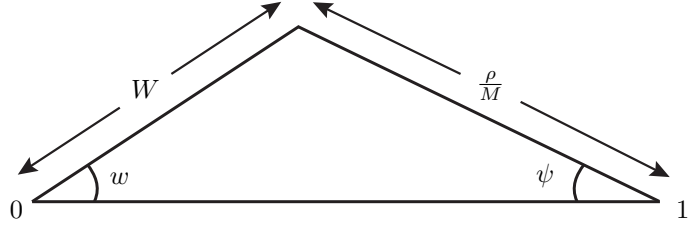
\begin{figure}
\bigskip
 \begin{overpic}[width=.7\textwidth]{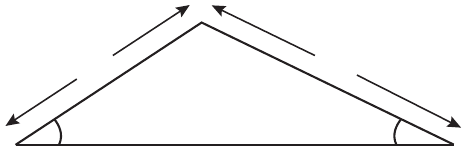}
 \put(0,0){$0$}
 \put(100,0){$1$}
 \put(15,5){$w$}
 \put(18.2,18){$W$}
 \put(79,5.3){$\psi$}
 \put(70,19){$\frac{\rho}{M}$}
   \end{overpic}
   \caption{The triangle $D_2$.}\label{triangleD2.pdf}
\end{figure}
$$ (0,1), \quad \left(1, 1 + \frac{\rho}{M} e^{i(\pi-\psi)}\right), \quad \left( 0,We^{iw}\right).$$
The cosine rule for the triangle $D_2$ implies the following identities:
\begin{align*}
W^2 & = 1 + \frac{\rho^2}{M^2} - 2 \frac{\rho}{M} \cos \psi
	\\
& = \left( 1 - \frac{\rho}{M}\right)^2 + \frac{2\rho}{M}(1-\cos\psi).
\end{align*}
Hence,
$$W = \left( 1 - \frac{\rho}{M}\right) \sqrt{ 1 + \frac{2\frac{\rho}{M}}{\left(1 - \frac{\rho}{M}\right)^2} (1-\cos\psi)}.$$

Thus, using the definition (\ref{6.5b}) of $\psi$, together with equation (\ref{6.6b}), we find
$$W = 1 - \frac{\rho}{M} + \frac{1}{2} \frac{\frac{\rho}{M}}{1 - \frac{\rho}{M}} \left( \frac{\epsilon^2\eta^2}{t^2} + O(\epsilon^4)\right), \qquad \epsilon \to 0.$$
Substituting in this equation the expression (\ref{6.6a}) for $M$, we find equation (\ref{6.8a}).

The sine rule for the triangle $D_2$ implies the identity
$$\frac{\sin w}{\rho/M} = \frac{\sin \psi}{W},$$
which, using the expressions for $\psi$ (equations (\ref{6.5b}) and (\ref{6.6b})), for $M$ (equation (\ref{6.6a})) and for $W$ (equation (\ref{6.8a})), yields equation (\ref{6.8b}).

Using the identity
$$ \left( 1 - \frac{\rho}{M} e^{-i\psi}\right)^{s-1} = W^{s-1} e^{i(s-1)w},$$
as well as the expressions for $M$, $\mu$, $W$ and $w$ obtained earlier, it is possible to compute all the expressions
appearing in the rhs of (\ref{6.5a}):
\begin{align*}
&M^{s-1} = e^{(it+\sigma-1)\ln M} = e^{- \frac{i\epsilon^2\eta}{2}}e^{O\big( \frac{\epsilon^2}{t}\big)}e^{O\big( \frac{\epsilon^2}{t^2}\big)},
	\\
& e^{-i(s-1)\mu} = e^{\mu t} e^{-i(\sigma-1)\mu} = e^{\epsilon t} e^{-\epsilon\eta} e^{O(\epsilon^2t)} e^{i(\sigma-1)\left(\frac{\epsilon \eta}{t}-\epsilon \right)} e^{O(\epsilon^2)},
	\\
& W^{s-1} = (1-\rho)^{s-1}e^{ (it+\sigma-1)\left[ O\big( \frac{\epsilon^2\eta}{t}\big) + O\big( \frac{\epsilon^2\eta^2}{t^2}\big)\right]} = (1-\rho)^{s-1} e^{O(\epsilon^2\eta)} e^{O\big( \frac{\epsilon^2\eta^2}{t}\big)},
	\\
& e^{iw(s-1)} = e^{-tw} e^{i(\sigma-1)w} = e^{-\epsilon\eta \frac{\rho}{1-\rho}} e^{\frac{\rho}{1-\rho}\left[ O(\epsilon) + O\left( \frac{\epsilon}{t}\right) \right]}.
\end{align*}
Substituting the above expressions in (\ref{6.5a}) and using the identities
$$ e^{-ite^{-i\epsilon}} = e^{-it(1-i\epsilon + O(\epsilon^2))} = e^{-it -\epsilon t + O(\epsilon^2t)},$$
in order to simplify the integrand of the integral in equation (\ref{6.5a}), we find the equation
\begin{align}\label{6.9}
e^{-\frac{i\pi s}{2}} \int_{C_\eta^t} \frac{e^{-z}z^{s-1}dz}{1-e^{-z}}=t^s  \lim_{\substack{\epsilon \to 0 \\ \epsilon^2 t \to 0}} \int_0^{1-\frac{\eta}{t}}\frac{e^{it(\rho-1)e^{-i\epsilon}}(1-\rho)^{s-1}d\rho}{1 - e^{it(\rho-1)e^{-i\epsilon}}}, 
\qquad  0<\eta<t. 
\end{align}
The transformation $\rho=1-\frac{u}{t}$ maps this equation to equation (\ref{6.1}).

We note that a detailed analysis in the neighborhood of $\rho=1$ shows that similar estimates are also valid near $\rho=1$.

\appendix
\chapter{The Asymptotics of $\Gamma(1-s)$ and $\chi(s)$} \label{appA}
The asymptotic formulae given in theorems \ref{zetath2} and \ref{ZETATH} involve the functions $\chi(s)$ and $\Gamma(1-s)$. In this appendix, we derive the asymptotics of these functions as $t \to \infty$.

We use the following well-known formula for the asymptotics of the Gamma function (see for example Olver \cite{O1974}, page 294):
\begin{align} \nonumber
& \Gamma(s) = e^{-s}e^{s\ln s} \left( \frac{2\pi}{s}\right)^{\frac{1}{2}} \left[ 1 + \frac{1}{12s} + O\left( \frac{1}{s^2}\right)\right], 
\qquad s\to \infty, \quad |\arg s| \leq \pi - \delta,
\end{align}
where $\delta >0$. Hence, as $t \to \infty$,
\begin{align}\nonumber
\Gamma(s) & = \sqrt{2\pi} e^{-\sigma-it} e^{s\ln [(it)\left( 1 + \frac{\sigma}{it}\right)]} \frac{1 + \frac{1}{12it\left( 1 + \frac{\sigma}{it}\right)}
+ O\left( \frac{1}{s^2}\right)}{(it)^{\frac{1}{2}}\left(1 + \frac{\sigma}{it}\right)^{\frac{1}{2}}}
	\\ \nonumber
& = \frac{\sqrt{2\pi}}{\sqrt{t}} \left( e^{\frac{i\pi}{2}}t\right)^s e^{-\sigma-it} e^{(\sigma+it)\left[ \frac{\sigma}{it} + \frac{\sigma^2}{2t^2} + O\left( \frac{1}{t^3}\right)\right]} e^{- \frac{i\pi}{4}}
	\\ \label{Gammaasymptotics}
& \qquad \times \left[ 1 + \frac{1}{12it} + O\left( \frac{1}{t^2}\right)\right] \left[ 1 - \frac{\sigma}{2it} + O\left( \frac{1}{t^2}\right) \right]
	\\ \nonumber
& = \sqrt{2\pi} t^{s-\frac{1}{2}} e^{\frac{i\pi s}{2}} e^{- \frac{i\pi}{4}} e^{-it} e^{-\frac{i\sigma^2}{2t} + O\left( \frac{1}{t^2}\right)}
 \left[ 1 + \frac{1 - 6\sigma}{12 i t} + O\left( \frac{1}{t^2} \right)\right] 
	\\ \nonumber
& = \sqrt{2\pi} t^{s-\frac{1}{2}} e^{\frac{i\pi s}{2}} e^{- \frac{i\pi}{4}} e^{-it} \left[ 1 - \frac{i\sigma^2}{2t} + O\left( \frac{1}{t^2}\right)\right]  \left[ 1 + \frac{1 - 6\sigma}{12 i t} + O\left( \frac{1}{t^2} \right)\right] 
	\\ \nonumber
& =\sqrt{2\pi} t^{s-\frac{1}{2}} e^{\frac{i\pi s}{2}} e^{- \frac{i\pi}{4}} e^{-it} \left[ 1 + \frac{ic(\sigma)}{t} + O\left( \frac{1}{t^2}\right) \right], \qquad 0 \leq \sigma \leq 1, \quad t \to \infty,
\end{align}
where $c(\sigma)$ is defined by
$$c(\sigma) = \frac{\sigma}{2}(1-\sigma) - \frac{1}{12}.$$
Replacing $\sigma$ by $1-\sigma$ in (\ref{Gammaasymptotics}) and taking the complex conjugate of the resulting equation, we find
\begin{align}\nonumber
  &\Gamma(1-s) = \sqrt{2\pi} t^{\frac{1}{2}-s} e^{\frac{-i\pi (1-s)}{2}} e^{\frac{i\pi}{4}} e^{it} \left[ 1 - \frac{ic(\sigma)}{t} + O\left( \frac{1}{t^2}\right) \right], 
  	\\ \label{Gamma1minuss}
& \hspace{6cm}   0 \leq \sigma \leq 1, \quad t \to \infty.
\end{align}
On the other hand, the definition (\ref{2.2}) of $\chi(s)$ implies
$$ (2\pi)^s\chi(1-s) = \frac{(2\pi)^s}{2i} \frac{(2\pi)^{1-s}}{\pi} \Gamma(s) \left[ e^{\frac{i\pi}{2}(1-s)} - e^{- \frac{i\pi}{2}(1-s)}\right] = \left( e^{-\frac{i\pi s}{2}} + e^{\frac{i\pi s}{2}}\right) \Gamma(s).$$
Replacing in this equation $\Gamma(s)$ by the rhs of (\ref{Gammaasymptotics}) we find 
\begin{align}\label{chi1minussasymptotics}
\chi(1-s) = (2\pi)^{\frac{1}{2}-s} t^{s-\frac{1}{2}} e^{- \frac{i\pi}{4}} e^{-it} \left[ 1 + \frac{ic(\sigma)}{t} + O\left( \frac{1}{t^2}\right)\right], \qquad 0 \leq \sigma \leq 1, \quad t \to \infty.
\end{align}
Equation (\ref{chi1minussasymptotics}) together with the identity (\ref{2.15}) imply the following asymptotic expression for $\chi(s)$:
$$\chi(s)=(2\pi)^{s-\frac{1}{2}} t^{\frac{1}{2}-s}e^{\frac{i \pi}{4}}e^{it} \left[1-\frac{i c(\sigma)}{t}+O \left( \frac{1}{t^2}\right) \right ], \qquad 0 \leq \sigma \leq 1, \quad t \to \infty.$$
It is of course straightforward to extend the above asymptotic formulae to higher order.

\chapter{Numerical Verifications} \label{appB}

\section{Verification of Theorem  \ref{th3.1}}
Letting $\sigma = 1/2$, the error term in equation (\ref{zetaformula2.1b}) is given by

\begin{center}
\begin{tabular}{|c|c|c|c|}
 & $t = 10$ & $t = 10^2$ & $t = 10^3$  \\
  \hline
$\eta = t^2$ &  $-(10.4 + 5.22i)\times 10^{-5}$ & $(-10.2 + 2.97i)\times10^{-9}$ & $(15.1 - 4.46i)\times 10^{-13}$ \\
	\hline
$\eta = t^{\frac{3}{2}} $ &  $-(4.00 + 4.19i)\times 10^{-3}$ & $-(1.40 + 1.11i)\times10^{-5}$ & $-(7.80+ 9.81 i)\times 10^{-8}$ 
\\ \hline
\end{tabular}
\medskip
\end{center}

Note that the error is proportional to $t^3/\eta^{3+\sigma}$ as expected.

In order to demonstrate the effect of the higher order terms in (\ref{zetaformula2.1b}), we also consider the difference between $\zeta(s)$ and the first term on the rhs of (\ref{zetaformula2.1b}), i.e., the difference
$$\zeta(s) -  \sum_{n =1}^{[\frac{\eta}{2\pi}]} n^{-s}.$$
For $\sigma = 1/2$, this difference is given by

\begin{center}
\begin{tabular}{|c|c|c|c|}
 & $t = 10$ & $t = 10^2$ & $t = 10^3$  \\
  \hline
$\eta = t^2$ &  $-0.291 + 0.274i$ & $-0.341 + 0.207 i$ & $-0.380 + 0.121i$ \\
	\hline
$\eta = t^{\frac{3}{2}} $ &  $0.266 + 0.0471i$ & $0.127 + 0.020 i$ & $0.0360 + 0.0612 i$ 
\\ \hline
\end{tabular}
\medskip
\end{center}

Similarly, the difference between $\zeta(s)$ and the first two terms on the rhs of (\ref{zetaformula2.1b}), i.e., the difference
$$\zeta(s) - \Biggl\{\sum_{n =1}^{[\frac{\eta}{2\pi}]} n^{-s} -  \frac{1}{1-s} \left(\frac{\eta}{2\pi}\right)^{1-s}\Biggr\},$$
is given by

\begin{center}
\begin{tabular}{|c|c|c|c|}
 & $t = 10$ & $t = 10^2$ & $t = 10^3$  \\
  \hline
$\eta = t^2$ &  $-(8.27 + 6.52i)\times 10^{-2}$ & $-(6.95 + 10.3i)\times10^{-4}$ & $-(3.40 + 10.6i)\times 10^{-4}$ \\
	\hline
$\eta = t^{\frac{3}{2}} $ &  $(1.58 - 1.50i)\times 10^{-1}$ & $(9.37 - 26.0i)\times 10^{-3}$ & $(-4.93 + 3.32i)\times 10^{-3}$ 
\\ \hline
\end{tabular}
\medskip
\end{center}

The above example illustrates that the inclusion of the higher order terms in (\ref{zetaformula2.1b}) improves the convergence of the asymptotic series  considerably.

\section{Verification of Theorem  \ref{th3.2}}
The error term in equation (\ref{zetaformula2.2b}) is given by

\begin{center}
\begin{tabular}{|c|c|c|c|}
 & $t = 10$ & $t = 10^2$ & $t = 10^3$  \\
  \hline
$\sigma = 0$ &  $(1.97 - 3.81i)\times 10^{-3}$ & $(-7.76 + 65.1i)\times10^{-7}$ & $(5.62 - 3.40i)\times 10^{-9}$ \\
	\hline
$\sigma = 1/2$ &  $(2.23 - 4.34i)\times 10^{-3}$ & $(-2.74 + 23.5i)\times10^{-7}$ & $(6.46-3.82i)\times 10^{-10}$ \\
	\hline
$\sigma = 1$ &  $(2.42 - 4.78i)\times 10^{-3}$ & $(-9.41 + 82.5i)\times10^{-8}$ & $(7.13 - 4.22i)\times 10^{-11}$ 
\\ \hline
\end{tabular}
\medskip
\end{center}

Note that the error is proportional to $t^{-\sigma-3}$ as expected.

\section{Verification of Theorem \ref{zetath2}}
Letting $\sigma = 1/2$, the error term in equation (\ref{zetaformula2b}), i.e., the term
$$e^{-i\pi s}\Gamma(1-s) e^{-\frac{\pi t}{2}}\eta^{\sigma-1}
\times \begin{cases}	
O\bigl(\frac{\eta}{t}\bigr), & 1 < \eta < t^{\frac{1}{3}} < \infty, \\
  O\bigl(e^{- \frac{At}{\eta^2}} + \frac{\eta^4}{t^2}\bigr), & t^{\frac{1}{3}} < \eta < \sqrt{t} < \infty, \end{cases} $$
   is given by

\begin{center}
\begin{tabular}{|c|c|c|c|}
 & $t = 10^2$ & $t = 10^4$ & $t = 10^6$  \\
  \hline
$\eta = 10$ &  $(3.04 + 7.27i)\times 10^{-3}$ & $(8.05 - 2.53i)\times10^{-6}$ & $(84.1 - 5.12i) \times 10^{-10}$ \\
	\hline
$\eta = t^{1/4}$ &  $(45.4 - 6.75i)\times 10^{-5}$ & $(8.05 - 2.53i) \times10^{-6}$ & $(-443.6 + 6.91i) \times 10^{-7}$ \\
	\hline
$\eta = t^{5/12}$ &  $7.27 - 1.45i) \times 10^{-1}$ & $(-8.69 + 9.53i)\times10^{-5}$ & $(4.38 - 13.7i)\times 10^{-6}$ 
\\ \hline
\end{tabular}
\medskip
\end{center}

In order to illustrate the increased accuracy achieved by including the higher order terms in (\ref{zetaformula2b}), we also display the corresponding table when only the first two sums on the rhs of (\ref{zetaformula2b}) are included: For $\sigma = 1/2$, the difference\footnote{This difference is exactly the error in the ``approximate functional equation'' of Hardy and Littlewood cf. \cite{HL1929}.}
$$\zeta(s) -   \sum_{n=1}^{[\frac{t}{\eta}]} \frac{1}{n^s} 
- \chi(s) \sum_{n=1}^{[\frac{\eta}{2\pi}]} \frac{1}{n^{1-s}}$$
is given by
\begin{center}
\begin{tabular}{|c|c|c|c|}
 & $t = 10^2$ & $t = 10^4$ & $t = 10^6$  \\
  \hline
$\eta = 10$ &  $(5.55 - 14.1i)\times 10^{-2}$ & $(-14.4 + 7.92i) \times10^{-3}$ & $(-15.9+4.40i)\times 10^{-4}$ \\
	\hline
$\eta = t^{1/4}$ &  $(4.76 - 7.53i) \times 10^{-2}$ & $(-14.4 + 7.92i)\times10^{-3}$ & $-(26.3 + 6.92i)\times 10^{-3}$ \\
	\hline
$\eta = t^{5/12}$ &  $-(2.14 + 2.15i)\times 10^{-1}$ & $(3.09 - 1.91i) \times10^{-2}$ & $(-5.70 + 10.4i)\times 10^{-3}$ 
\\ \hline
\end{tabular}
\medskip
\end{center}

\section{Verification of Corollary \ref{zetacor2}}
Letting $\sigma = 1/2$ and $N = 2$, we find that the error term in equation (\ref{zetacor2formula}) is given by

\begin{center}
\begin{tabular}{|c|c|c|c|}
 & $t = 10$ & $t = 10^2$ & $t = 10^3$  \\
  \hline
$\eta = t^{\frac{7}{12}}$ &  $-(4.08 + 2.47i)\times 10^{-1}$ & $(-2.39 + 3.50i)\times10^{-1}$ & $(-6.42 + 15.9i)\times 10^{-2}$ \\
	\hline
$\eta = t^{\frac{3}{4}}$ &  $(6.92 + 55.7i) \times 10^{-2}$ & $(9.93 + 25.7i)\times10^{-2}$ & $-(8.75+12.2i)\times 10^{-3}$ \\
	\hline
$\eta = 10t$ &  $(1.40 - 1.12 i)\times 10^{-2}$ & $(2.24 + 5.28i)\times10^{-4}$ & $(2.45 + 186.1i)\times 10^{-7}$ 
\\ \hline
\end{tabular}
\medskip
\end{center}


\section{Verification of Theorem \ref{ZETATH}}
Letting $\sigma = 1/2$, the error term in equation (\ref{zetaformulab}), i.e., the term
$$e^{-i\pi s}\Gamma(1-s) e^{-\frac{\pi t}{2}}\eta^{\sigma-1}
\times \begin{cases}	
O\bigl(\frac{\eta}{t}\bigr), & 1 < \eta < t^{\frac{1}{3}} < \infty, \\
  O\bigl(e^{- \frac{At}{\eta^2}} + \frac{\eta^4}{t^2}\bigr), & t^{\frac{1}{3}} < \eta < \sqrt{t} < \infty, \end{cases} $$
   is given by
\begin{center}
\begin{tabular}{|c|c|c|c|}
 & $t = 10^2$ & $t = 10^4$ & $t = 10^6$  \\
  \hline
$\eta = \sqrt{\frac{2\pi t}{100}}$ &  $(-5.96 + 8.83i)\times 10^{-4}$ & $(6.44 + 1.30i)\times10^{-4}$ & $(4.64 + 69.7i) \times 10^{-7}$ \\
	\hline
$\eta = \sqrt{2\pi t}$ &  $(3.59 - 10.8i)\times 10^{-3}$ & $(2.72 - 28.1i) \times10^{-5}$ & $(104.6 + 8.97i) \times 10^{-7}$ \\
	\hline
$\eta = \sqrt{200\pi t}$ &  $(13.7 - 5.21i) \times 10^{-6}$ & $(-9.02 + 2.27i)\times10^{-4}$ & $-(14.2 + 6.74i)\times 10^{-7}$ 
\\ \hline
\end{tabular}
\medskip
\end{center}

\begin{remark}\upshape
In the particular case of $\eta = \sqrt{2\pi t}$, the formula of Siegel (Eq. (32) of \cite{S1932}) and the corresponding formula of Titchmarsh (Theorem 4.16 of \cite{T1986}) are equivalent to the formula of our theorem \ref{ZETATH}. This equivalence can easily be checked numerically---for example, for $N = 3$ all three formulas yield the same numerical values for the error terms. In this regard, we note that Titchmarsh states Siegel's formula in terms of the function $\Psi(a)$ defined by
$$\Psi(a) = \frac{e^{-i\pi(\frac{a^2}{2} - \frac{5}{8})}}{2\pi} \int_L \frac{e^{\frac{iw^2}{4\pi}} + aw}{e^w -1} dw = \frac{\cos{\pi(\frac{a^2}{2} - a - \frac{1}{8})}}{\cos{\pi a}},$$
which is related to our function $\Phi(\tau, u)$ defined in (\ref{Phidef}) by
$$\Psi(a) = -i\Phi\Bigl(-1, a- \frac{1}{2}\Bigr) e^{-i\pi(\frac{a^2}{2} - \frac{5}{8})}, \qquad a \in \C.$$
\end{remark}

\backmatter
\bibliographystyle{amsalpha}
\bibliography{is}

\printindex

\end{document}